\documentclass[12pt,a4paper]{article}
\usepackage[utf8]{inputenc}
\usepackage[english]{babel} 
\usepackage{amsfonts, amsthm}
\usepackage{mathrsfs}
\usepackage{setspace}
\usepackage{extsizes}
\usepackage{colortbl}
\usepackage[left=2cm,right=2cm,top=2cm,bottom=1.5cm,bindingoffset=0cm]{geometry}
\usepackage{graphicx, color}
\graphicspath{ {./images/} }

\theoremstyle{definition}
\newtheorem{definition}{Definition}

\newtheorem{remark}{Remark}

\theoremstyle{plain}
\newtheorem{theorem}{Theorem}
\newtheorem{prop}{Proposition}

\usepackage{hyperref}

\begin{document}
\begin{center}
    {\large\bf Concrete examples of the rate of convergence of Chernoff approximations: numerical results for the heat semigroup and open questions on them (with appendix: full list of pictures and Python code)}


\medskip
{\bf K.A.~Katalova~(Dragunova),$^{1}$ N.~Nikbakht,$^{2}$  I.D.~Remizov$^{1,3}$} 
\medskip

\end{center}

\footnotesize

$^{1}$ National Research University Higher School of Economics. Address: 25/12 Bol. Pecherskaya Ulitsa, Room 224, Nizhny Novgorod, 603155, Russia.

$^{2}$ Department of Mathematics, Faculty of Science, University of Auckland. Address: 23 Symonds Street, Auckland 1010, New Zealand.

$^{3}$ Institute for Information Transmission Problems of the Russian Academy of Sciences. Address: Bolshoy Karetny per. 19, build.1, Moscow 127051 Russia. 

\textbf{Contact information:}

Ksenia Aleksandrovna Katalova (Dragunova). Email: k.dragunova13@mail.ru

Nasrin Nikbakht. Email: nasrin.nikbakht@gmail.com

Corresponding author: Ivan Dmitrievich Remizov. Email: ivremizov@yandex.ru

\normalsize

\medskip



\medskip
 
{\bf Abstract.} The article is devoted to the construction of examples that illustrate (using computer calculations) the rate of convergence of Chernoff approximations to the solution of the Cauchy problem for the heat equation. We are interested in the Chernoff theorem in general and select the heat semigroup as a model case because this semigroup (and solutions of the heat equations) are known, so it is easy to measure the distance between the exact solution and its Chernoff approximations. Two Chernoff functions (of the first and second order of Chernoff tangency to the generator of the heat semigroup, i.e. to the operator of taking the second derivative) and several initial conditions of different smoothness are considered. From the numerically plotted graphs, visually, it is determined that the approximations are close to the solution. For each of the two Chernoff functions, for several initial conditions of different smoothness and for approximation numbers up to 11 inclusive, the error (i.e. the supremum of the absolute value of the difference between the exact solution and the approximating function) corresponding to each approximation was numerically found. As it turned out, in all the cases studied, the dependence of the error on the number of the approximation has an approximately power-law form (we call this power the order of convergence). This follows from the fact that, as we discovered, the dependence of the logarithm of the error on the logarithm of the approximation number is approximately linear. Using the considered family of initial conditions, an empirical dependence of the order of convergence on the smoothness class of the initial condition is found. The orders of convergence for all the initial conditions studied are collected in a table. We also found some behavior of the error that we summarize in the conclusion and that we do not know how to explain with the existing theory, which opens a field of work for further theoretical research.

\newpage
\tableofcontents
\normalsize

\newpage
\section{Introduction}

$C_0$-semigroup of operators is a generalization of the notion of an exponent. For a bounded linear operator $L$ in the Banach space and real number $t$ one can define the exponent $e^{tL}$ via standard power series $\sum_{k=0}^\infty (tL)^k/k!$, which converges with respect to the operator norm to the bounded linear operator which is called $e^{tL}$. For an unbounded operator $L$ this procedure fails, a reasonable meaning of exponent $e^{tL}$ can still be introduced, and the collection of exponents $(e^{tL})_{t\geq 0}$ is called \cite{EN2, EN1} the $C_0$-semigroup of operators with infinitesimal generator $L$. The formal definition will follow below; for now let us mention only that for the $C_0$ -semigroup we have standard exponent properties: $e^{(t_1+t_2)L}=e^{t_1L}e^{t_2L}$, $e^{0L}=I$, $\frac{d}{dt}e^{tL}=Le^{tL}$. 

This is why the solution of the Cauchy problem $[U'(t)=LU(t), U(0)=u_0]$ in the Banach space is given by the equality $U(t)=e^{tL}u_0$. If $L$ is a differential operator, then this ordinary differential equation $U'(t)=LU(t)$ for the Banach-valued function $U$ can be interpreted as a partial differential equation $u_t(t,x)=Lu(t,x)$ for the number-valued function $u$. In particular, if $L=\partial$ is the differentiation operator, then $U'(t)=LU(t)$ can be interpreted as the transport equation $u_t(t,x)=u_x(t,x)$. If $L=\partial^2$ is the operator of taking the second derivative, then $U'(t)=LU(t)$ can be interpreted as the heat equation $u_t(t,x)=u_{xx}(t,x)$. For more profound operators $L$, equation $U'(t)=LU(t)$ can be interpreted \cite{EN2, EN1} as a parabolic equation with variable coefficients, Schr\"odinger-type equation, etc. 

For unbounded operators $L$ the definition of $e^{tL}$ does not give a rule for calculation of $e^{tL}$ using $L$ explicitly, e.g. as in power series. However, $e^{tL}$ can be calculated approximately with an arbitrary small error if the so-called \cite{Butko2020, R-2017} operator-valued Chernoff function $C$ for operator $L$ is found: we have $e^{tL}=\lim_{n\to\infty}C(t/n)^n$ in the strong operator topology. Joining this remark with the previous remarks, we come to the conclusion that Chernoff approximations are a flexible and powerful tool of functional analysis \cite{Chernoff, EN2, EN1}, which can be used, in particular, to find numerically approximate solutions of some differential equations with variable coefficients. 

For a given Cauchy problem for the linear evolution equation $[u_t(t,x)=Lu(t,x), u(0,x)=u_0(x)$, the Chernoff approximation method generates a sequence of functions $u_n(t,x)=(C(t/n)^nu_0)(x)$ that converges to the exact solution $u(t,x)=(e^{tL}u_0)(x)=U(t)(x)$ of the equation studied. For an arbitrary fixed moment of time $t$ functions $x\longmapsto u(t,x)$ and $x\longmapsto u_n(t,x)$ are elements of some Banach space, and Chernoff's theorem guarantees that $\|u(t,\cdot)-u_n(t,\cdot)\|\to 0$ as $n\to\infty$. 

Chernoff approximations of $C_0$-semigroups are used not only in the approximate solving of differential equations. They have other applications, both theoretical and practical; some of them are as follows:

\begin{itemize}
\item Definition and calculation of Feynman path integrals. In fact, if $C(t)$ is an integral operator, then $C(t/n)^n$ is an $n$-tuple integral operator. The limits of multiple integrals as multiplicity tends to infinity (such expressions are called Feynman formulas \cite{R-2017}) are one of the ways (which was originally used by Richard Feynman \cite{F1,F2}) to define the Feynman path integral \cite{Maz}.

\item Obtaining Feynman-Kac formulas as a corollary of Feynman formulas which were proven with the use of the Chernoff theorem \cite{BGS2010}. It is also interesting to compare this with an approach to the Feynman-Kac formula without Chernoff approximations \cite{Fee}. 

\item Statistical modeling. See the paper \cite{K2023} that numerically investigates the efficiency of the Monte Carlo method based on the application of the Chernoff theorem, and the papers \cite{KOS2022, KOS2022} that mathematically substantiate such an approach to the Chernoff approximations.

\item The law of large numbers for one-parameter semigroups \cite{OSS2019}.

\item Averaging of random Hamiltonians \cite{OSS2014} and averaging of random affine transformations \cite{KOS2022, KOS2023}.

\item Construction of measures on infinite-dimensional spaces \cite{SO2021}

\item Parabolic equations (in $\mathbb{R}^1$ \cite{R4,Vedenin2026}, in $\mathbb{R}^n$ \cite{Rjmp}, in Hilbert space \cite{R2018}, on manifolds \cite{MMRS2023})

\item Schr\"odinger equation  \cite{R-2017,RS2018,R-PotAn2020, GP2021}

\item Approximate calculation of resolvents of generators of $C_0$-semigroups, and of solutions to ODEs and elliptic PDEs \cite{Rarx2023}

\end{itemize}

A brief history and an overview of the results obtained up to 2017 in the construction of Chernoff approximations of $e^{tL}$ for several classes of operators $L$ can be found in \cite{Butko2020}. Several papers on the topic showing the diversity of cases studied are \cite{R2018, R-2017, RS2018, R-PotAn2020, R5, R1, Rjmp}, see also \cite{OSS2019, ST2020}.

To our current knowledge, all contributions to a rather young ``theory of rates of convergence in Chernoff's theorem'' can be found in \cite{Prudnikov2020, Bent2004, Bent2009,NSZ2018, GSK2019,Vedenin2020, Zagrebnov2020, GR2021, Galkin2022, VVGKR, Dragunova2023, GR2024} and references therein. These papers provide estimates for the rate of convergence under some conditions, but if these conditions are not satisfied, then one can say nothing about the quality of Chernoff approximations. There are also very few ``practical'' research papers \cite{OSS2012, Prudnikov2020} that measure the speed of convergence in particular cases obtained by numerical simulations. In our research, we continue to contribute to this field of study. Our long-term goal is to study the convergence rate in the Chernoff theorem. That is, our aim is to gather information on the convergence of Chernoff approximations for all possible semigroups, all possible Chernoff functions, and all possible vectors in a general setting. However, theoretical knowledge in this field is limited, so we collected data numerically for some particular semigroup (heat semigroup), two concrete Chernoff functions, and several vectors - functions of different smoothness. We do not have a clear explanation of each fact that we have found, however we have some guesses based on our observations. We believe that collecting reports on numerical measurements of the rate of convergence in the Chernoff theorem in concrete cases will help to state some conjectures and prove some theorems in the future.

We consider initial value problem for the heat equation
\begin{equation}\label{ZK}
\left\{ \begin{array}{ll}
   u_{t}(t,x)=u_{xx}(t,x) \textrm{ for }t>0, x\in\mathbb{R}^1  \\
   u(0,x)=u_0(x) \textrm{ for }x\in\mathbb{R}^1
\end{array} \right.
\end{equation}
which is a good model example because its bounded solution $u(t,x)$ is already known for each uniformly continuous bounded initial condition $u_0$ and is given by the formula
\begin{equation}\label{exact}
u(t,x)=\int_{\mathbb{R}}\Phi(x-y,t)u_0(y)dy, \textrm{ where }\Phi(x,t)=(2\sqrt{\pi t})^{-1}\exp\left(\frac{-x^2}{4t} \right).
\end{equation}

Then we use Chernoff functions $G$ and $S$ defined below in the paper to obtain Chernoff approximations $u_n(t,x)=(C(t/n)^nu_0)(x)$, $C\in\{G,S\}$ to the exact solution $u(t,x)=(e^{t\partial^2}u_0)(x)$ for $n=4,5,\dots,11$ and fixed time $t=1/2$. We are interested in large $n$ and agree that $n=1,2,3$ is not large. For brevity, we use the same letter $u_n$ to denote approximations generated by Chernoff functions $G$ and $S$, however, these approximations are different but never meet in the same expression, so this does not create ambiguity. Chernoff function $G$ has the first and the Chernoff function $S$ has the second order of the Chernoff tangency (the meaning of these words will be given in the main text of the paper) to the operator $\partial^2$, and we aim to find the relationship of this fact with the rate of convergence. One should expect that $S$ will provide a smaller error than $G$, but we consider the initial conditions $u_0$ with properties which are not covered by the theory that exists in 2025, e.g. $u_0\notin D(L)$.

\textbf{Main result of the present paper.} Numerically we find the approximation error $\sup_{x\in\mathbb{R}}|u(t,x)-u_n(t,x)|$ which is the norm of the difference between the Chernoff approximations and the heat semigroup. Then we use linear regression (ordinary least squares method) and discover that 
$$
\ln\sup_{x\in\mathbb{R}}|u(t,x)-u_n(t,x)|\approx \alpha \ln n + \beta
$$
which implies
$$
\sup_{x\in\mathbb{R}}|u(t,x)-u_n(t,x)|\approx n^\alpha\cdot e^\beta
$$
with a reasonable accuracy (coefficient of determinance $R^2>0.98$). We see that coefficients $\alpha$ and $\beta$ depend on the smoothness of the initial condition $u_0$ and on the Chernoff function used in the construction of Chernoff approximations. However, the dependence of $\alpha$ and $\beta$ on the initial condition has some other features that we do not know how to explain; they are listed in the end of the paper.

P.S.~Prudnikov in 2020 studied \cite{Prudnikov2020} this question in a similar setting, but his approach does not allow for a direct generalization. Meanwhile, the simulation method that we use allows us to study not only the heat equation, but also the equations with variable coefficients. In addition, we consider more initial conditions than were studied in \cite{Prudnikov2020}.

Now let us provide the necessary background on the topic to explain the notion of Chernoff tangency and the Chernoff operator-valued function that are important to understand how we obtain Chernoff approximations $u_n(t,x)$.

\section{Preliminaries}

\subsection{$C_0$-semigroups and their Chernoff approximations}

\begin{definition}\label{C0sgdef}
Let $\mathcal{F}$ be a Banach space. Let $\mathscr{L}(\mathcal{F})$ be the set of all bounded linear operators in $\mathcal{F}$. Suppose that we have a mapping $V\colon [0,+\infty)\to \mathscr{L}(\mathcal{F})$, that is, $V(t)$ is a linear bounded operator $V(t)\colon \mathcal{F}\to \mathcal{F}$ for each $t\geq 0.$ The mapping $V$ is called \cite{EN1} a \textit{$C_0$-semigroup}, or \textit{a strongly continuous one-parameter semigroup of operators} iff it satisfies the following conditions: 
	
	1) $V(0)$ is the identity operator $I$, i.e. $\forall \varphi\in \mathcal{F}: V(0)\varphi=\varphi;$ 
	
	2) $V$ maps the addition of numbers in $[0,+\infty)$ to the composition of operators in $\mathscr{L}(\mathcal{F})$, that is, $\forall t\geq 0,\forall s\geq 0: V(t+s)=V(t)\circ V(s),$ where for each $\varphi\in\mathcal{F}$ the notation $(A\circ B)(\varphi)=A(B(\varphi))=AB\varphi$ is used;
	
	3) $V$ is continuous with respect to the strong operator topology in $\mathscr{L}(\mathcal{F})$, i.e. $\forall \varphi\in \mathcal{F}$ function $t\longmapsto V(t)\varphi$ is continuous as a mapping $[0,+\infty)\to \mathcal{F}.$
	
The definition of a \textit{$C_0$-group} is obtained by the substitution of $[0,+\infty)$ by $\mathbb{R}$ in the paragraph above.

\end{definition}

It is known \cite{EN1} that if $(V(t))_{t\geq 0}$ is a $C_0$-semigroup in Banach space $\mathcal{F}$, then the set $$\left\{\varphi\in \mathcal{F}: \exists \lim_{t\to +0}\frac{V(t)\varphi-\varphi}{t}\right\}\stackrel{denote}{=}D(L)
$$ 
is a dense linear subspace in $\mathcal{F}$. The operator $L$ defined on the domain $D(L)$ by the equality $$L\varphi=\lim_{t\to +0}\frac{V(t)\varphi-\varphi}{t}$$ is called \textit{an infinitesimal generator} (or just \textit{generator} to make it shorter) of the $C_0$-semigroup $(V(t))_{t\geq 0}$, and notation $V(t)=e^{tL}$ is widely used. 

One of the reasons for the study of $C_0$-semigroups is their connection with differential equations. If $Q$ is a set, then the function $u\colon [0,+\infty)\times Q\to \mathbb{R}$, $u\colon (t,x)\longmapsto u(t,x)$ of two variables $(t,x)$ can be considered as a function $u\colon t\longmapsto [x\longmapsto u(t,x)]$ of one variable $t$
with values in the space of functions of the variable $x$. If $U(t)=u(t,\cdot)\in\mathcal{F}$ for all $t\geq0$ then one can define $Lu(t,x)=(Lu(t,\cdot))(x).$ If there exists a $C_0$-semigroup $(e^{tL})_{t\geq 0}$ then the Cauchy problem for the linear evolution equation
$$
\left\{ \begin{array}{ll}
U'(t)=LU(t) \ \mathrm{ for }\ t>0,\\
U(0)=u_0
\end{array} \right.
$$
has a unique solution $U\in C^1([0,+\infty),\mathcal{F})$ given by the formula $U(t)=e^{tL}$. Moreover, the Cauchy problem for the linear evolution equation
\begin{equation}\label{generalCP}
\left\{ \begin{array}{ll}
u'_t(t,x)=Lu(t,x) \ \mathrm{ for }\ t>0, x\in Q\\
u(0,x)=u_0(x)\ \mathrm{ for } \ x\in Q
\end{array} \right.
\end{equation} 
has a (unique in the sense of $\mathcal{F}$) solution $u$  given by the formula $$u(t,x)=(e^{tL}u_0)(x).$$ The solution depends on $u_0$ continuously and is understood in the sense of $\mathcal{F}$, where $u(t,\cdot)\in\mathcal{F}$ for every $t\geq 0$. Compare also different meanings of the solution \cite{EN1}, including mild solution which is given by the same formula $u(t,\cdot)=e^{tL}=U(t)$ and solves the corresponding integral equation $U(t)=u_0+L\int_0^tU(t)dt$. Note that if there exists a strongly continuous group $(e^{tL})_{t\in\mathbb{R}}$ then in (\ref{generalCP}) the equation $u'_t(t,x)=Lu(t,x)$ can be considered not only for $t>0$, but also for $t\in\mathbb{R}$, and the solution is provided by the same formula $u(t,x)=(e^{tL}u_0)(x)$.

\begin{definition}\label{defCT} (\textit{Introduced in \cite{R1}}). Let us say that $C$ is \textit{Chernoff-tangent} to $L$ iff the following conditions of Chernoff tangency (CT) hold: 

(CT0). Let $\mathcal{F}$ be a Banach space, and let $\mathscr{L}(\mathcal{F})$ be the space of all linear bounded operators in $\mathcal{F}$. Suppose that we have an operator-valued function $C\colon [0, +\infty) \to \mathscr{L}(\mathcal{F})$, or, in other words, we have a family $(C(t))_{t\geq 0}$ of linear bounded operators in $\mathcal{F}$. We also have a closed linear operator $L\colon D(L) \to \mathcal{F}$ defined on the linear subspace $D(L)\subset\mathcal{F}$ that is dense in $\mathcal{F}$.

(CT1) The function $t\longmapsto C(t)f\in\mathcal{F}$ is continuous for each $f\in\mathcal{F}$. 

(CT2) $C(0)=I$, that is, $C(0)f=f$ for each $f\in\mathcal{F}$.

(CT3) There exists such a dense subspace $\mathcal{D}\subset D(L)\subset \mathcal{F}$ that for each $f\in \mathcal{D}$ there exists a limit $$C'(0)f=\lim_{t\to 0}\frac{C(t)f-f}{t}.$$ 

(CT4) The closure of the operator $(C'(0),\mathcal{D})$ is equal to $(L,D(L)).$
\end{definition}

\begin{remark} Let us consider one-dimensional example $\mathcal{F}=\mathscr{L}(\mathcal{F})=\mathbb{R}$. Then $g\colon [0,+\infty)\to\mathbb{R}$ is Chernoff-tangent to $l\in\mathbb{R}$ iff $g$ is continuous and $g(t)=1+tl+o(t)$ as $t\to+0$.
\end{remark}

\begin{theorem}\label{chth} (\textsc{P.\,R.~Chernoff (1968)}, see \cite{EN1, Chernoff}). Let $\mathcal{F}$ be a Banach space, and $\mathscr{L}(\mathcal{F})$ be the space of all linear bounded operators in $\mathcal{F}$. Suppose that the operator $L\colon \mathcal{F}\supset D(L)\to \mathcal{F}$ is linear and closed, and function $C$ is defined on $[0,+\infty)$ and takes values in $\mathscr{L}(\mathcal{F})$. Suppose that these assumptions are fulfilled:

(E) There exists a $C_0$-semigroup $(e^{tL})_{t\geq 0}$ with infinitesimal generator $(L,D(L))$.

(CT) $C$ is Chernoff tangent to $(L,D(L)).$

(N) There exists a number $\omega\in\mathbb{R}$, such that $\|C(t)\|\leq e^{\omega t}$ for all $t\geq 0$.

Then for each $f\in \mathcal{F}$  we have $(C(t/n))^nf\to e^{tL}f$ as $n\to \infty$ with respect to norm in $\mathcal{F}$ uniformly with respect to $t\in[0,T]$ for each $T>0$, i.e.
$$\lim_{n\to\infty}\sup_{t\in[0,T]}\left\|e^{tL}f - (C(t/n))^nf \right\|=0.$$
\end{theorem}

\begin{remark} 
In the one-dimensional example ($\mathcal{F}=\mathscr{L}(\mathcal{F})=\mathbb{R}$) the Chernoff theorem says $e^{tl}=\lim_{n\to\infty}g(t/n)^n=\lim_{n\to\infty}(1+tl/n+o(t/n))^n$, which is a simple fact of calculus.
\end{remark}

\begin{definition} 
Let $\mathcal{F}, \mathscr{L}(\mathcal{F}), L$ be as above in the definition \ref{defCT}. If $C$ is the Chernoff tangent to $L$ and equality $\lim_{n\to\infty}\sup_{t\in[0,T]}\left\|e^{tL}f - (C(t/n))^nf \right\|=0$ holds, then $C$ is called \textit{a Chernoff function} for the operator $L$, and $(C(t/n))^nf$ is called \textit{a Chernoff approximation expression} to $e^{tL}f$. 
\end{definition}

\begin{remark}\label{noteasy} 
If $L$ is a linear bounded operator in $\mathcal{F}$, then $e^{tL}=\sum_{k=0}^{\infty}(tL)^k/k!$ where the series converges in the standard norm topology in $\mathscr{L}(\mathcal{F})$. If $L$ is not bounded (such as the Laplacian and many other differential operators), especially if $L$ has variable coefficients, then expressing $(e^{tL})_{t\geq 0}$ in terms of $L$ usually is a difficult problem which is equivalent to the problem of finding (for each $u_0\in D(L)$) the $\mathcal{F}$-valued function $U$ that solves the Cauchy problem $[U'(t)=LU(t); U(0)=u_0]$. If we find this solution, then $e^{tL}$ is obtained for each $u_0\in D(L)$ and each $t\geq 0$ in the form $e^{tL}u_0=U(t)$. Using the fact that $D(L)$ is dense in $\mathcal{F}$ because $L$ is a generator of a $C_0$-semigroup, and the fact that $e^{tL}$ is a bounded operator for each $t\geq0$, we find $e^{tL}u_0$ for all $u_0\in\mathcal{F}$. Alternatively, we can solve the integral equation $U(t)=u_0+L\int_0^tU(t)dt$ for all $u_0\in\mathcal{F}$, and find $e^{tL}u_0=U(t)$ for all $u_0\in\mathcal{F}$. These two tasks are not trivial, so expressing $e^{tL}$ in terms of $L$ directly is also not easy. 
\end{remark}

\begin{remark} 
In the definition of the Chernoff tangency, the family $(C(t))_{t\geq 0}$ usually does not have a semigroup composition property, that is, $C(t_1+t_2)\neq C(t_1)C(t_2)$, while $(e^{tL})_{t\geq 0}$ has it: $e^{t_1L}e^{t_2L}=e^{(t_1+t_2)L}$. However, each $C_0$-semigroup $(e^{tL})_{t\geq 0}$ is Chernoff-tangent to its generator $L$ and appears to be its Chernoff function. If the coefficients of the operator $L$ are variable and $L$ is unbounded, there is usually no simple formula for $e^{tL}$ due to Remark \ref{noteasy}. On the other hand, even in this case one can find a rather simple formula to construct a Chernoff function $C$ for the operator $L$, because there is no need to worry about the composition property, and then obtain $e^{tL}$ in the form $e^{tL}=\lim_{n\to\infty}C(t/n)^n$ via the Chernoff theorem. 
\end{remark}

\subsection{Theoretical bounds for the speed of convergence}

Paul Chernoff published his theorem on approximations of operator semigroups in 1968; however, in 2025 the rate of convergence is not completely studied yet. The known results can be found in the papers \cite{Prudnikov2020, Bent2004, Bent2009,NSZ2018, GSK2019,Vedenin2020, Zagrebnov2020, GR2021, Galkin2022, VVGKR, Dragunova2023, GR2024, Vedenin2022} and the references therein. Let us sketch a few of them. The first fact to mention is that all rates of convergence in the Chernoff theorem are possible.

\begin{theorem}\label{allrates} (O.E.Galkin, I.D.Remizov \cite{GR2024})
Suppose that $(h_n)_{n\in\mathbb{N}}$ is an arbitrary sequence of non-negative real numbers and $\lim_{n\to\infty}h_n=0$. Suppose that $(e^{tL})_{t\geq0}$ is a $C_0$-semigroup in the Banach space $\mathcal{F}$ and $e^{tL}\neq 0$ for all $t>0$. Then for each $\tau>0$ there exists a Chernoff function $C_\tau\colon[0,+\infty)\to\mathscr{L}(\mathcal{F})$ such that
		$$
		\|C_\tau(\tau/n)^n-e^{\tau L}\|=h_n+O(h_n^2) \textrm{ as }n\to\infty.
		$$
\end{theorem}

	\begin{remark}
		In the above theorem we can set, e.g. $h(n)=\frac{1}{2^{2^n}}$ and obtain very fast convergence, but also the case $h(n)=\frac{1}{\ln \ln n}$ with very slow convergence is possible. So, Chernoff approximation can be an extremely (in)effective tool for finding operator exponents.
	\end{remark}

Another example can be found in \cite{GR2021} for a translation semigroup. See also the example in \cite{GSK2019}, Proposition 4.8 for $\alpha=0$. In \cite{GSK2019} one considers
the case $\alpha\in (0,2]$, but the proof works for $\alpha=0$ as well, see also \cite{GT2014}. Here, the symbol $\alpha$ is used in the notation of \cite{GSK2019}, having nothing in common with the usage of this symbol in the present paper.
    
The second fact to mention is that Chernoff approximations can converge on every vector but do not converge in norm. More precisely, the following fact is known:

\begin{theorem} (O.E.Galkin, I.D.Remizov \cite{GR2021}) There exists a
Banach space $\mathcal{F}$, $C_0$-semigroup $(e^{tL})_{t\geq 0}$ in $\mathcal{F}$ with generator $(L,D(L))$ and Chernoff function $C$ for operator $(L,D(L))$ such that for each $t>0$: 

1. $\lim_{n\to\infty}\|C(t/n)^nf-e^{tL}f\|=0$ for all $f\in\mathcal{F}$,

2. $\|e^{tL}\|=\|C(t)\|=1$,

3.  $\|C(t/n)^n-e^{tL}\|\geq 1 \not\to 0$ as $n\to\infty$.
\end{theorem}

The final fact to mention is the following. The numerical experiment \cite{OSS2012, Prudnikov2020, Dragunova2023} shows that the rate of convergence $C(t/n)^nf\to e^{tA}f$ as $n\to\infty$ can depend on both $C$ and $f$ even if $\|C(t)\|=\|f\|=1$. So natural questions appear: Is it possible to find conditions that guarantee for fixed $k\in\mathbb{N}$ that $\|C(t/n)^nf-e^{tA}f\|=O(1/n^k)$ as $n\to\infty$? The following theorem answers this question.

\begin{theorem}\label{grth}(O.E.Galkin, I.D.Remizov \cite{GR2024, Galkin2022})
Suppose that the following three conditions are met:
\begin{enumerate}
\item 
			$C_0$-semigroup $(e^{tL})_{t\ge 0}$ in Banach space $\mathcal{F}$ is given, and there exists a number $T>0$ 
			such that for some $M_1\geq 1$ and $w\geq 0$ the inequality $\|e^{tL}\| \le M_1e^{wt}$ holds 
			for all $t\in[0,T]$. 
			
			\item 
			A mapping $C\colon (0,T]\to\mathscr{L}(\mathcal{F})$ is given and for some constant $M_2\geq 1$ the inequality $\|C(t)^k\| \le M_2e^{kwt}$ holds
for all $t\in (0,T]$ and all $k\in\mathbb{N}$.
			
			\item 
The numbers $m\in\mathbb{N}_0=\mathbb{N}\cup\{0\}$ and $p\in\mathbb{N}$ are fixed. There exist a $(e^{tL})_{t\ge 0}$-invariant subspace $\mathcal{D}\subset D(L^{m+p})\subset \mathcal{F}$  and functions $K_j\colon (0,T]\to [0,+\infty)$, $j=0,1,\dots,m+p$ such that for all $t\in(0,T]$ and all $f\in \mathcal{D}$ we have 
			\begin{equation} \label{ocdiffmain}
				\bigg\|C(t) f - \sum_{k=0}^{m} \frac{t^kL^k f}{k!}\bigg\| \le t^{m+1}\sum_{j=0}^{m+p}K_j(t)\|L^{j}f\|.
			\end{equation}
			
		\end{enumerate}
		
		Then the following two statements hold:
		\begin{enumerate}
			\item 
			For all $t>0$, all integer $n\geq t/T$ and all $f\in \mathcal{D}$ the estimate is true:
			\begin{equation} \label{ocresdiff1main}
				\|C(t/n)^n f - e^{tL}f\|\leq \frac{M_1M_2t^{m+1}e^{wt}}{n^{m}}\sum_{j=0}^{m+p}c_j(t/n)\|L^{j}f\|,
			\end{equation}
			where $c_{m+1}(t)=K_{m+1}(t)e^{-wt}+M_1/(m+1)!$ and $c_j(t)=K_j(t)e^{-wt}$ 
			for all such $j\in\{0,1,\ldots,m+p\}$, that $j\neq m+1$.
			
			\item 
			If $\mathcal{D}$ is dense in $\mathcal{F}$ and for all $j=0,1,\dots, m+p$ we have 
			$K_j(t)=o(t^{-m})$ as $t\to+0$,
			then for all $g\in \mathcal{F}$ and all $\mathcal{T}>0$
			the following equality is true:
			\begin{equation} \label{eqanChernoff}
				\lim_{\mathcal{T}/T\leq n\to\infty}\sup_{t\in(0,\mathcal{T}]}\left\|C(t/n)^ng-e^{tL}g\right\| = 0.
			\end{equation}
		\end{enumerate}
	\end{theorem}

\begin{remark}\label{remord}
In the case of the heat semigroup in $\mathcal{F}=UC_b(\mathbb{R})$, the generator $L=\partial^2$ is equal to the operator of taking the second derivative, so the domain of $L, L^2, L^3,\dots$ are functions that have $2,4,6,\dots$ uniformly bounded and continuous derivatives. The above theorem \ref{grth} says that for a high convergence rate we should have $f\in D(L^q)$ with $q$ large enough, that is, for the case of the heat semigroup function $f$ should have bounded derivative of high order. Also we see from theorem \ref{grth} that Chernoff functions with higher order of Chernoff tangency to $L$ provide a better rate of convergence. Here we come to the main question that we answer in this paper.

\textbf{The main question attacked in the paper.} What happens if conditions of the above theorem \ref{grth} do not hold? For example, what the convergence rate will be if $f$ is differentiable only once but not twice, and $f'$ is e.g. H\"older continuous with H\"older exponent equal to $1/2$?  Then $f\notin D(\partial^2)$ and the above theorem say nothing about the rate of convergence. In the present paper, we study this situation numerically for several functions $f$ with different smoothness.
\end{remark}

\subsection{Open problem: do superfast converging approximations exist?} 

Recall that Chernoff approximations are called fast converging if the convergence rate is higher than $\mathrm{const}/n$, and superfast converging if the convergence rate is higher than $\mathrm{const}/n^k$ for each $k\in\mathbb{N}$. Due to Theorem \ref{allrates} all rates of convergence are possible, but in a setting that can be called artificial. In Theorem \ref{allrates} we start from the given semigroup and the given rate of convergence, and then construct a Chernoff function using them. Meanwhile, in applications, we never know the semigroup and use Chernoff approximations to find it. Most examples deal with generator equal to a differential operator with variable coefficients, and for that case natural examples of super-fast convergent approximations are not known.

In the List of open problems in the One-parameter semigroups theory \cite{Open} there is the following problem: \textit{Provide a concrete example of a super-fast (at least exponentially) converging Chernoff approximation to a semigroup generated by a first-order differential operator with variable coefficients, with Chernoff function explicitly expressed in terms of these coefficients and without circulus vitiosus of any kind, e.g.\ in the form of using the values of the semigroup for creating a Chernoff function for it.} If this is possible then one can expect rapidly grown interest to the Chernoff approximation method.

\textit{Formal statement.} Consider Banach space $UC_b(\mathbb{R},\mathbb{R})$ of all bounded and uniformly continuous real-valued functions of one real variable, with the uniform norm $\|f\|=\sup_{x\in\mathbb{R}}|f(x)|$. Recall that the space $C_b^\infty(\mathbb{R},\mathbb{R})$ of all functions that are bounded and have bounded derivatives of all orders is a dense linear subspace in $UC_b(\mathbb{R},\mathbb{R})$. Suppose that arbitrary functions $a,b,c\in UC_b(\mathbb{R},\mathbb{R})$ are given, and inequality $a(x)\geq a_0\equiv\mathrm{const}>0$ holds for all $x\in\mathbb{R}$ . Define an operator $L$ by the equality 
\[
(Lf)(x)=a(x)f''(x)+b(x)f'(x)+c(x)f(x) \textrm{ for all }  f\in C_b^\infty(\mathbb{R},\mathbb{R}), x\in\mathbb{R}
\] and recall (see e.g. lemma 2 in \cite{Rjmp}, theorem 4 in \cite{GR2024}) that the closure of $L$ generates a $C_0$-semigroup $(e^{tL})_{t\geq 0}$ in $UC_b(\mathbb{R},\mathbb{R})$, satisfying 
\[
\|e^{tL}\|\leq e^{t\max(0, \sup\limits_{x\in\mathbb{R}} c(x))} \quad \textrm{ for all } \quad  t\geq 0.
\]
\textbf{The open problem} is to find a family $(C(t))_{t\geq 0}$ of linear bounded operators defined everywhere in $UC_b(\mathbb{R},\mathbb{R})$ which has the following three properties:

\begin{itemize}
\item $C$ and $L$ satisfy the conditions of at least one version of the Chernoff theorem on approximation of operator semigroups, i.e.\ $C$ is Chernoff-equivalent to the semigroup $(e^{tL})_{t\geq 0}$, i.e.\ $C$ is a Chernoff function for the operator $L$;

\item There exists a dense linear subspace $D$ in $UC_b(\mathbb{R},\mathbb{R})$ and constants $M>0, q>1, n_0>0$ such that 
\[
\|e^{tL}f-C(t/n)^nf\|\leq M/q^{n} \quad  \textrm{ for all }  \quad n\geq n_0 \textrm{ and all } f\in D;
\]

\item $C$ is defined by an explicitly given formula $(C(t)f)(x)=F(a,b,c,f,x)$ and the calculation of $e^{tL}$ is not included into this formula, even in a hidden way, so expression $e^{tL}=\lim\limits_{n\to\infty}C(t/n)^n$ is free from circulus vitiosus of any kind and can be used as a practical method of calculation of $e^{tL}$.
\end{itemize}

After the Chernoff theorem was published in 1968, it has found many applications in random processes theory, partial differential equations, functional integration and other fields closely related to operator semigroups theory, see \cite{Butko2020}. 

Chernoff's original theorem states that the fact of convergence $\|e^{tL}f-C(t/n)^nf\|\to 0$ for all $f\in\mathcal{F}$, however, says nothing about the rate of convergence. It was shown (see papers \cite{Zagrebnov2020, Zagrebnov2022, Gomilko2019, CACHIA2001176, IcTa1997, NeZa1998}) that as $n$ tends to infinity, in some cases $\|e^{tL}-C(t/n)^n\|$ tends to zero with a rate of $const/n$, and under some conditions may reach $const/n^2$. Meanwhile, in other cases \cite{GR2024} we have $\|e^{tL}f-C(t/n)^nf\|\to 0$ for all $f\in\mathcal{F}$ but $\|e^{tL}-C(t/n)^n\|\geq 1$ and it does not tend to zero. The rate of tending $\|e^{tL}f-C(t/n)^nf\|$ to zero depends heavily on $f$ (see \cite{Prudnikov2020}, \cite{Dragunova2023}) and may be slower than $const/n$. Model examples show that $\|e^{tL}f-C(t/n)^nf\|$ can tend to zero with an arbitrary high or arbitrary low rate, and this is also possible in a general setting for an arbitrary never-zero semigroup (see  \cite{Galkin2022}). If we already know the semigroup $e^{tL}$ then we can construct a Chernoff function $C$ that provides Chernoff approximations $C(t/n)^n$ that converge to $e^{tL}$ with arbitrarily high or arbitrary slow rate (even with infinitely high rate if we take $C(t)=e^{tL}$) but they are not useful in practice because a circulus vitiosus arises. In these constructions, we employ the semigroup to build a Chernoff function with desired properties and then use this Chernoff function to approximate the semigroup. 

In summary, we see that the rate of convergence of Chernoff approximations depends on all three elements: $e^{tL}$, $C$, $f$. If $e^{tL}$ and $C$ have the same Taylor polynomial of order $k$ and $C$ is close to its Taylor polynomial in some sense (see  \cite{Galkin2022}), and $f$ is good enough in terms of $L$, then $\|e^{tL}f-C(t/n)^nf\|\leq const/n^k$. If the rate of convergence of Chernoff approximations is higher than $const/n$, then we say that they are fast-converging Chernoff approximations, and if the rate is higher than $const/n^k$ for all $k\in\mathbb{N}$ then we say that they are super-fast-converging Chernoff approximations. The Chernoff functions $C$ that provide super-fast converging approximations are expected to satisfy $\frac{d^k}{dt^k}C(t)\big|_{t=0}=\frac{d^k}{dt^k}e^{tL}\big|_{t=0}=L^k$ for each $k\in\mathbb{N}$. 

Examples of fast converging approximations for the operator $L$ given by $Lf=af''+bf'+cf$ were constructed for $k=2$ and it is clear (see  \cite{Vedenin2020}) how to extend this method to arbitrary positive integer values of $k$; however, a significant amount of work should be done to complete this task. Effective methods for constructing superfast converging Chernoff approximations for this operator $L$ are not known at present. However, there are no fundamental objections to their existence either, which presents a challenge to either find such approximations or to prove that they are impossible for some reason.

\textbf{Theoretical value of superfast converging approximations.} Examples of fast-converging Chernoff approximations known at present all involve exponential growth in computational complexity: to calculate $C(t/n)^n$ we need to perform $const\cdot p^n$ arithmetic operations where $p$ ranges from 3 to 9 for different variants of $C$ known at present. Even worse is the fact that high values of $k$ in the estimate $\|e^{tL}f-C(t/n)^nf\|\leq const/n^k$ result in high values of $p$. However, if one can find the Chernoff function $C$ that provides $\|e^{tL}f-C(t/n)^nf\|\leq const/q^{n}$ for $q>1$, then this exponential decay of the approximation error may balance the (hoped to be still not worse than exponential) growth of computational complexity. If it happens then these super fast converging Chernoff approximations can be a good alternative to existing numerical methods of solving various linear differential equations with variable coefficients: evolution PDEs (Schr\"odinger type equations with complex coefficients, parabolic equations with real coefficients) as well as ODEs and elliptic PDEs. 

\textbf{Connection with the research presented in the present paper.} Numerical study of the rate of convergence on model cases helps to gather information on the Chernoff approximation method, and in the future may give a clue how to construct these super-fast converging approximations, or to prove that it is not possible.

\subsection{$C_0$-semigroups and dynamical systems}

Being one of the backbone structures in mathematics, $C_0$-semigroups are connected with many other main theories and structures, including PDEs, stochastic equations and random processes, infinite-dimensional analysis, measure theory, quantum mechanics, quantum information, control theory, and some other, see \cite{EN1, EN2}.

The connection with dynamical system theory is the following. On the one hand, $C_0$-semigroup is an example of linear semiflow in a Banach space. On the other hand, each uniformly continuous semiflow in metric space (e.g. on manifold) naturally generates a $C_0$-semigroup in the Banach space of all real-valued uniformly continuous functions on this metric space -- using the shift through the trajectories of the semiflow. This shift generates a composition operator that is commonly referred to as the Koopman operator, named after Bernard Koopman. It is the left-adjoint of the Frobenius-Perron transfer operator. We state these two results (each $C_0$-semigroup is a specific type of semiflow, and each semiflow naturally defines a specific type of $C_0$-semigroup) and provide short proofs that are better than references. We start from the definition of a dynamical system.

\begin{definition}\label{defds} Let $(M,d)$ be metric space and $T$ be the time interval $T\in\{\mathbb{R},[0,+\infty),\mathbb{Z},\mathbb{N}_0=\{0,1,2,3,\dots\}\}$. One says that dynamical system in $M$ is a tuple $(M,T,\Phi)$, where $\Phi\colon T\times M\to M$ a continuous (in some sense) function, which satisfies two properties:

1. $\Phi(0,x)=x$ for all $x\in M$.

2. $\Phi(t_1,\Phi(t_2,x))=\Phi(t_1+t_2,x)$ for all $x\in M$ and all $t_1,t_2\in T$.\\
Function $\Phi$ is called the evolution function of the dynamical system, but often $\Phi$ itself is called the dynamical system for brevity. Depending on the notion of continuity that $\Phi$ satisfies, we have different types of dynamical systems. Depending on the choice of $T$, dynamical systems are divided into two classes, each with two subclasses:

i. Continuous time dynamical system: flow in $M$ for $T=\mathbb{R}$, semiflow in $M$ for $T=[0,+\infty)$.
	
ii. Discrete time dynamical system: cascade in $M$ for $T=\mathbb{Z}$, semicascade in $M$ for $T=\mathbb{N}_0$.
\end{definition}

\begin{prop}
Suppose that $(e^{tA})_{t\geq0}$ is a $C_0$-semigroup in the Banach space $\mathcal{F}$, denote $\Phi(t,f)=e^{tA}f$ for all $t\geq 0$, $f\in\mathcal{F}$. Then $\Phi$ is a semiflow in $\mathcal{F}$ with the following notion of continuity: functions $t\mapsto\Phi(t,f)$ and $f\mapsto\Phi(t,f)$ are continuous for each $t\geq 0, f\in\mathcal{F}$. The semiflow $\Phi$ is linear, i.e. $\Phi(t,f_1+f_2)=\Phi(t,f_1) + \Phi(t,f_2)$ for all $f_1,f_2\in\mathcal{F}$ and all $t\geq0$. 

Moreover, if $(e^{tA})_{t\geq0}$ is a contraction $C_0$-semigroup (i.e. $\|e^{tA}\|\leq 1$ for all $t\geq 0$) then function $f\mapsto\Phi(t,f)$ is Lipschiz-continuous (with Lipschiz constant equal to 1) uniformly with respect to $t$, i.e. $\|\Phi(t,f_1)-\Phi(t,f_2)\|\leq \|f_1-f_2\|$ for all $f_1,f_2\in\mathcal{F}$ and all $t\geq0$.
\end{prop}

\begin{proof}
Conditions 1 and 2 of the definition \ref{defds} follow directly from conditions 1 and 2 of the definition \ref{C0sgdef}. Condition 3 of the definition \ref{C0sgdef} states that for each $f\in\mathcal{F}$ we have $0=\lim_{t\to t_0}\|e^{tA}f-e^{t_0A}f\|=\lim_{t\to t_0}\|\Phi(t,f)-\Phi(t_0,f) \|$ which means that $\Phi(t,f)$ is continuous in $t\geq 0$ for each fixed $f\in \mathcal{F}$. Due to definition \ref{C0sgdef} operator $e^{tA}$ is linear and bounded (hence continuous), so $\Phi(t,f)$ is continuous in $f\in \mathcal{F}$ for each fixed $t\geq0$. If $\|e^{tA}\|\leq 1$ then $\|\Phi(t,f_1)-\Phi(t,f_2)\|=\|e^{tA}f_1-e^{tA}f_2\|=\|e^{tA}(f_1-f_2)\|\leq \|e^{tA}\|\cdot \|f_1-f_2\|\leq \|f_1-f_2\|$.
\end{proof}

\begin{prop}
Let $M$ be a metric space, and $\mathcal{F}=UC_b(M)$ be the a Banach space of all uniformly continuous and bounded functions $f\colon M\to\mathbb{R}$ with the uniform norm $\|f\|=\sup_{x\in M}|f(x)|$. Suppose that $\phi\colon [0,+\infty)\times M\to M$ is a semiflow in $M$ with the following notion of continuity: a) $x\mapsto \phi(t,x)$ is uniformly continuous in $x$ for each $t\geq 0$, and b) $t\mapsto \phi(t,x)$ is continuous in $t$ uniformly with respect to $x\in M$, i.e. $\lim_{\tau\to0}\sup_{x\in M}|\phi(t+\tau,x)-\phi(t,x)|=0$ for all $t\geq 0, x\in M$. 
Let us define $\Phi\colon [0,+\infty)\times\mathcal{F}\to\mathcal{F}$ by the equality
\begin{equation}\label{fidef}
(\Phi(t,f))(x)=f(\phi(t,x))\textrm{ for all }t\geq0, x\in M, f\in\mathcal{F}; \quad Q(t)f=\Phi(t,f).
\end{equation}
Then $Q(t)\in\mathscr{L}(\mathcal{F})$ for each $t\geq0$, and $(Q(t))_{t\geq0}$ is a contraction $C_0$-semigroup in $\mathcal{F}$. 
\end{prop}

\begin{proof} Let us first prove that $f\in\mathcal{F}$ implies $Q(t)f\in\mathcal{F}$. Indeed, for each $t\geq 0$ and each $f\in UC_b(M)$ function $x\mapsto f(\phi(t,x))$ is bounded and uniformly continuous because $f$ is bounded, and $x\mapsto \phi(t,x)$ is uniformly continuous due to condition a). Conditions 1 and 2 of the definition \ref{C0sgdef} follow directly from (\ref{fidef}) and conditions 1, 2 of the definition \ref{defds}. The linearity of $Q(t)f$ in $f$ for fixed $t$ follows from (\ref{fidef}) and from the standard definition of the linear structure in $\mathcal{F}=UC_b(M)$. To prove the continuity of $Q(t)f$ in $t$ for fixed $f$ we start from the equality $\lim_{\tau\to0}\|Q(t+\tau)f-Q(t)\|=\lim_{\tau\to0}\sup_{x\in M}|f(\phi(t+\tau,x))-f(\phi(t,x)|))$, and then mention that $\lim_{\tau\to0}\sup_{x\in M}|\phi(t+\tau,x)-\phi(t,x)|=0$ due to condition b), and that function $f$ is uniformly continuous and bounded. Hence, the condition 3 of the definition \ref{C0sgdef} holds. Also $\|Q(t)f\|=\sup_{x\in M}|f(\phi(t,x))|\leq \sup_{x\in M}|f(x)|=\|f\|$ so the operator $Q(t)$ is bounded and $\|Q(t)\|\leq 1$ for each $t\geq 0$. The condition $\|Q(t)\|\leq 1$ means that $(Q(t))_{t\geq0}$ is not just a $C_0$-semigroup, but a contraction $C_0$-semigroup.
\end{proof}

With the above propositions, we see that Koopman operators corresponding to nonlinear dynamical system (a semiflow in the metric space $M$) appear to be a linear dynamical system ($C_0$-semigroup in $UC_b(M)$). This procedure is sometimes called koopmanization in analogy with linearization and quantization, and its wild surge in popularity is sometimes joked about as ``Koopmania''. Let us also remark that depending on the problem, not only the space $UC_b(M)$ is used, one can also employ $L_p(M)$ or even the space of all measurable functions on $M$ -- which are not continuous; this setting often appears in the ergodic theory. All of this shows a deep connection between nonlinear dynamical systems and $C_0$-semigroups. 

\begin{remark}
Recall that the heat semigroup (i.e. semigroup generated by the heat equation) that we use in the present paper as a model example is a contraction $C_0$-semigroup. Polynomial dynamical systems and ordinary differential equations associated with the heat equation were studied by V.M.~Buchstaber and E.Yu.~Bunkova \cite{BB}. See also A.V.~Vinogradov's paper \cite{Vin} on the Painlev\'e test for ordinary differential equations associated with the heat equation.
\end{remark}

\section{Numerical simulation results}

\subsection{Problem setting}

\begin{definition}
We say that operator-valued function $C$ is \textit{Chernoff-tangent of order $k$} to operator $L$ iff $C$ is Chernoff-tangent to $L$ in the sense of definition \ref{defCT} and the following condition (CT3-k) holds:

There exists such a dense subspace $\mathcal{D}\subset \mathcal{F}$ that for each $f\in \mathcal{D}$ we have
$$
C(t)f=\left(I+tL+\frac{1}{2}t^2L^2+\dots+\frac{1}{k!}t^kL^k\right)f+o(t^{k})\textrm{ as }t\to 0.\eqno(CT3-k)
$$

\end{definition}

\begin{remark} It is clear that for $k=1$ condition (CT3-k) becomes just (CT3). For the semigroup $C(t)=e^{tL}$ condition (CT3-k) holds for all $k=1,2,3,\dots$ So one can expect that the bigger $k$ is the better rate of convergence $C(t/n)^nf\to e^{tL}f$ as $n\to\infty$ will be, if $f$ belongs to the space $D$. This idea was proposed in  \cite{R2}, where two conjectures about the convergence speed were formulated explicitly, and one of them were recently proved in \cite{GR2021, Galkin2022}. For initial conditions that are good enough and $t$ fixed, Chernoff function with Chernoff tangency of order $k$ (by conjecture from \cite{R2} and by theorem \ref{grth} above) should provide $\|u(t,\cdot)-u_n(t,\cdot)\|=O(1/n^k)$ as $n\to\infty$. However, if $f\not\in D$ then nothing is known on the rate of convergence. In the present paper we are starting to fill this gap for operator $L$ given by $(Lf)(x)=f''(x)$ for all $x\in\mathbb{R}$ and all bounded, infinitely smooth functions $f\colon \mathbb{R}\to\mathbb{R}$, and $k=1,2$. 
\end{remark}

\textbf{Problem setting.} In the initial value problem (\ref{generalCP}) consider $Q=\mathbb{R}$, and the Banach space $\mathcal{F}=UC_b(\mathbb{R})$ of all bounded, uniformly continuous functions $f\colon\mathbb{R}\to\mathbb{R}$ endowed with the uniform norm $\|f\|=\sup_{x\in\mathbb{R}}|f(x)|$. Consider the operator $L$ given by $(Lf)(x)=f''(x)$ for all $x\in\mathbb{R}$ and all $f\in D=C^\infty_b(\mathbb{R})$ of all infinitely smooth functions $\mathbb{R}\to\mathbb{R}$ that are bounded with all the derivatives. Then (\ref{generalCP}) reads as (\ref{ZK}). Cauchy problem (\ref{ZK}) is the constant (one, zero, zero) coefficients particular case of the Cauchy problem considered in \cite{R4}, and the corresponding Chernoff function was found in \cite{R4}. The particular case of this Chernoff function reads as
$$
(G(t)f)(x)=\frac{1}{2}f(x)+\frac{1}{4}f(x+2\sqrt{t})+\frac{1}{4}f(x-2\sqrt{t})
$$
where we write $G(t)$ instead of $C(t)$ to show that $C(t)$ is a general abstract Chernoff function for some operator $L$, while $G(t)$ is this particular Chernoff function given above for operator $d^2/dx^2$. In \cite{R4} it was proved that $G(t)$ is the first-order Chernoff tangent to $d^2/dx^2$.

A.Vedenin (see \cite{Vedenin2020, Vedenin2022}) proposed another Chernoff function for the operator $L$ considered in \cite{R4}, and the particular case of constant coefficient of this operator is $d^2/dx^2$. The particular case of the Chernoff function obtained by A.Vedenin reads as 
$$
(S(t)f)(x)=\frac{2}{3}f(x)+\frac{1}{6}f(x+\sqrt{6t})+\frac{1}{6}f(x-\sqrt{6t}),
$$
and it was proved by A.Vedenin that $S(t)$ is second order Chernoff tangent to $d^2/dx^2$.

In the paper, we study how $\sup_{x\in\mathbb{R}}|u(t,x)-u_n(t,x)|$ depends on $n$ while $t=1/2$ is fixed and $u_n(t,x)$ is given by
$$
u_n(t,x)=(C(t/n)^nu_0)(x)
$$
where $C\in\{G,S\}$, $C(t/n)$ is obtained by substitution of $t$ by $t/n$ in the formula that defines $C(t)$, and $C(t/n)^n=C(t/n)C(t/n)\dots C(t/n)$ is a composition of $n$ copies of linear bounded operator $C(t/n)$. We consider several initial conditions $u_0$ that are all H\"older continuous (hence all belong to the $UC_b(\mathbb{R})$ space) but have different H\"older exponents. Then we remark how the rate of tending of $\sup_{x\in\mathbb{R}}|u(t,x)-u_n(t,x)|$ to zero depends on these H\"older exponents and the order of Chernoff tangency (which is 1 for $G(t)$, and 2 for $S(t)$).

\textbf{Comments on computational techniques.} Calculations were performed in the Python 3 environment using a program we wrote and which is available in the Appendix. All measurements, for the sake of reducing computational complexity, for each value of n (varying from 1 to 11) were carried out for 1000 points uniformly dividing the segment $[a,b]=[-\pi,\pi]$ for trigonometric initial conditions. Initial conditions of the form $u_0(x)=|\sin x|^{\xi}$ for various $\xi\in$\{9/2, 7/2, 5/2, 3/2, 1, 3/4, 1/2, 1/4\}, like any of Chernoff's approximations based on them, are periodic functions. So, the standard norm in $UC_b(\mathbb {R})$, namely 
$$
\|u_n(t,\cdot)-u(t,\cdot)\|=\sup_{x\in\mathbb{R}}|u_n(t,x)-u(t,x)|=\sup_{x\in[a,b]}|u_n(t,x)-u(t,x)|,
$$ 
where $u$ is the exact solution of (\ref{ZK}) and $u_n$ is the Chernoff approximation, is reached at the interval $[a,b]$ corresponding to the period. So we have
$$
d=\max_{k=1,\dots,1000} \left|u_n\left(t,a+\frac{k}{1000}(b-a)\right)-u\left(t,a+\frac{k}{1000}(b-a)\right)\right| \approx  \sup_{x\in[a,b]}|u_n(t,x)-u(t,x)| 
$$ 
and this value $d$ is what we measure numerically. We use the value $t=1/2$.

For exponential initial conditions, we used the segment $[a,b]=[-10,10]$, because these conditions tend to zero very fast when $|x|$ increases, so effectively they are zero for $|x|\geq 10$.

The program code was written with the possibility to set any operator and any initial condition, i.e. without simplifying Chernoff functions and using binomial coefficients, in contrast to the work \cite{Prudnikov2020} published earlier. Moreover, the initial condition does not necessarily have to be a smooth function. The number of iterations is not limited to 11, the value $n$ can be changed, both upward and downward. 

We have chosen the optimal value $n$ since the program is rather time consuming: using Jupyter Notebook 6.1.4 Anaconda 3 Python 3.8.3 set on a personal computer with Windows 10, CPU Intel Core i5-1035G1, 1.0-3.6 GHz, 8 Gb RAM, it takes about 20 minutes to complete the program for all initial conditions with the construction of graphs for them. At the research stage of the new method (Chernoff approximations), this is acceptable, but in the future, of course, the code will be optimized for a better speed since this is important in practice. After reasonable research in this direction, it will be possible in the future to write a library that allows us to solve partial derivative equations in this way.

\begin{remark} We are interested in the study of the asymptotic behavior, i.e. in the large values of $n$, so we exclude $n=1,2,3$ because these numbers are not large. However, the values $n>11$ take too much computation time. So, the coefficients of the straight lines were calculated by running linear regression (least squares method), using points corresponding to $4\leq n\leq 11$. 
\end{remark}

\subsection{Graphs of approximations for several trigonometric initial conditions}
\subsubsection{Approximations for initial condition $u_0(x)=\sin(x)$ as the first model case}

Let us first analyze the approximations for the initial condition $ u_0(x) =\sin x$. The exact solution $u(1/2,x)=e^{-1/2}\sin(x)$ was obtained by setting $t=1/2$ in the formula $u(t,x)=e^{-t}\sin(x)$ that satisfies $u'_t=u''_{xx}$ and $u(0,x)=\sin x$.

\begin{center}
\includegraphics[scale=0.6]{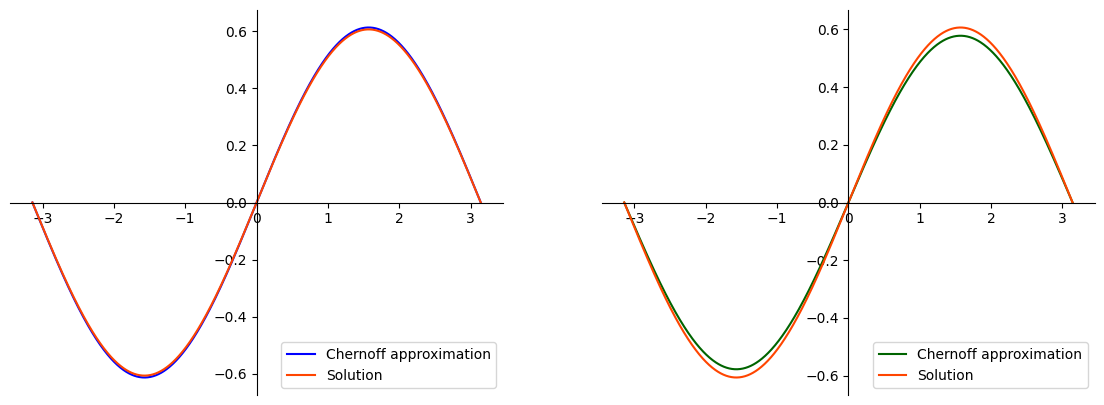}
   {fig. 1.1, $n=1$, $u_0(x)=\sin x$, $t=\frac{1}{2}$} 
\end{center} 

Figure~1.1 shows the exact solution $u(1/2,x)=e^{-1/2}\sin(x)$ (red) for the initial condition $u_0(x)=\sin x$, and approximate solutions constructed with the use of Chernoff functions $S(t)$ (left) and $G(t)$ (right) for $n=1$ (blue). The blue and red curves almost coincide on the left plot because the difference is too small to be observed visually. The initial condition $u_0=\sin x$ is very good in the sense of Chernoff approximations, since its derivatives of any order exist, have no discontinuities, and are bounded. And already at $n=1$ the function $S(t)$ gives a good approximation. For $n=2$ the quality of the approximation is again so good that visually it is impossible to distinguish the approximation and the exact solution, so we do not present a plot. However, it is possible to measure the distance between the exact solution and the approximation anyway, because the picture is just an illustration and is not used to measure the quality of the approximation. 

Figure 1.2 below shows plots of the decreasing error of Chernoff approximations as a function of $n$, where $1\leq n\leq 11$. On the left are plots of decreasing error for Chernoff functions $S(t)$ (in blue) and $G(t)$ (in green) in regular scale, and on the right are the same plots in logarithmic scale. The graph on the logarithmic scale allows us to estimate how much the convergence rate for the function $G(t)$ is less than the convergence rate for the function $S(t)$. 

\begin{center}
\includegraphics[scale=1.1]{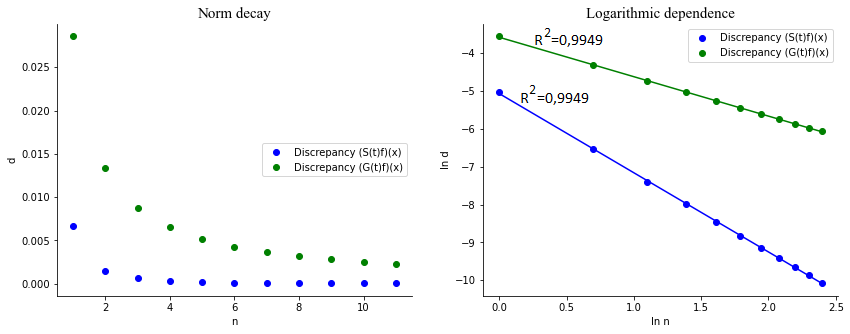}
    {fig. 1.2, $1\leq n\leq 11$, $u_0(x)=\sin x$, $t=\frac{1}{2}$} 
\end{center} 

One can see that the points on the right graph lie on straight lines with good accuracy. Using the least-squares method, we found the equations of these lines. For both lines, the least squares method gives the coefficient of determination $R^2>0.99$, which is numerical evidence of the visually obtained effect: the points on the right graph lie on straight lines with good accuracy. Rounding off the coefficients, we see that for the blue line (obtained via the Chernoff function $S(t)$) the equation is as follows:
$\ln(d)= -2.0416\ln(n) -5.1659$, i.e. $d=n^{-2.0416}e^{-5.1659}$. Similarly, for the green line (obtained via the Chernoff function $G(t)$), the equation $\ln(d)= -1.0178 \ln(n) -3.6262$, i.e. $d=n^{-1.0178}e^{-3.6262}$.

Using the same approach, we study the behavior of the error for other initial conditions.

\subsubsection{Approximations for initial condition $u_0(x)=|\sin(x)|$}
\begin{center}
\includegraphics[scale=0.5]{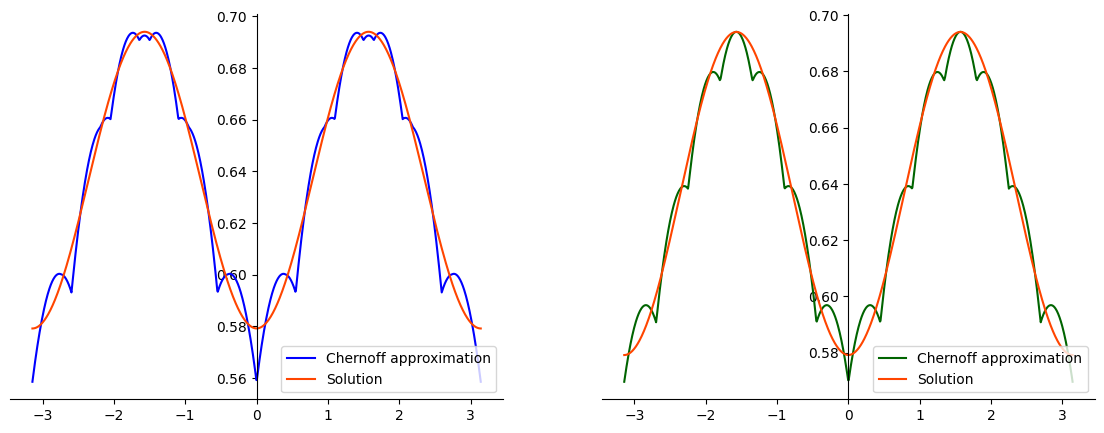}
    {\\fig. 2.1, $n=10$, $u_0(x)=|\sin x|$, $t=\frac{1}{2}$} 
\end{center}
Figure~2.1 shows two graphs of the approximate solution for the functions we study in $n=10$ and the exact solution under the initial condition $u_0(x)= |{\sin x}|$. The exact solution was obtained by numerical calculation of the integral in formula (\ref{exact}).

\begin{center}
\includegraphics[scale=1.0]{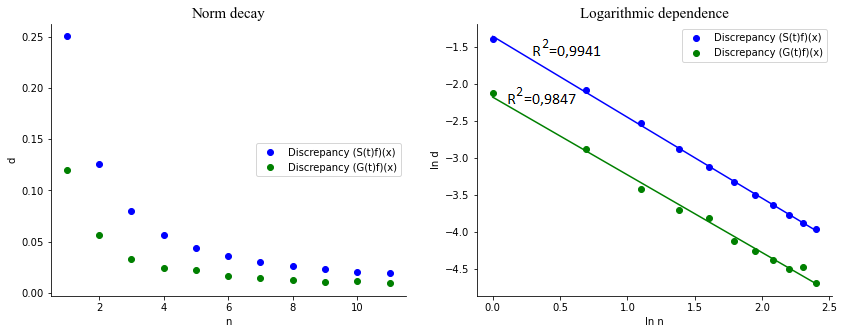}
    {\\fig. 2.2, $1\leq n\leq 11$, $u_0(x)=|\sin x|$, $t=\frac{1}{2}$} 
\end{center}

The blue color corresponds to the Chernoff function $S(t)$, and the green color corresponds to the Chernoff function $G(t)$. Rounding off the coefficients, we see that for the blue line (see Fig. 2.2) the equation looks as follows: $\ln(d)= -1.0818 \ln(n) -1.3846$, i.e. $d=n^{-1.0818}e^{-1.3846}$. Similarly, for the green line (see Figure 2.2), the equation $\ln(d)= -0.9655 \ln(n) -2.3452$, i.e. $d=n^{-0.9655}e^{-2.3452}$. The coefficients of the straight lines were calculated by running linear regression (least squares method), using points corresponding to $4\leq n\leq 11$.

\subsubsection{Approximations for initial condition $\sqrt[4]{|\sin x|}$}

Note that all special cases $|\sin x|^\xi$, where $0<\xi<1$, are similar to those already considered.
Consider $\xi=1/4$. The blue color corresponds to the Chernoff function $S(t)$, and the green color corresponds to the Chernoff function $G(t)$. The exact solution was obtained by numerical calculation of the integral in formula (\ref{exact}).

\begin{center}
\includegraphics[scale=0.5]{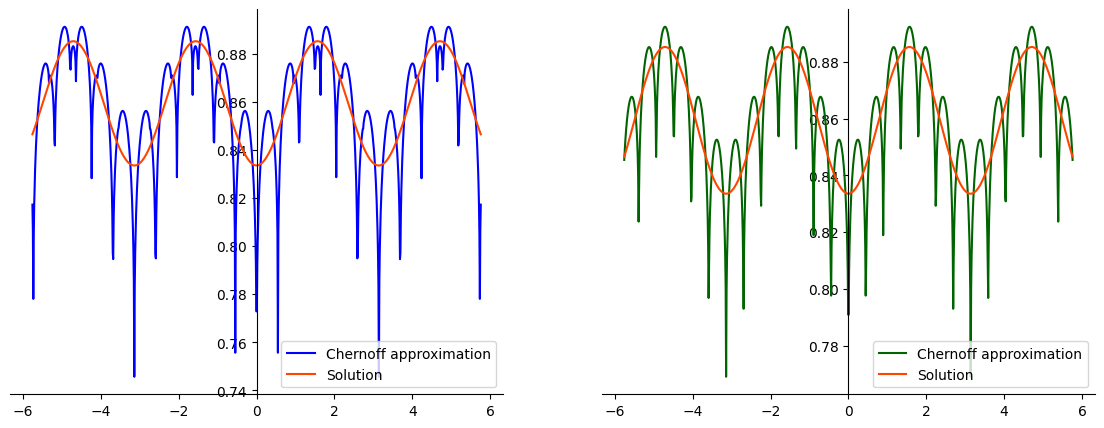}
    {\\fig. 3.1, $n=10$, $u_0(x) = \sqrt[4]{|\sin x|}$, $t=\frac{1}{2}$} 
\end{center}
\begin{center}
\includegraphics[scale=1.0]{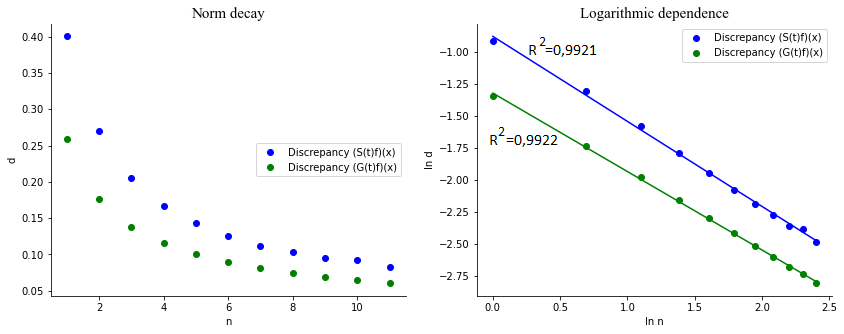}
    {\\fig. 3.2, $1\leq n\leq 11$, $u_0(x) = \sqrt[4]{|\sin x|}$, $t=\frac{1}{2}$} 
\end{center}

Rounding off the coefficients, we see that for the blue line (see Figure 3.2, right), the equation looks as follows: $\ln(d)= -0.6725 \ln(n) -0.8677$, $d=n^{-0.6725}e^{-0.8677}$. Similarly, for the green line (see Fig. 3.2, right), the equation $\ln(d)= -0.6441 \ln(n) -1.2630$, i.e. $d=n^{-0.6441}e^{-1.2630}$. The coefficients of the straight lines were calculated by running linear regression (least squares method), using points corresponding to $4\leq n\leq 11$.

\subsubsection{Approximations for initial condition $u_0(x)=|\sin(x)|^{3/2}$}
\begin{center}
\includegraphics[scale=0.6]{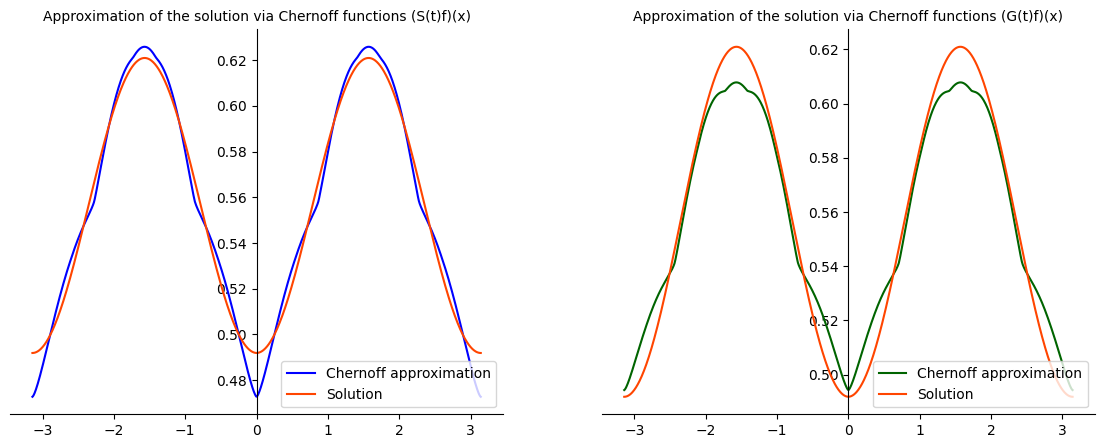}
    {fig. 4.1, $n=4$, $u_0(x)=|\sin(x)|^{3/2}$, $t=\frac{1}{2}$} 
\end{center}

\begin{center}
\includegraphics[scale=1.1]{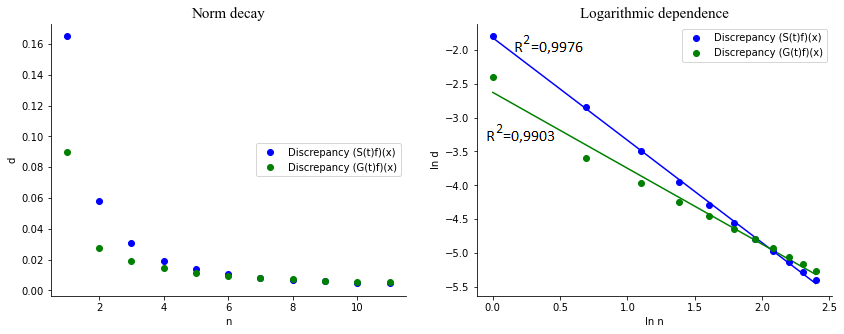}
    {fig. 4.2, $1\leq n\leq 11$, $u_0(x)=|\sin(x)|^{3/2}$, $t=\frac{1}{2}$} 
\end{center}

The exact solution was obtained by numerical calculation of the integral in formula (\ref{exact}).
 The blue color corresponds to the Chernoff function $S(t)$, and the green color corresponds to the Chernoff function $G(t)$. For the green line (see Fig. 4.2, right) the equation $\ln(d)= -0.9997 \ln(n) -2.8538$, i.e. $d=n^{-0.9997}e^{-2.8538}$. Similarly, for the blue line (see Figure 4.2), the equation is as follows: $\ln(d)= -1.4330 \ln(n) -1.9811$, i.e. $d=n^{-1.4330}e^{-1.9811}$. The coefficients of the straight lines were calculated by running linear regression (least squares method), using points corresponding to $4\leq n\leq 11$.

As can be seen from Figure~4.2, the difference between the error decay rates using Chernoff functions $S(t)$ and $G(t)$ for $u_0(x)=|\sin(x)|^{3/2}$ is greater than for $u_0(x) =|\sin x|$. This is due to the greater smoothness of $u_0(x)=|\sin(x)|^{3/2}$. 

\subsubsection{Approximations for initial condition $u_0(x)=|\sin(x)|^{13/2}$}
\begin{center}
\includegraphics[scale=0.6]{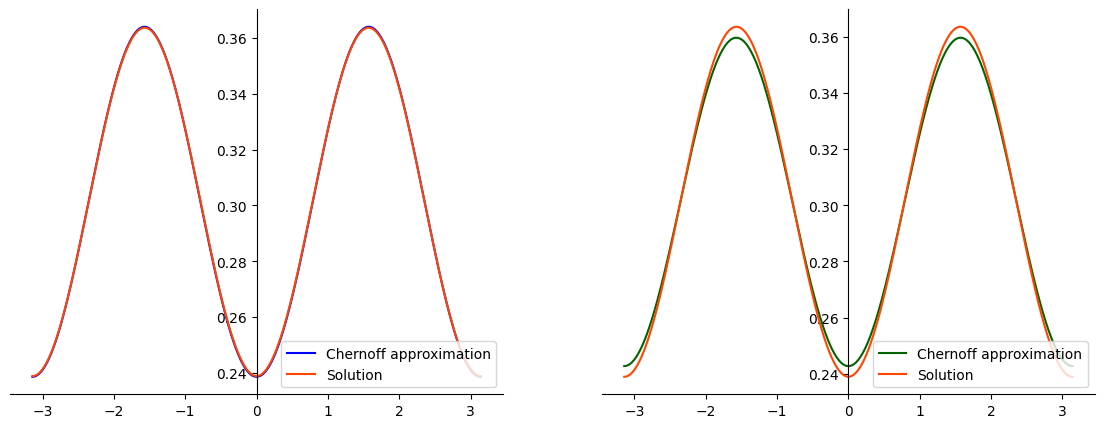}
    {fig. 5.1, $n=11$, $u_0(x)=|\sin(x)|^{13/2}$, $t=\frac{1}{2}$} 
\end{center}

\begin{center}
\includegraphics[scale=0.55]{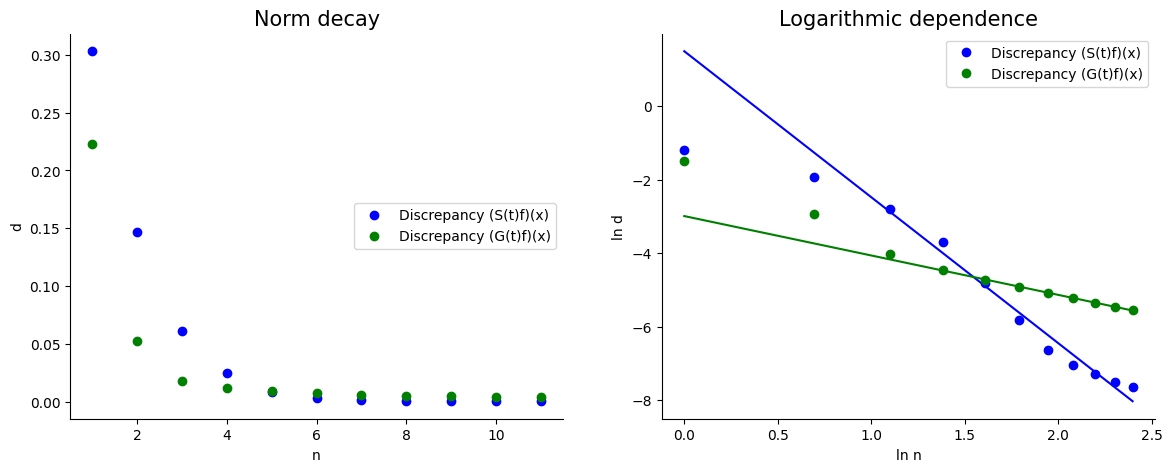}
    {\\fig. 5.2, $1\leq n\leq 11$, $u_0(x)=|\sin(x)|^{13/2}$, $t=\frac{1}{2}$} 
\end{center}

The exact solution was obtained by numerical calculation of the integral in formula (\ref{exact}).
 The blue color corresponds to the Chernoff function $S(t)$, and the green color corresponds to the Chernoff function $G(t)$. For the green line (see Fig. 5.2, right) the equation $\ln(d)= -1.0719 \ln(n) -2.9886$, i.e. $d=n^{-1.0719}e^{-2.9886}$. Similarly, for the blue line (see Figure 5.2), the equation is as follows: $\ln(d)= -3.9717 \ln(n) -1.4969$, i.e. $d=n^{-3.9717}e^{-1.4969}$. The coefficients of the straight lines were calculated by running linear regression (least squares method), using points corresponding to $4\leq n\leq 11$.

For the blue line, the determination coefficient is equal to $R_S^{2}=0.9613$, and for the green line it is equal to $R_G^{2}=0.9989$. However, we visually observe that for the blue line the approximation is not good, so the values of coefficients of the straight lines should not be trusted completely.

\subsection{Graphs of approximations for several exponential initial conditions}

\subsubsection{Approximations for initial condition $u_0(x)=e^{-|x|}$}
Let us consider a non-smooth and non-periodic function $e^{-|x|}$ as an initial condition. The exact solution was obtained by numerical calculation of the integral in formula (\ref{exact}).

\begin{center}
\includegraphics[scale=0.6]{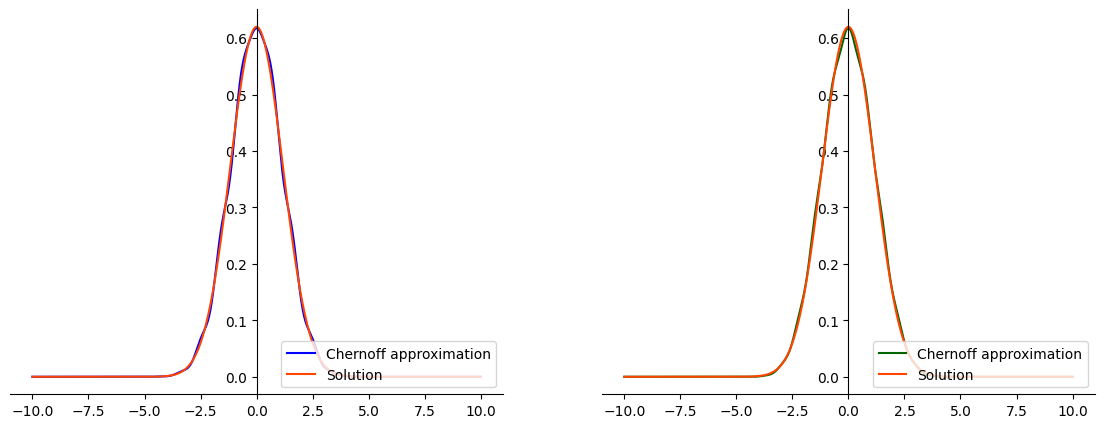}
    {fig. 6.1, $n=4$, $u_0(x) = e^{-|x|}$, $t=\frac{1}{2}$} 
\end{center}
\begin{center}
\includegraphics[scale=1.1]{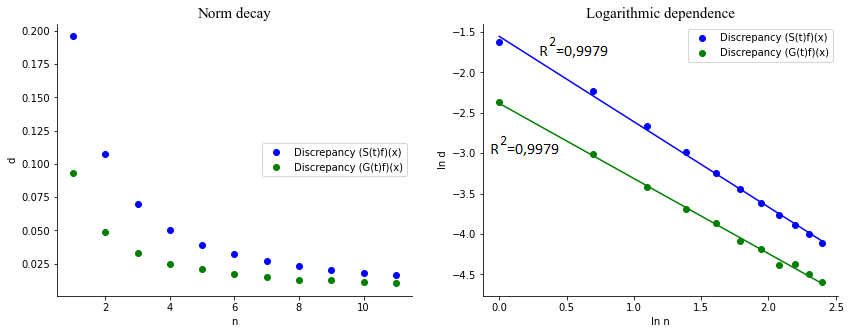}
    {fig. 6.2, $1\leq n\leq 11$, $u_0(x) = e^{-|x|}$, $t=\frac{1}{2}$} 
\end{center}

Figures 6.1 and 6.2 show plots of the exact solution, the approximation to the solution, and the rate of convergence of the error to zero. The blue color corresponds to the Chernoff function $S(t)$, and the green color corresponds to the Chernoff function $G(t)$. As can be seen, the convergence rate for the function $S(t)$ is higher than for $G(t)$.

For the green line (see Figure 6.2, right), the equation is as follows: $\ln(d)= -1.1202 \ln(n) -2.7179$, i.e., $d=n^{-1.1202}e^{-2.7179}$. Similarly, for the blue line (see Figure 6.2), the equation looks like this: $\ln(d)= -1.7803 \ln(n) -1.4788$, that is, $d=n^{-1.7803}e^{-1.4788}$. The coefficients of the straight lines were calculated by running linear regression (least squares method), using points corresponding to $4\leq n\leq 11$.

\subsubsection{Approximations for initial condition $u_0(x)=e^{-|x|^{1/2}}$}

\begin{center}
\includegraphics[scale=0.6]{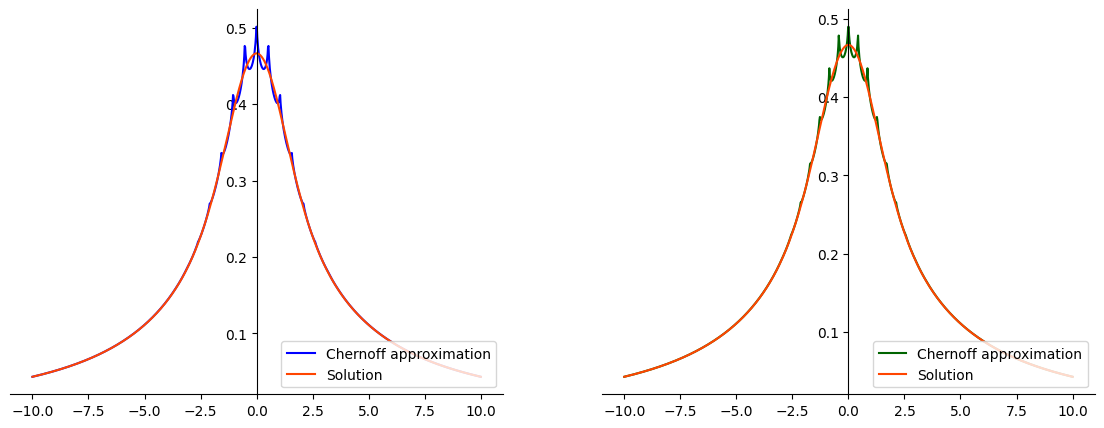}
    {fig. 7.1, $n=4$, $u_0(x) = e^{-|x|^{1/2}}$, $t=\frac{1}{2}$} 
\end{center}
\begin{center}
\includegraphics[scale=0.6]{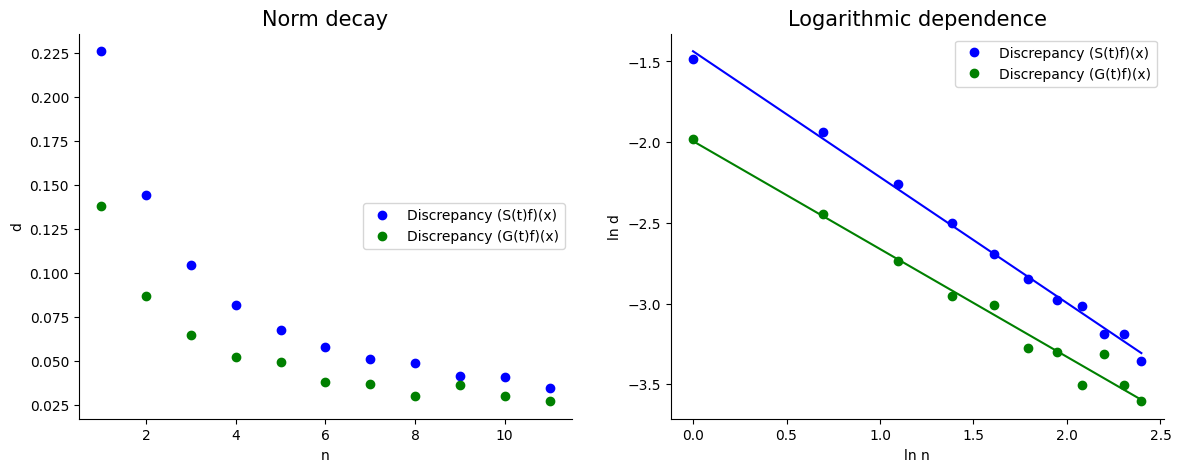}
    {fig. 7.2, $1\leq n\leq 11$, $u_0(x) = e^{-|x|^{1/2}}$, $t=\frac{1}{2}$} 
\end{center}
The exact solution was obtained by numerical calculation of the integral in formula (\ref{exact}).
 The blue color corresponds to the Chernoff function $S(t)$, and the green color corresponds to the Chernoff function $G(t)$. For the green line (see Fig.7.2, right), the equation is as follows: $\ln(d)= -0.6255 \ln(n) -2.0796$, i.e. $d=n^{-0.6255}e^{-2.0796}$. Similarly, for the blue line (see Figure 7.2) the equation is as follows: $\ln(d)= -0.8003\ln(n) -1.3993$, i.e. $d=n^{-0.8003}e^{-1.3993}$. The coefficients of the straight lines were calculated by running linear regression (least squares method), using points corresponding to $4\leq n\leq 11$.

\subsubsection{Approximations for initial condition $u_0(x)=e^{-|x|^{3/2}}$}

\begin{center}
\includegraphics[scale=0.6]{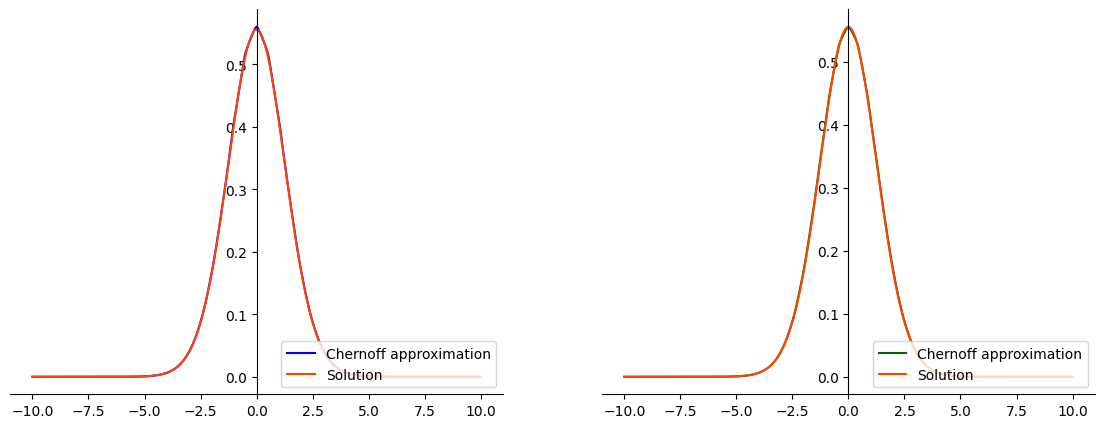}
    {fig. 8.1, $n=11$, $u_0(x) = e^{-|x|^{3/2}}$, $t=\frac{1}{2}$} 
\end{center}
\begin{center}
\includegraphics[scale=0.6]{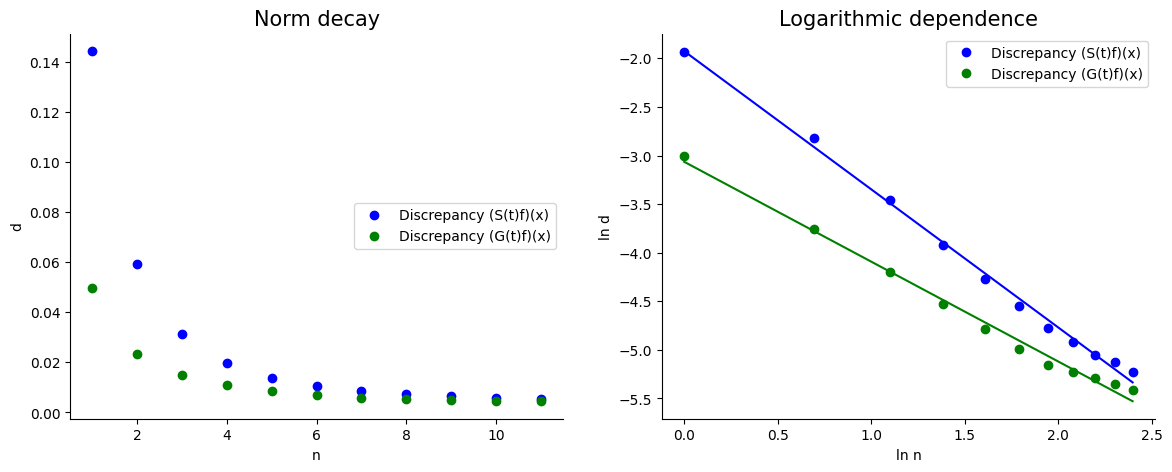}
    {fig. 8.2, $1\leq n\leq 11$, $u_0(x) = e^{-|x|^{3/2}}$, $t=\frac{1}{2}$} 
\end{center}
The exact solution was obtained by numerical calculation of the integral in formula (\ref{exact}).
 The blue color corresponds to the Chernoff function $S(t)$, and the green color corresponds to the Chernoff function $G(t)$. For the green line (see Fig. 8.2, right), the equation is as follows: $\ln(d)= -0.8578 \ln(n) -3.4081$, i.e. $d=n^{-0.8578}e^{-3.4081}$. Similarly, for the blue line (see Figure 8.2) the equation is as follows: $\ln(d)= -1.2859\ln(n) -2.2059$, i.e. $d=n^{-1.2859}e^{-2.2059}$. The coefficients of the straight lines were calculated by running linear regression (least squares method), using points corresponding to $4\leq n\leq 11$.

\subsection{Tables of rates and free terms of the approximation error for all initial conditions}

Experimentally (using numerical simulation in Python 3) we calculated the error of Chernoff approximations for many initial conditions, not only for those that are presented above with graphs. The tables below show the orders of decreasing of the error depending on the smoothness class of the initial condition and the Chernoff function.

\begin{remark}\label{notrem}
\textbf{We use the following notation.} In the expression 
$$
\ln\sup_{x\in\mathbb{R}}|u(t,x)-u_n(t,x)|\approx \alpha \ln n + \beta
$$ 
we call $\alpha$ the order of decreasing error, and we call $\beta$ the free term of decreasing error. We use the letter $\xi$ to denote the smoothness class of the initial condition. We use letters $\alpha_S$ and $\beta_S$ for the Chernoff function $S$, and $\alpha_G$ and $\beta_G$ for the Chernoff function $G$. Remark \ref{remord} provides the necessary background.
\end{remark}

\subsubsection{Trigonometric initial conditions}

\scriptsize

\begin{tabular}{|m{6cm}|m{4cm}|m{4cm}|} 
\hline
\textbf{$\xi$, the smoothness class of the initial condition $u_0$}
& \textbf{$\alpha_G$, the order of decreasing error on the Chernoff function $G(t)$, which has the 1st order of the Chernoff tangency to the operator $L=\frac{d^2}{dx^2}$}
& \textbf{$\alpha_S$, the order of decreasing error on the Chernoff function $S(t)$, which has the 2nd order of tangency by Chernoff to the operator $L=\frac{d^2}{dx^2}$}\\ 
\hline

$H^{6\frac{1}{2}}$,  derivatives of orders 1-6 exist and are bounded, and the 6th derivative is H\"older with a H\"older exponent 1/2, $u_0(x)=|\sin(x)|^{13/2}$ & -1.0719 & -3.9717 (but the quality of the regression is not perfect)\\
\hline

$H^{5\frac{1}{2}}$,  derivatives of orders 1-5 exist and are bounded, and the 5th derivative is H\"older with a H\"older exponent 1/2, $u_0(x)=|\sin(x)|^{11/2}$ & -1.0451 & -3.3196 \\
\hline

$H^{4\frac{1}{2}}$,  the first, second, third, and fourth derivatives exist and are bounded, and the fourth is H\"older with a H\"older exponent 1/2, $u_0(x)=|\sin(x)|^{9/2}$ & -1.0456 & -2.5275\\
\hline

$H^{3\frac{1}{2}}$,  the first, second, and third derivatives exist and are bounded, and the third is H\"older with H\"older exponent 1/2, $u_0(x)=|\sin(x)|^{7/2}$ & -1.0948 & -2.0340\\
\hline
$H^{3}$, the first and second derivative exist and are bounded, while the second derivative is H\"older continuous with H\"older exponent 1/2, $u_0(x)=\sin(x)\cdot|\sin(x)|^{2}$ & -1.1358 & -1.9600\\ 
\hline
$H^{2\frac{1}{2}}$, the first and second derivative exist and are bounded, while the second derivative is H\"older continuous with H\"older exponent 1/2, $u_0(x)=|\sin(x)|^{5/2}$ & -1.1259 & -1.7464\\ 
\hline
$H^{2}$, the first and second derivative exist and are bounded, while the second derivative is H\"older continuous with H\"older exponent 1/2, $u_0(x)=|\sin(x)|\cdot\sin(x)$ & -1.1082& -1.2573\\ 
\hline

$H^{1\frac{1}{2}}$, the first derivative: exists, is bounded and H\"older continuous with H\"older exponent 1/2, $u_0(x)=|\sin(x)|^{3/2}$ & -0.9997 & -1.4330\\ 
\hline

$H^{1}$, H\"older continuous with the H\"older exponent 1, $u_0(x)=|\sin(x)|$ & -0.9655 & -1.0818\\ 
\hline

$H^{3/4}$, H\"older continuous with the H\"older exponent 3/4, $u_0(x)=|\sin(x)|^{3/4}$ & -0.8409 & -0.9286\\
\hline

$H^{1/2}$, H\"older continuous with the H\"older exponent 1/2, $u_0(x)=|\sin(x)|^{1/2}$ & -0.7170 & -0.7835\\
\hline

$H^{1/4}$, H\"older continuous with the H\"older exponent 1/4, $u_0(x)=|\sin(x)|^{1/4}$ & -0.6441 & -0.6725\\
\hline
\end{tabular}

\begin{flushleft}
Table 1. Orders of decreasing error for trigonometric initial conditions.
\end{flushleft}
\normalsize

\vskip2mm

Keeping in mind Remark \ref{remord} we see that on the initial condition with high smoothness (first line in the table), the first order of Chernoff tangency corresponds to a decreasing error rate of about $1/n$, and the second order corresponds to a decreasing error rate of about $1/n^2$. This is in accordance with theorem \ref{grth}, the conjecture in \cite{R2} and the theorem from \cite{Galkin2022}. 

As smoothness decreases (other lines in the table), the theory from \cite{Galkin2022} stops working, and the experimental evidence is the following: the convergence speed gradually decreases and the advantages of the Chernoff function with the second order of the Chernoff tangency gradually vanish. Let us present the results from the table graphically:

\begin{center}
\includegraphics[scale=0.6]{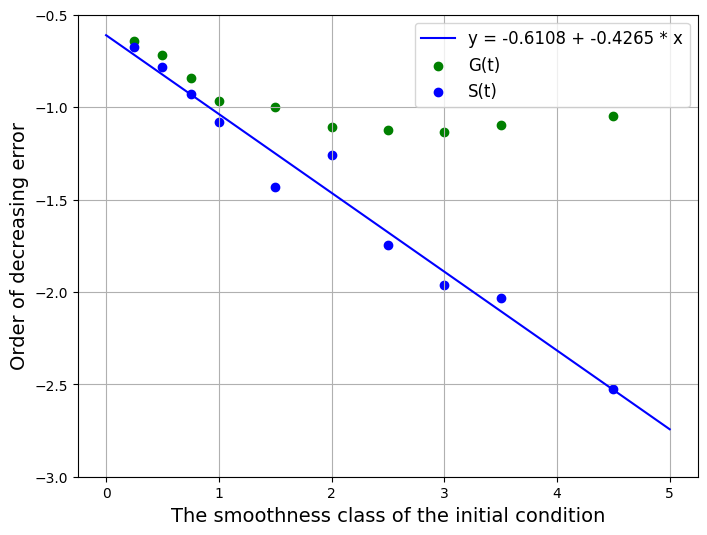}\\
    {fig. 9} 
\end{center}

We see that first several blue points visually lie on a straight line, so we used linear regression (ordinary least squares method) to find the equation of this line.  The equation of the approximating line: $y=-0.6108-0.4265x$. This can be interpreted as follows: When the smoothness class $\xi$ of the initial condition $u_0$ is not greater than the order of Chernoff tangency then$$
d=\|u_n(t,\cdot)-u(t,\cdot)\|=\sup_{x\in\mathbb{R}}|u_n(t,x)-u(t,x)|\approx\mathrm{const}\cdot\left(\frac{1}{n}\right)^{0.61\xi+0.42}.
$$ 
Meanwhile when the smoothness class $\xi$ of the initial condition $u_0$ is greater than the order of Chernoff tangency then there is no such easy-to-state dependence, but the Chernoff function $S(t)$ with the second-order Chernoff tangency provides better approximations than the Chernoff function $G(t)$ with the first-order Chernoff tangency. 

\scriptsize

\begin{tabular}{|m{6cm}|m{4cm}|m{4cm}|} 
\hline
\textbf{$\xi$, the smoothness class of the initial condition $u_0$}
& \textbf{$\beta_G$, the free term of the decreasing error on the Chernoff function $G(t)$, which has the 1st order of the Chernoff tangency to the operator $L=\frac{d^2}{dx^2}$}
& \textbf{$\beta_S$, the free term of the decreasing error on the Chernoff function $S(t)$, which has the 2nd order of tangency by Chernoff to the operator $L=\frac{d^2}{dx^2}$}\\ 
\hline

$H^{6\frac{1}{2}}$,  derivatives of orders 1-6 exist and are bounded, and the 6th derivative is H\"older with a H\"older exponent 1/2, $u_0(x)=|\sin(x)|^{13/2}$ & -2.9886 & 1.4969 (but the quality of the regression is not perfect)\\
\hline

$H^{5\frac{1}{2}}$,  derivatives of orders 1-5 exist and are bounded, and the 5th derivative is H\"older with a H\"older exponent 1/2, $u_0(x)=|\sin(x)|^{11/2}$ & -3.0152 & -0.0066\\
\hline

$H^{4\frac{1}{2}}$,  the first, second, third, and fourth derivatives exist and are bounded, and the fourth is H\"older with a H\"older exponent 1/2, $u_0(x)=|\sin(x)|^{9/2}$ & -2.9836 & -1.6601\\
\hline

$H^{3\frac{1}{2}}$,  the first, second, and third derivatives exist and are bounded, and the third is H\"older with H\"older exponent 1/2, $u_0(x)=|\sin(x)|^{7/2}$ & -2.8339 &  -2.4633\\
\hline
$H^{3}$, the first and second derivative exist and are bounded, while the second derivative is H\"older continuous with H\"older exponent 1/2, $u_0(x)=\sin(x)\cdot|\sin(x)|^{2}$ & -2.6948 & -2.2677\\ 
\hline
$H^{2\frac{1}{2}}$, the first and second derivative exist and are bounded, while the second derivative is H\"older continuous with H\"older exponent 1/2, $u_0(x)=|\sin(x)|^{5/2}$ & -2.6883 & -2.6341\\ 
\hline

$H^{2}$, the first and second derivative exist and are bounded, while the second derivative is H\"older continuous with H\"older exponent 1/2, $u_0(x)=|\sin(x)|\cdot\sin(x)$ & -3.1800& -2.6904\\ 

\hline
$H^{1\frac{1}{2}}$, the first derivative: exists, is bounded and H\"older continuous with H\"older exponent 1/2, $u_0(x)=|\sin(x)|^{3/2}$ & -2.8538 & -1.9811\\ 
\hline

$H^{1}$,  H\"older continuous with the H\"older exponent 1, $u_0(x)=|\sin(x)|$ & -2.3452 & -1.3846\\ 
\hline

$H^{3/4}$, H\"older continuous with the H\"older exponent 3/4, $u_0(x)=|\sin(x)|^{3/4}$ & -1.8579 & -1.1302\\
\hline

$H^{1/2}$, H\"older continuous with the H\"older exponent 1/2, $u_0(x)=|\sin(x)|^{1/2}$ & -1.4186 & -0.8820\\
\hline

$H^{1/4}$, H\"older continuous with the H\"older exponent 1/4, $u_0(x)=|\sin(x)|^{1/4}$ & -1.2630&-0.8677 \\
\hline
\end{tabular}
\begin{flushleft}
Table 2. Free terms of decreasing error for trigonometric initial conditions.
\end{flushleft}
\normalsize
\normalsize

\begin{center}
\includegraphics[scale=0.6]{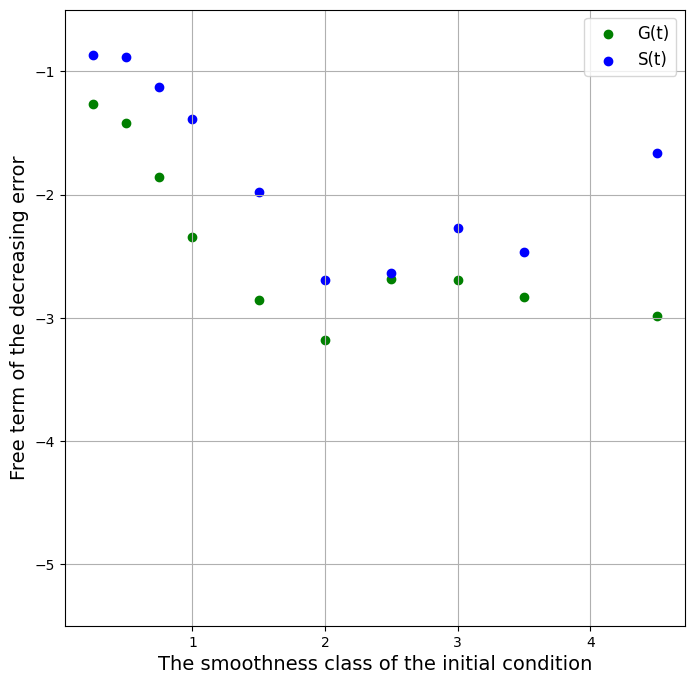}\\
    {fig. 10} 
\end{center}

In the following, we present tables for initial conditions in the form of an exponent of a power-law function.

\subsubsection{Exponential initial conditions}

\scriptsize

\begin{tabular}{|m{6cm}|m{4cm}|m{4cm}|} 
\hline
\textbf{$\xi$, the smoothness class of the initial condition $u_0$}
& \textbf{$\alpha_G$, the order of decreasing error on the Chernoff function $G(t)$, which has the 1st order of the Chernoff tangency to the operator $L=\frac{d^2}{dx^2}$}
& \textbf{$\alpha_S$, the order of decreasing error on the Chernoff function $S(t)$, which has the 2nd order of tangency by Chernoff to the operator $L=\frac{d^2}{dx^2}$}\\ 
\hline

$H^{6\frac{1}{2}}$,  derivatives of orders 1-6 exist and are bounded, and the 6th derivative is H\"older with a H\"older exponent 1/2, $u_0(x)=e^{-|x|^{13/2}}$ & -1.5215 & -1.3880\\
\hline

$H^{5\frac{1}{2}}$,  derivatives of orders 1-5 exist and are bounded, and the 5th derivative is H\"older with a H\"older exponent 1/2, $u_0(x)=e^{-|x|^{11/2}}$ & -1.5545 & -1.6191\\
\hline

$H^{4\frac{1}{2}}$,  the first, second, third, and fourth derivatives exist and are bounded, and the fourth is H\"older with a H\"older exponent 1/2, $u_0(x)=e^{-|x|^{9/2}} $ & -1.3374 & -1.8367 \\
\hline

$H^{3\frac{1}{2}}$,  the first, second, and third derivatives exist and are bounded, and the third is H\"older with H\"older exponent 1/2, $u_0(x)=e^{-|x|^{7/2}}$ & -0.8595 & -1.4800 \\
\hline
$H^{3}$, the first and second derivative exist and are bounded, while the second derivative is H\"older continuous with H\"older exponent 1/2, $u_0(x)= e^{-|x|\cdot x^{2}}$ & -0.8144 & -1.1437\\ 
\hline
$H^{2\frac{1}{2}}$, the first and second derivative exist and are bounded, while the second derivative is H\"older continuous with H\"older exponent 1/2, $u_0(x)=e^{-|x|^{5/2}}$ & -0.9374 & -1.0803\\ 
\hline
$H^{2}$, the first and second derivative exist and are bounded, while the second derivative is H\"older continuous with H\"older exponent 1/2, $u_0(x)=e^{x|x|} \cdot e^{x^4}$ & -1.6107 & -2.3204\\ 
\hline

$H^{1\frac{1}{2}}$, the first derivative: exists, is bounded and H\"older continuous with H\"older exponent 1/2, $u_0(x)=e^{-|x|^{3/2}}$ & -0.8578 & -1.2859\\ 
\hline

$H^{1}$, H\"older continuous with the H\"older exponent 1, $u_0(x)=e^{-|x|}$ & -1.1202 & -1.7803\\
\hline

$H^{3/4}$, H\"older continuous with the H\"older exponent 3/4, $u_0(x)=e^{-|x|^{3/4}}$ & -0.7637 & -0.9526\\
\hline

$H^{1/2}$, H\"older continuous with the H\"older exponent 1/2, $u_0(x)=e^{-|x|^{1/2}}$ & -0.6255 & -0.8003\\
\hline

$H^{1/4}$, H\"older continuous with the H\"older exponent 1/4, $u_0(x)=e^{-|x|^{1/4}}$ & -0.5302 & -0.6543\\
\hline
\end{tabular}
\begin{flushleft}
Table 3. Orders of decreasing error for exponential initial conditions.
\end{flushleft}
\normalsize

Figure 11 shows the data from Table 3 graphically.

\begin{center}
\includegraphics[scale=0.6]{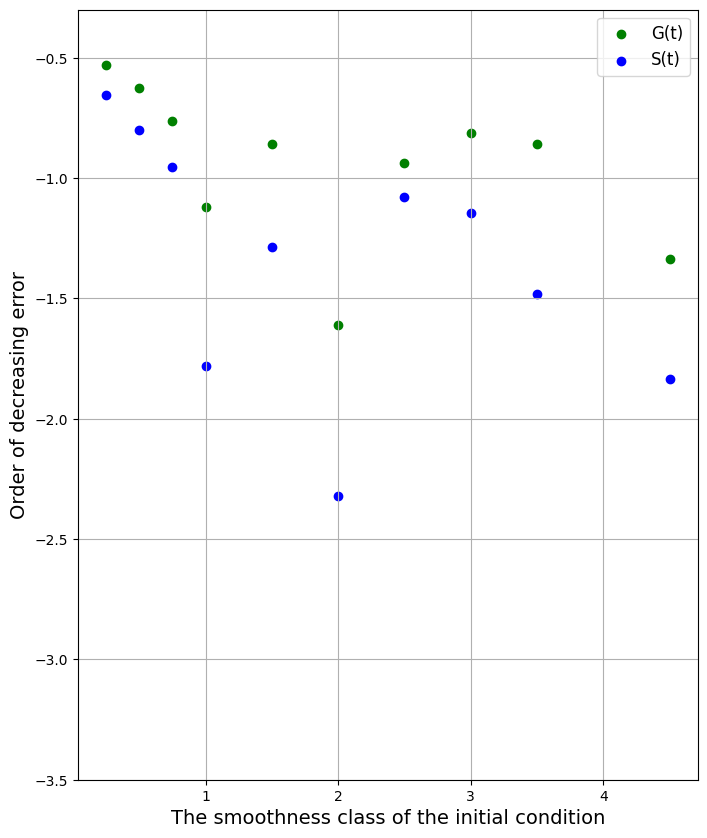}\\
    {fig. 11} 
\end{center}


\scriptsize

\begin{tabular}{|m{6cm}|m{4cm}|m{4cm}|} 
\hline
\textbf{$\xi$, the smoothness class of the initial condition $u_0$}
& \textbf{$\beta_G$, the free term of the decreasing error  on the Chernoff function $G(t)$, which has the 1st order of the Chernoff tangency to the operator $L=\frac{d^2}{dx^2}$}
& \textbf{$\beta_S$, the free term of the decreasing error  on the Chernoff function $S(t)$, which has the 2nd order of tangency by Chernoff to the operator $L=\frac{d^2}{dx^2}$}\\ 
\hline

$H^{6\frac{1}{2}}$,  derivatives of orders 1-6 exist and are bounded, and the 6th derivative is H\"older with a H\"older exponent 1/2, $u_0(x)=e^{-|x|^{13/2}}$ & -1.0662 & -0.9040\\
\hline

$H^{5\frac{1}{2}}$,  derivatives of orders 1-5 exist and are bounded, and the 5th derivative is H\"older with a H\"older exponent 1/2, $u_0(x)=e^{-|x|^{11/2}}$ & -1.3350 & -0.8259\\
\hline

$H^{4\frac{1}{2}}$,  the first, second, third, and fourth derivatives exist and are bounded, and the fourth is H\"older with a H\"older exponent 1/2, $u_0(x)=e^{-|x|^{9/2}} $ & -2.1309 & -0.9857 \\
\hline

$H^{3\frac{1}{2}}$,  the first, second, and third derivatives exist and are bounded, and the third is H\"older with H\"older exponent 1/2, $u_0(x)=e^{-|x|^{7/2}}$ & -3.3532 & -2.4567 \\
\hline
$H^{3}$, the first and second derivative exist and are bounded, while the second derivative is H\"older continuous with H\"older exponent 1/2, $u_0(x)= e^{-|x|\cdot x^{2}}$ & -3.4844 & -3.3309\\ 
\hline
$H^{2\frac{1}{2}}$, the first and second derivative exist and are bounded, while the second derivative is H\"older continuous with H\"older exponent 1/2, $u_0(x)=e^{-|x|^{5/2}}$ & -3.2500& -3.4135\\ 
\hline
 $H^{2}$, the first and second derivative exist and are bounded, while the second derivative is H\"older continuous with H\"older exponent 1/2, $u_0(x)=e^{x|x|}\cdot e^{x^4}$ & -1.5613 & 0.0706\\ 
\hline 

$H^{1\frac{1}{2}}$, the first derivative: exists, is bounded and H\"older continuous with H\"older exponent 1/2, $u_0(x)=e^{-|x|^{3/2}}$ & -3.4081 & -2.2059\\ 
\hline

$H^{1}$, H\"older continuous with the H\"older exponent 1, $u_0(x)=e^{-|x|}$ & -2.7179 & -1.4788\\
\hline

$H^{3/4}$, H\"older continuous with the H\"older exponent 3/4, $u_0(x)=e^{-|x|^{3/4}}$ & -2.1767 & -1.3566\\
\hline

$H^{1/2}$, H\"older continuous with the H\"older exponent 1/2, $u_0(x)=e^{-|x|^{1/2}}$ & -2.0796 & -1.3993\\
\hline

$H^{1/4}$, H\"older continuous with the H\"older exponent 1/4, $u_0(x)=e^{-|x|^{1/4}}$ &  -2.2717 & -1.7918\\
\hline
\end{tabular}
\begin{flushleft}
Table 4. Free terms of decreasing error for exponential initial conditions.
\end{flushleft}
\normalsize

Figure 12 shows the data from Table 4 graphically.

\begin{center}
\includegraphics[scale=0.6]{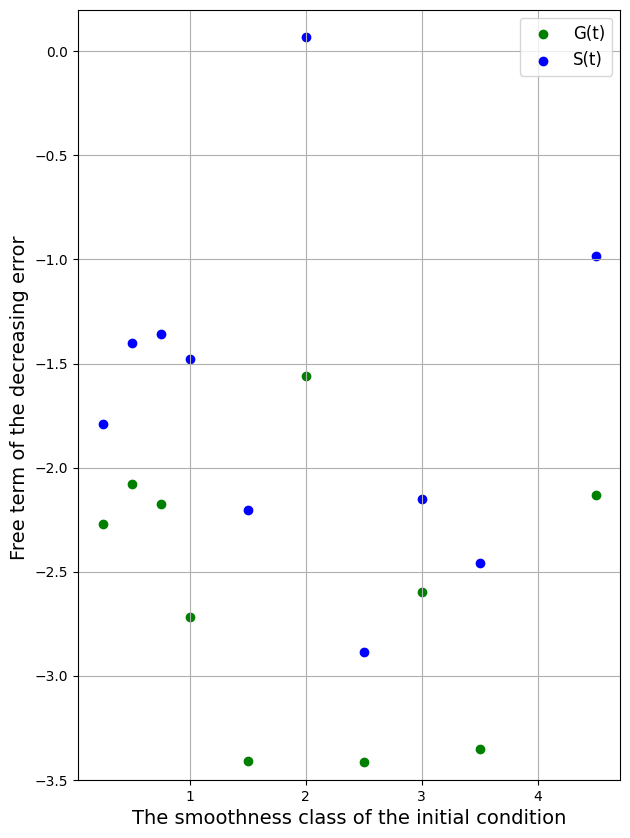}\\
    {fig. 12} 
\end{center}

\normalsize

\subsubsection{Infinitely smooth initial conditions}

Now let us consider only infinitely smooth initial conditions. The graphs, as well as the values of $\alpha$ and $\beta$ are given below. The exact solution in each case was obtained by numerical calculation of the integral in formula (\ref{exact}).

\begin{center}
\includegraphics[scale=0.6]{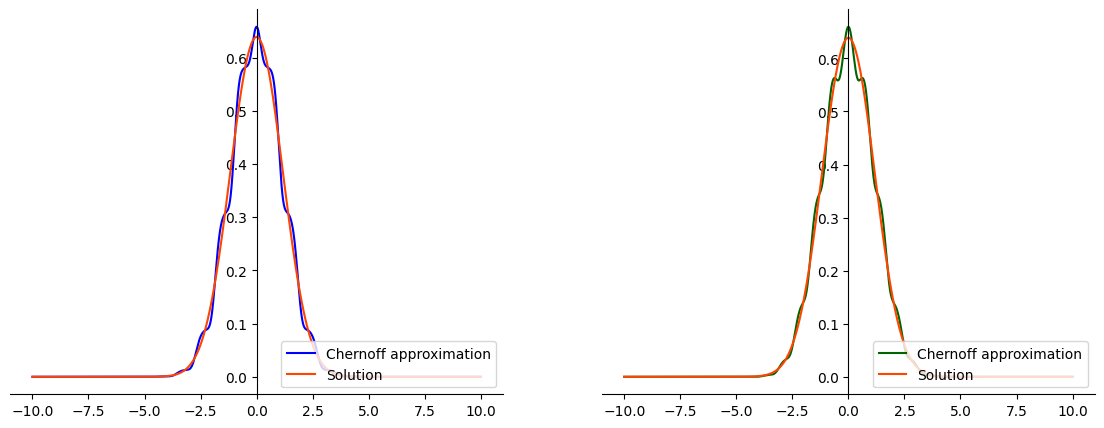}
    {fig. 13.1, $n=4$, $u_0(x) = e^{-x^{6}}$, $t=\frac{1}{2}$} 
\end{center}
\begin{center}
\includegraphics[scale=0.6]{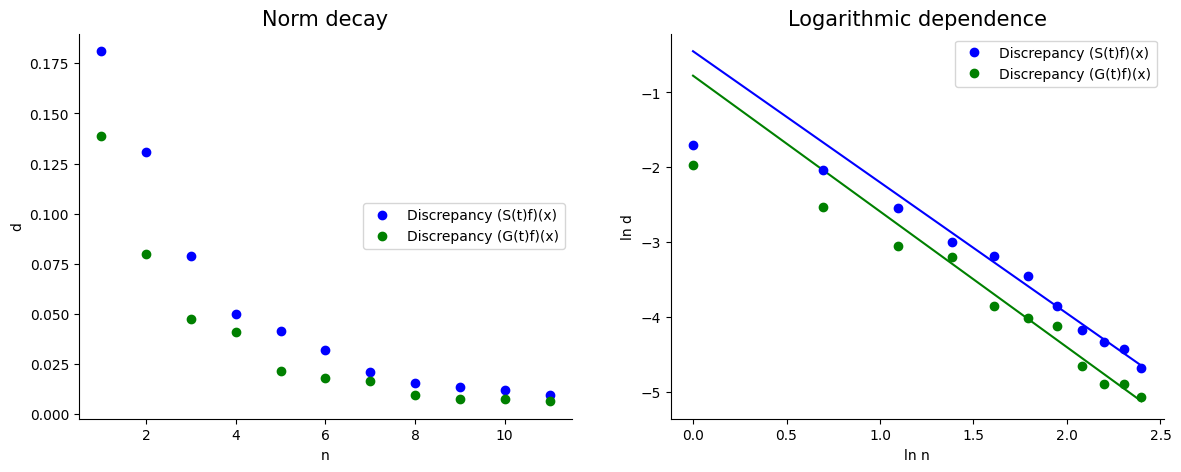}
    {fig. 13.2, $1\leq n\leq 11$, $u_0(x) = e^{-x^{6}}$, $t=\frac{1}{2}$} 
\end{center}
The blue color corresponds to the Chernoff function $S(t)$, and the green color corresponds to the Chernoff function $G(t)$. For the green line (see Fig. 13.2, right), the equation is as follows: $\ln(d)= -1.8139 \ln(n) -0.7760$, i.e. $d=n^{-1.8139}e^{-0.7760}$. Similarly, for the blue line (see Figure 13.2) the equation is as follows: $\ln(d)= -1.7517\ln(n) -0.4493$, i.e. $d=n^{-1.7517}e^{-0.4493}$. The coefficients of the straight lines were calculated by running linear regression (least squares method), using points corresponding to $4\leq n\leq 11$.

\begin{center}
\includegraphics[scale=0.6]{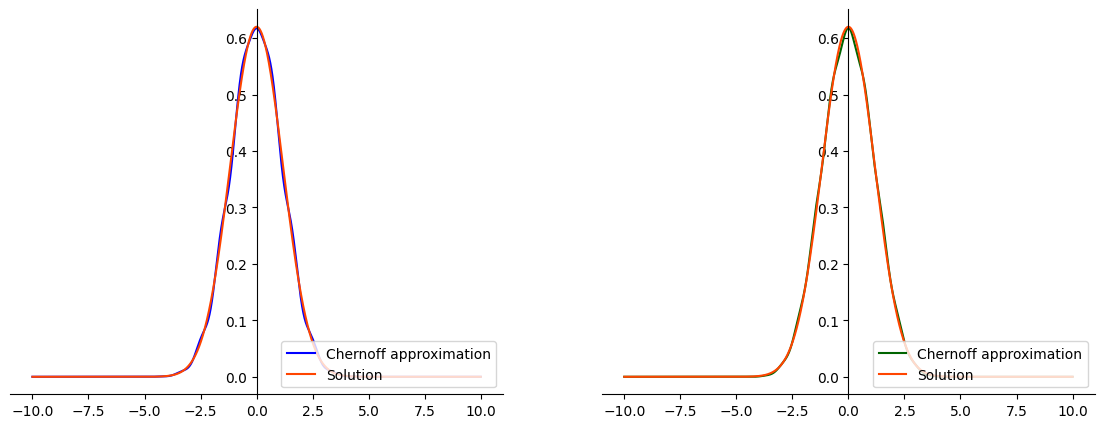}
    {fig. 14.1, $n=4$, $u_0(x) = e^{-x^4}$, $t=\frac{1}{2}$} 
\end{center}

\begin{center}
\includegraphics[scale=0.6]{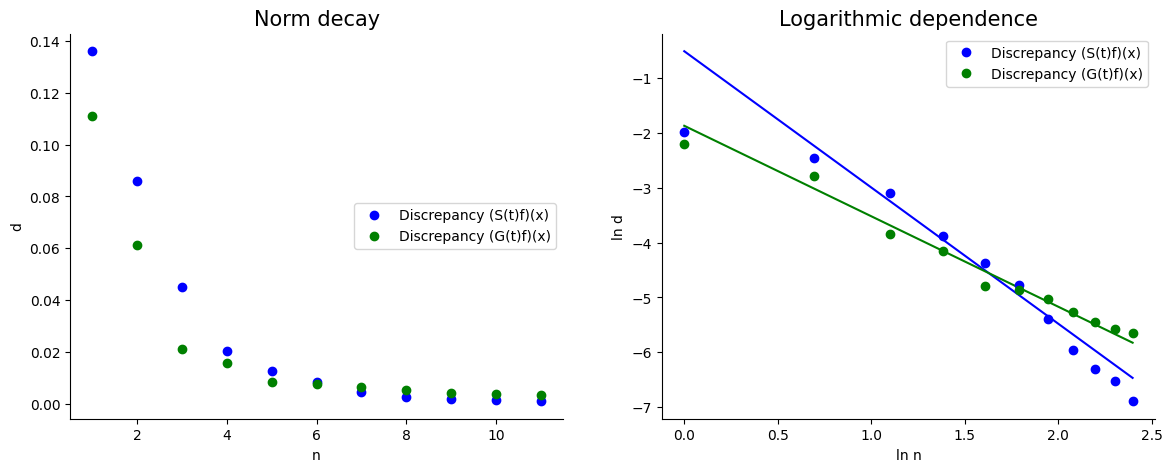}
    {fig. 14.2, $1\leq n\leq 11$, $u_0(x) = e^{-x^4}$, $t=\frac{1}{2}$} 
\end{center}
The blue color corresponds to the Chernoff function $S(t)$, and the green color corresponds to the Chernoff function $G(t)$. Figures 14.1 and 14.2 show plots of the exact solution, the approximation to the solution, and the rate of convergence of the error to zero. As can be seen, the rate of convergence of the function $S(t)$ is higher than that of $G(t)$.

For the green line (see Figure 14.2, right), the equation is as follows: $\ln(d)= -1.3948 \ln(n) -2.3548$, i.e., $d=n^{-1.3948}e^{-2.3548}$. Similarly, for the blue line (see Figure 14.2), the equation looks like this: $\ln(d)= -3.0768 \ln(n) -0.5226$, i.e., $d=n^{-3.0768}e^{-0.5326}$. The coefficients of the straight lines were calculated by running linear regression (least squares method), using points corresponding to $4\leq n\leq 11$.

\begin{center}
\includegraphics[scale=0.6]{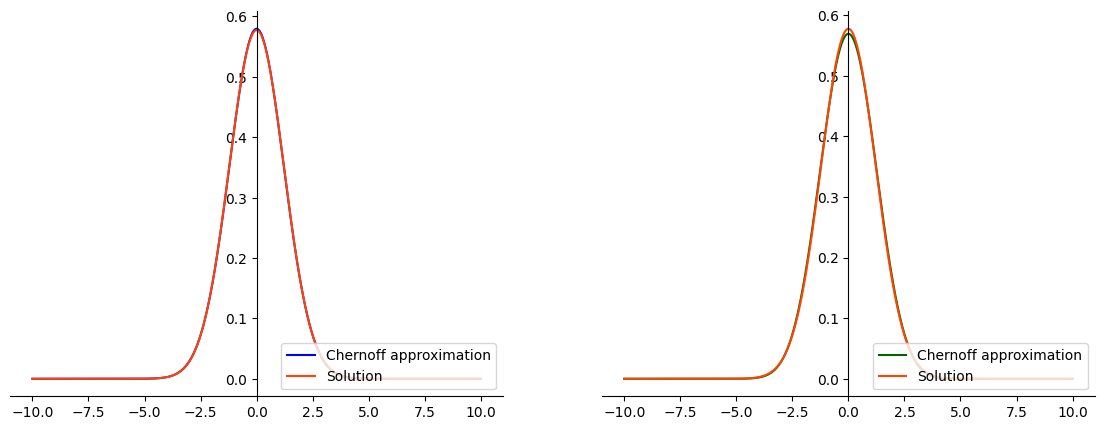}\\
    {fig. 15.1, $n=4$, $u_0(x) = e^{-x^2}$, $t=\frac{1}{2}$}  
\end{center}

\begin{center}
\includegraphics[scale=0.6]{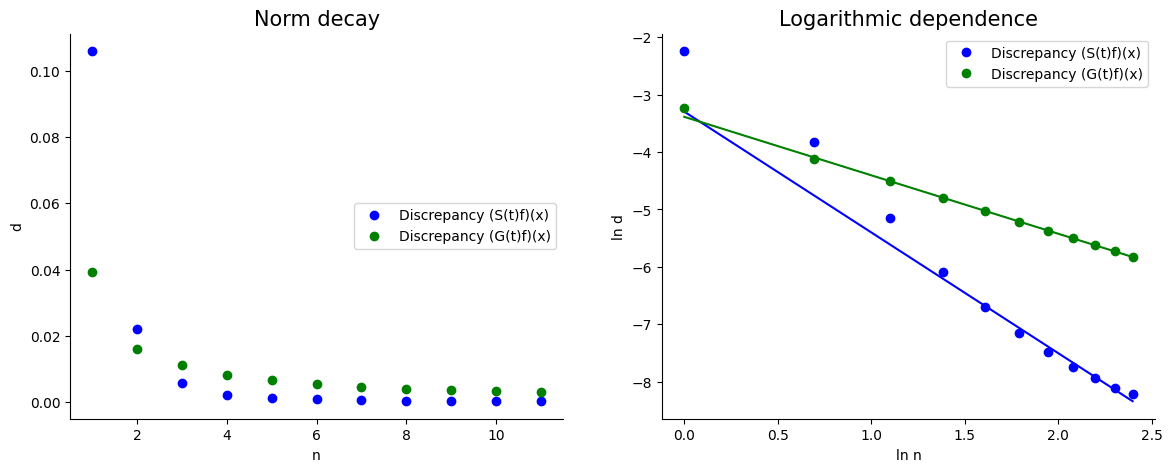}\\
    {fig. 15.2, $1\leq n\leq 11$, $u_0(x) = e^{-x^2}$, $t=\frac{1}{2}$} 
\end{center}

For the green line (see Figure 15.2, right), the equation is as follows: $\ln(d)= -1.0177 \ln(n) -3.3874$, i.e., $d=n^{-1.0177}e^{-3.3874}$. Similarly, for the blue line (see Figure 15.2), the equation looks like this: $\ln(d)= -2.1028 \ln(n) -3.2967$, i.e., $d=n^{-2.1028}e^{-3.2967}$. The coefficients of the straight lines were calculated by running linear regression (least squares method), using points corresponding to $4\leq n\leq 11$.

\begin{center}
\includegraphics[scale=0.4]{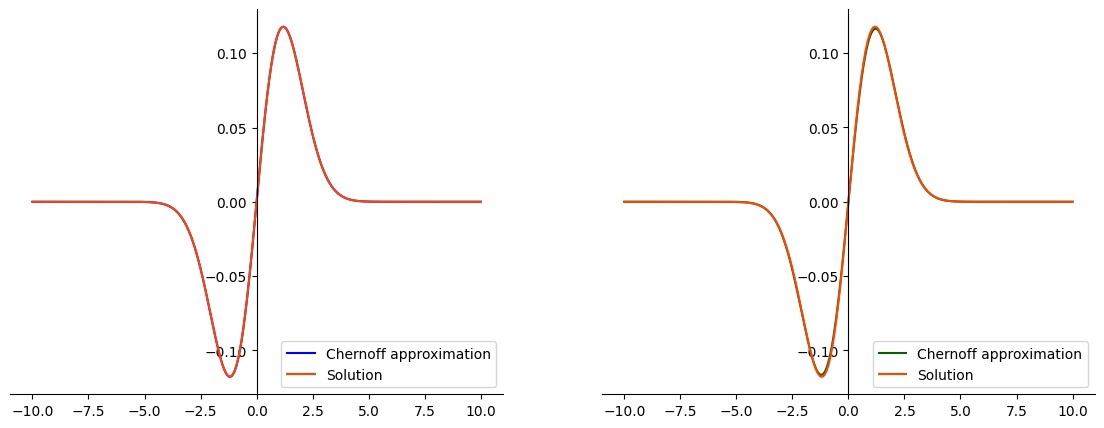}\\
    {fig. 16.1, $n=11$, $u_0(x) = e^{-x^2}\cdot\sin(x)$, $t=\frac{1}{2}$}  
\end{center}

\begin{center}
\includegraphics[scale=0.4]{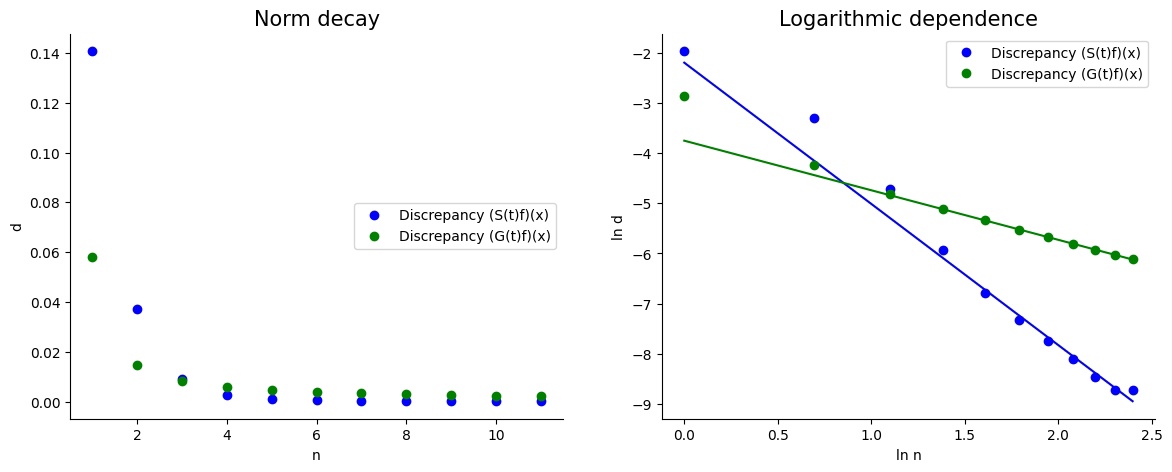}\\
    {fig. 16.2, $1\leq n\leq 11$, $u_0(x) = e^{-x^2}\cdot\sin(x)$, $t=\frac{1}{2}$} 
\end{center}

For the green line (see Figure 16.2, right), the equation is as follows: $\ln(d)= -0.9897 \ln(n) -3.7478$, i.e., $d=n^{-0.9897}e^{-3.7478}$. Similarly, for the blue line (see Figure 16.2), the equation looks like this: $\ln(d)= -2.8194 \ln(n) -2.1899$, i.e., $d=n^{-2.8194}e^{-2.1899}$. The coefficients of the straight lines were calculated by running linear regression (least squares method), using points corresponding to $4\leq n\leq 11$.

\begin{center}
\includegraphics[scale=0.5]{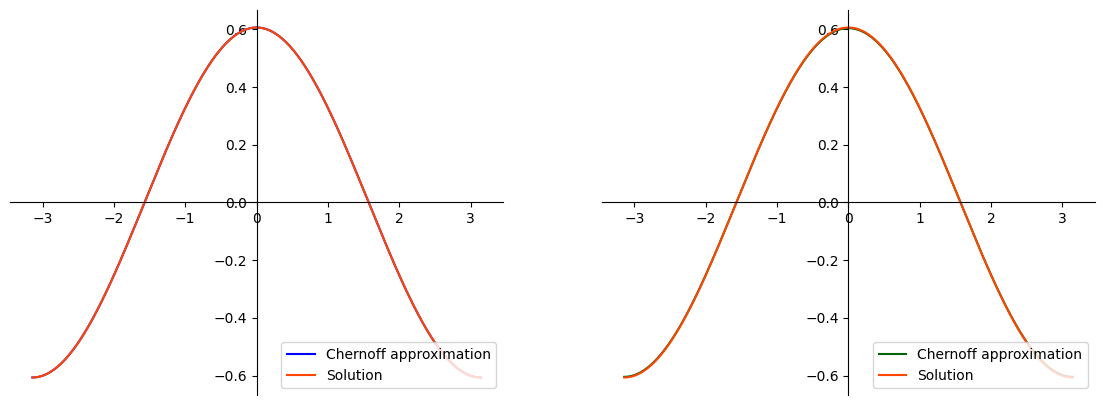}\\
    {fig. 17.1, $n=11$, $u_0(x) = \cos(x)$, $t=\frac{1}{2}$}  
\end{center}

\begin{center}
\includegraphics[scale=0.5]{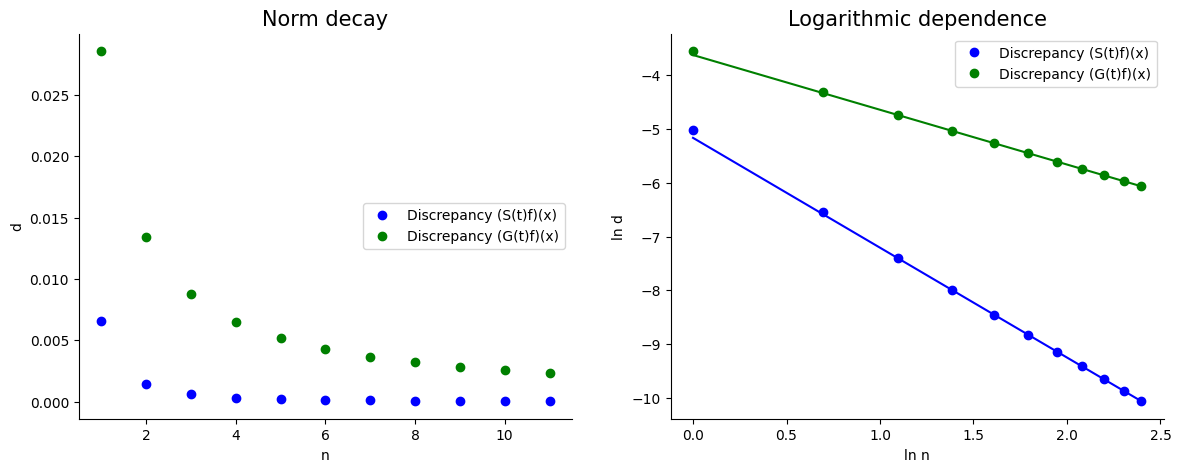}\\
    {fig. 17.2, $1\leq n\leq 11$, $u_0(x) = \cos(x)$, $t=\frac{1}{2}$} 
\end{center}

For the green line (see Figure 17.2, right), the equation is as follows: $\ln(d)= -1.0178 \ln(n) -3.6261 $, i.e., $d=n^{-1.0178}e^{-3.6261 }$. Similarly, for the blue line (see Figure 17.2), the equation looks like this: $\ln(d)= -2.0415 \ln(n) -5.1659$, i.e., $d=n^{-2.0415}e^{-5.16599}$. The coefficients of the straight lines were calculated by running linear regression (least squares method), using points corresponding to $4\leq n\leq 11$. 

In the Tables 5 and 6 below, we summarize the data provided above. 

\footnotesize

\begin{tabular}{|m{6cm}|m{4cm}|m{4cm}|} 
\hline
\textbf{Initial condition $u_0$}
& \textbf{$\alpha_G$, the order of decreasing error on the Chernoff function $G(t)$, which has the 1st order of the Chernoff tangency to the operator $L=\frac{d^2}{dx^2}$}
& \textbf{$\alpha_S$, the order of decreasing error on the Chernoff function $S(t)$, which has the 2nd order of tangency by Chernoff to the operator $L=\frac{d^2}{dx^2}$}\\ 
\hline

$u_0(x)=e^{-x^6}$ & -1.8139 & -1.7517\\
\hline

$u_0(x)=e^{-x^4}$ & -1.3948 & -3.0768\\
\hline

$u_0(x)=e^{-x^2}$ & -1.0177 & -2.1028\\
\hline

$u_0(x)=\sin(x)$ & -1.0178 & -2.0415\\
\hline

$u_0(x)=\cos(x)$ & -1.0178 & -2.0415\\
\hline

$u_0(x)=e^{-x^2}\cdot\sin(x)$ & -0.9897 & -2.8194\\
\hline
\end{tabular}
\begin{flushleft}
Table 5. Orders of decreasing error for infinitely smooth initial conditions.
\end{flushleft}
\normalsize

\normalsize

\vskip5mm

Below are free terms of decreasing error that we found.

\vskip5mm

\footnotesize

\begin{tabular}{|m{6cm}|m{4cm}|m{4cm}|} 
\hline
\textbf{Initial condition $u_0$}
& \textbf{$\beta_G$, the free term of the decreasing error  on the Chernoff function $G(t)$, which has the 1st order of the Chernoff tangency to the operator $L=\frac{d^2}{dx^2}$}
& \textbf{$\beta_S$, the free term of the decreasing error  on the Chernoff function $S(t)$, which has the 2nd order of tangency by Chernoff to the operator $L=\frac{d^2}{dx^2}$}\\ 
\hline

$u_0(x)=e^{-x^6}$ & -0.7760 & -0.4493\\
\hline

$u_0(x)=e^{-x^4}$ & -2.3548 & -0.5326\\
\hline

$u_0(x)=e^{-x^2}$ & -3.3874 & -3.2967\\
\hline

$u_0(x)=\sin(x)$ & -3.6262 & -5.1659\\
\hline

$u_0(x)=\cos(x)$ & -3.6261 & -5.1659\\
\hline

$u_0(x)=e^{-x^2}\cdot\sin(x)$ & -3.7478 & -2.1899\\
\hline
\end{tabular}
\begin{flushleft}
Table 6. Free terms of decreasing error for infinitely smooth initial conditions.
\end{flushleft}
\normalsize

\normalsize

\section{Conclusion}

Let us recall the context of what we did and point out the key findings. $C_0$-semigroups have various applications in mathematics and mathematical physics \cite{EN1, EN2}. In the case where the infinitesimal generator $L$ of the semigroup $(e^{tL})_{t\geq0}$ is a differential operator with variable coefficients, often Chernoff's theorem is the most clear way to express $e^{tL}$ explicitly in terms of $L$, e.g., in terms of $a,b,c$ if $(Lf)(x)=a(x)f''(x)+b(x)f'(x)+c(x)f(x)$, $x\in\mathbb{R}$. In this approach we have  $e^{tL}=\lim_{n\to\infty}C(t/n)^n$ where $C$ is the so-called Chernoff function for the operator $L$, and $C(t/n)^n$ is the composition of $n$ copies of linear bounded operator $C(t/n)$. Expressions $C(t/n)^n$ are called the Chernoff approximations to $e^{tL}$ provided by Chernoff function $C$, and it is possible  \cite{BGS2010,R4,Rjmp, Butko2020} to express $C$ as a short formula containing $a,b,c$ as parameters. One should keep in mind that for a given operator $L$ there exists only one semigroup $e^{tL}$ but many Chernoff functions $C$, and they provide different rate of convergence $C(t/n)^n\to e^{tL}$ as $n\to\infty$.

The knowledge of $e^{tL}=\lim_{n\to\infty}C(t/n)^n$ allows to represent the bounded solution of the Cauchy problem $u_t(t,x)=a(x)u_{xx}+b(x)u_x+c(x)u$, $u(0,x)=u_0(x)$ as a formula in terms of $a,b,c,u_0$, because due to the general theory of $C_0$-semigroups \cite{EN1, EN2} we have $u(t,x)=(e^{tL}u_0)(x)$. Also, one can represent (see \cite{Rarx2023}) the bounded solution of the equation $a(x)f''(x)+b(x)f'(x)+(c(x)-\lambda)f(x)=-g(x)$ as a formula in terms of $a,b,c,g,\lambda$ because (again due to the general theory of $C_0$-semigroups \cite{EN1, EN2}) we have $f(x)=\int_0^\infty e^{-\lambda t}(e^{tL}g)(x)dt$. Chernoff's theorem has other applications, but these two are the most easy to explain in a few words. 

Our goal is to study the rate of convergence in the Chernoff theorem numerically (theoretical knowledge in this area is still very limited), and we selected the heat semigroup on the real line as a simple yet informative model case. We rewrite the initial value problem  
$$
\left\{ \begin{array}{ll}
	u_{t}(t,x)=u_{xx}(t,x) \textrm{ for }t>0, x\in\mathbb{R}^1  \\
	u(0,x)=u_0(x) \textrm{ for }x\in\mathbb{R}^1
\end{array} \right.
$$

as $[U'(t)=LU(t),\, U(0)=u_0]$ with $Lf=f''$, $u(t,\cdot)=U(t)=e^{tL}u_0$. We know that $(e^{tL})_{t\geq0}$ is a $C_0$-semigroup in the Banach space $\mathcal{F}=UC_b(\mathbb{R})$ of all uniformly continuous and bounded functions $f\colon\mathbb{R}\to\mathbb{R}$, $\|f\|=\sup_{x\in\mathbb{R}}|f(x)|$. For each $t\geq0$ we introduce linear bounded operators $G(t)$ and $S(t)$ in $\mathcal{F}$ as 
$$(G(t)f)(x)=\frac{1}{2}f(x)+\frac{1}{4}f(x+2\sqrt{t})+\frac{1}{4}f(x-2\sqrt{t}),$$ 
$$(S(t)f)(x)=\frac{2}{3}f(x)+\frac{1}{6}f(x+\sqrt{6t})+\frac{1}{6}f(x-\sqrt{6t}).$$

It is known that $G$ and $S$ are Chernoff functions for the operator $L=\partial^2$. Paul Chernoff's theorem \cite{Chernoff} implies that for each $t>0$, each $u_0\in\mathcal{F}$ sequences 
$$
d_n^G=\|e^{tL}u_0-G(t/n)^nu_0\|,\quad d_n^S=\|e^{tL}u_0-S(t/n)^nu_0\|$$ converge to 0 as $n\to\infty$, but it says \textbf{nothing} about the rate of convergence. We try to say at least something via numerical simulation. 

We considered 25 initial conditions of different smoothness: trigonometric ($u_0(x)=|\sin x|^\xi$), exponential ($u_0(x)=\exp (-|x|^\xi$), and also several infinitely smooth functions ($u_0(x)=e^{-x^{2k}}$, $u_0(x)=e^{-x^2}\sin(x)$). For all of them, we numerically calculated
$$\ln d_n^G\approx \alpha_G \ln n + \beta_G,\quad \ln d_n^S\approx \alpha_S \ln n + \beta_S.$$

We call $\alpha$ the order of decreasing error, $\beta$ the free term of decreasing error. Function $G$ has the first order of Chernoff tangency to $\partial^2$, and function $S$ has the second order of Chernoff tangency to $\partial^2$. So, one can expect $\alpha_G\approx -1$, $\alpha_S\approx -2$, but the numerical experiment shows that the reality is not that simple.

\textbf{These are not surprising findings} that were predicted by the existing theory: 

\begin{itemize}
	
	\item We see that $\alpha_S<\alpha_G$ for all initial conditions $u_0$ that we have tested. This means that a higher order of the Chernoff tangency of the Chernoff function to the generator of the semigroup provides a better order of decreasing error. Indeed, blue points are below green points on ($\xi, \alpha$) plane for trigonometric initial conditions (fig. 9) and for exponential initial conditions (fig. 11).
	
	\item The difference between the errors $d_n^G$ and $d_n^S$ is obvious if $\xi\in[2,4.5]$, i.e. when $u_0\in D(L)$.
	
	\item If $u_0\notin D(L)$, i.e. $\xi\in(0,2)$, then this difference is small. 
	
\end{itemize} 

\textbf{These are effects that we see numerically but do not know how to explain theoretically:}

\begin{enumerate}
	
	\item The linear approximation of error $\ln d_n$ vs. $\ln n$ is \textit{not equally perfect}, and it is often different for different initial conditions. But for $u_0(x)=\sin(x)$, $u_0(x)=\cos(x)$ the error is \textit{identical}. For smooth $u_0(x)=\sin(x)$ we see that points fit the line perfectly (fig. 1.2), as well as for non-smooth $u_0(x)=|\sin(x)|^{1/4}$, see fig. 3.2. But for smooth $u_0(x)=e^{-x^4}$ we see that points do not fit the line well (fig. 14.2). It seems that a higher smoothness provides a worse line fit, but for $u_0(x)=\sin(x)$ the fit is very good.
	
	\item The error \textit{oscillates}. In Figure 18 below in the right picture, one can see that the points, both green and blue, form ``waves'' that rise and fall, oscillating near the linear trend. This effect is more or less seen for all initial conditions that we tested. 
	\begin{center}
		\includegraphics[scale=0.5]{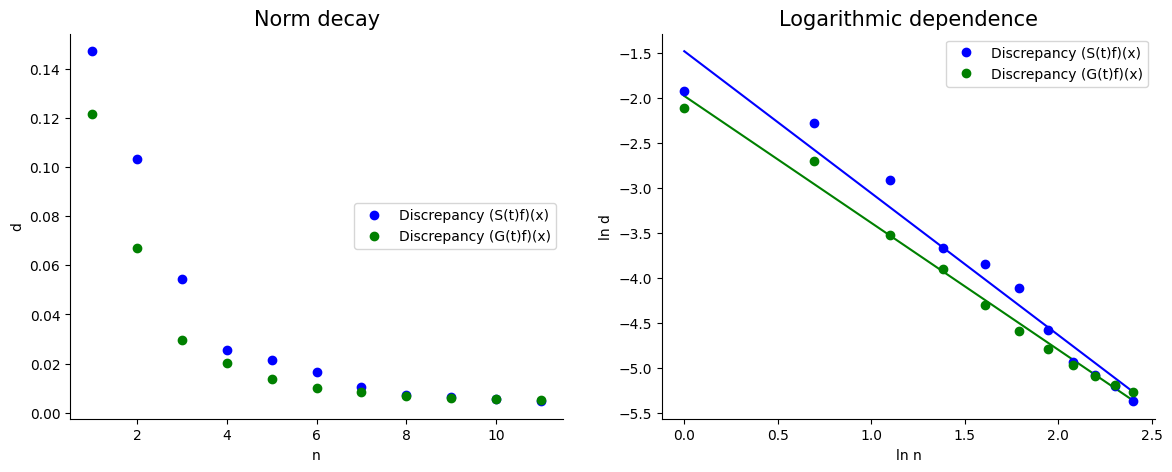}\\
		{fig. 18, $1\leq n\leq 11$, $u_0(x) = e^{-|x|^{9/2}}$, $t=\frac{1}{2}$} 
	\end{center}
	
	The blue color corresponds to the Chernoff function $S(t)$, and the green color corresponds to the Chernoff function $G(t)$. The coefficients of the straight lines were calculated by running linear regression (least squares method), using points corresponding to $4\leq n\leq 11$.

	\item The order of error behaves \textit{in the opposite} to the free term of the error. In other words: $\alpha_S<\alpha_G$ for all initial conditions, but $\beta_S>\beta_G$ (green points are below blue points) for almost all initial conditions, exeptions are $u_0(x)=\sin(x)$ (large difference in $\beta$) and $u_0(x)=e^{-|x|^{5/2}}$ (small difference in $\beta$). Cf. fig. 9 vs. fig. 10, fig. 11 vs. fig. 12.
	
	\item For $\xi\in(0,1)$ dependence of $\alpha_S$ and $\alpha_G$ on $\xi$ seems to be \textit{linear} for both trigonometric and exponential $u_0$.
	
	\textit{Smoothness matters, but it is clear that something else is also important. Indeed:}
	
	\item Infinitely smooth initial conditions can provide a different behavior of the error:
	
	$u_0(x)=\sin(x):\ \alpha_G\approx-1.0, \alpha_S\approx-2.0$
	
	$u_0(x)=e^{-x^2}:\ \alpha_G\approx-1.0,\alpha_S\approx-2.1$ 
	
	$u_0(x)=e^{-x^4}:\ \alpha_G\approx-1.4$, $\alpha_S\approx-3.1$
	
	$u_0(x)=e^{-x^6}:\ \alpha_G\approx-1.8,\alpha_S\approx-1.8$ 
	\item For trigonometric initial conditions $u_0$ with $\xi\in[1,4.5]$ the dependence of $\alpha_S$ on $\xi$ seems to be linear, \textit{but} it seems that there is no dependence of $\alpha_G$ on $\xi$. Meanwhile, for exponential $u_0$ with $\xi\in[1,4.5]$ it seems that there is no dependence of \textit{both} $\alpha_G$ and $\alpha_S$ on $\xi$.
	
\end{enumerate}

\textbf{Acknowledgements.} Research results (except item 3.4.3) was supported by the Russian Science Foundation (project 23-71-30008). Item 3.4.3 was obtained in the IITP (Dobrushin's Math. Lab.) within the state assignment of Ministry of Science and Higher Education of the Russian Federation for IITP RAS. The authors thank Prof. O.E.Galkin for discussions of the results presented in the paper. AI was not used in this research and in the preparation of this manuscript.

\section*{Appendix A: full list of pictures}
\addcontentsline{toc}{section}{Appendix A: full list of pictures}
\subsection{$u_0(x)=\sin(x)$}

n=1
\begin{center}
	\includegraphics[scale=0.5]{sin1.png}\\
\end{center}

\flushleft n=2
\begin{center}
	\includegraphics[scale=0.5]{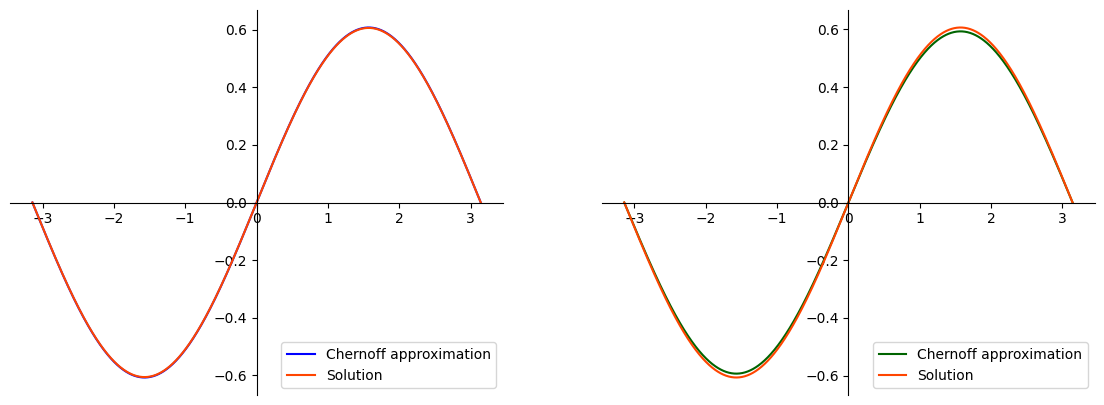}\\
\end{center}
n=3
\begin{center}
	\includegraphics[scale=0.5]{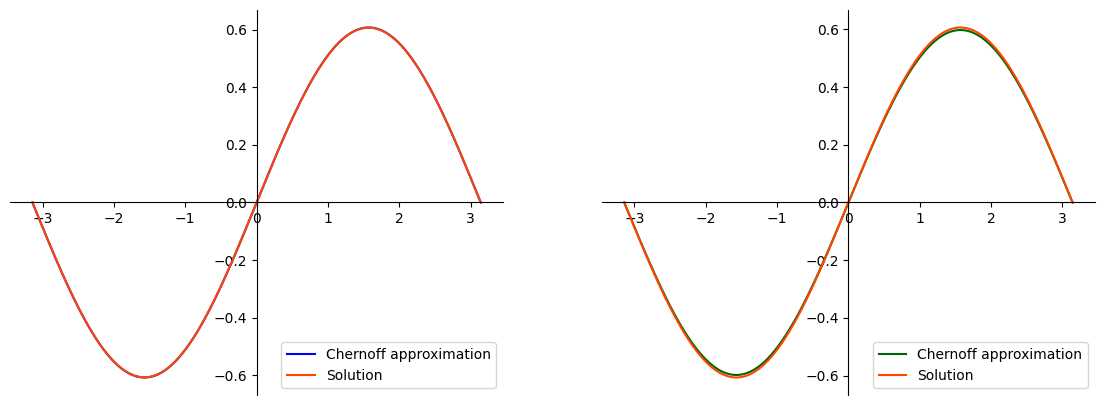}\\
\end{center}

\newpage
n=4

\begin{center}
	\includegraphics[scale=0.5]{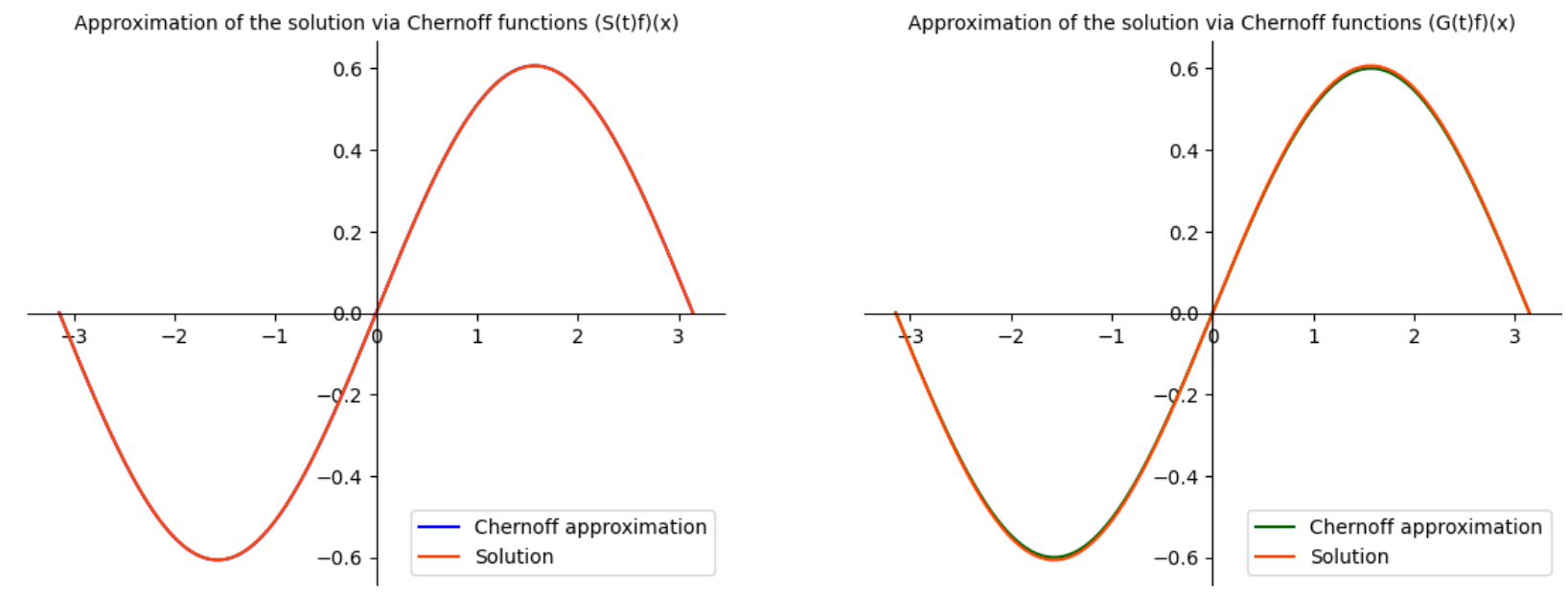}\\
\end{center} 
n=5
\begin{center}
	\includegraphics[scale=0.5]{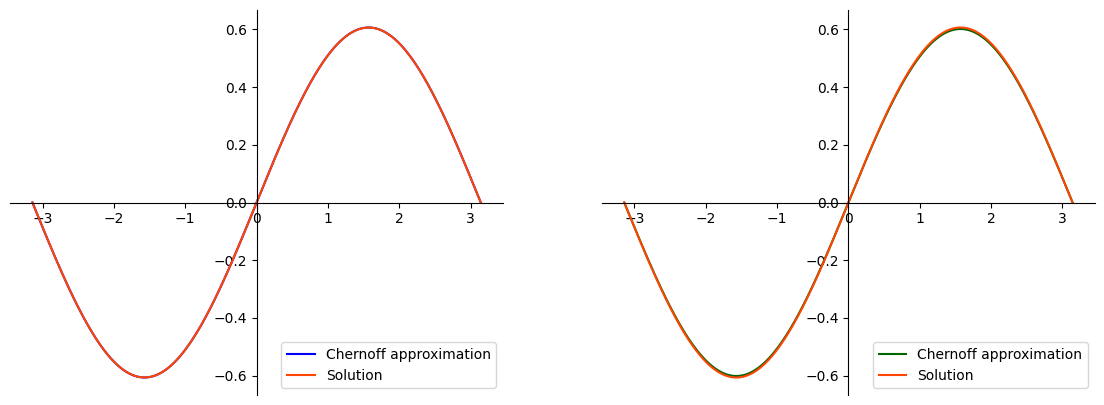}\\
\end{center}
n=6
\begin{center}
	\includegraphics[scale=0.5]{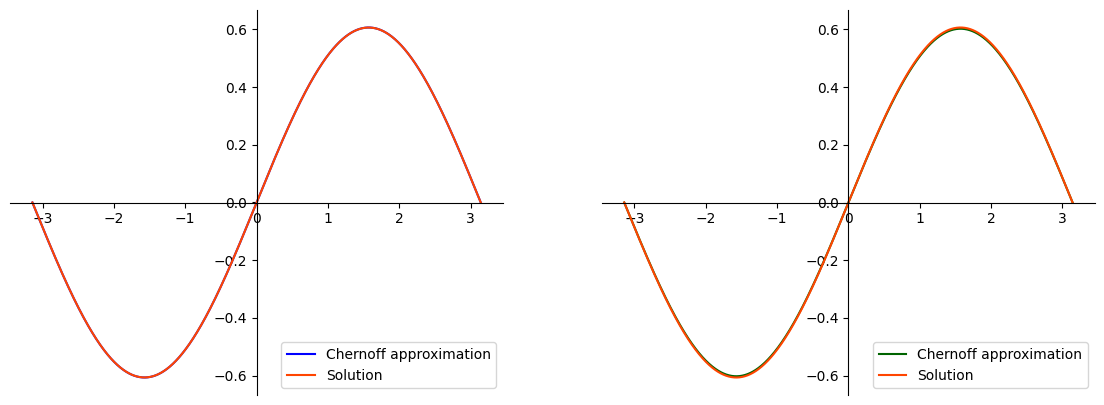}\\
\end{center}
n=7
\begin{center}
	\includegraphics[scale=0.5]{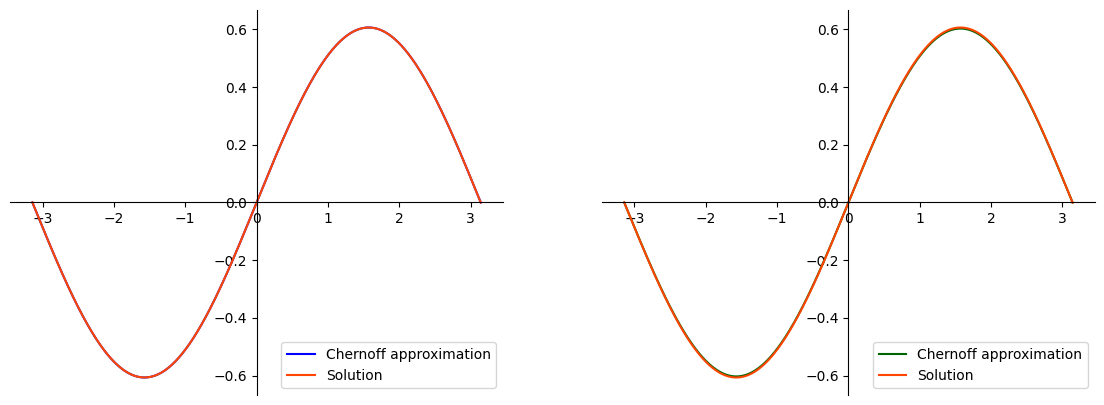}\\
\end{center}
n=8

\begin{center}
	\includegraphics[scale=0.5]{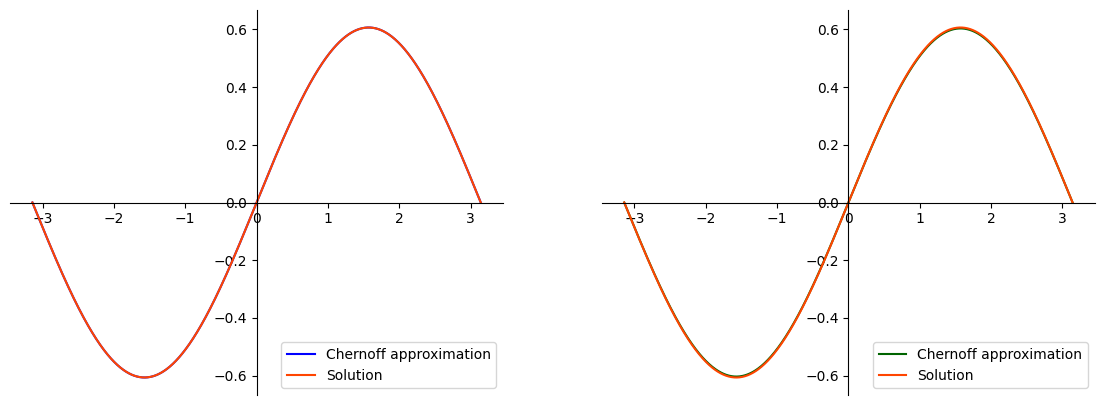}\\
\end{center}
n=9
\begin{center}
	\includegraphics[scale=0.5]{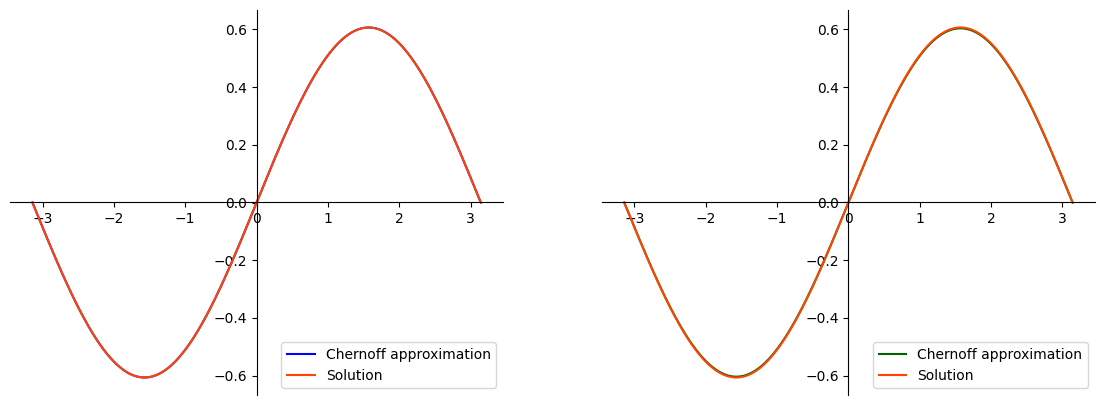}\\
\end{center}
n=10
\begin{center}
	\includegraphics[scale=0.5]{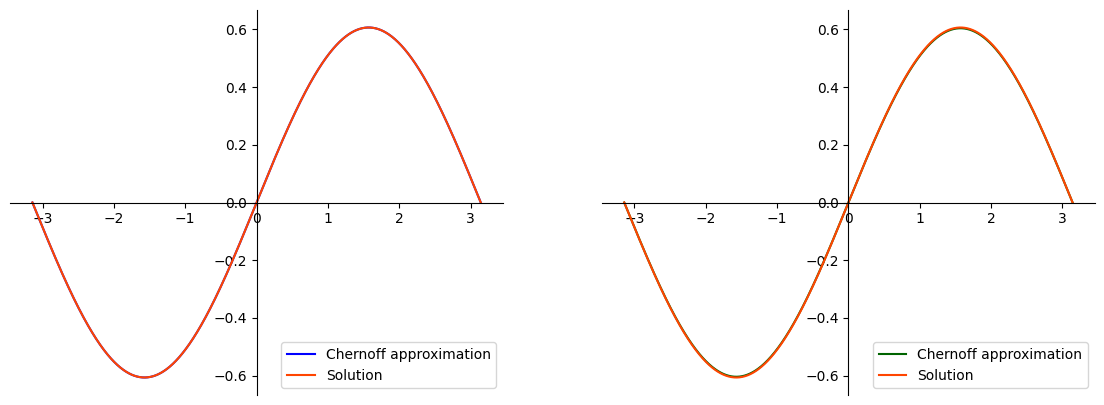}\\
\end{center}

\subsection{$u_0(x)=\cos(x)$}

n=1
\begin{center}
	\includegraphics[scale=0.5]{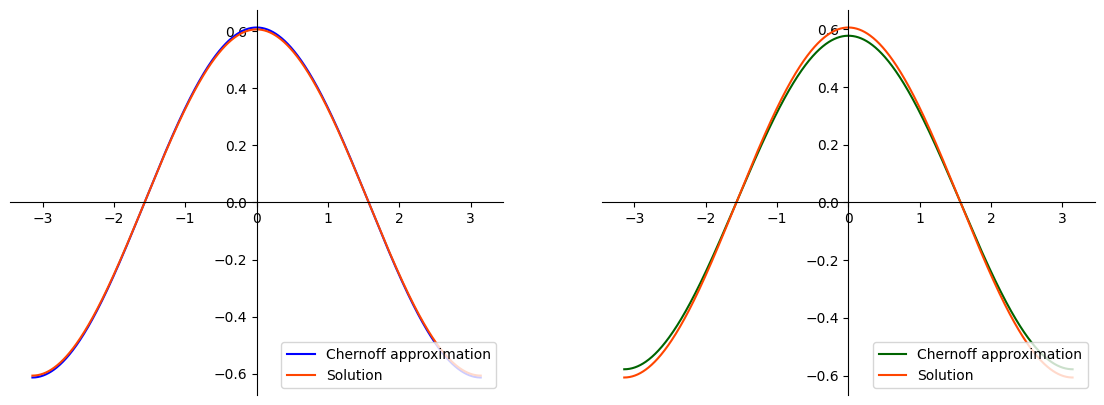}\\
\end{center}

\flushleft n=2
\begin{center}
	\includegraphics[scale=0.5]{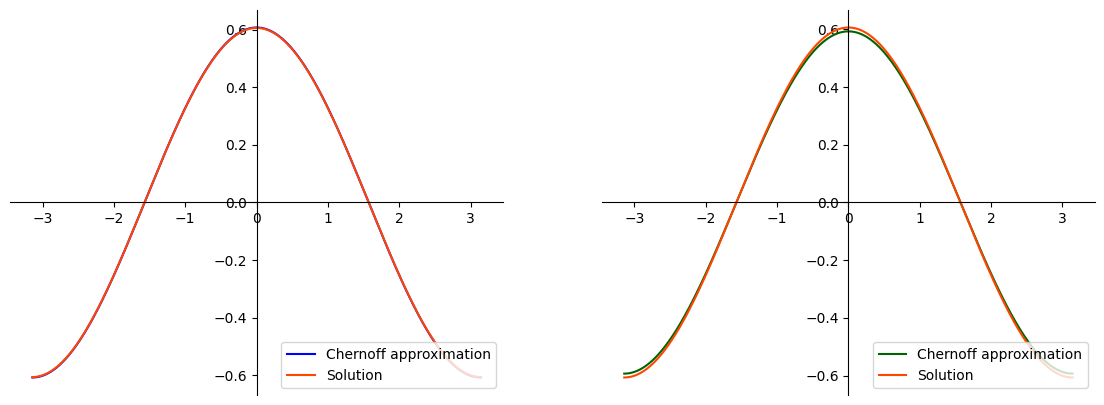}\\
\end{center}
n=3
\begin{center}
	\includegraphics[scale=0.5]{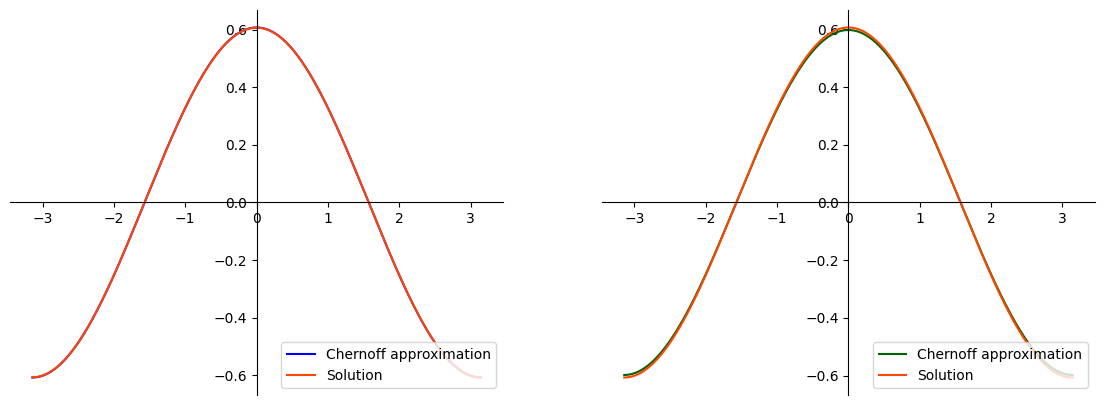}\\
\end{center}

\flushleft n=4

\begin{center}
	\includegraphics[scale=0.5]{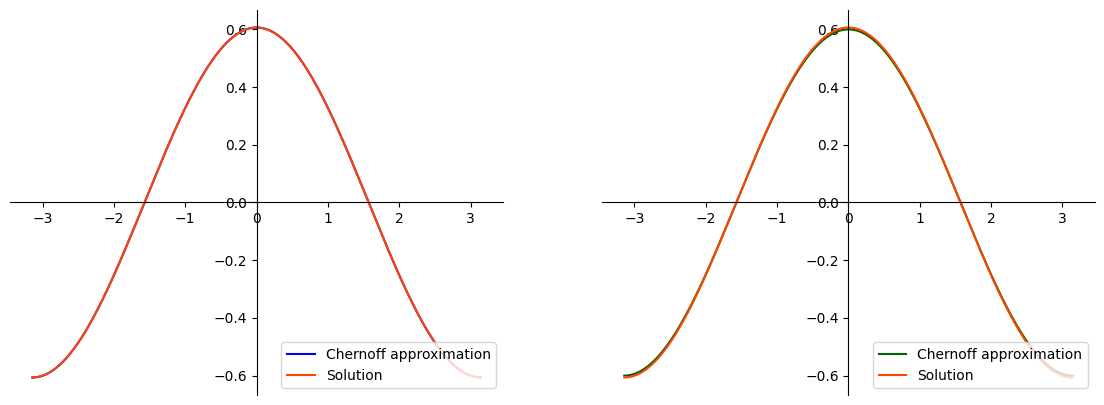}\\
\end{center} 
n=5
\begin{center}
	\includegraphics[scale=0.5]{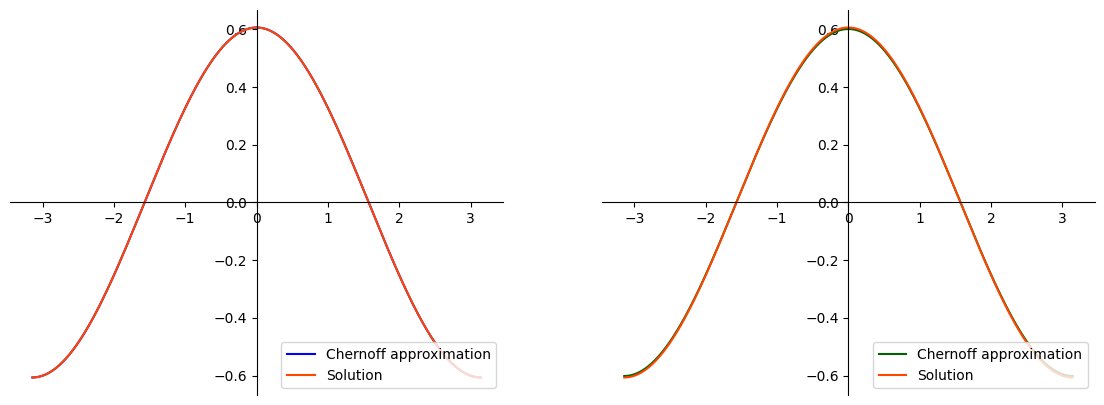}\\
\end{center}
n=6
\begin{center}
	\includegraphics[scale=0.5]{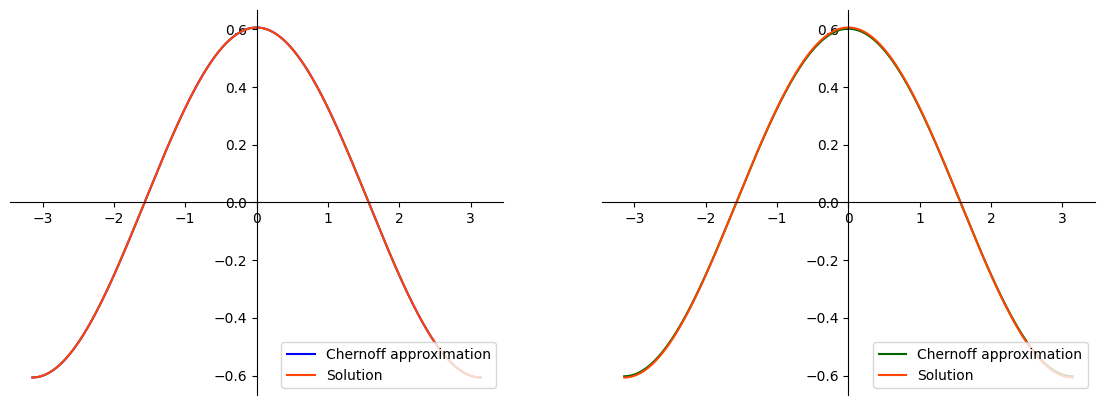}\\
\end{center}
n=7
\begin{center}
	\includegraphics[scale=0.5]{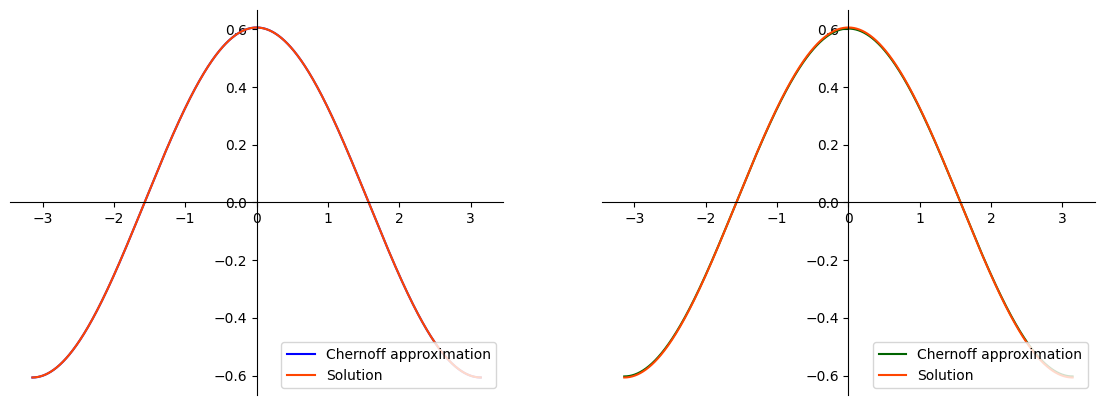}\\
\end{center}
n=8

\begin{center}
	\includegraphics[scale=0.5]{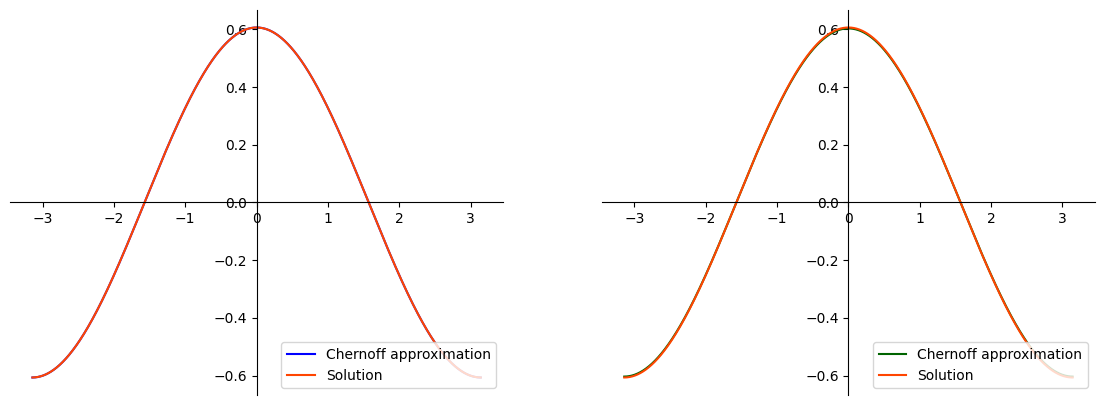}\\
\end{center}
n=9
\begin{center}
	\includegraphics[scale=0.5]{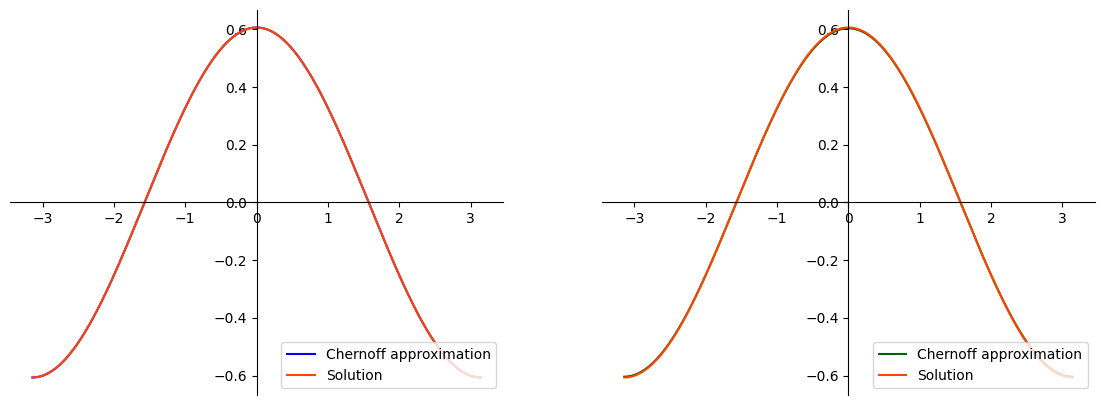}\\
\end{center}
n=10
\begin{center}
	\includegraphics[scale=0.5]{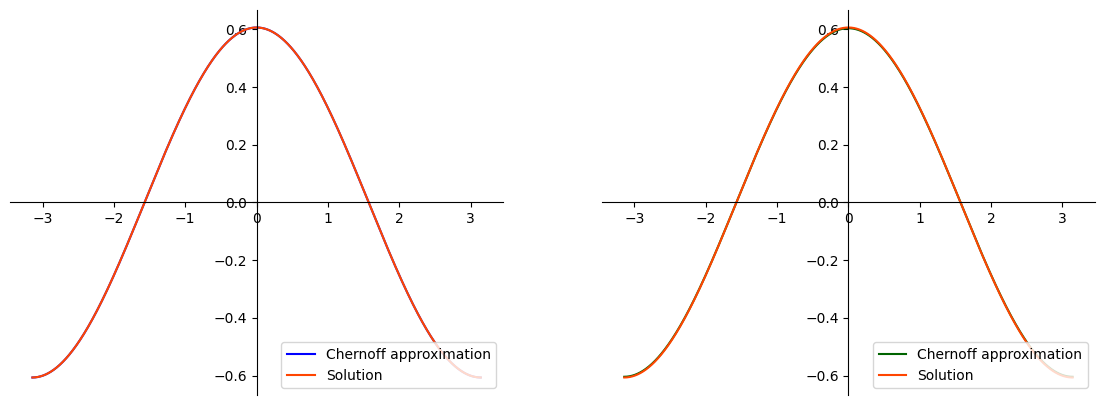}\\
\end{center}

\newpage
\subsection{$u_0(x)=|\sin(x)|$}
n=1
\begin{center}
	\includegraphics[scale=0.5]{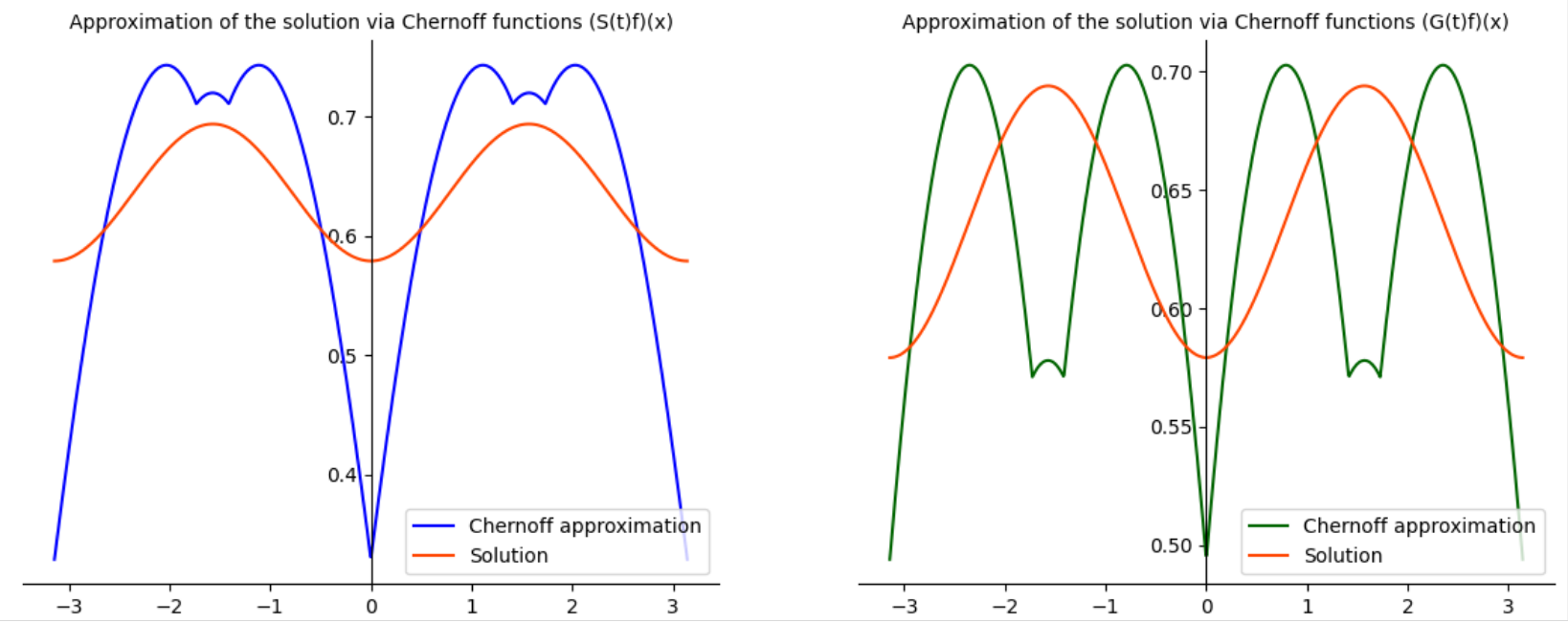}\\
\end{center}
n=2
\begin{center}
	\includegraphics[scale=0.5]{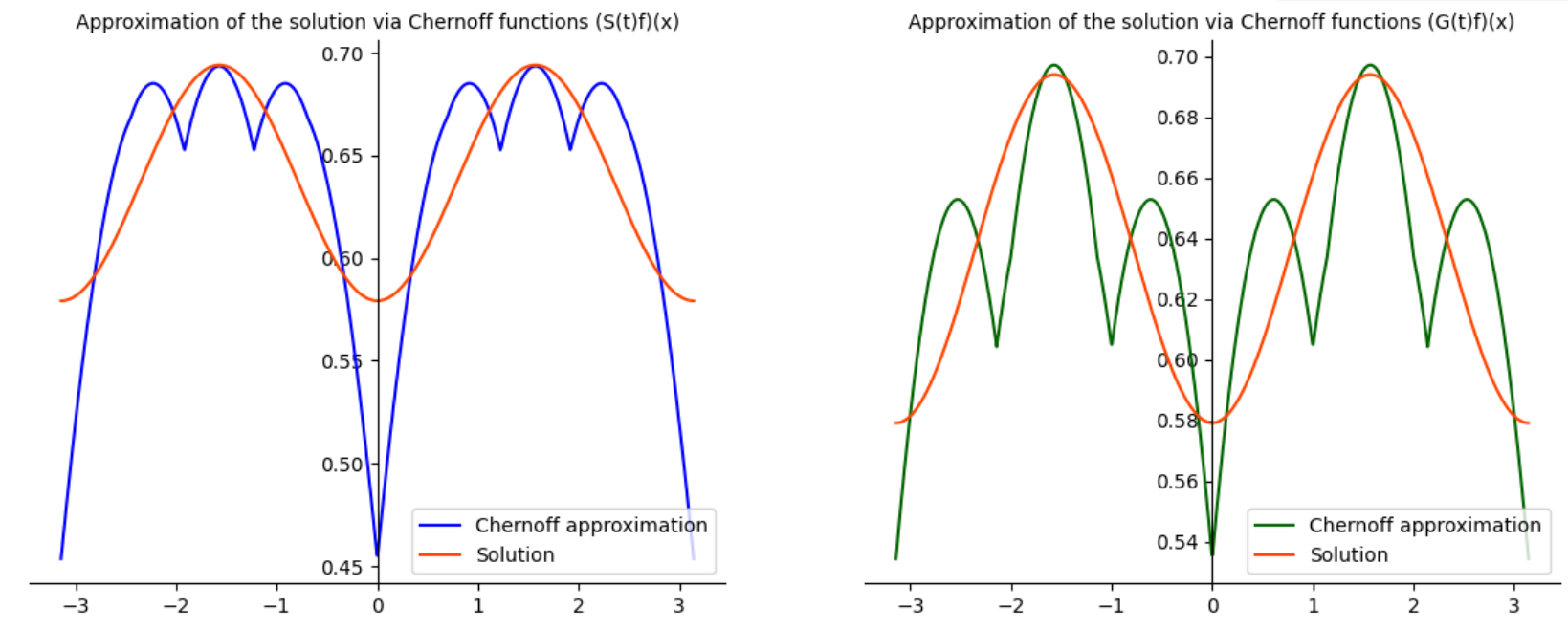}\\
\end{center}
n=3
\begin{center}
	\includegraphics[scale=0.5]{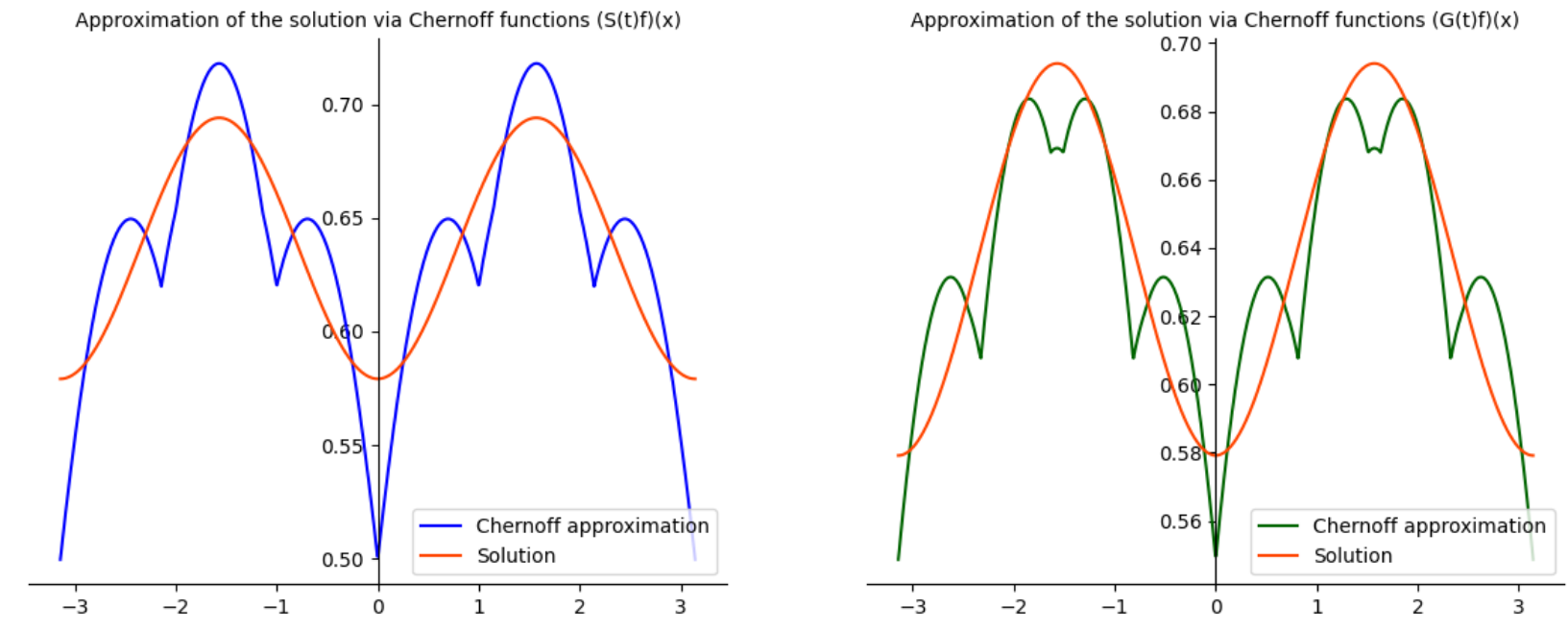}\\
\end{center}

\newpage
n=4
\begin{center}
	\includegraphics[scale=0.5]{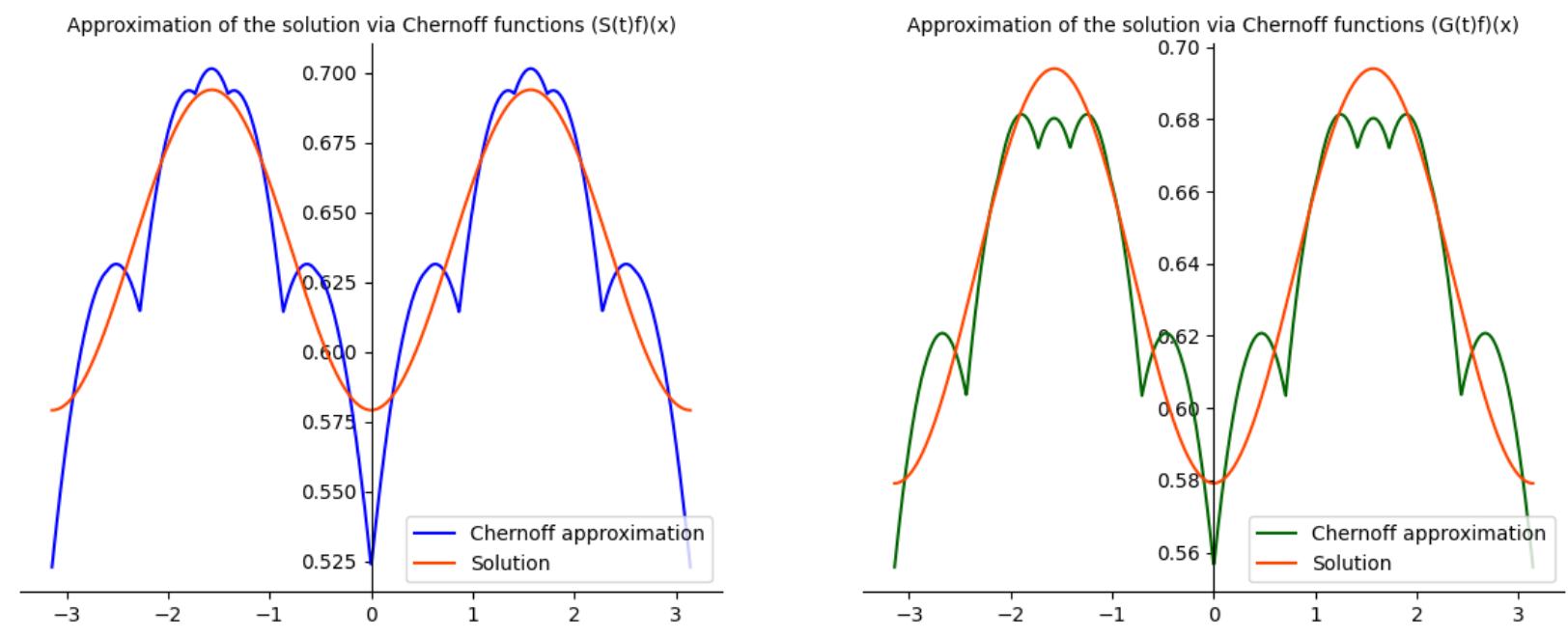}\\
\end{center}
n=5
\begin{center}
	\includegraphics[scale=0.5]{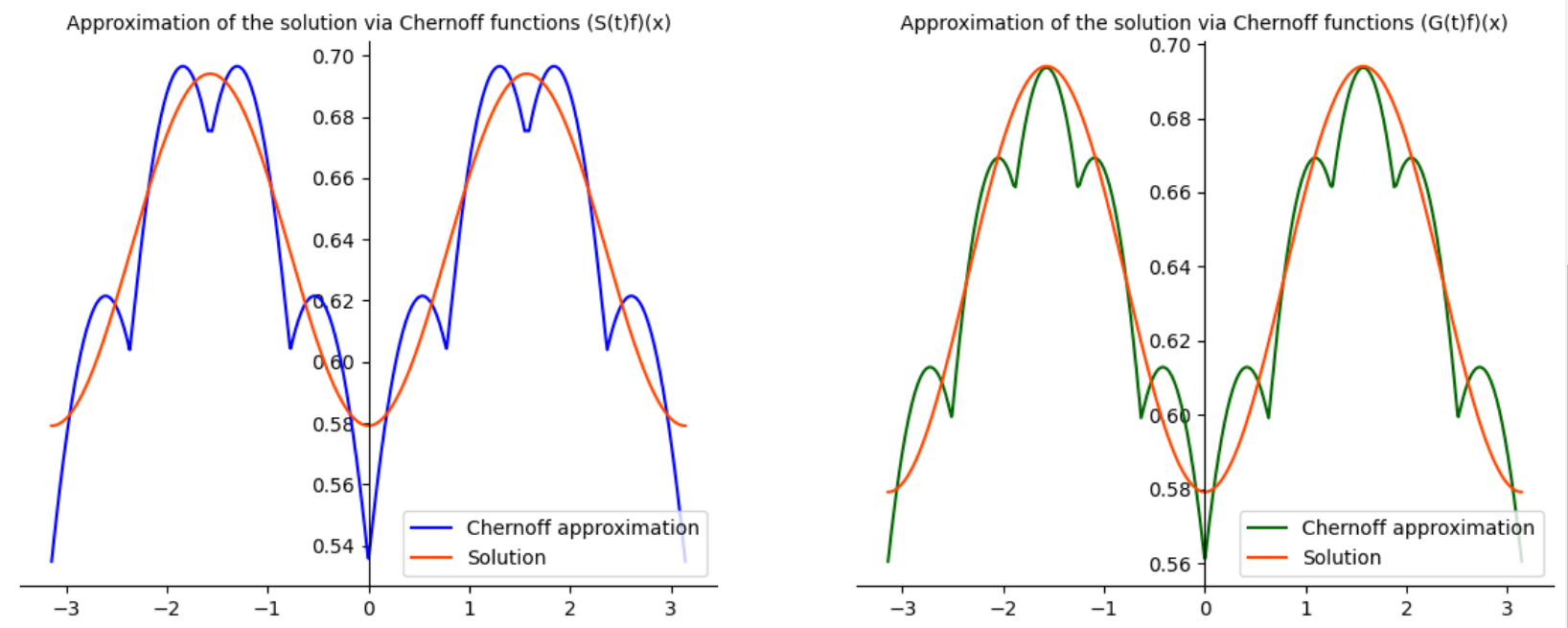}\\
\end{center}
n=6
\begin{center}
	\includegraphics[scale=0.5]{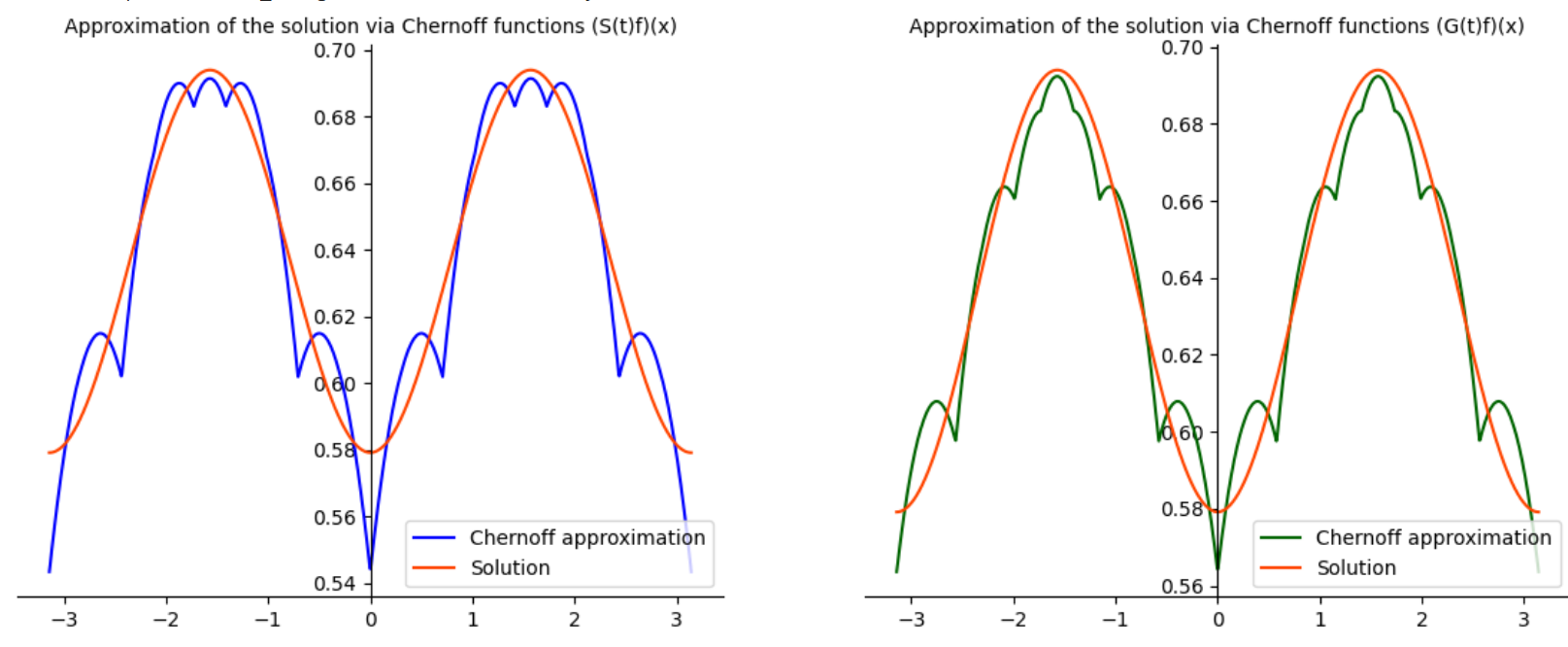}\\
\end{center}

\newpage
n=7
\begin{center}
	\includegraphics[scale=0.5]{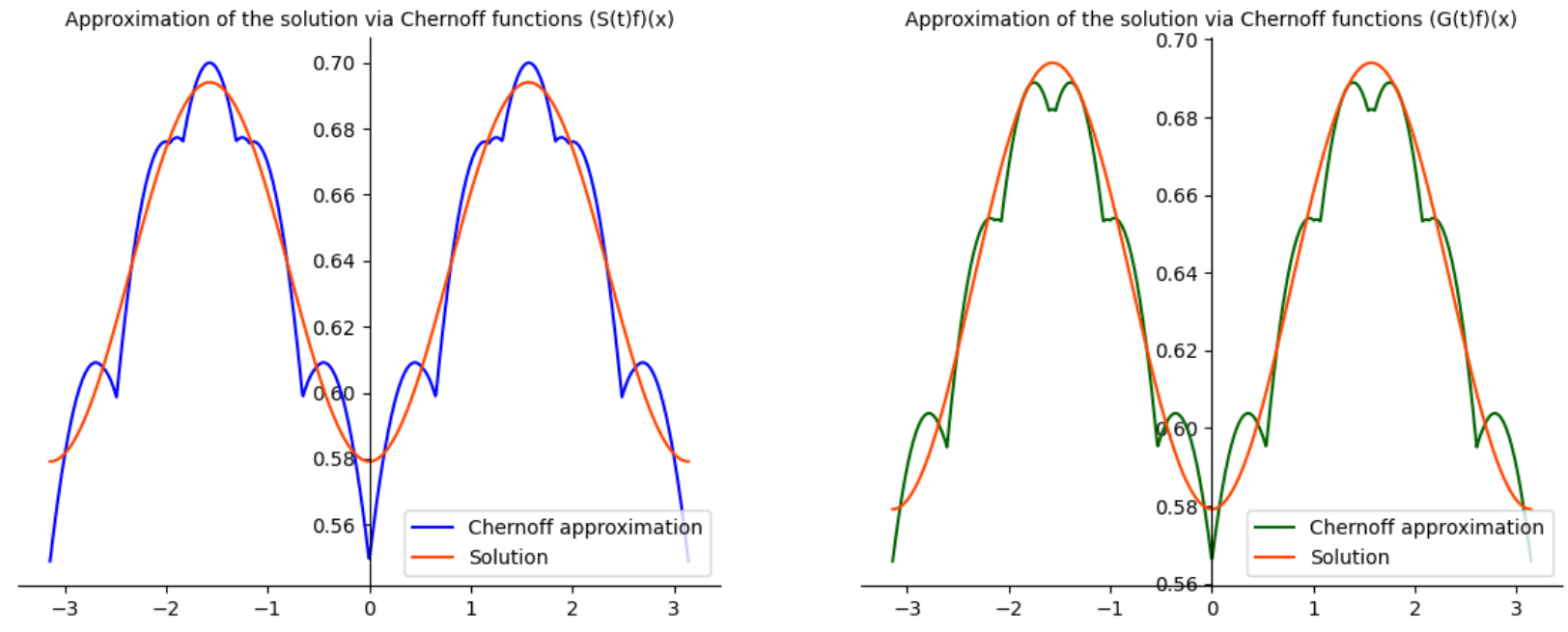}\\
\end{center}
n=8
\begin{center}
	\includegraphics[scale=0.5]{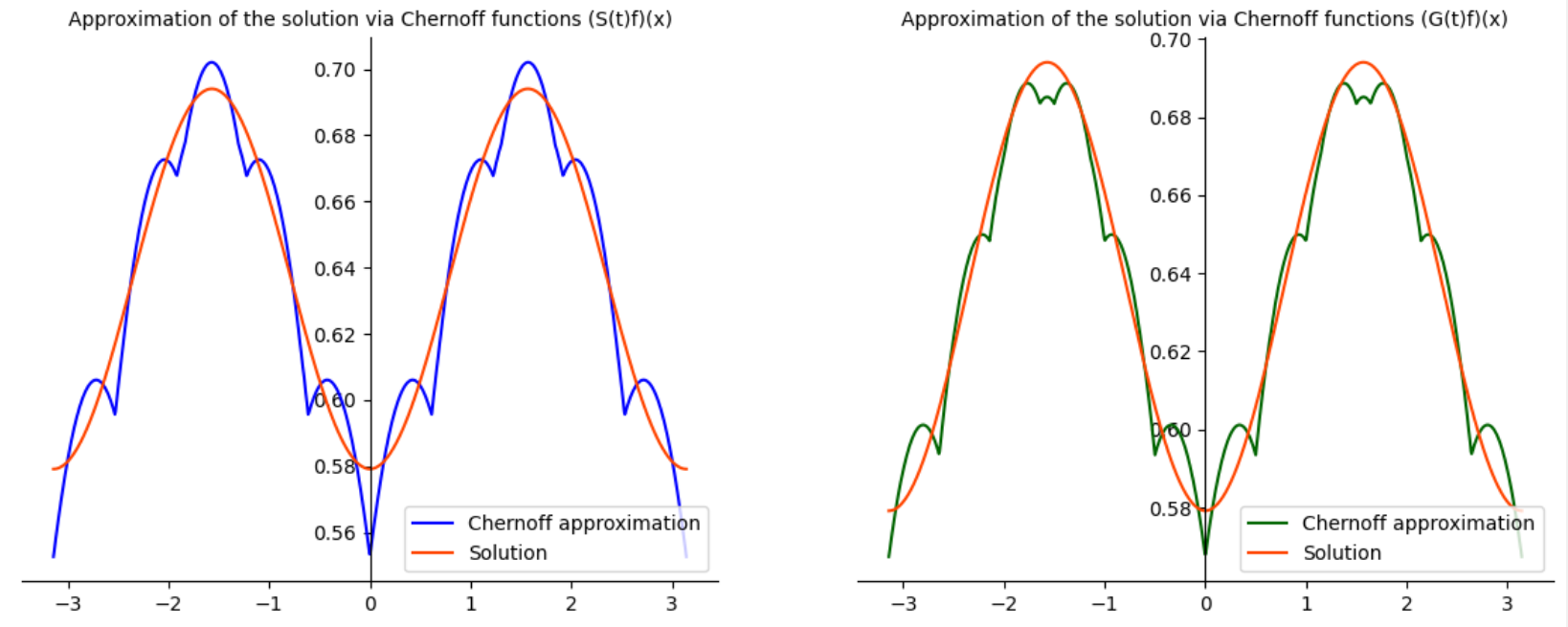}\\
\end{center}
n=9
\begin{center}
	\includegraphics[scale=0.5]{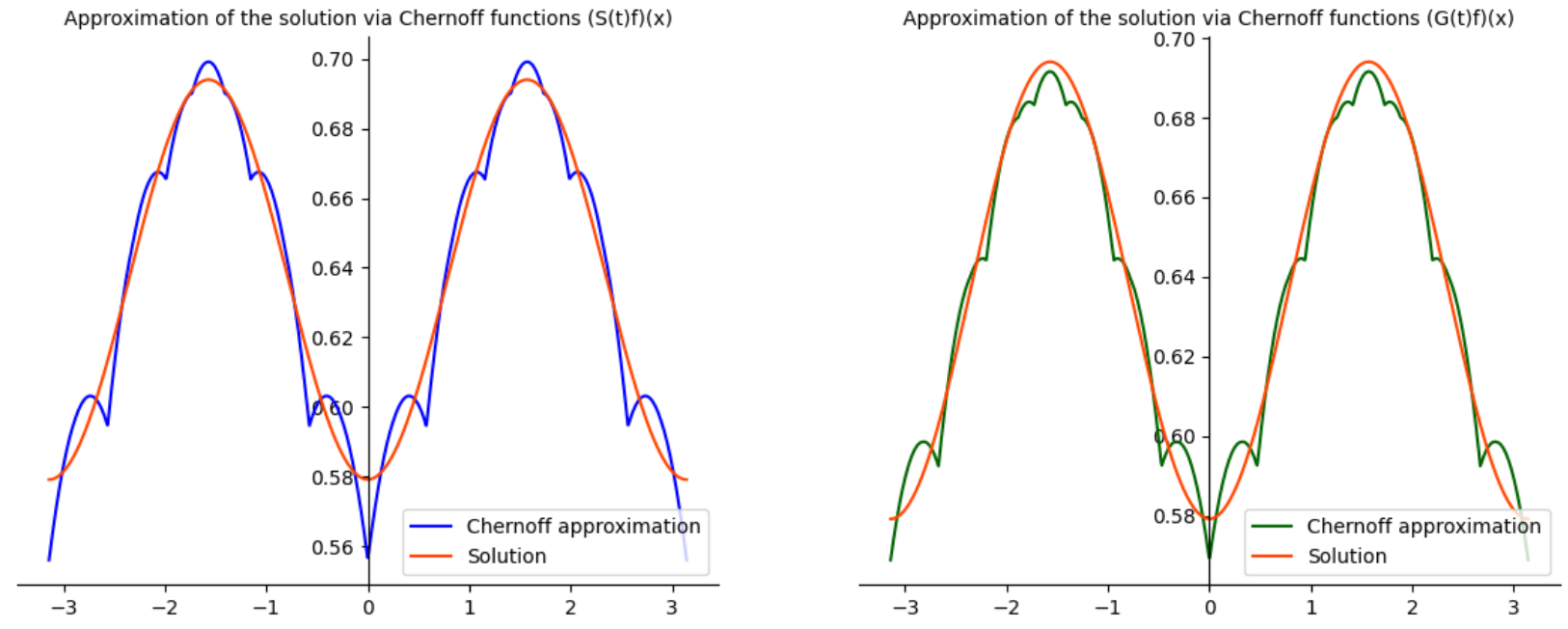}\\
\end{center}

\newpage
n=10
\begin{center}
	\includegraphics[scale=0.5]{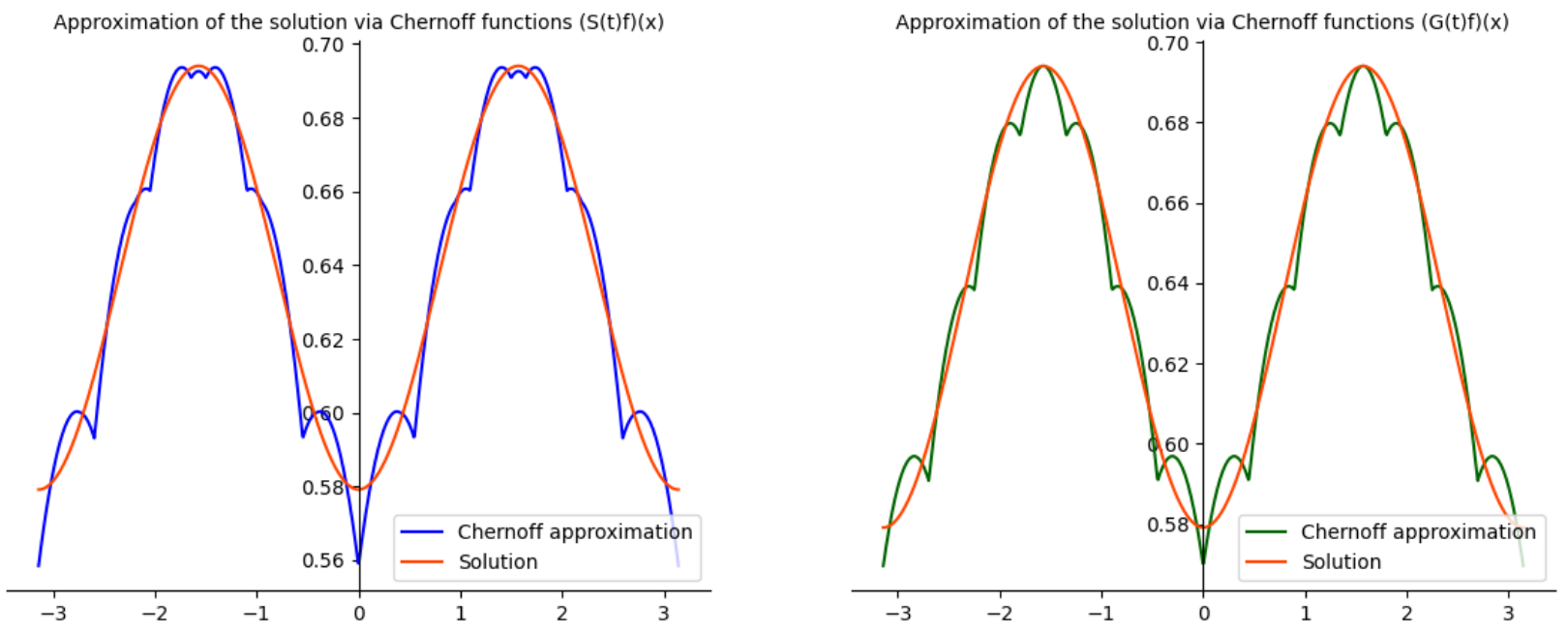}\\
\end{center}
\subsection{$u_0(x)=  \sqrt{|\sin x|}$}
n=1
\begin{center}
	\includegraphics[scale=0.5]{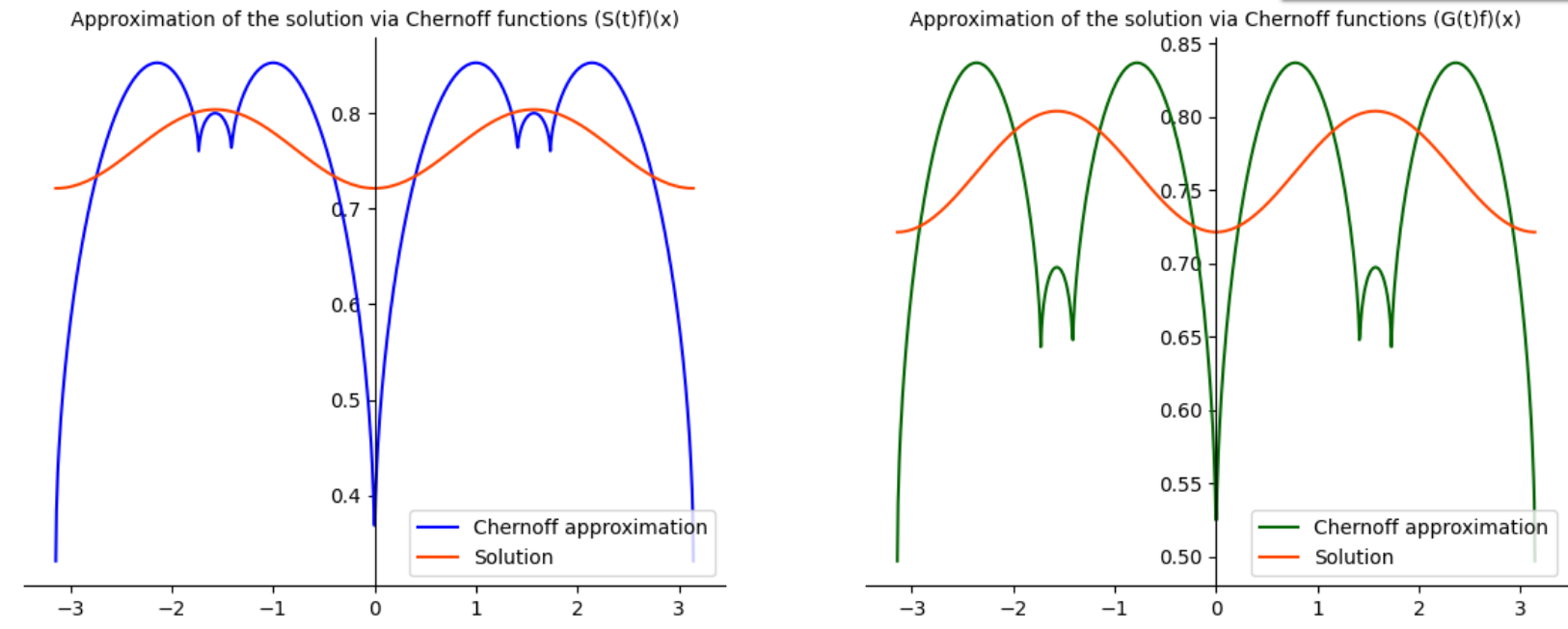}\\
\end{center}
n=2
\begin{center}
	\includegraphics[scale=0.5]{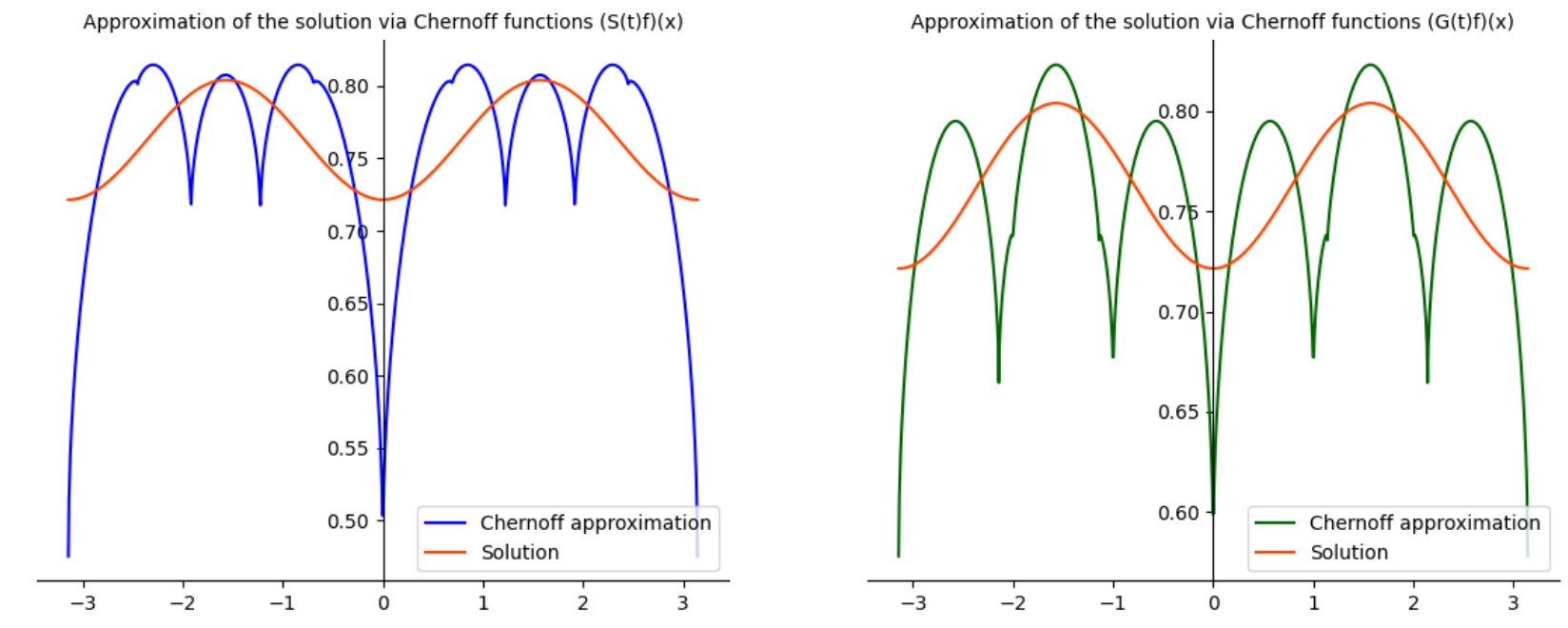}\\
\end{center}

\newpage
n=3
\begin{center}
	\includegraphics[scale=0.5]{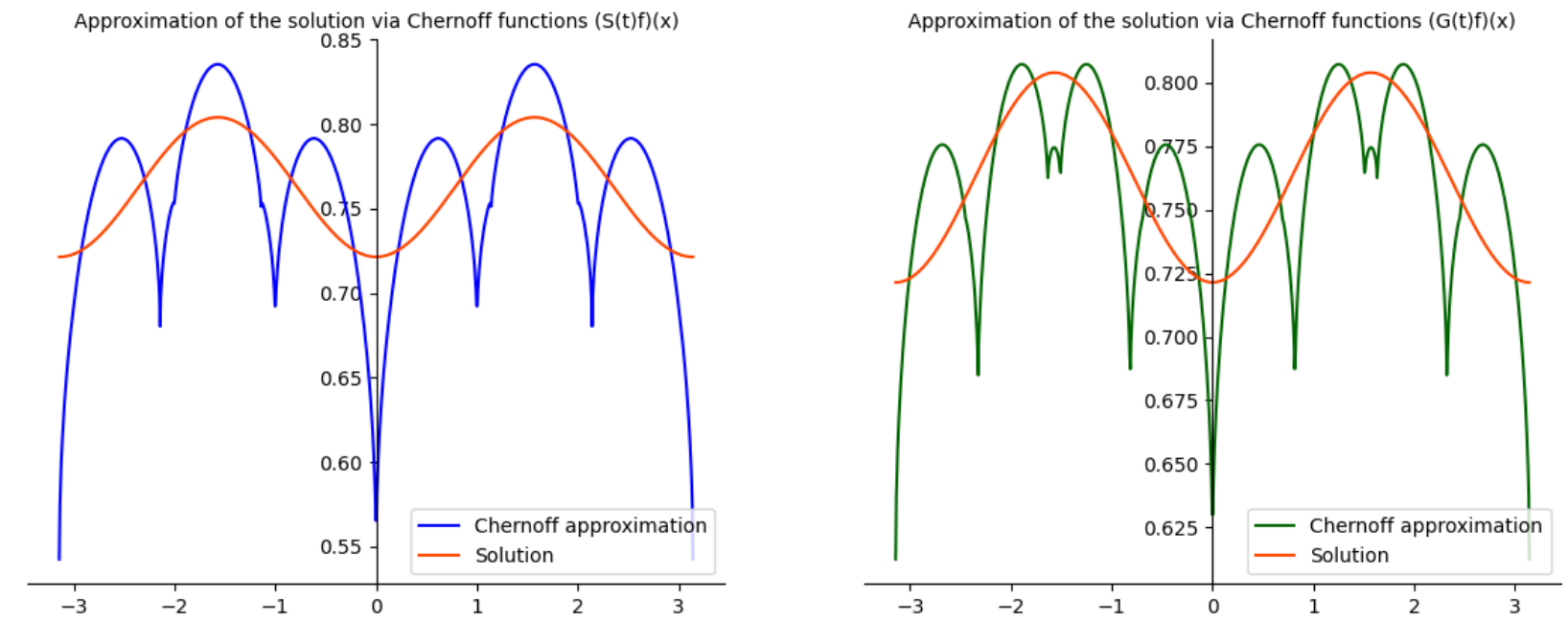}\\
\end{center}
n=4
\begin{center}
	\includegraphics[scale=0.5]{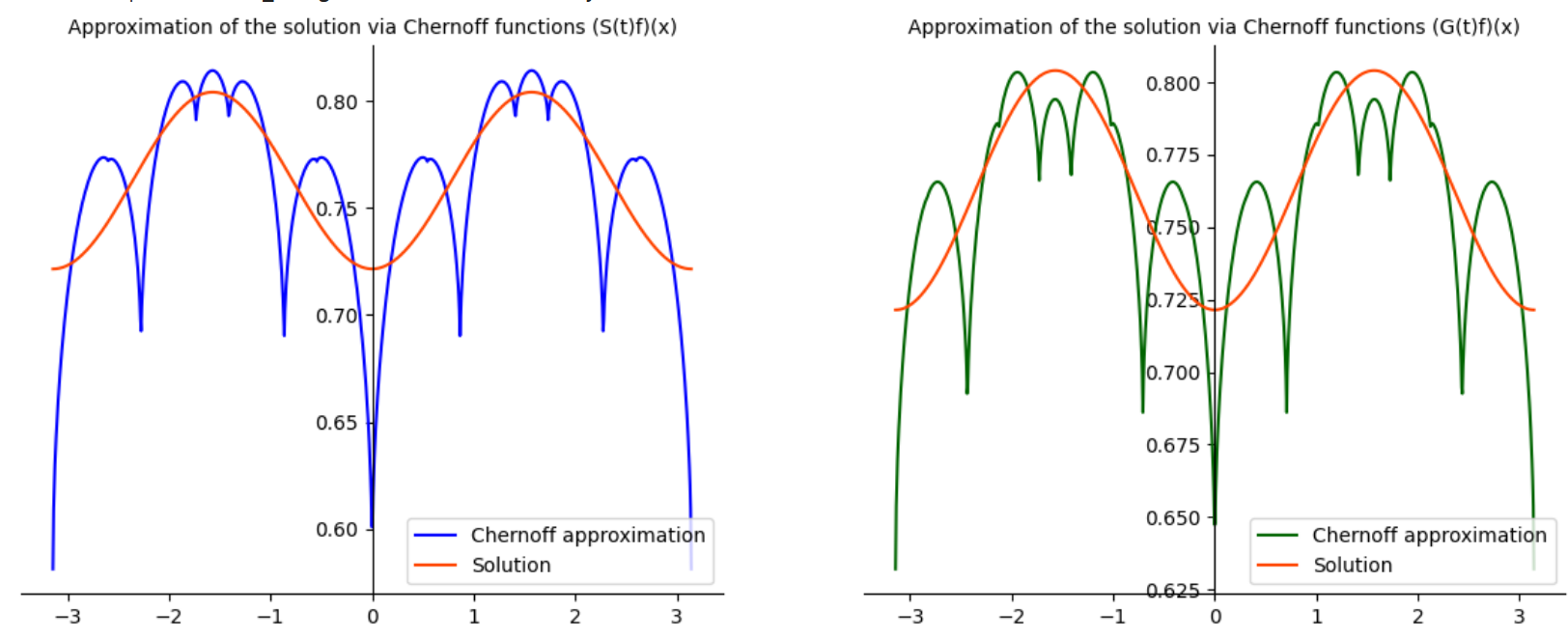}\\
\end{center}
n=5
\begin{center}
	\includegraphics[scale=0.5]{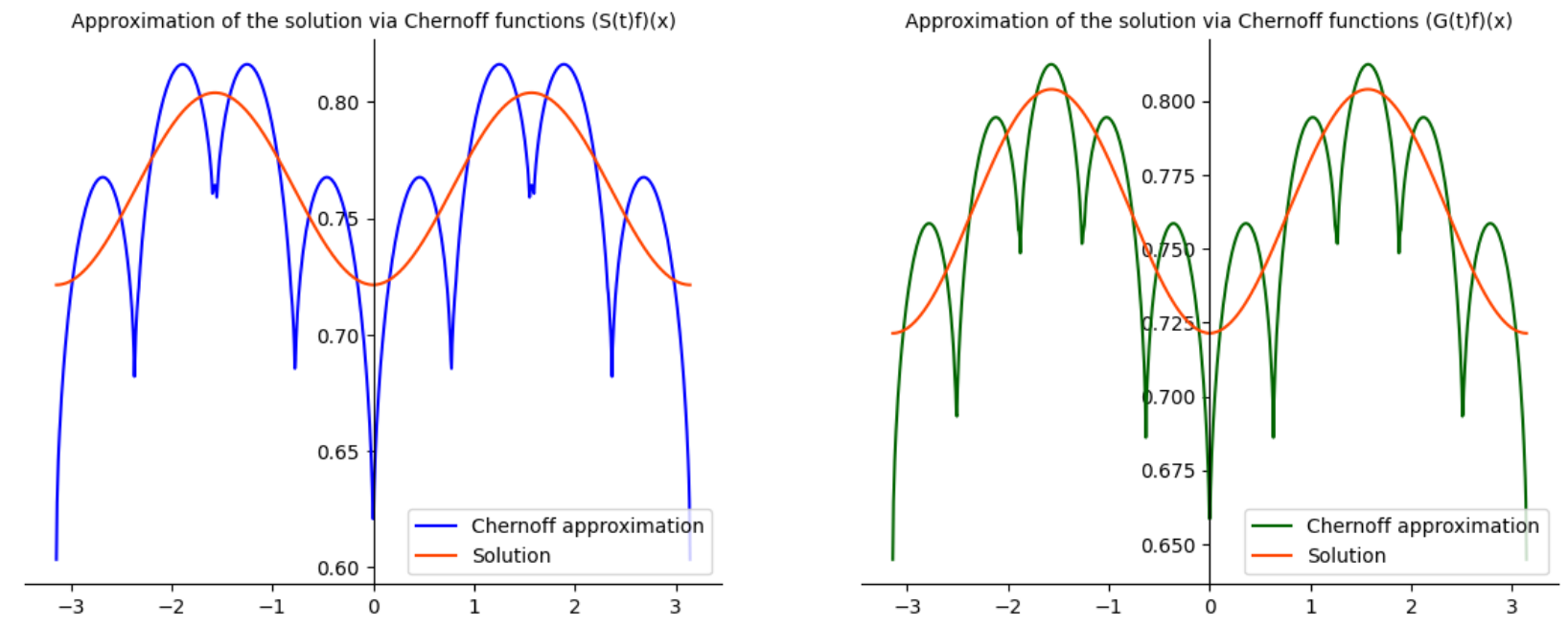}\\
\end{center}
\newpage
n=6
\begin{center}
	\includegraphics[scale=0.5]{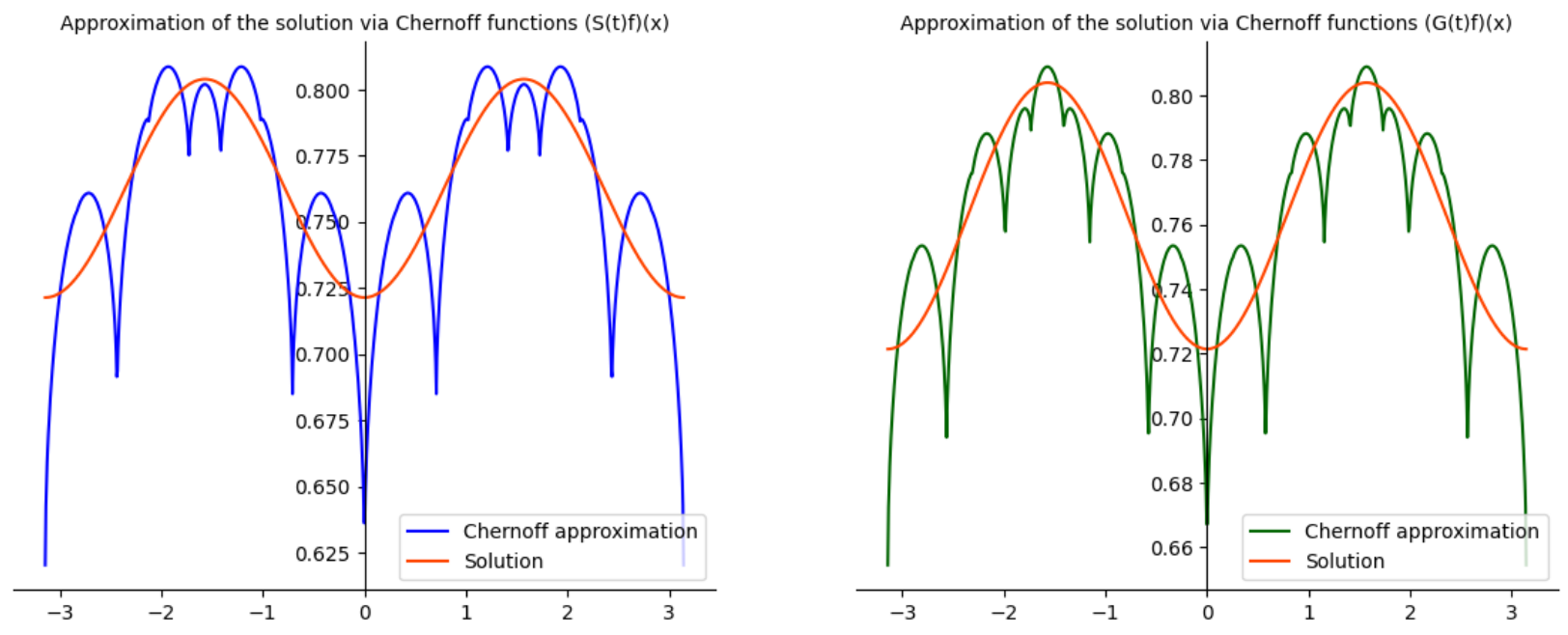}\\
\end{center}

n=7
\begin{center}
	\includegraphics[scale=0.5]{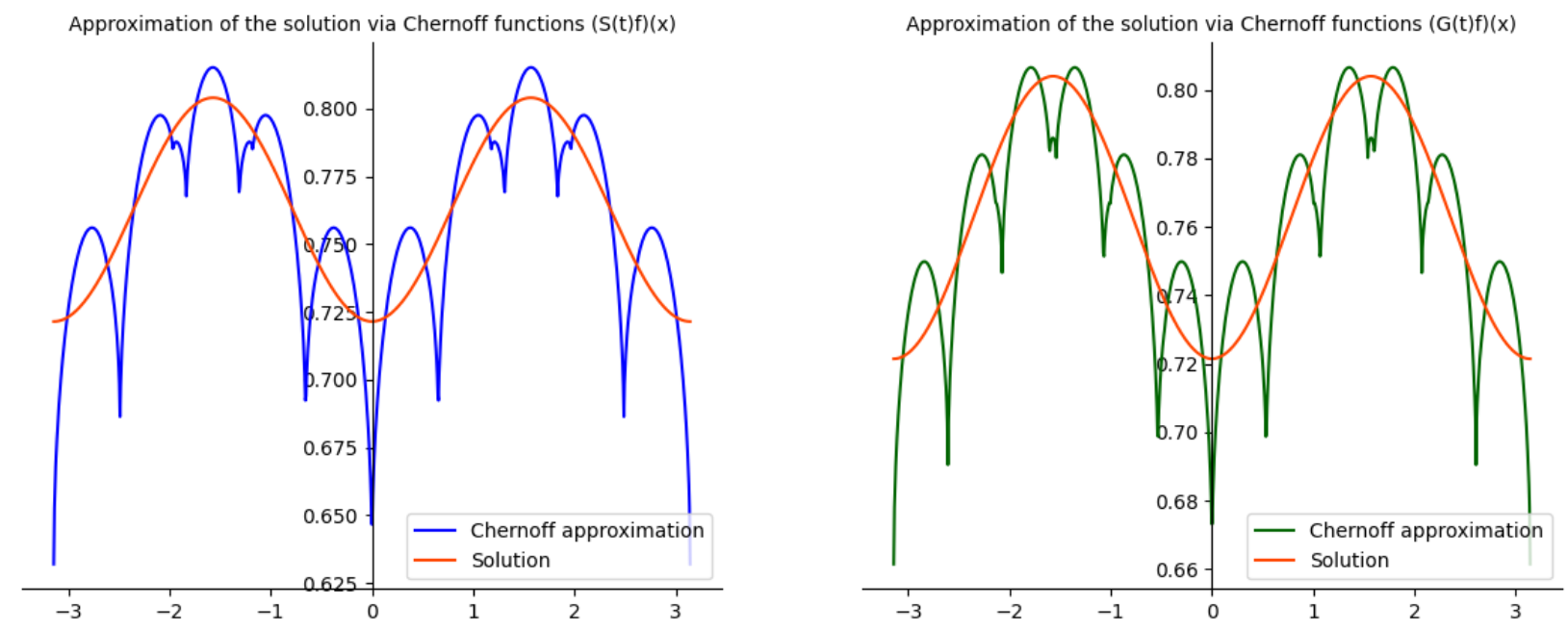}\\
\end{center}
n=8
\begin{center}
	\includegraphics[scale=0.5]{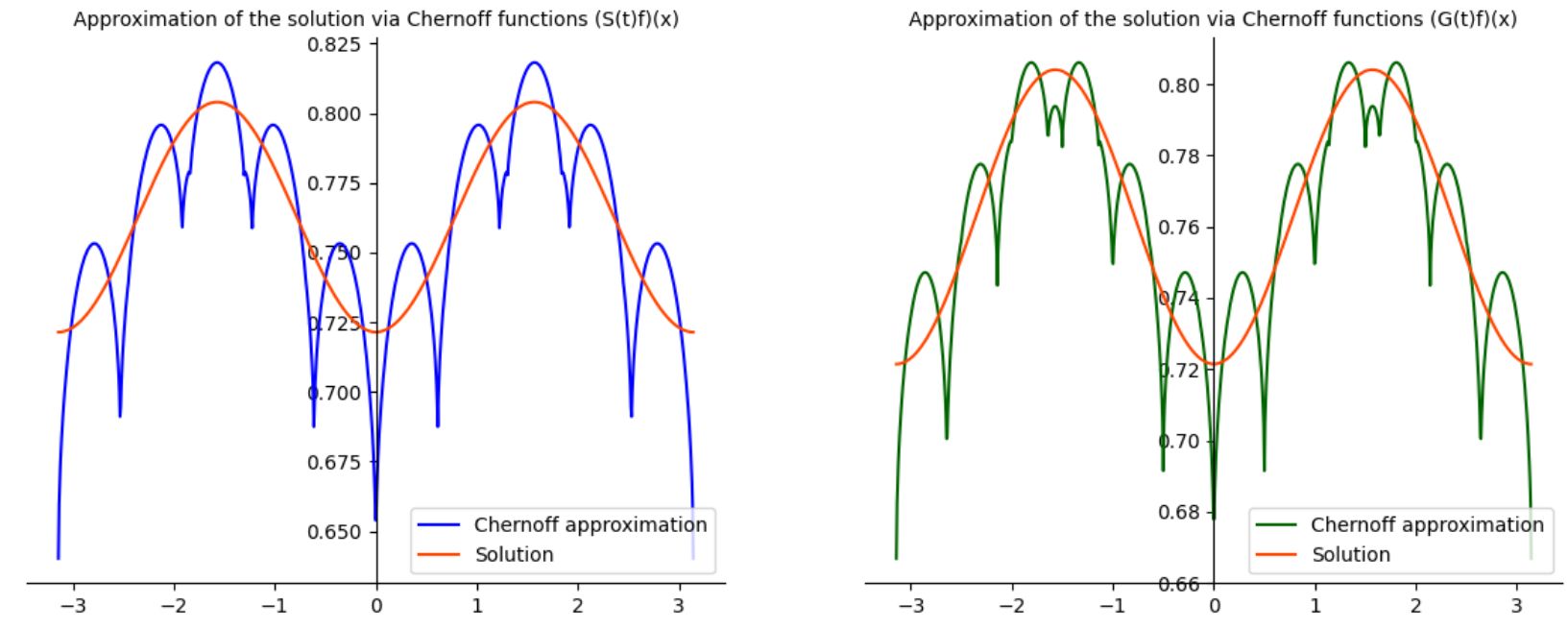}\\
\end{center}
n=9
\begin{center}
	\includegraphics[scale=0.5]{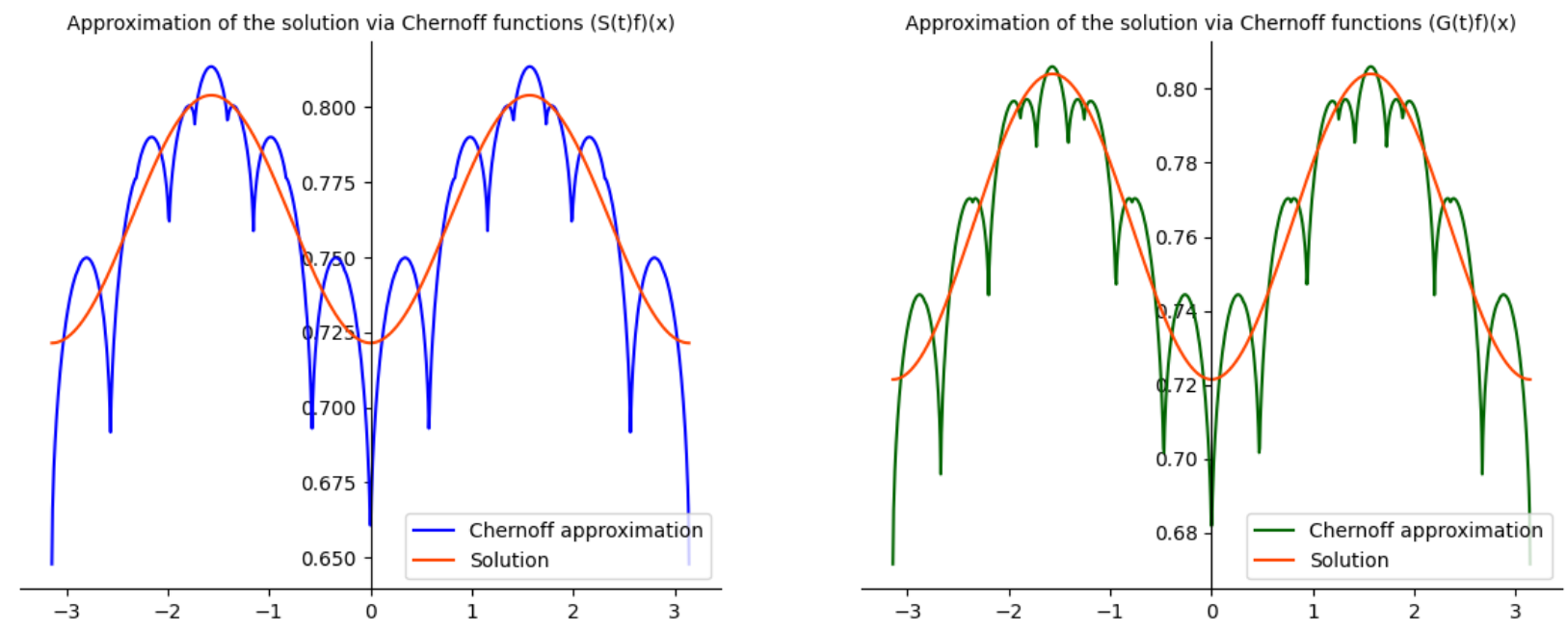}\\
\end{center}
n=10
\begin{center}
	\includegraphics[scale=0.5]{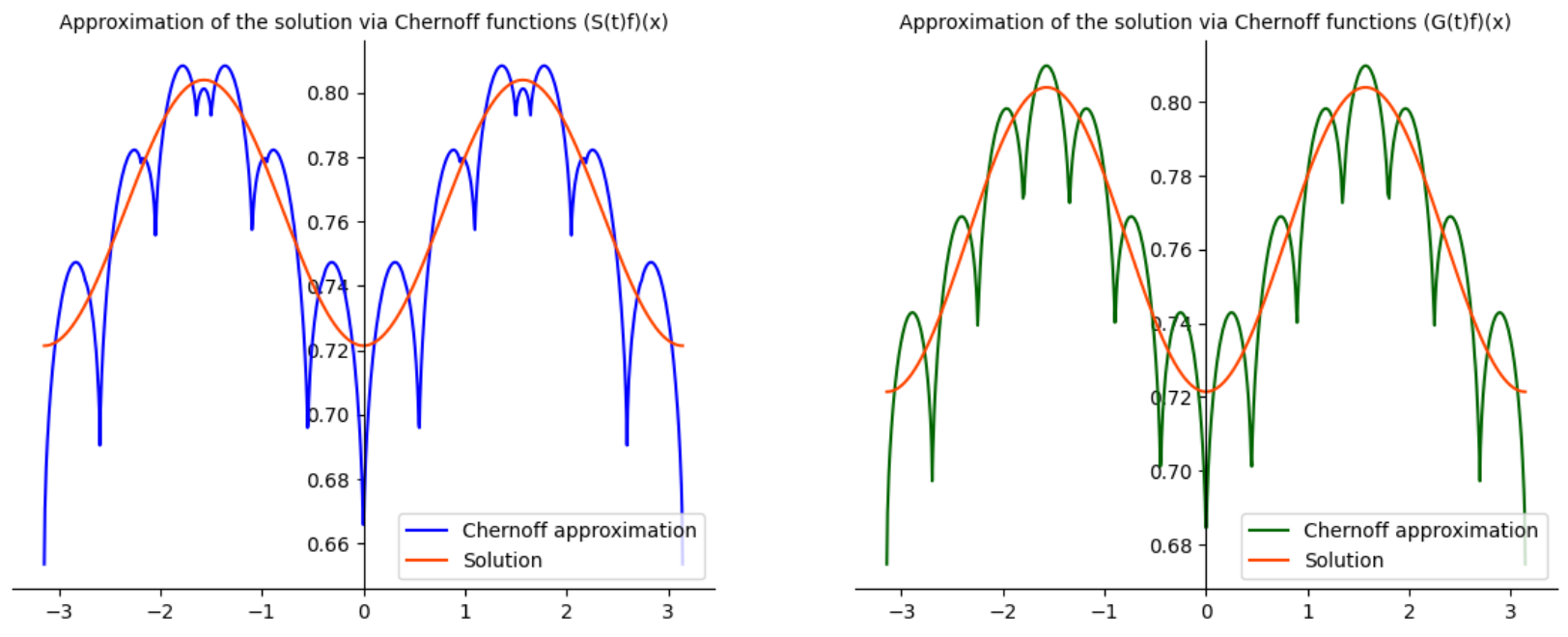}\\
\end{center}

\subsection{$u_0(x)=\sqrt[4]{|\sin x|}$}

n=1
\begin{center}
	\includegraphics[scale=0.5]{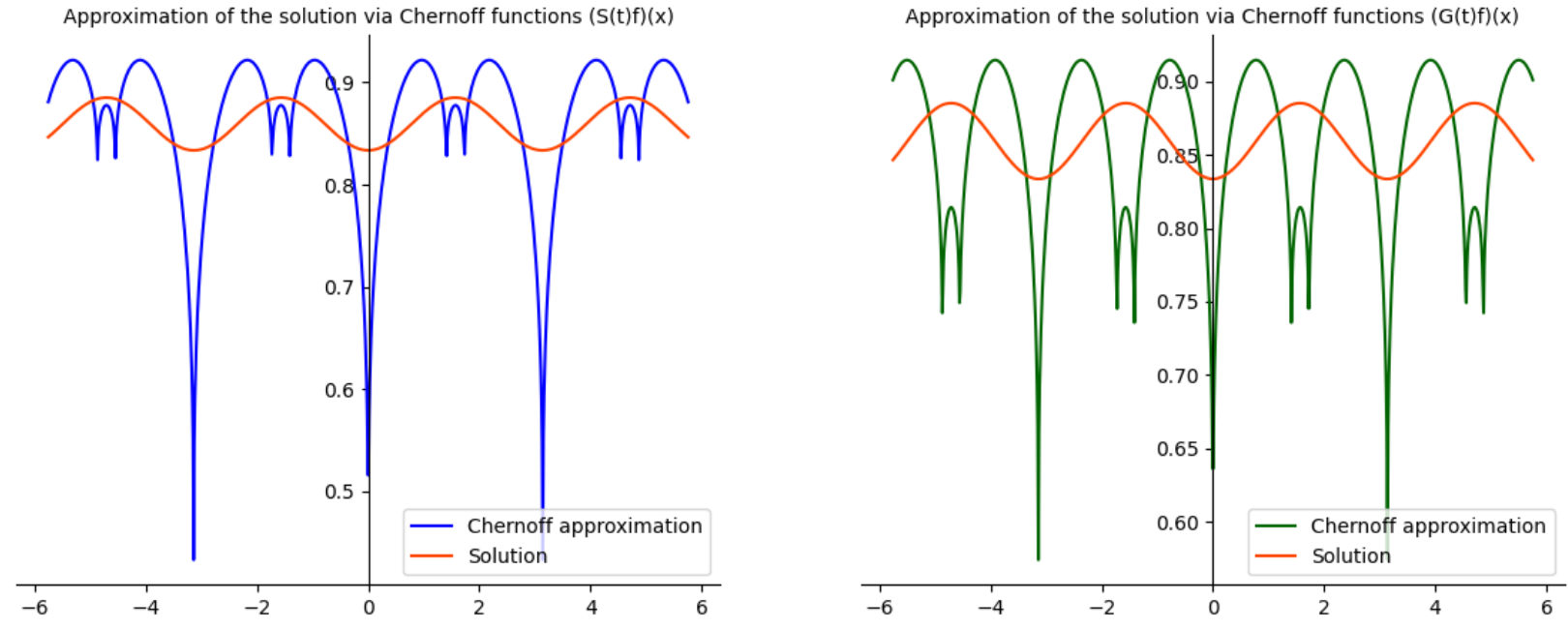}\\
\end{center}

\newpage
n=2
\begin{center}
	\includegraphics[scale=0.5]{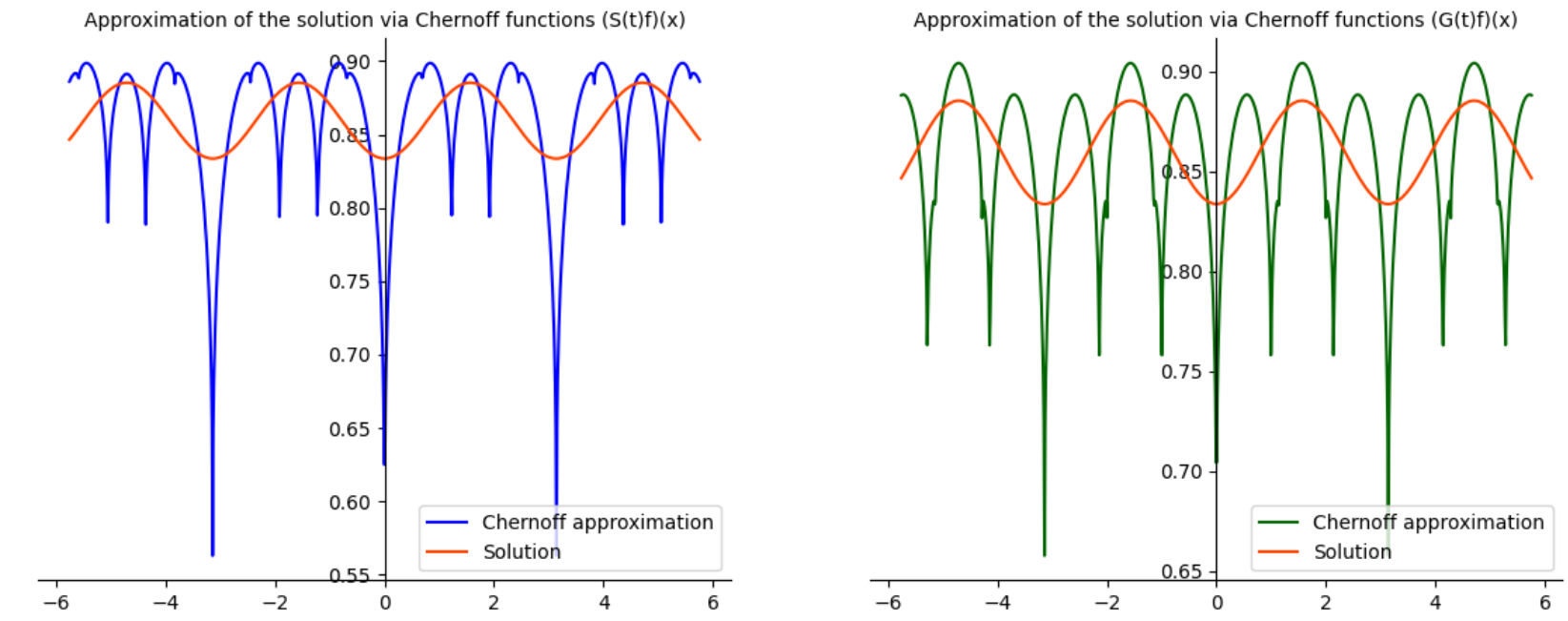}\\
\end{center}
n=3
\begin{center}
	\includegraphics[scale=0.5]{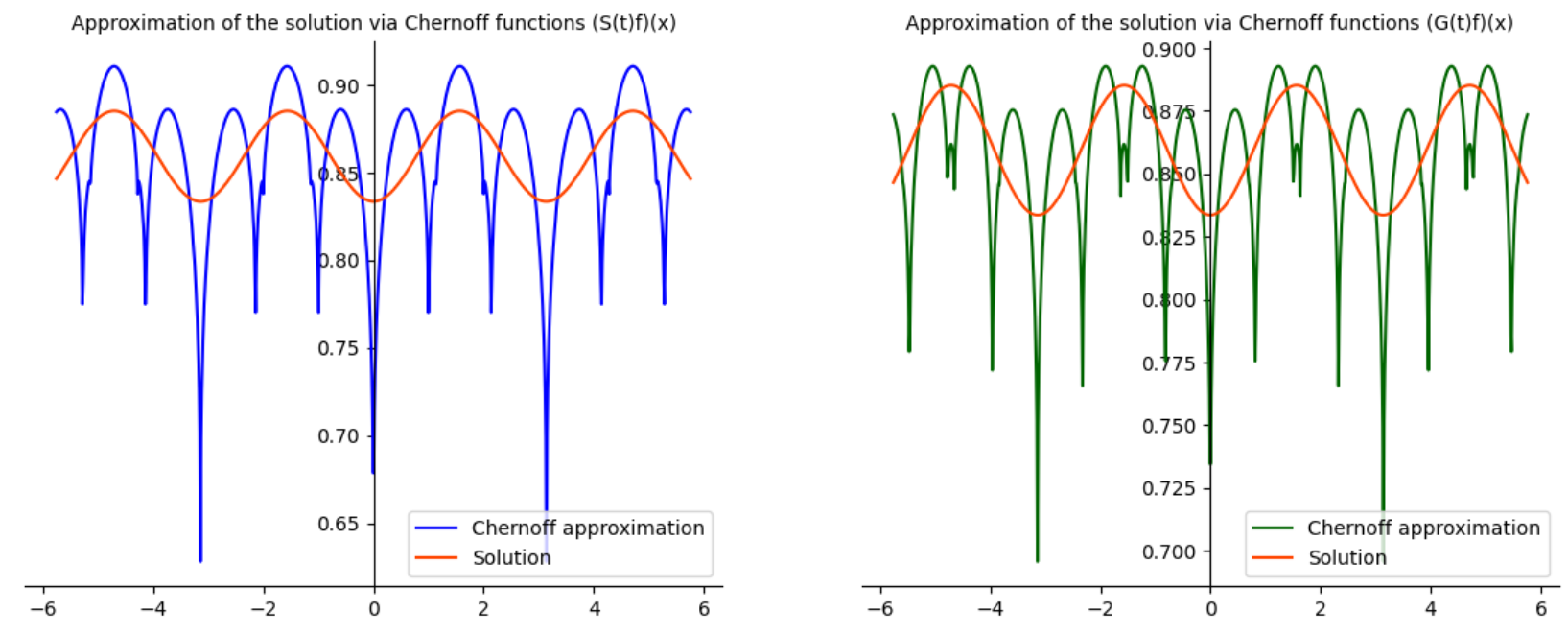}\\
\end{center}
n=4
\begin{center}
	\includegraphics[scale=0.5]{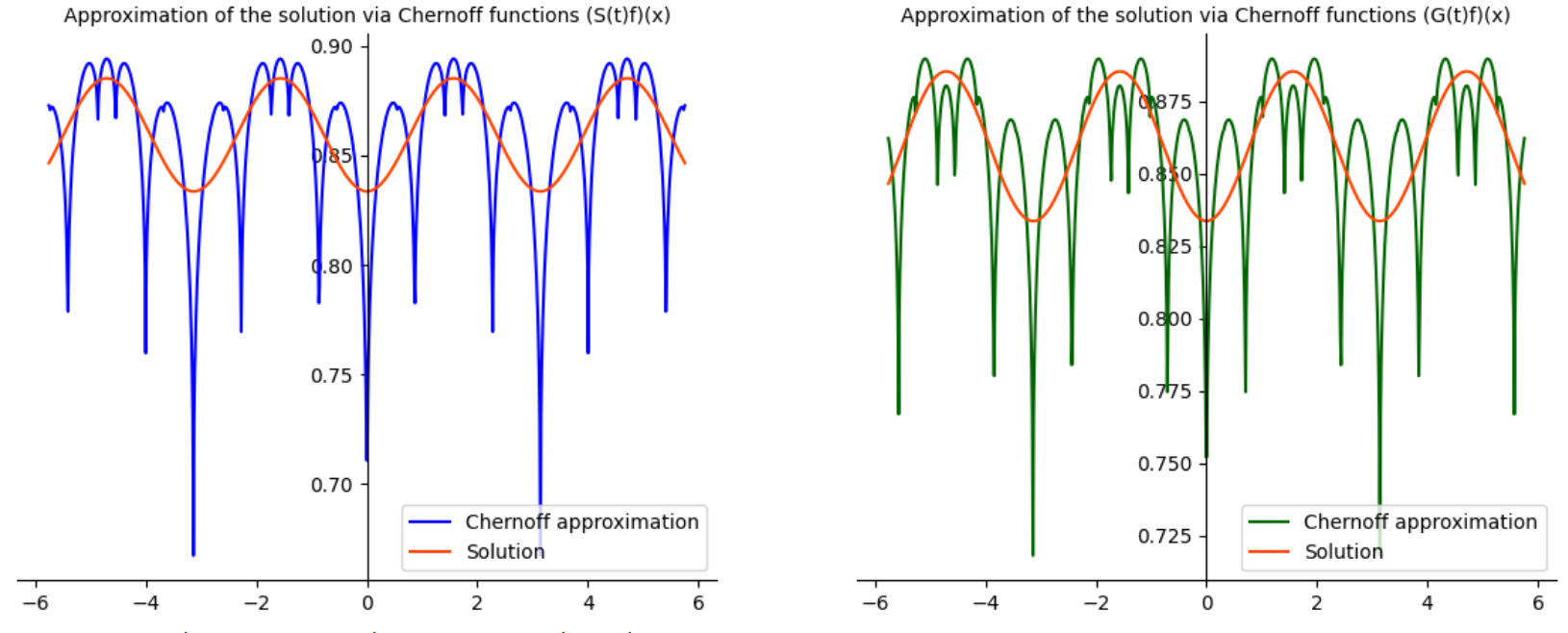}\\
\end{center}

\newpage
n=5
\begin{center}
	\includegraphics[scale=0.5]{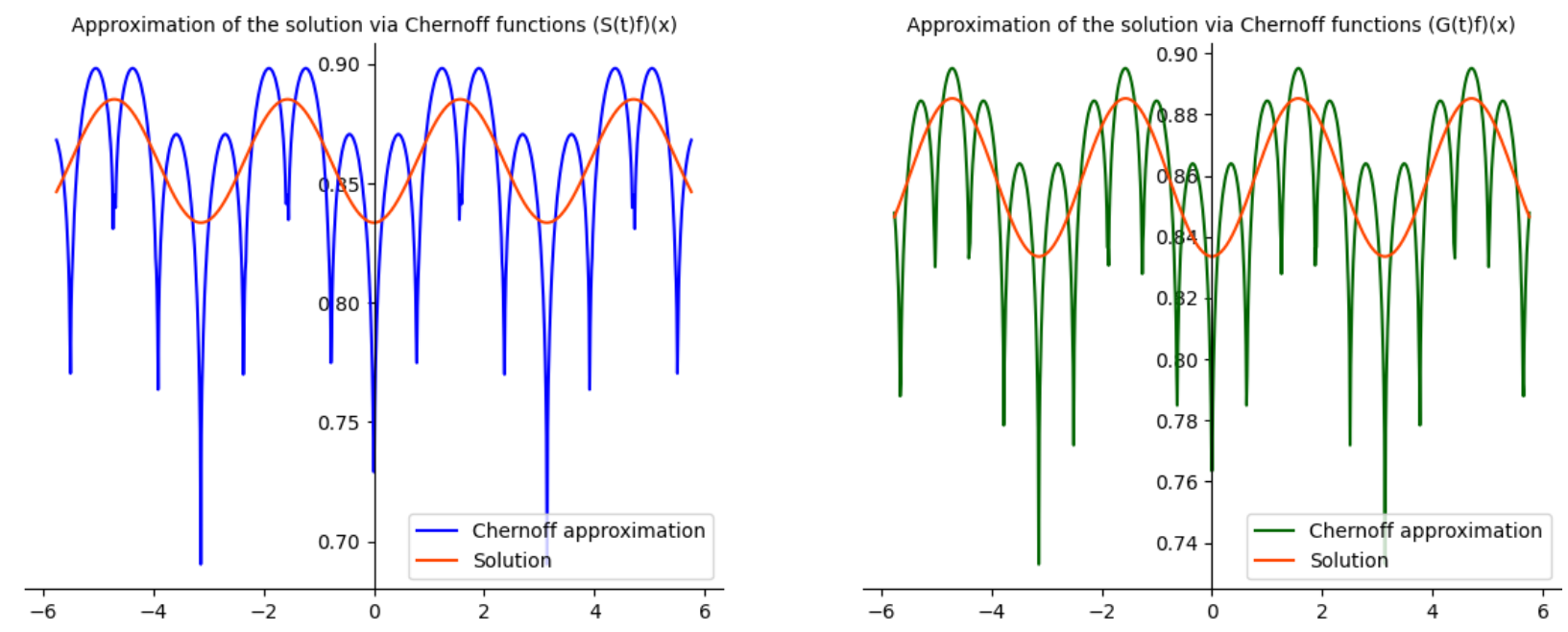}\\
\end{center}
n=6
\begin{center}
	\includegraphics[scale=0.5]{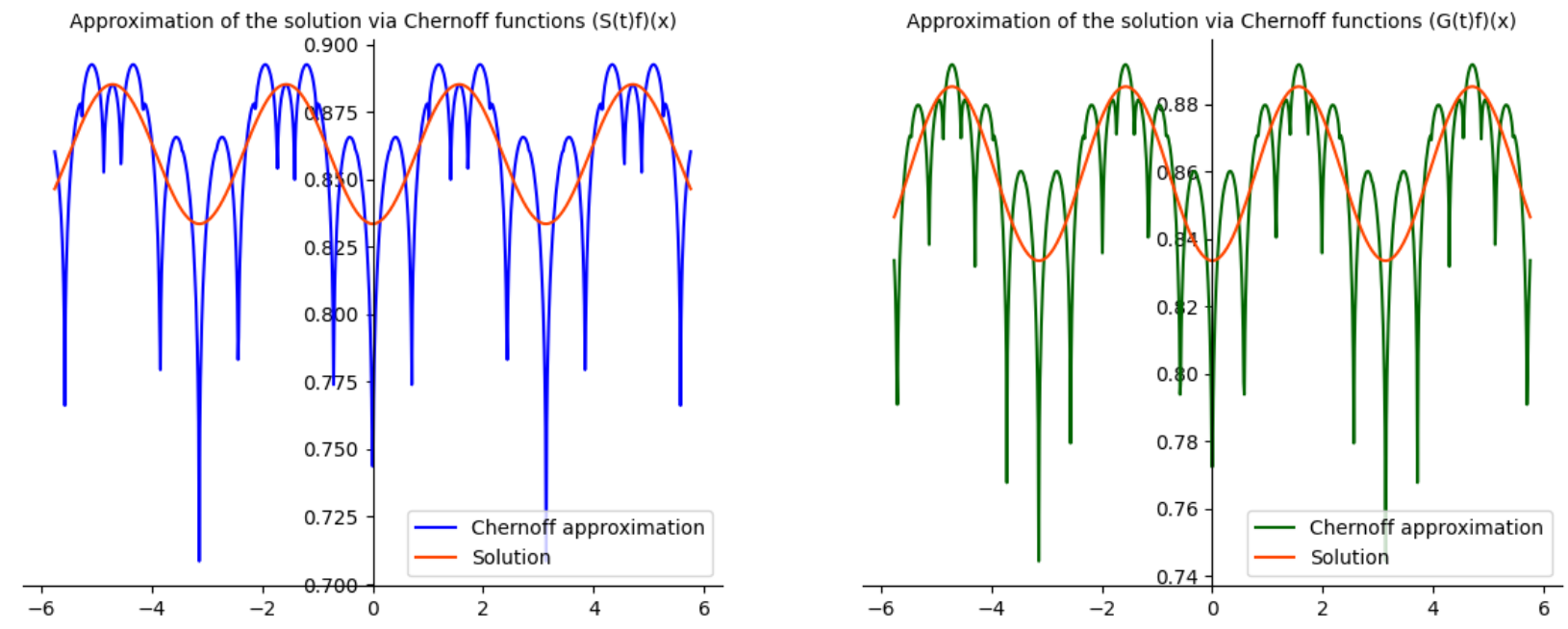}\\
\end{center}
n=7
\begin{center}
	\includegraphics[scale=0.5]{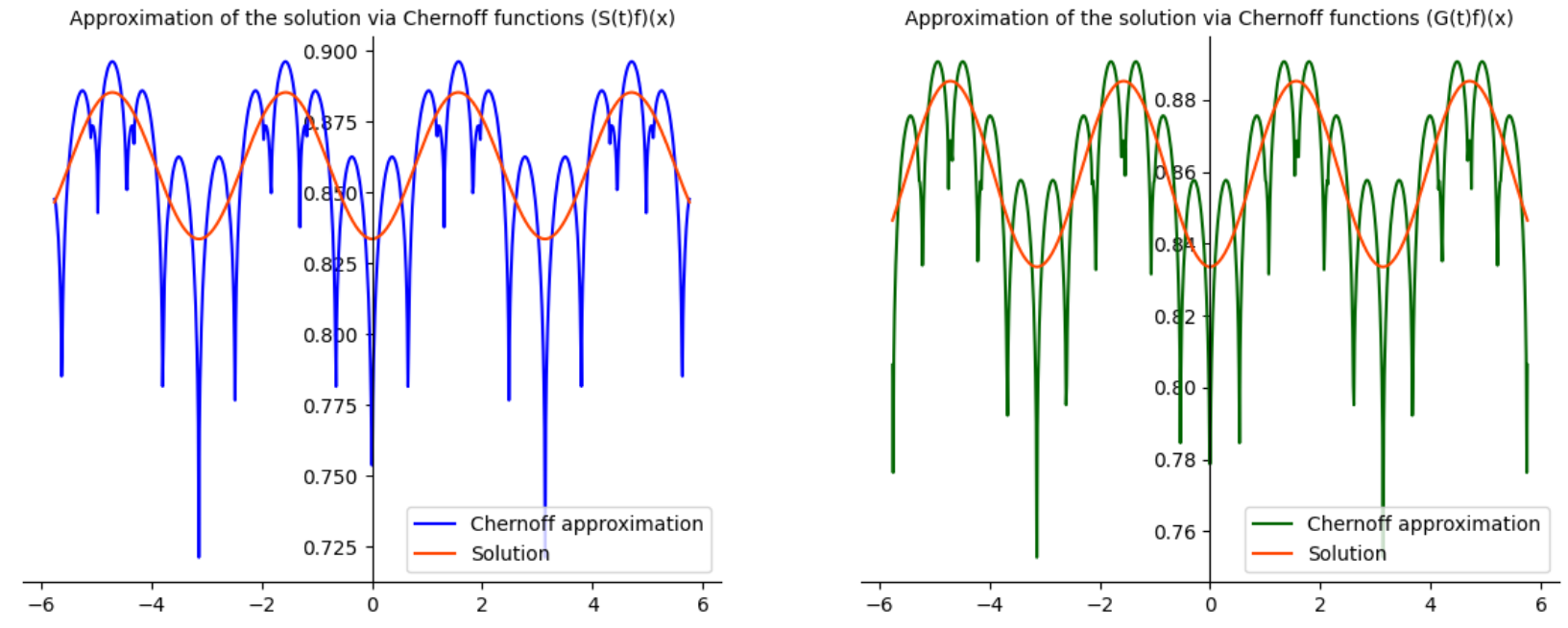}\\
\end{center}

\newpage
n=8
\begin{center}
	\includegraphics[scale=0.5]{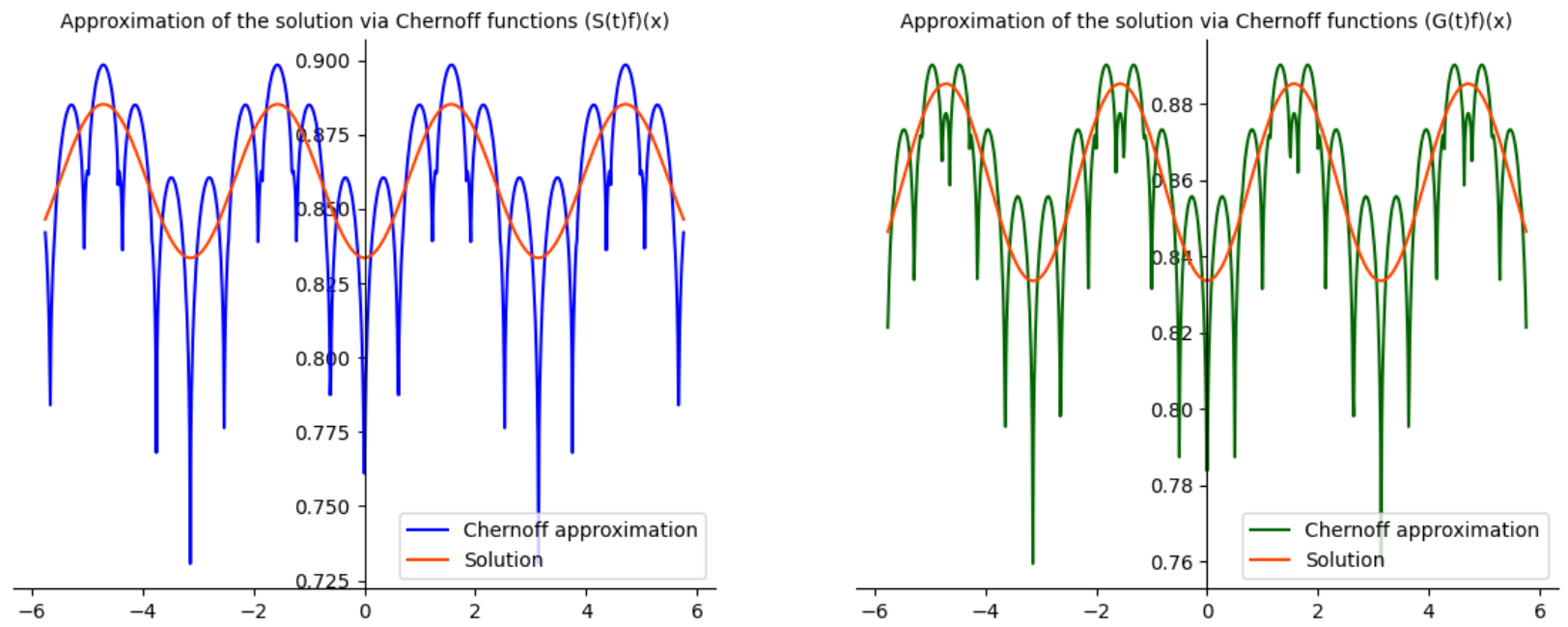}\\
\end{center}
n=9
\begin{center}
	\includegraphics[scale=0.5]{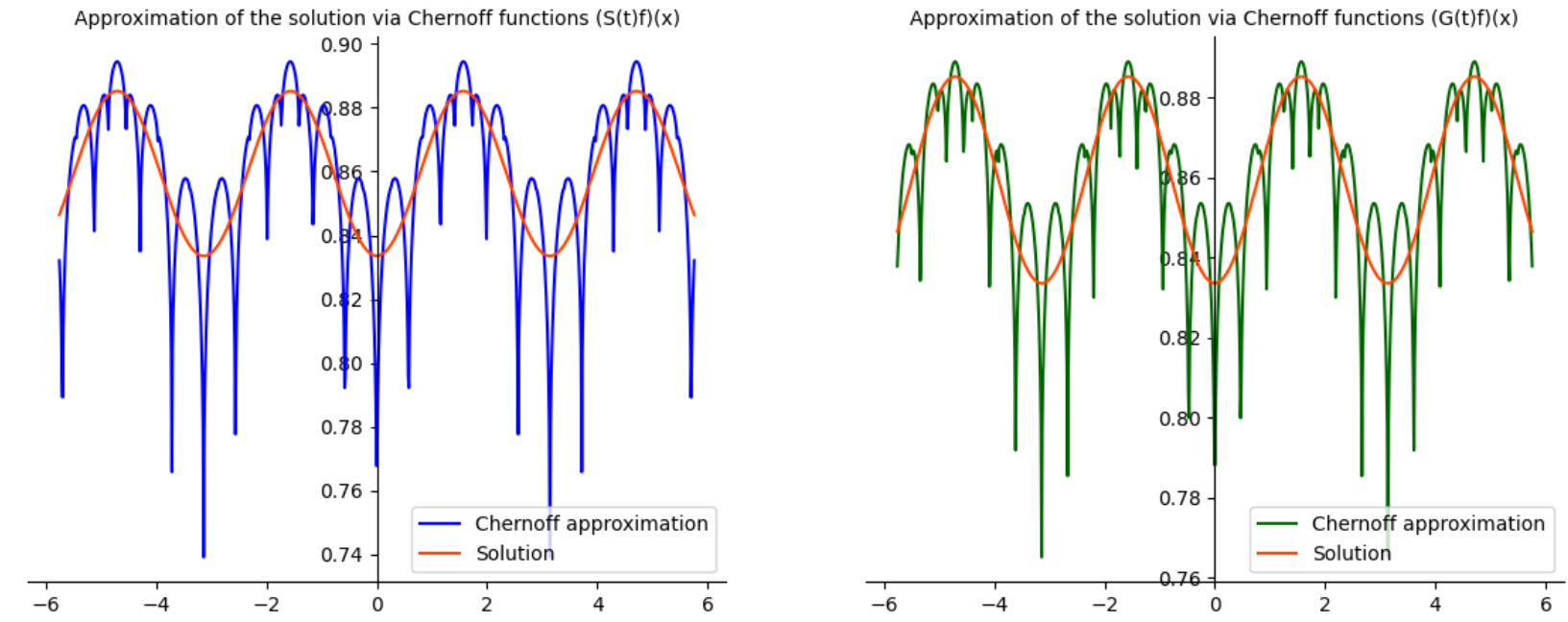}\\
\end{center}
n=10
\begin{center}
	\includegraphics[scale=0.5]{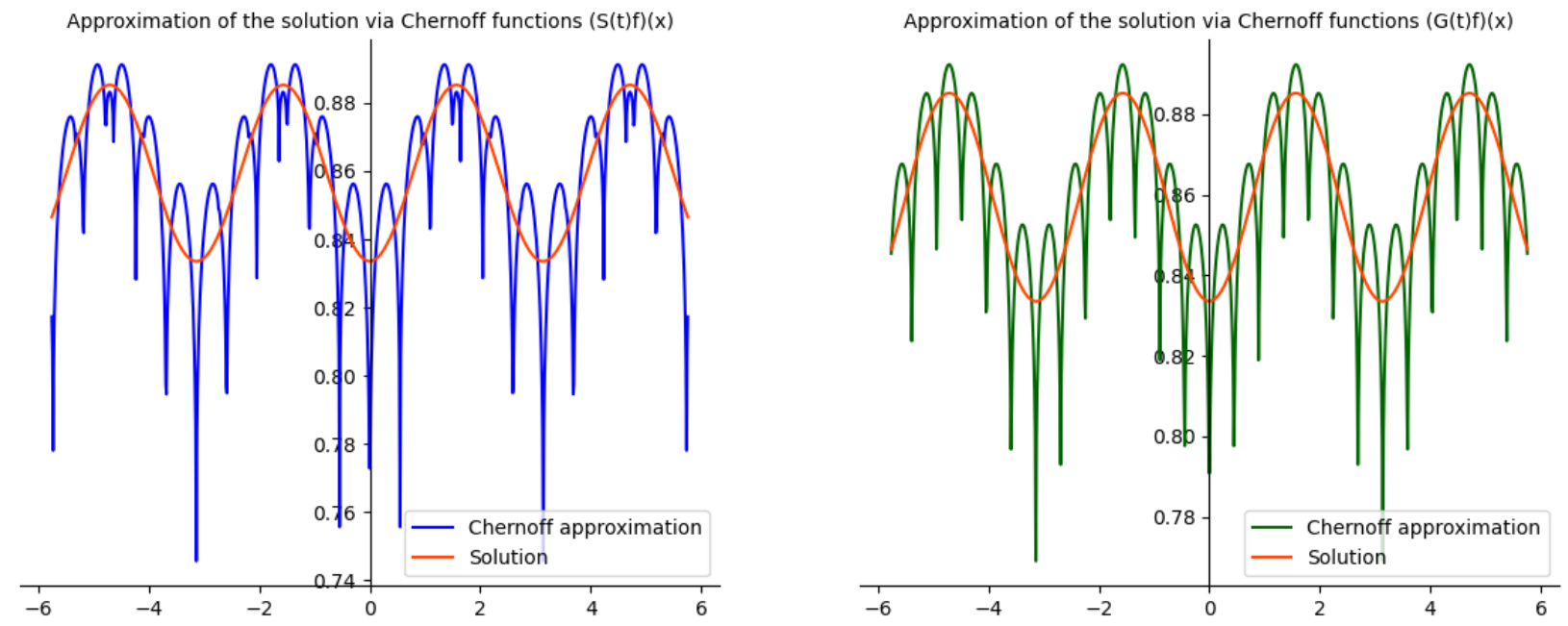}\\
\end{center}
\newpage
\subsection{$u_0(x)=|\sin(x)|^{3/2}$}
n=1
\begin{center}
	\includegraphics[scale=0.5]{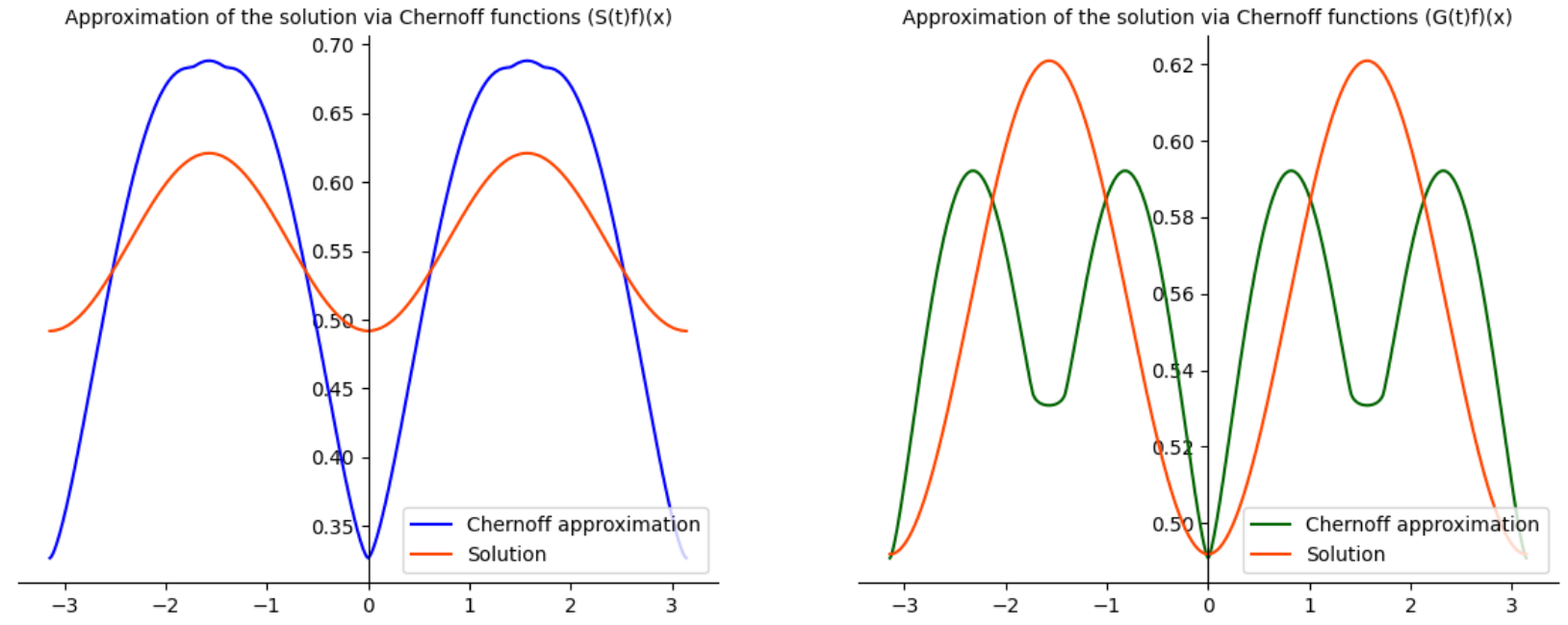}\\
\end{center}
n=2
\begin{center}
	\includegraphics[scale=0.5]{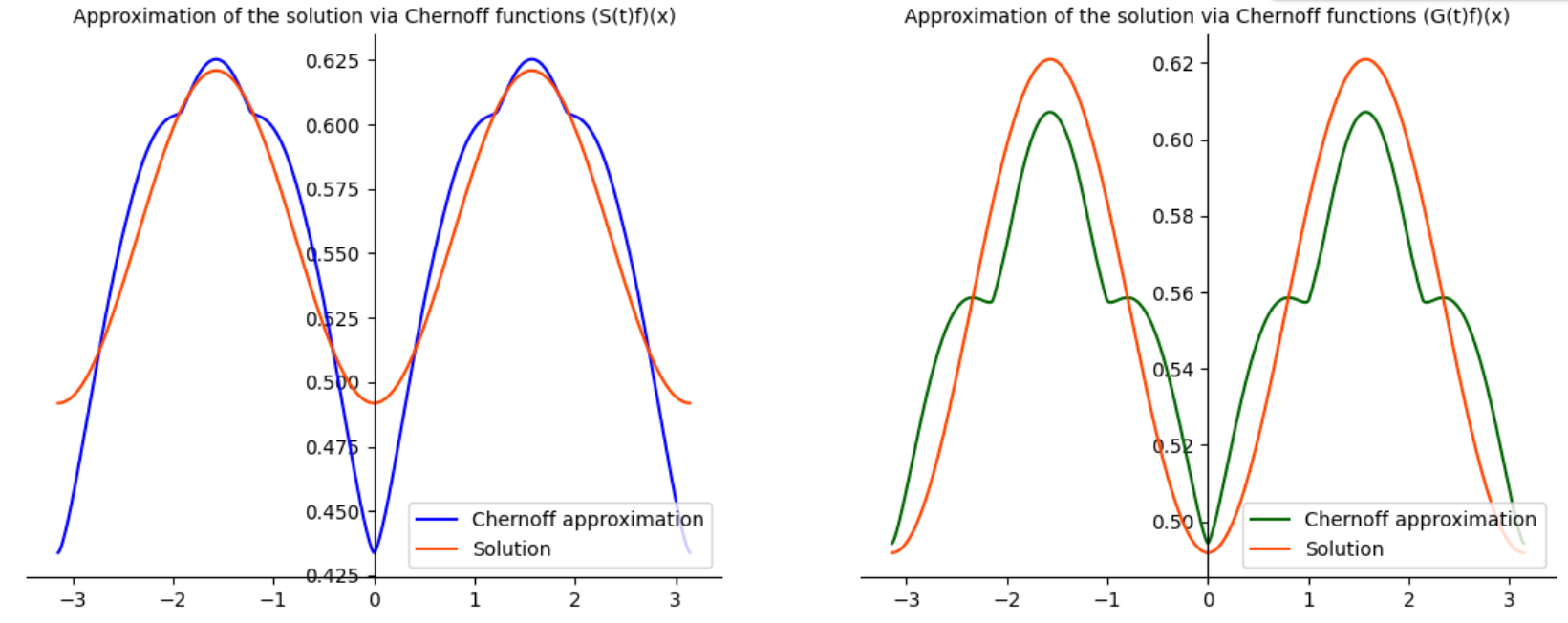}\\
\end{center}
n=3
\begin{center}
	\includegraphics[scale=0.5]{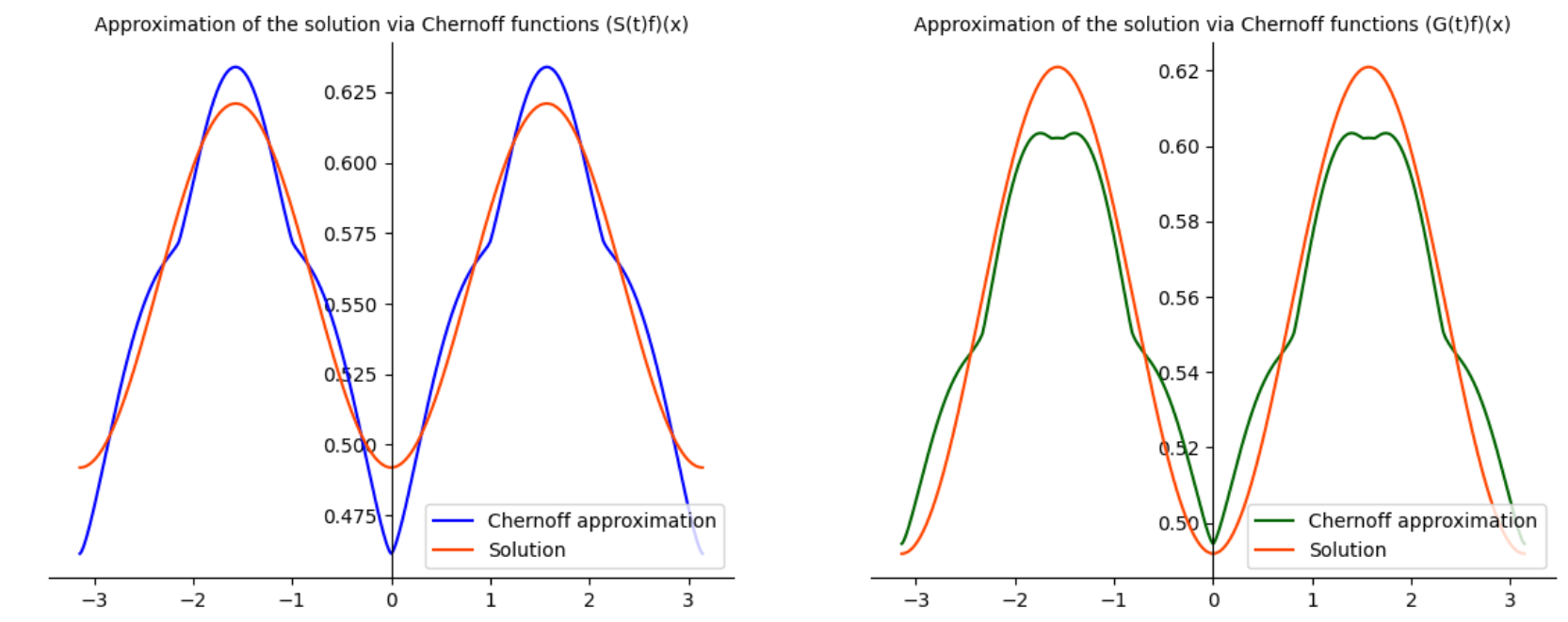}\\
\end{center}

\newpage
n=4
\begin{center}
	\includegraphics[scale=0.5]{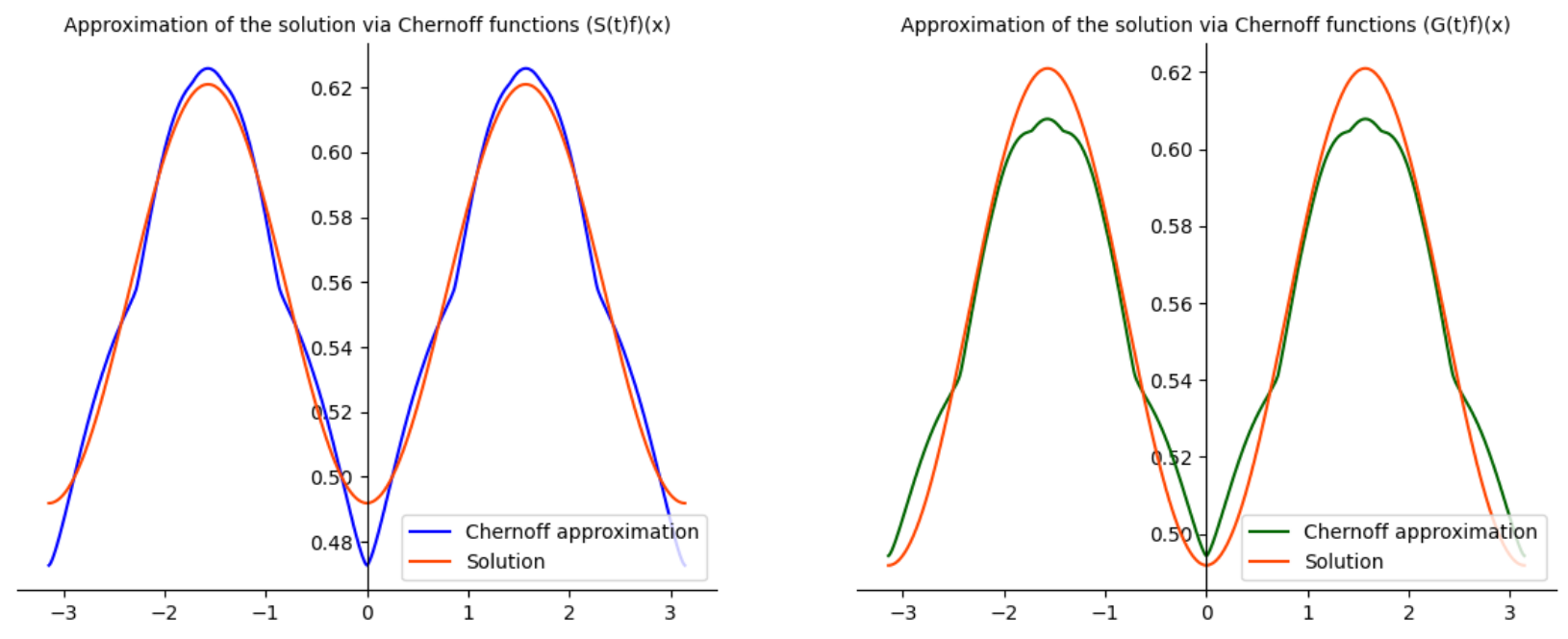}\\
\end{center}
n=5
\begin{center}
	\includegraphics[scale=0.5]{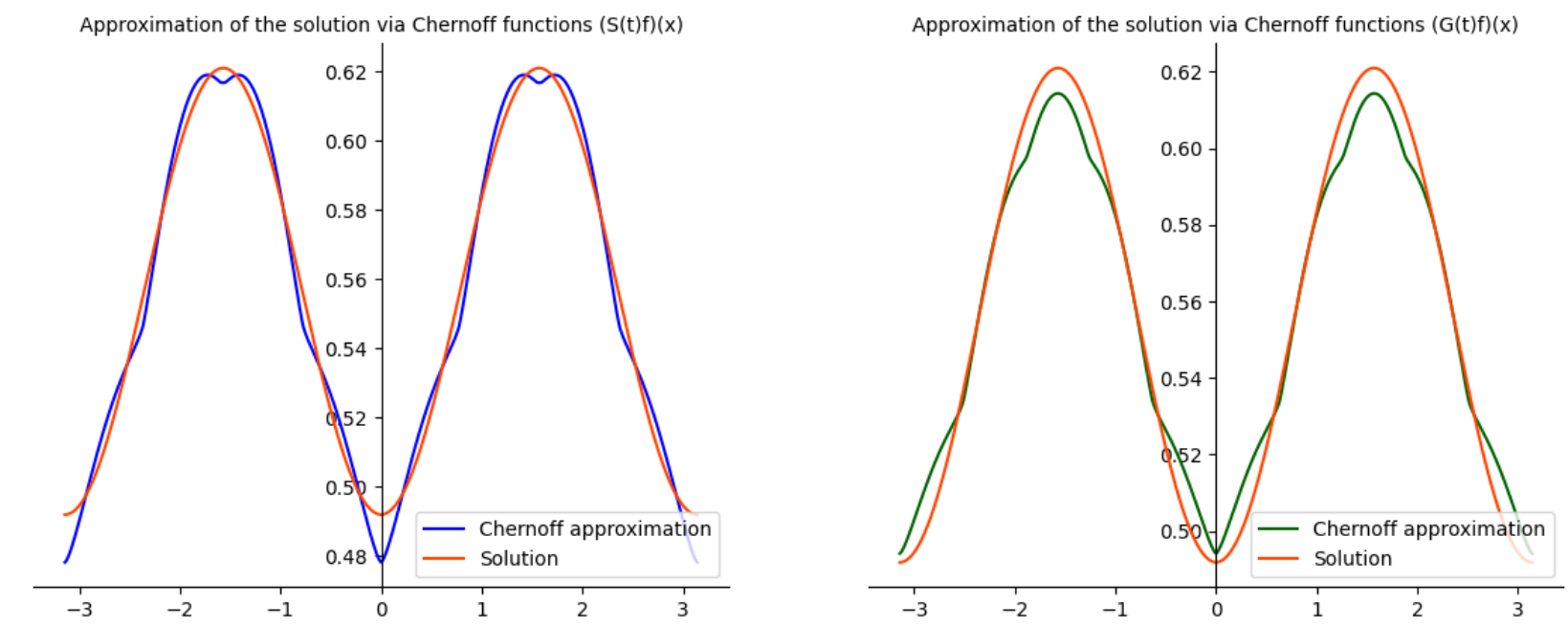}\\
\end{center}
n=6
\begin{center}
	\includegraphics[scale=0.5]{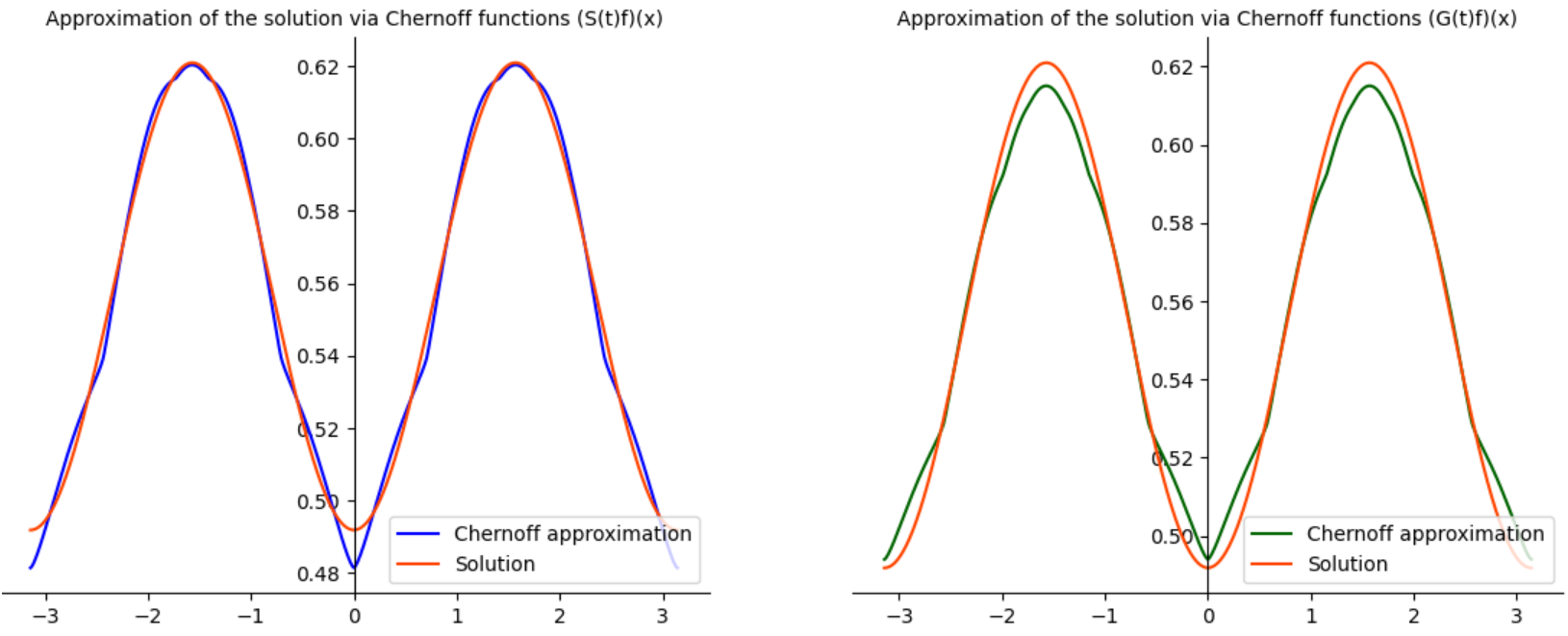}\\
\end{center}

\newpage
n=7
\begin{center}
	\includegraphics[scale=0.5]{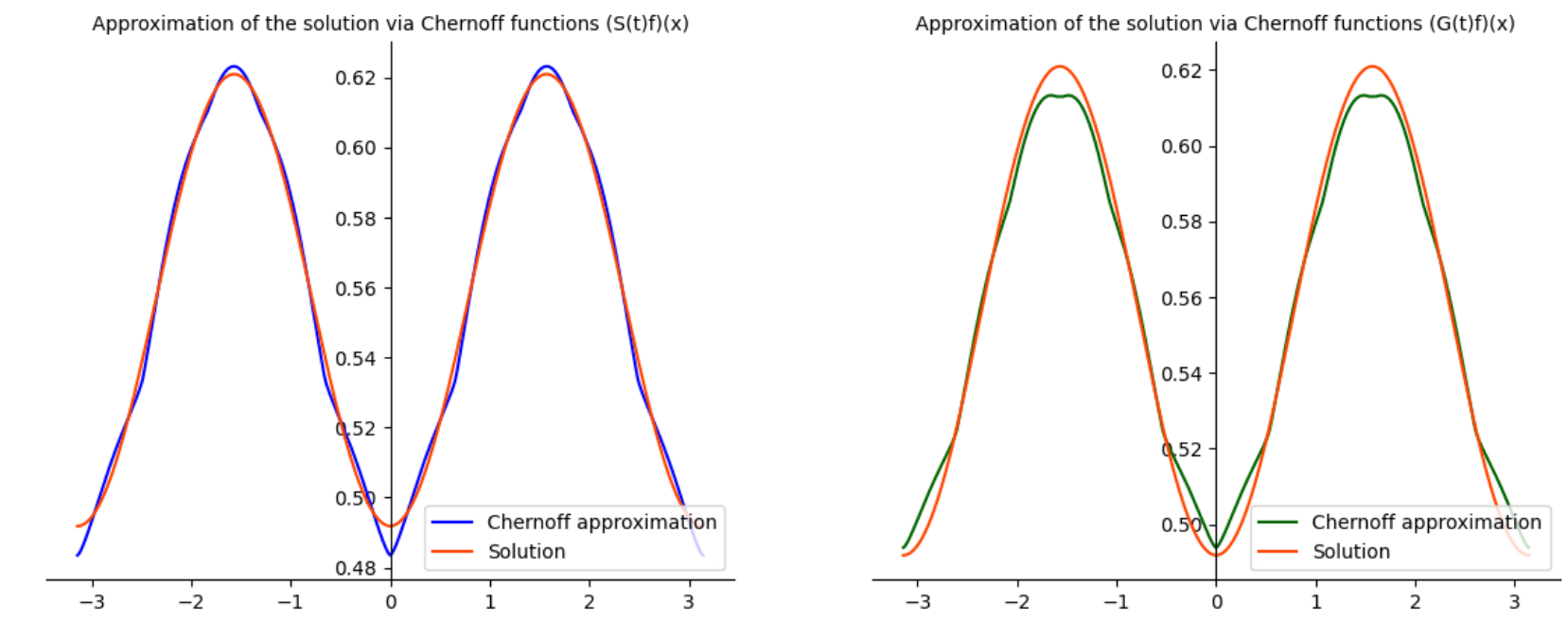}\\
\end{center}
n=8
\begin{center}
	\includegraphics[scale=0.5]{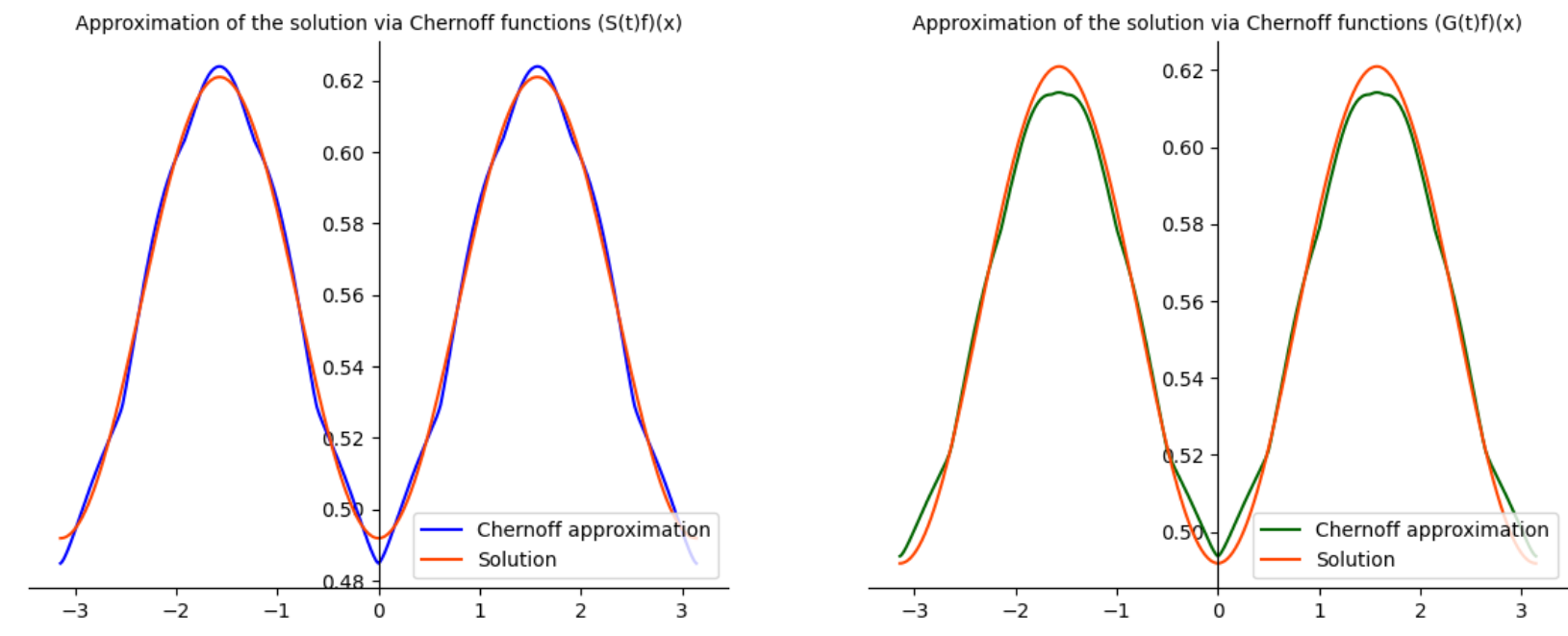}\\
\end{center}
n=9
\begin{center}
	\includegraphics[scale=0.5]{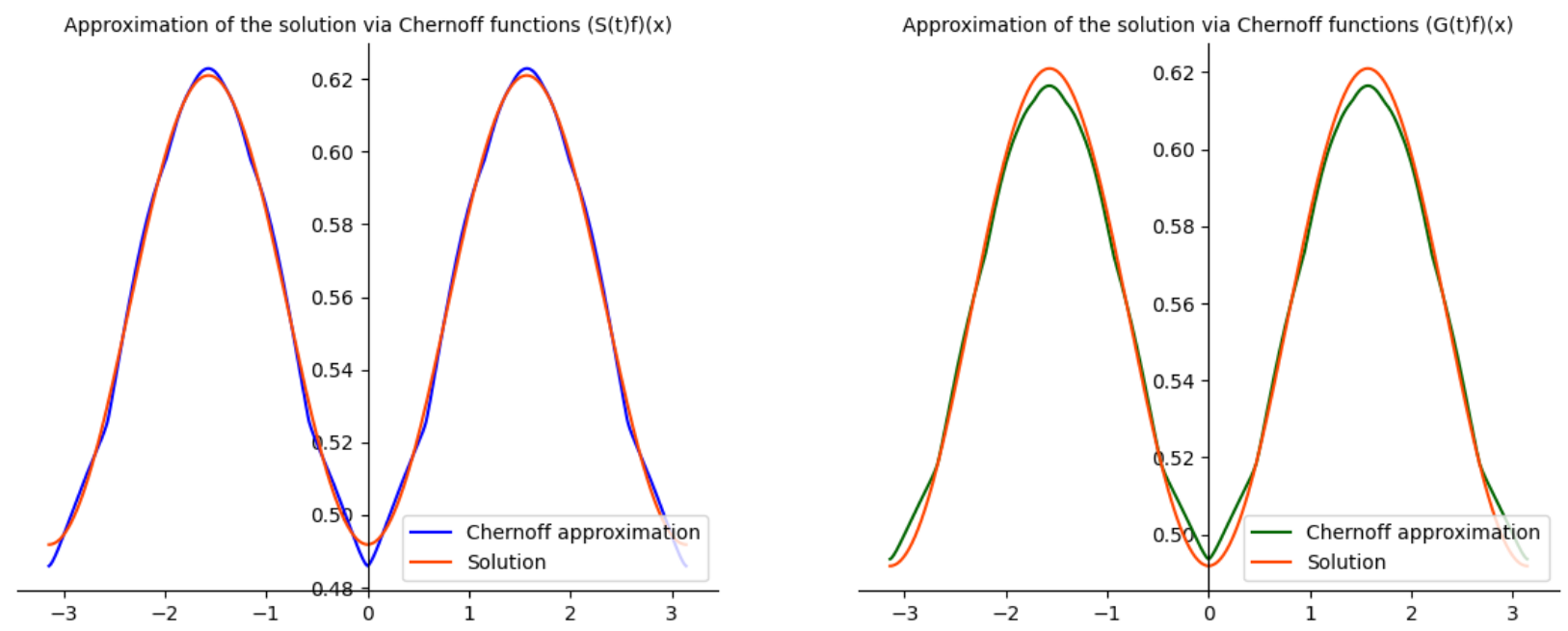}\\
\end{center}
\newpage
n=10
\begin{center}
	\includegraphics[scale=0.5]{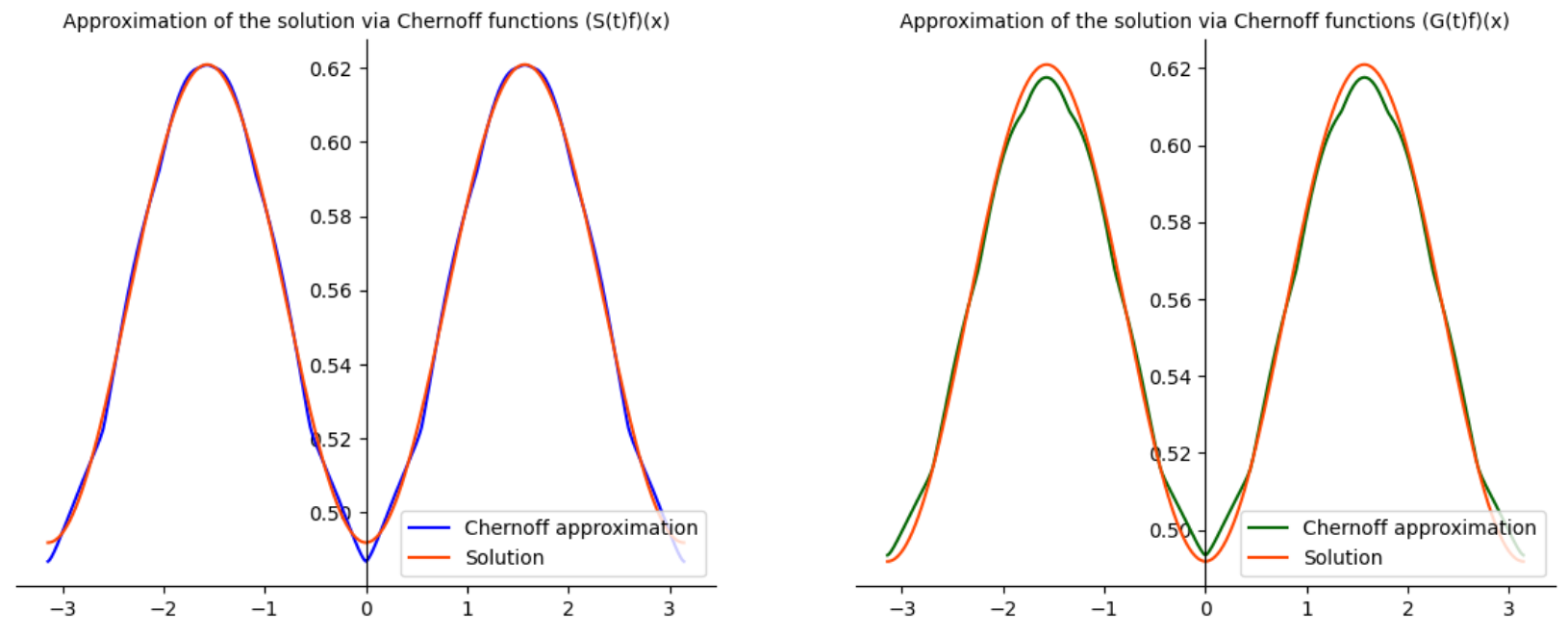}\\
\end{center}

\subsection{$u_0(x)=|\sin(x)|^{5/2}$}

n=1
\begin{center}
	\includegraphics[scale=0.5]{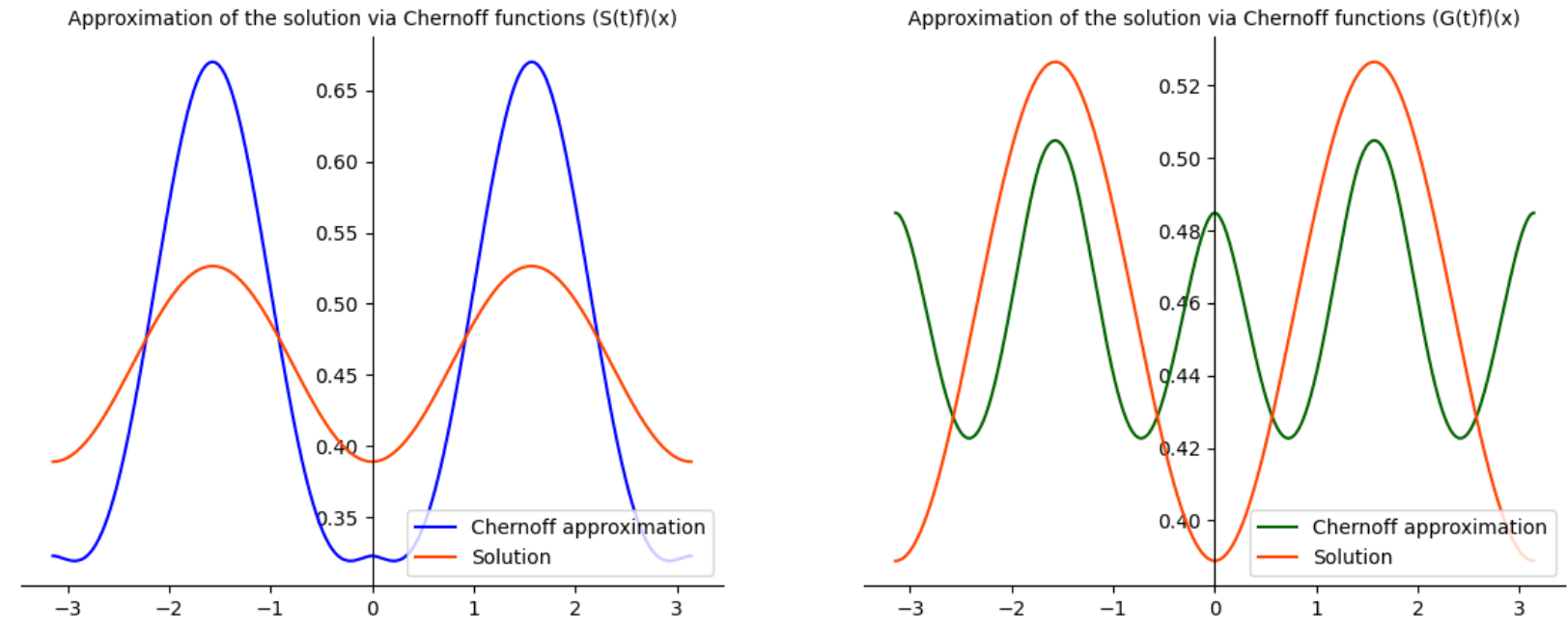}\\
\end{center}
\newpage
n=2
\begin{center}
	\includegraphics[scale=0.5]{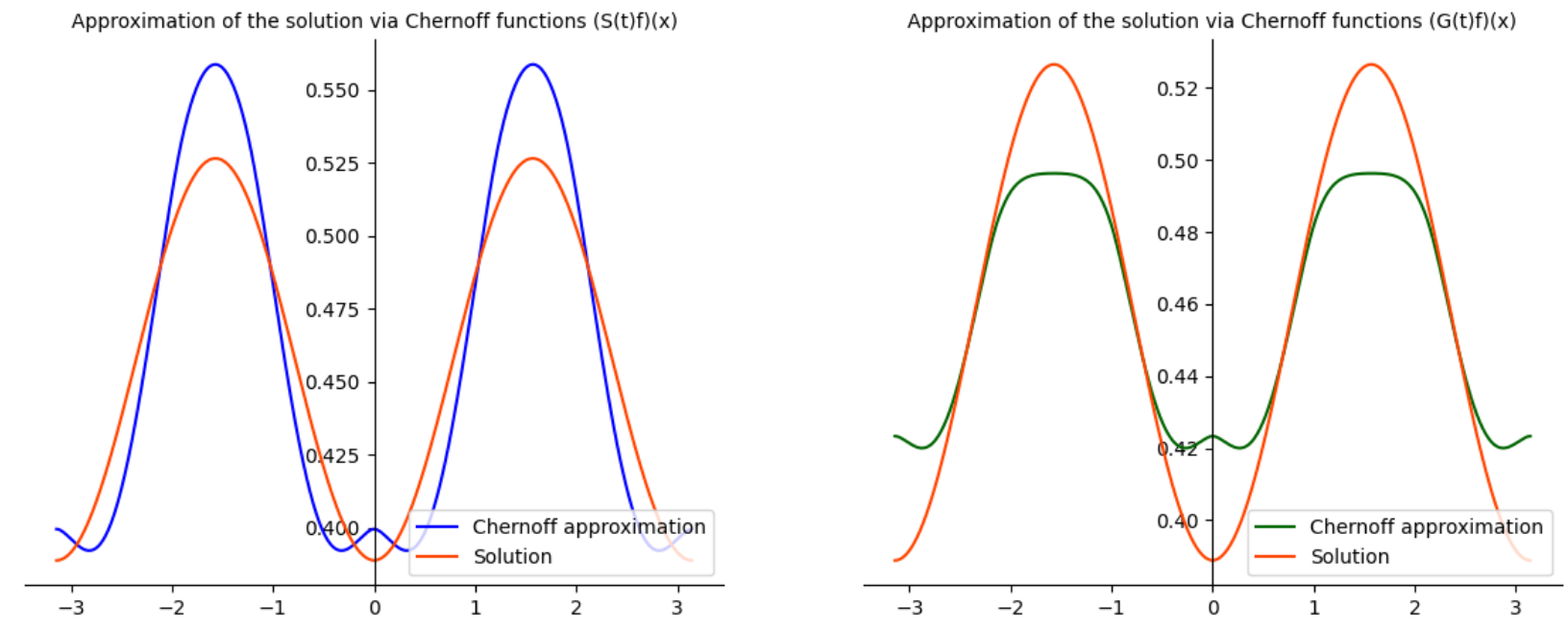}\\
\end{center}
n=3
\begin{center}
	\includegraphics[scale=0.5]{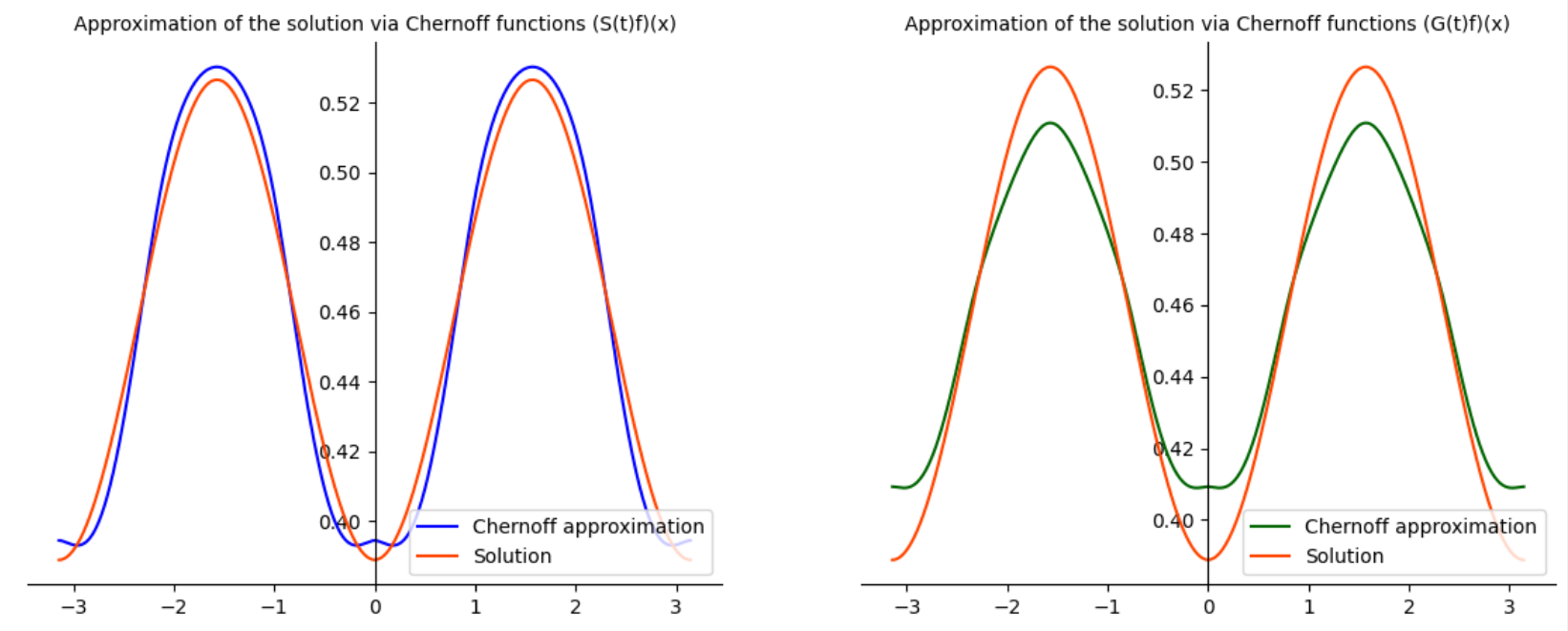}\\
\end{center}
n=4
\begin{center}
	\includegraphics[scale=0.5]{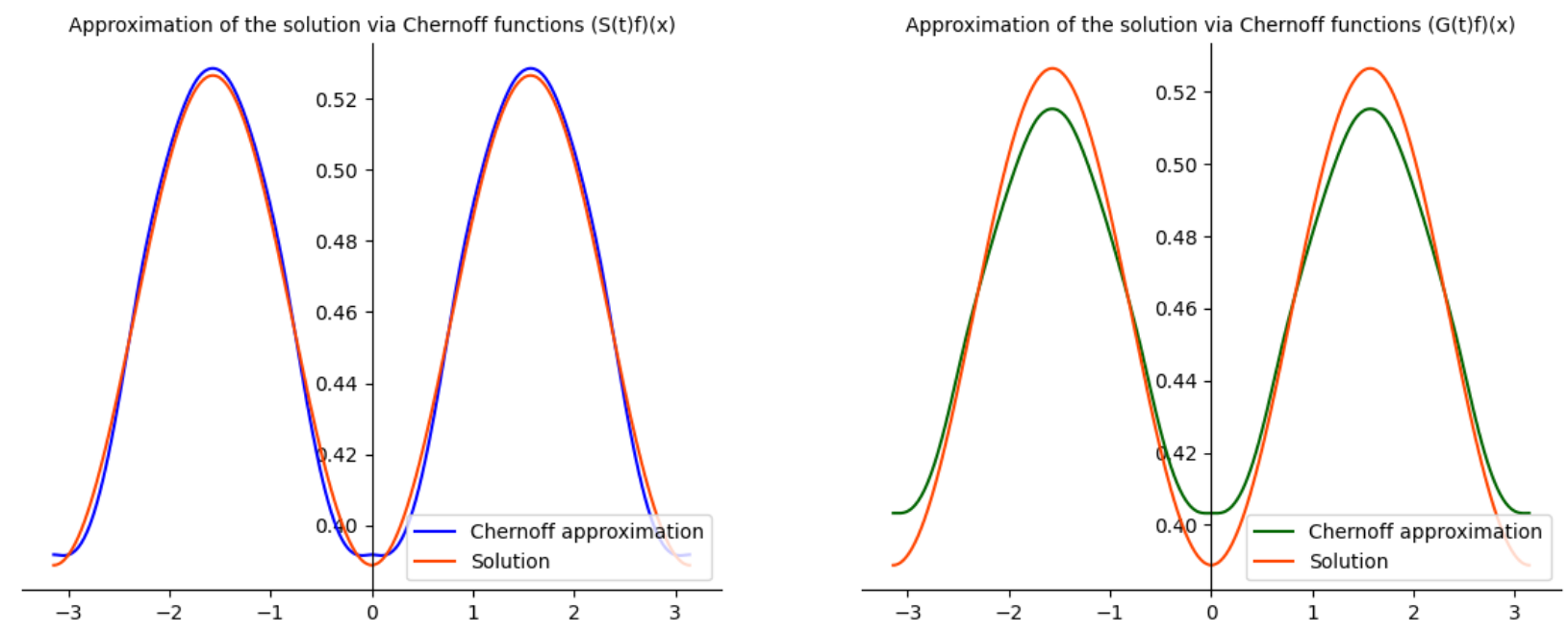}\\
\end{center}
\newpage

n=5
\begin{center}
	\includegraphics[scale=0.5]{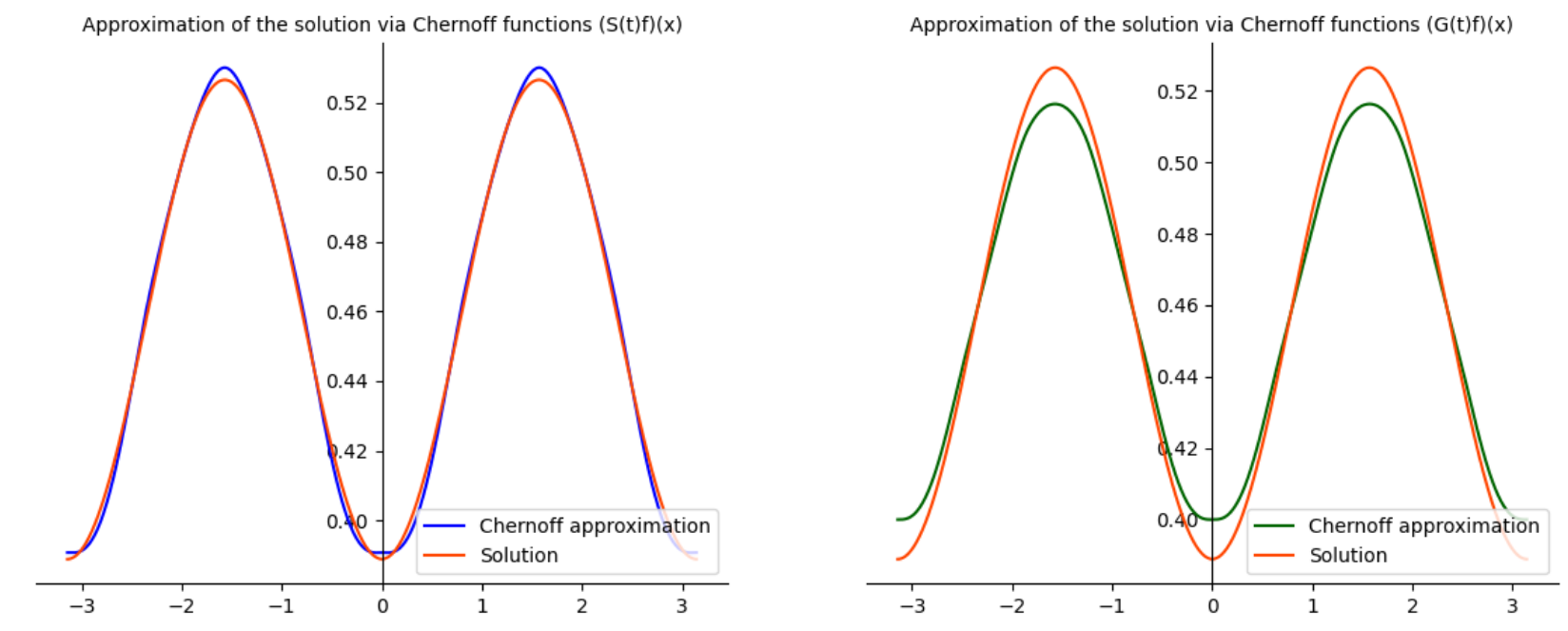}\\
\end{center}
n=6
\begin{center}
	\includegraphics[scale=0.5]{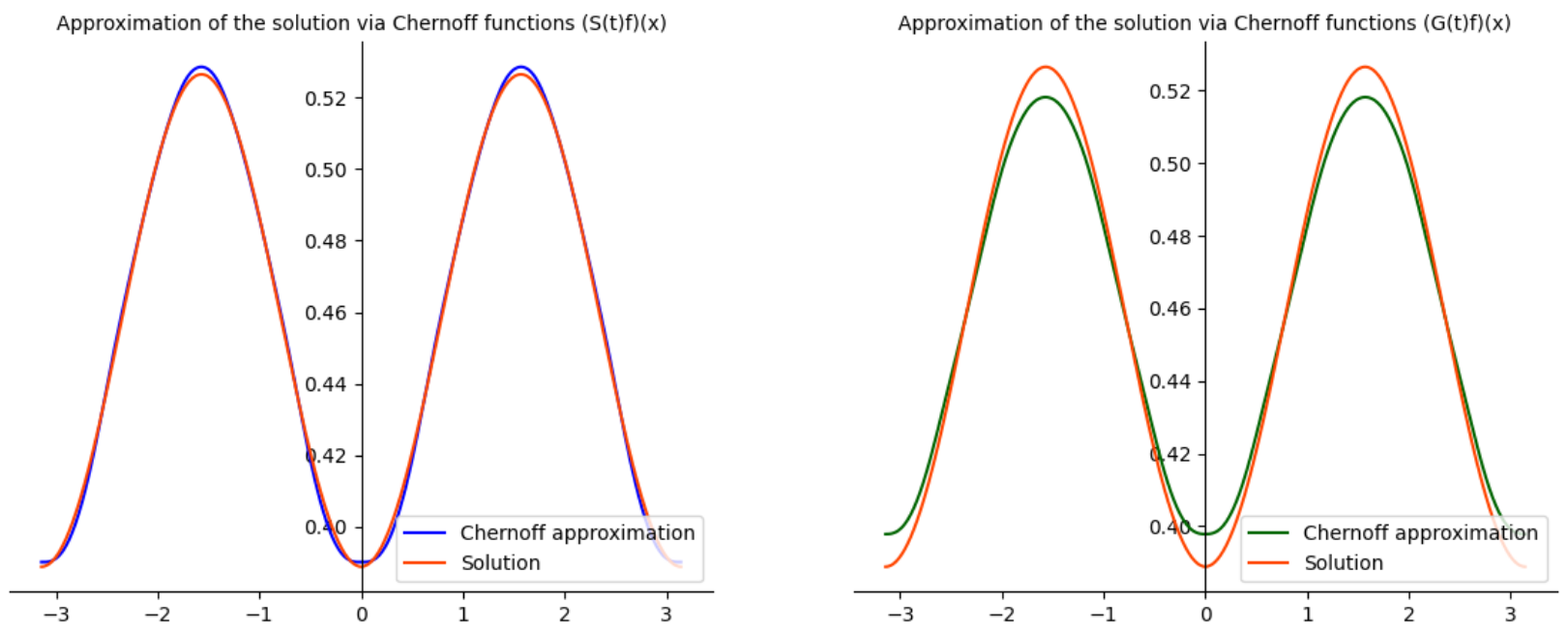}\\
\end{center}
n=7
\begin{center}
	\includegraphics[scale=0.5]{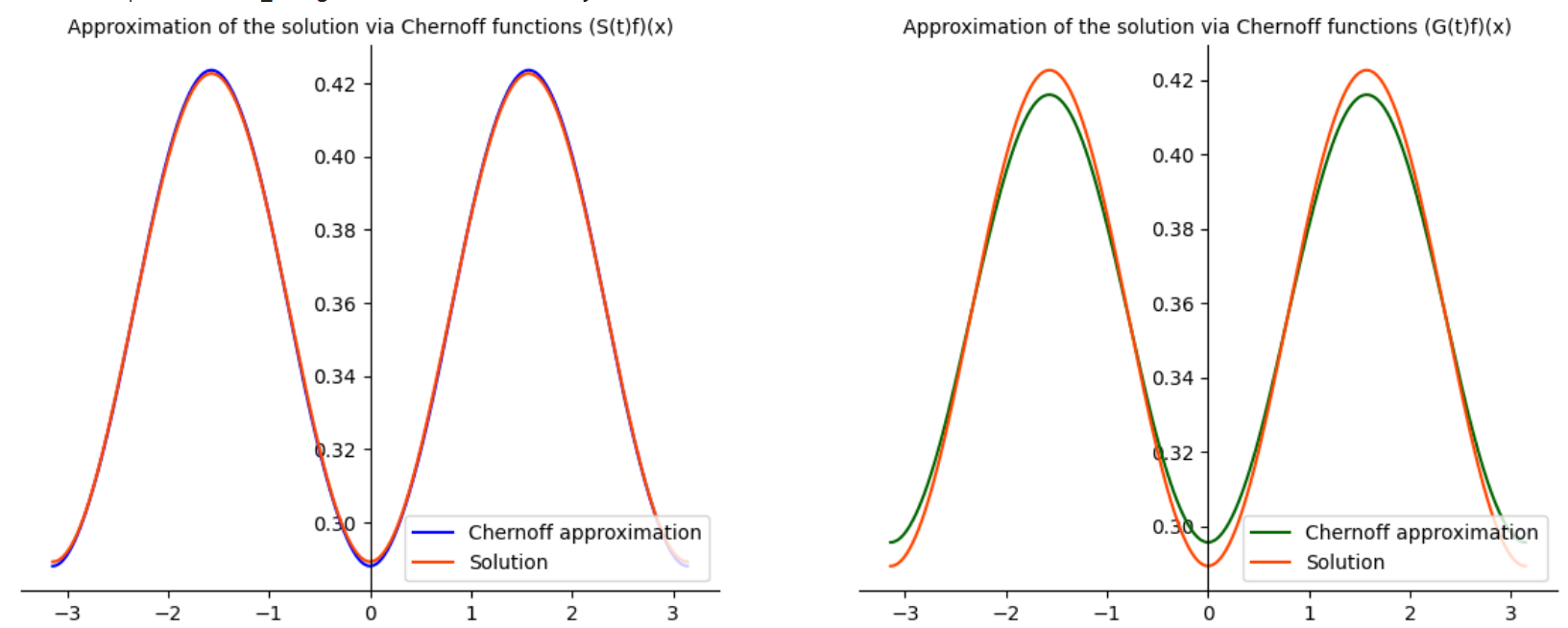}\\
\end{center}

\newpage
n=8
\begin{center}
	\includegraphics[scale=0.5]{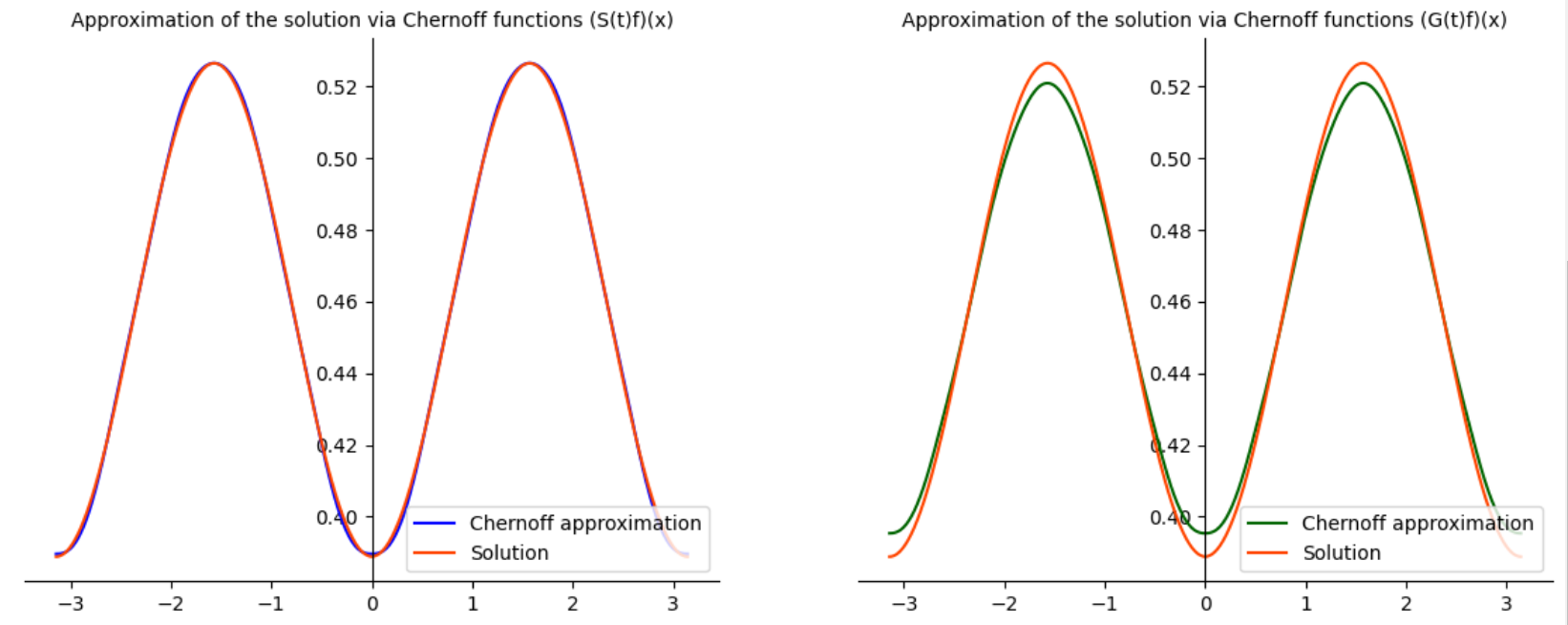}\\
\end{center}
n=9
\begin{center}
	\includegraphics[scale=0.5]{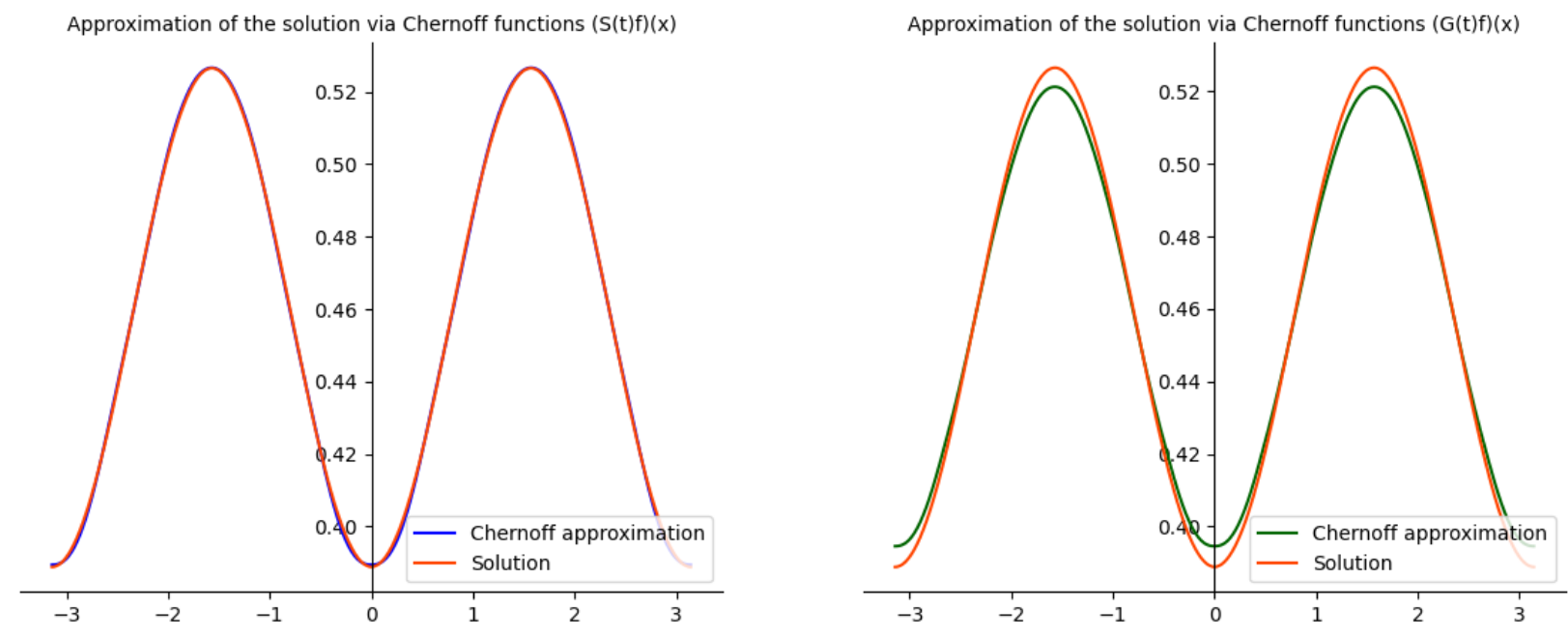}\\
\end{center}
n=10
\begin{center}
	\includegraphics[scale=0.5]{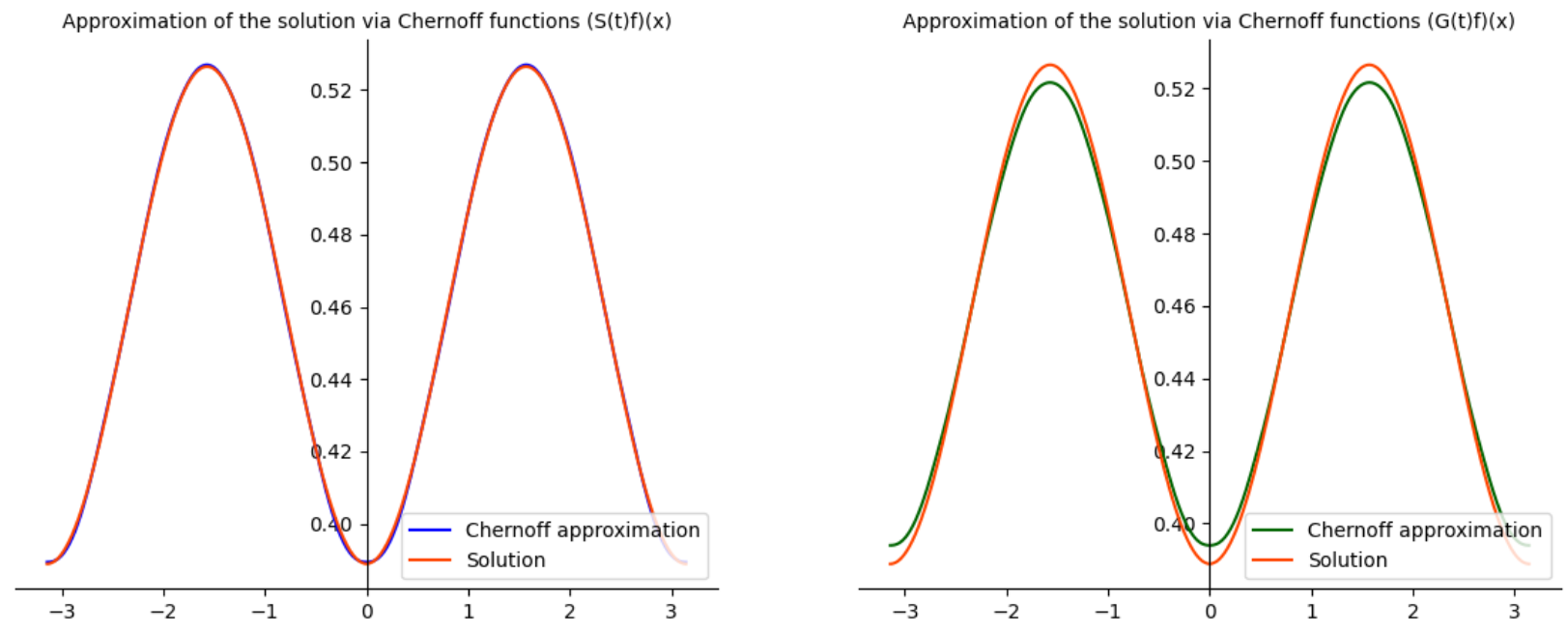}\\
\end{center}

\newpage
\subsection{$u_0(x)=|\sin(x)|^{7/2}$}
n=1
\begin{center}
	\includegraphics[scale=0.5]{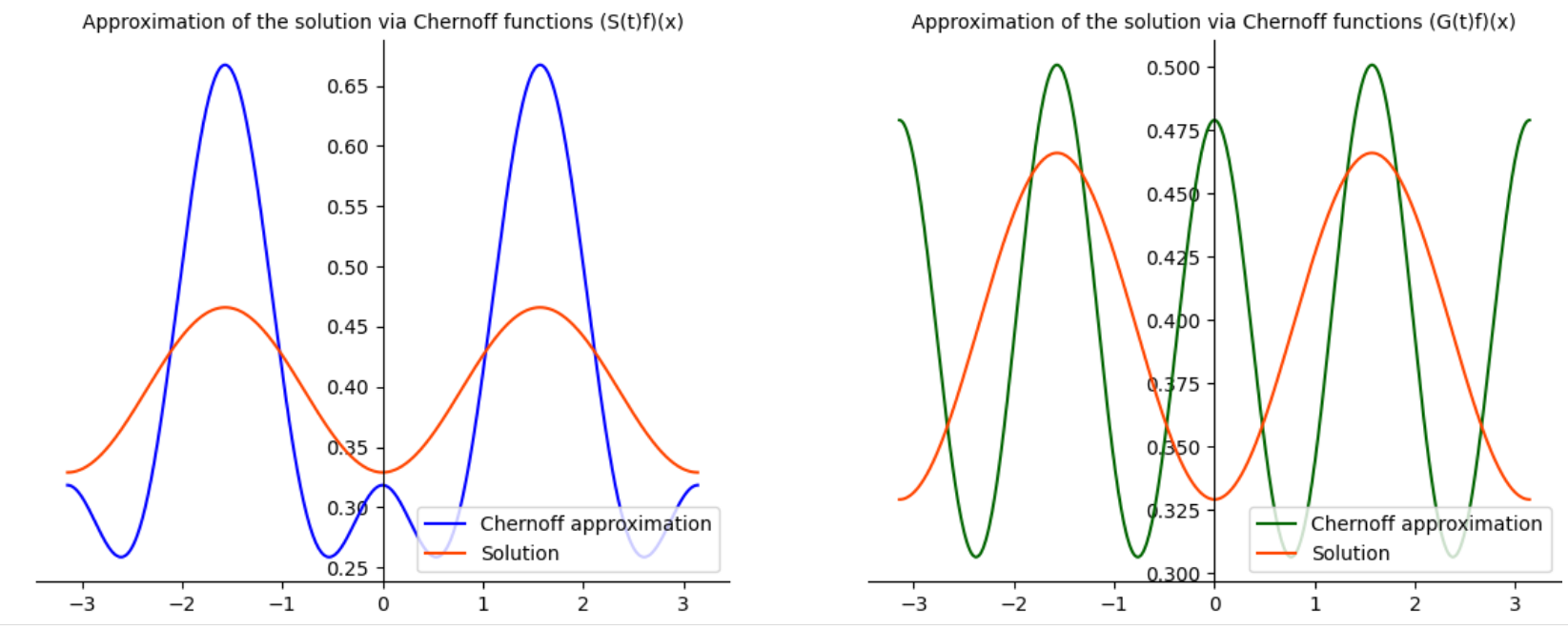}\\
\end{center}
n=2
\begin{center}
	\includegraphics[scale=0.5]{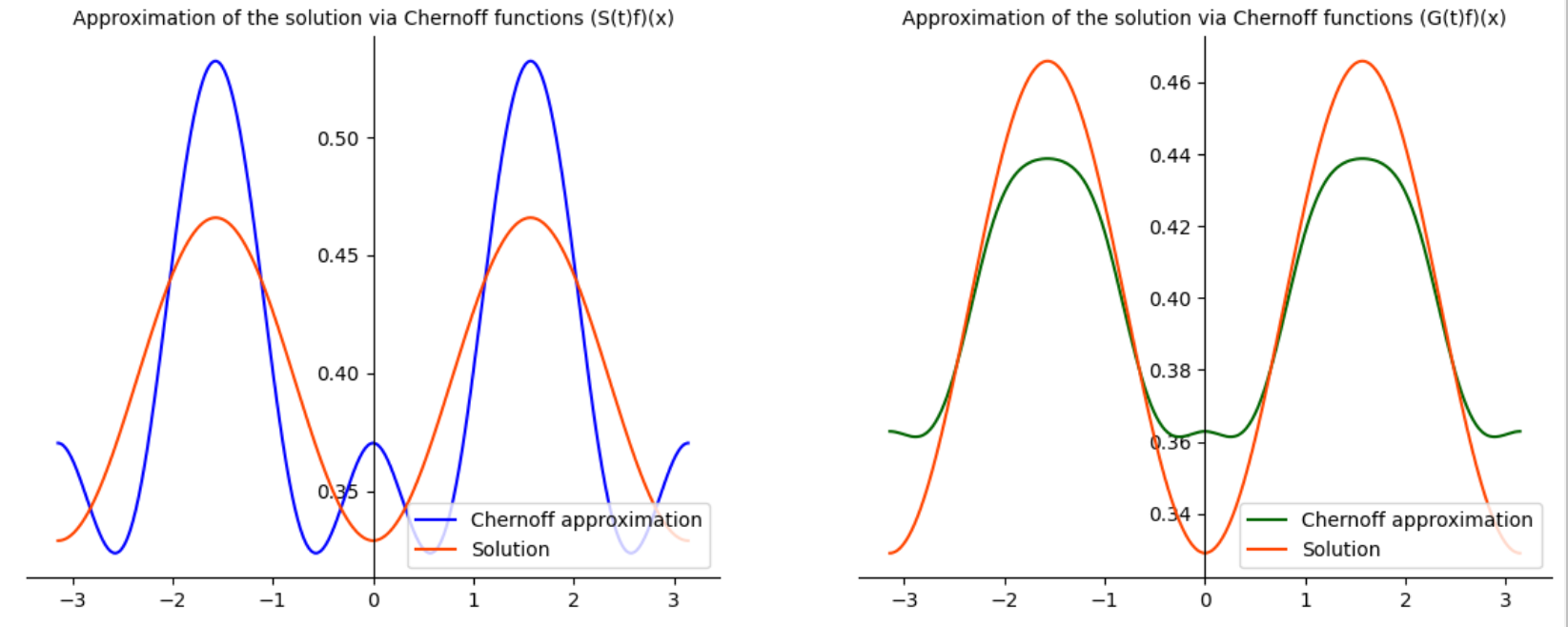}\\
\end{center}
n=3
\begin{center}
	\includegraphics[scale=0.5]{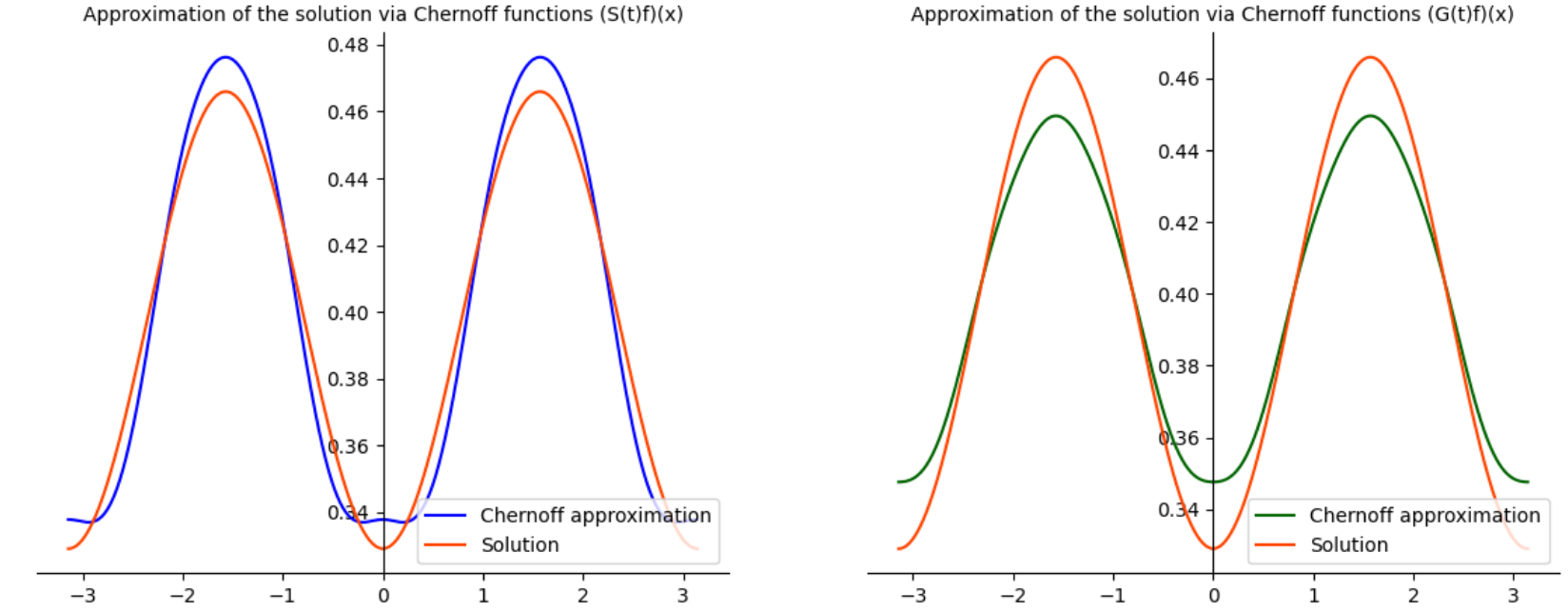}\\
\end{center}
\newpage
n=4
\begin{center}
	\includegraphics[scale=0.5]{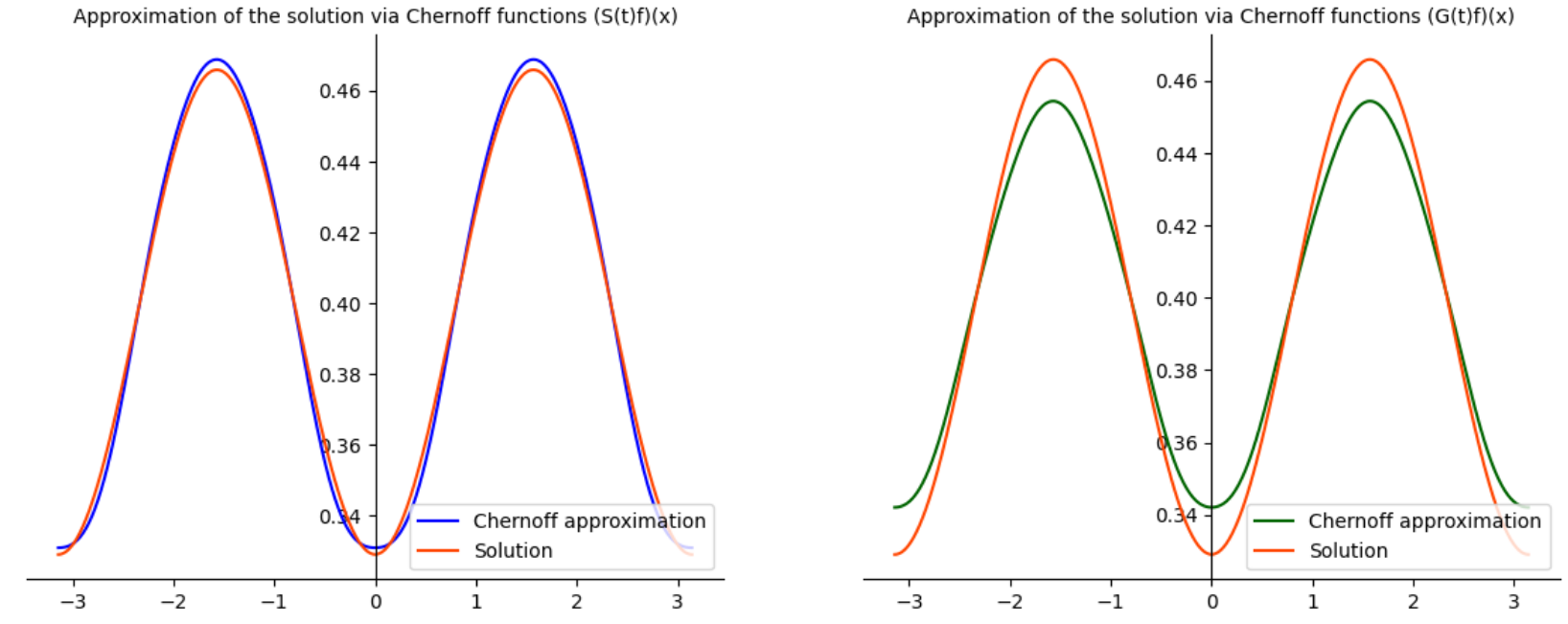}\\
\end{center}
n=5

\begin{center}
	\includegraphics[scale=0.5]{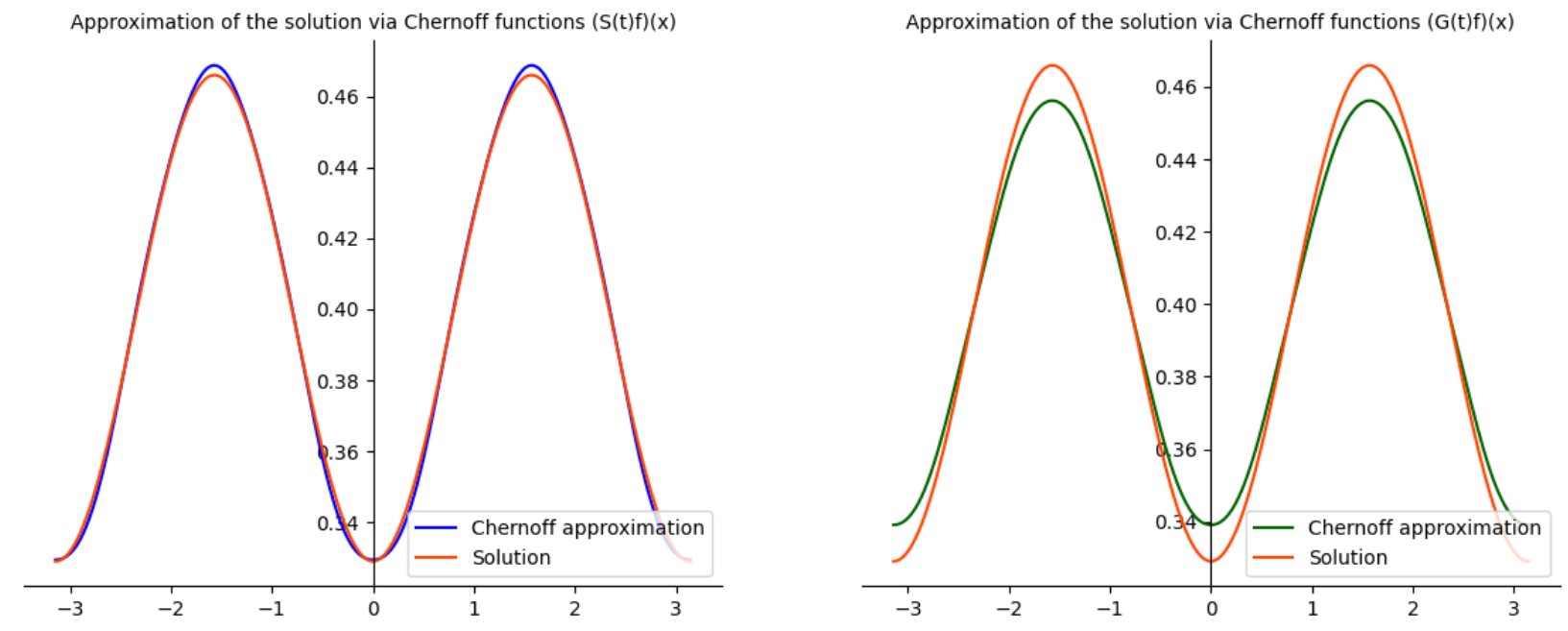}\\
\end{center}
n=6
\begin{center}
	\includegraphics[scale=0.5]{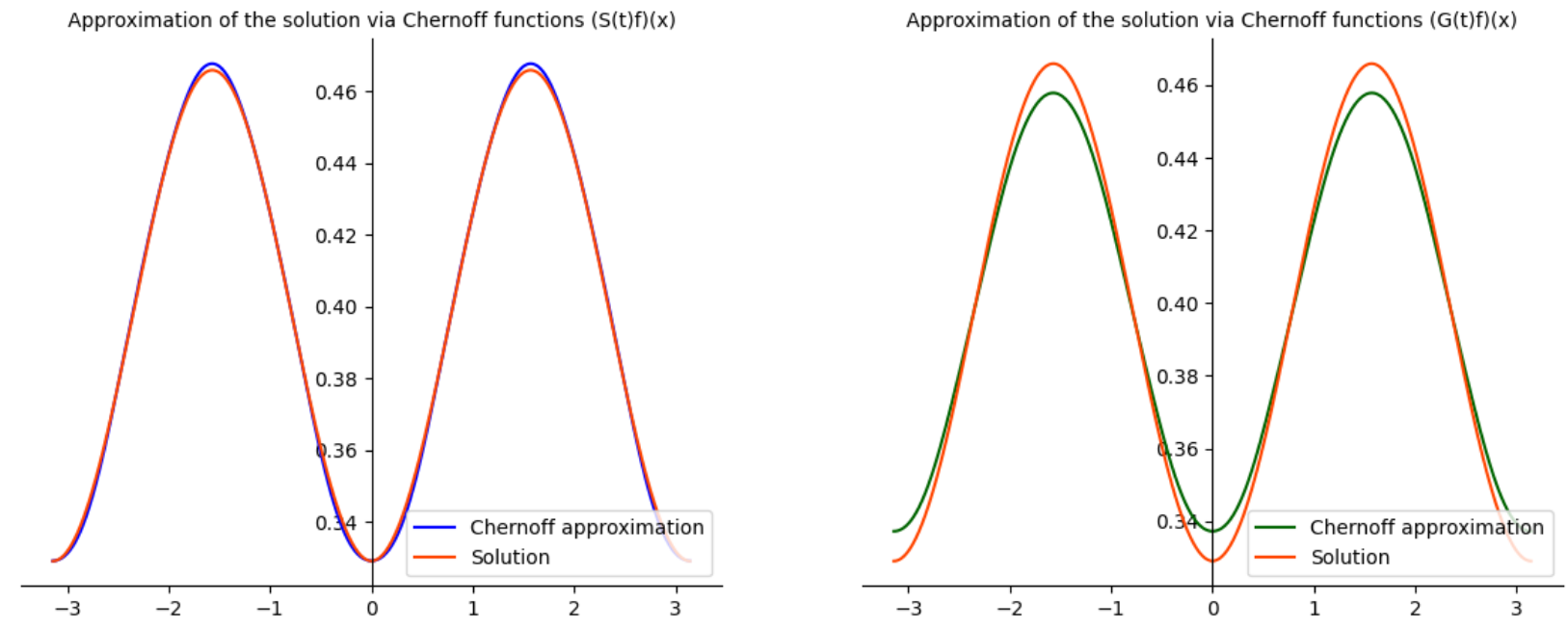}\\
\end{center}
\newpage

n=7
\begin{center}
	\includegraphics[scale=0.5]{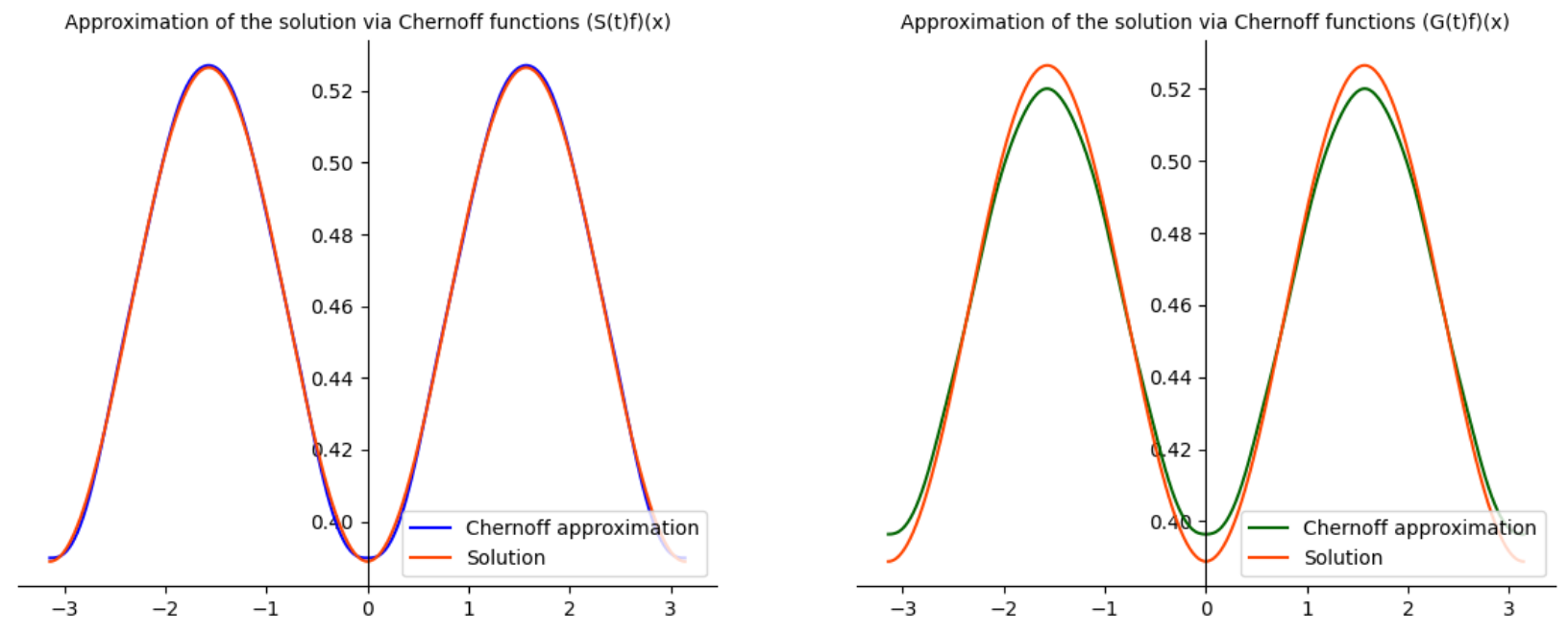}\\
\end{center}
n=8
\begin{center}
	\includegraphics[scale=0.5]{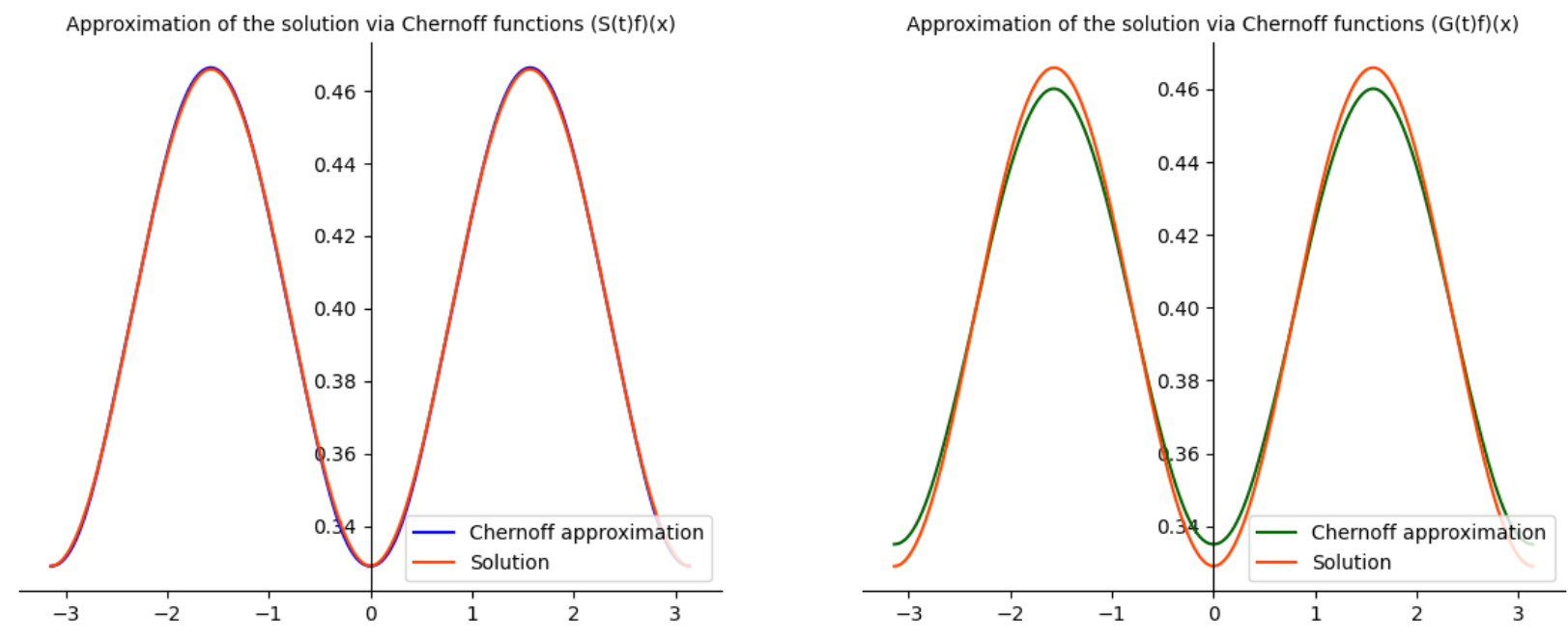}\\
\end{center}
n=9
\begin{center}
	\includegraphics[scale=0.5]{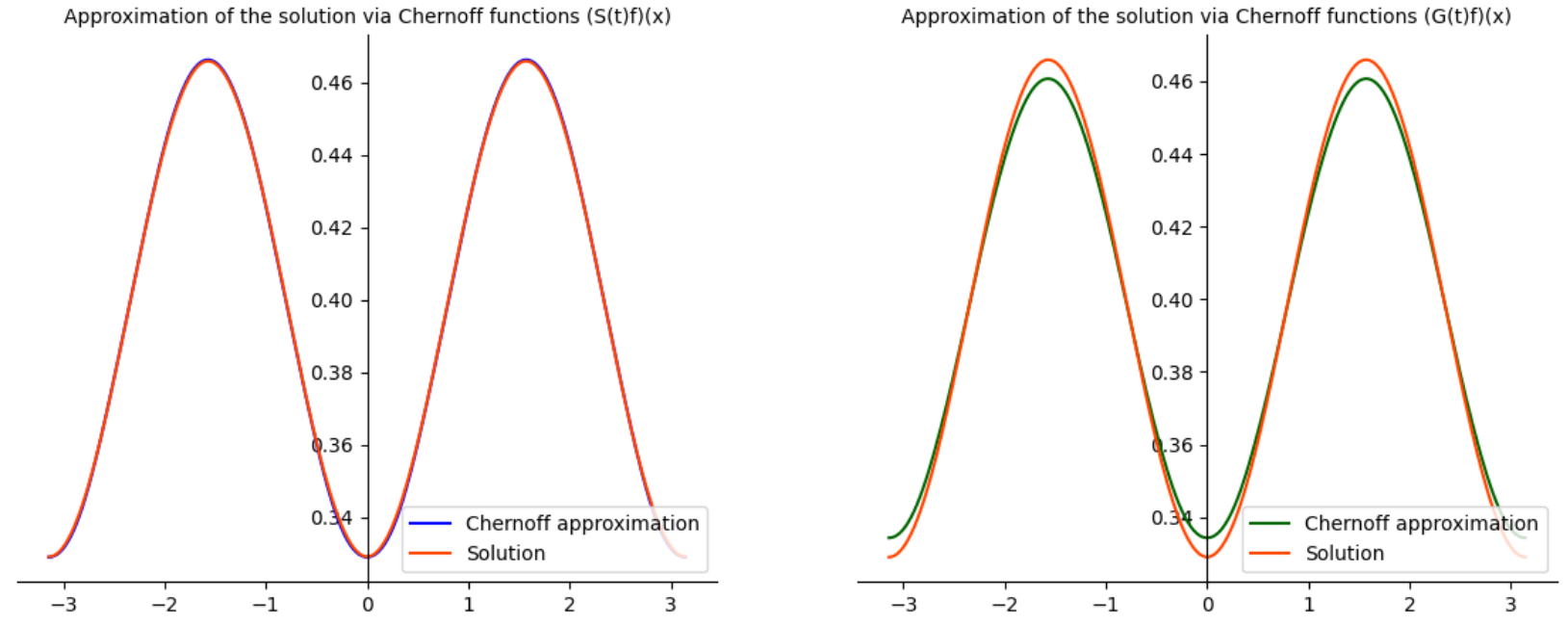}\\
\end{center}
\newpage
n=10
\begin{center}
	\includegraphics[scale=0.5]{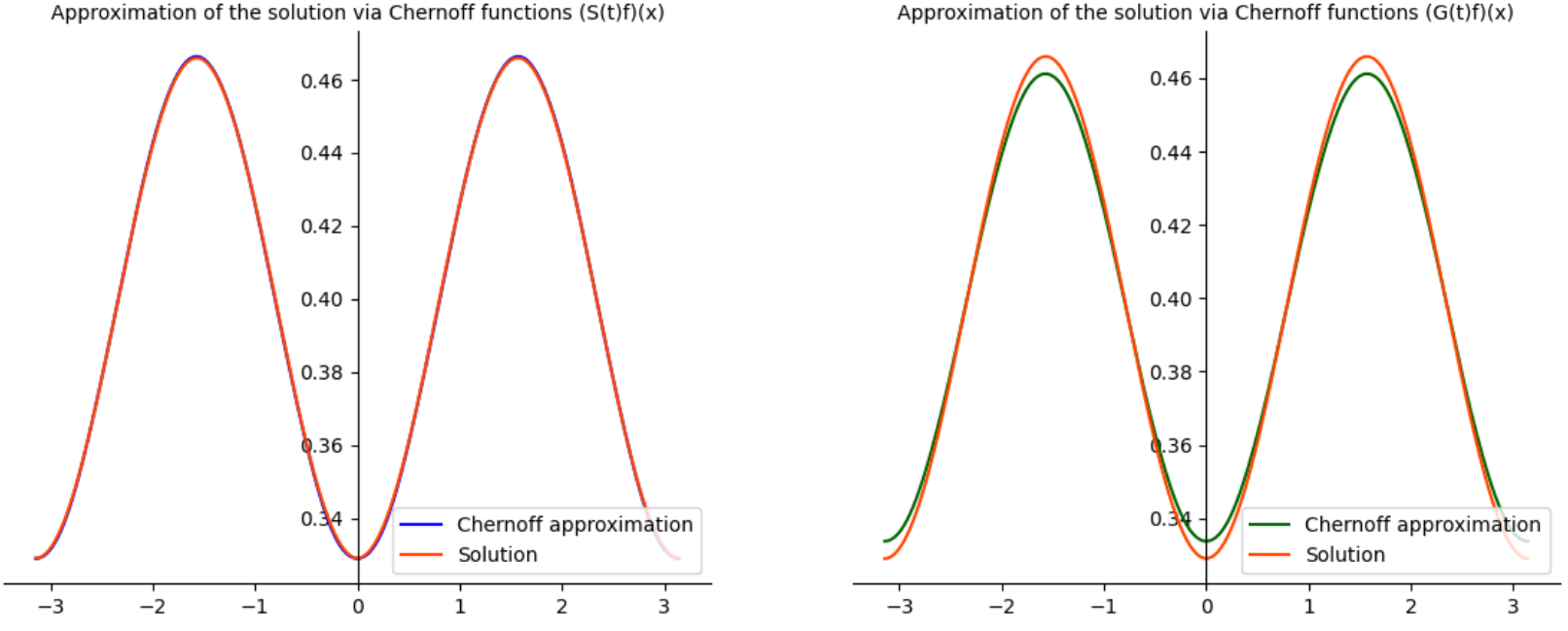}\\
\end{center}

\subsection{$u_0(x)=|\sin(x)|^{9/2}$}

n=1
\begin{center}
	\includegraphics[scale=0.5]{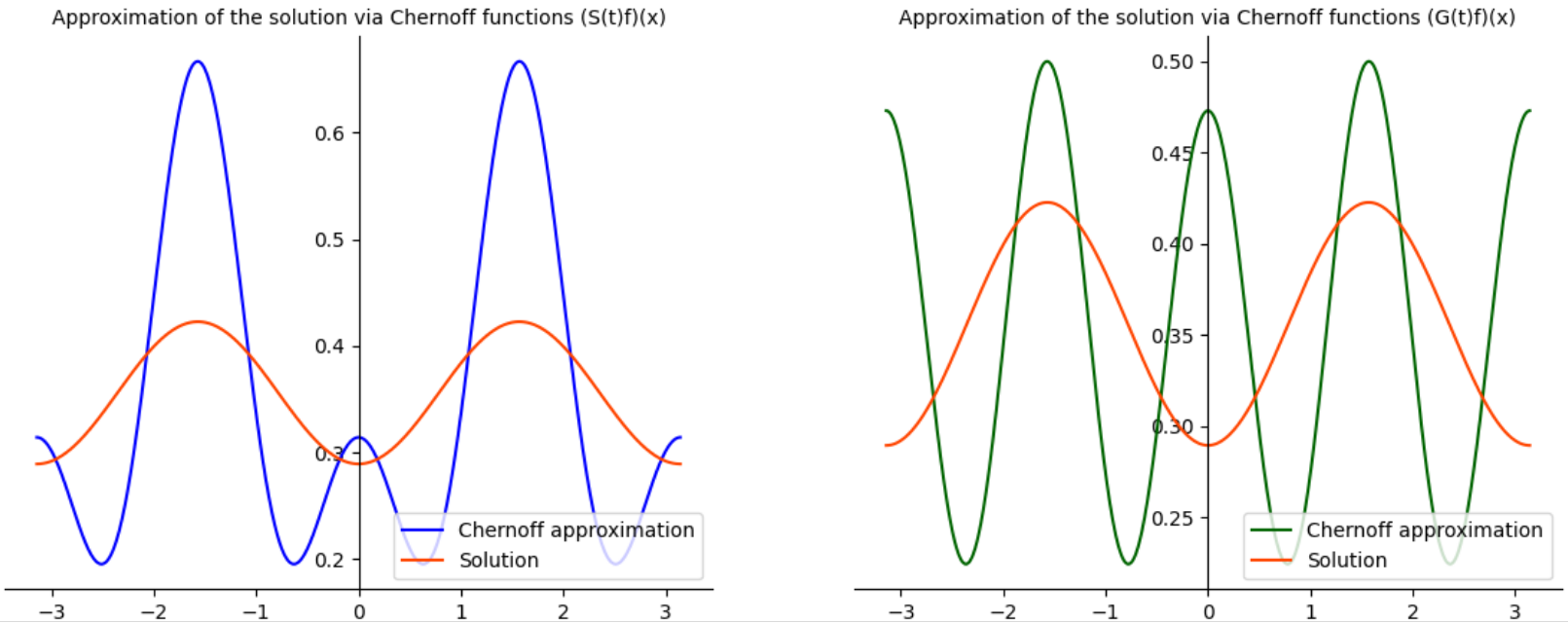}\\
\end{center}
n=2
\begin{center}
	\includegraphics[scale=0.5]{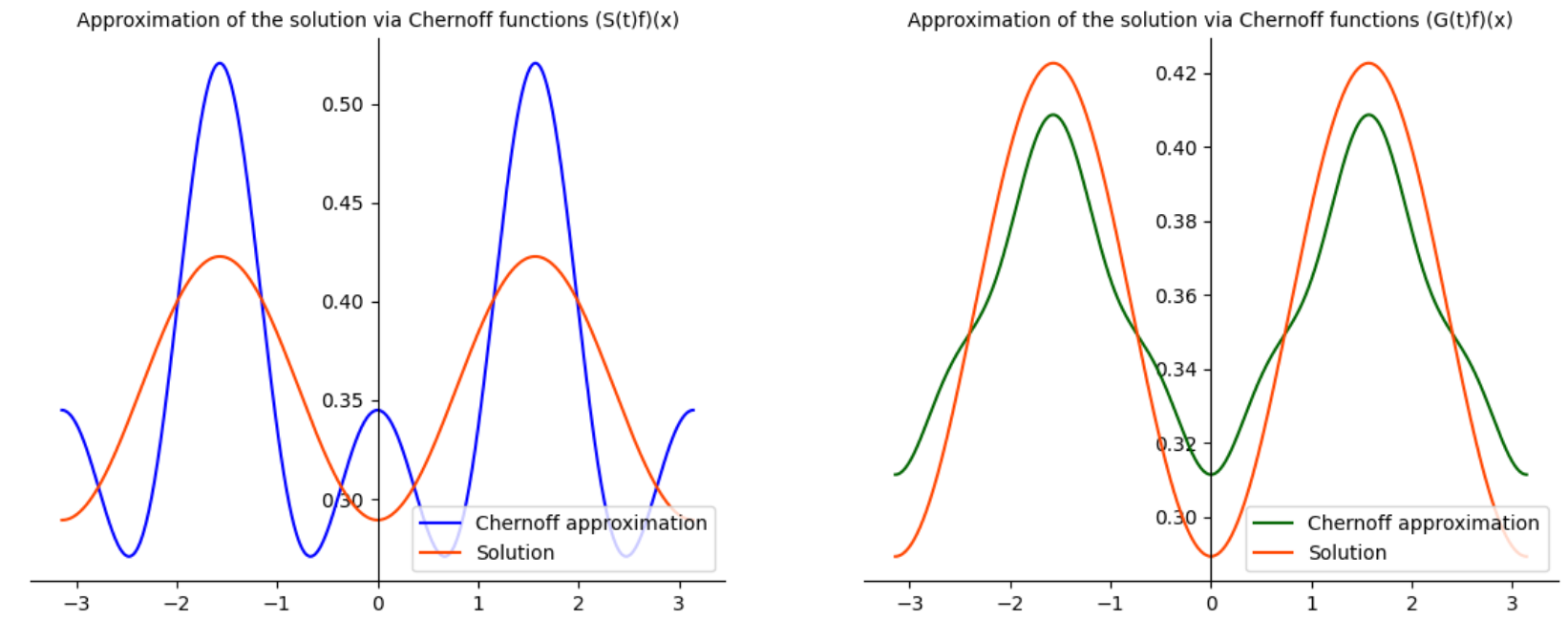}\\
\end{center}
\newpage
n=3
\begin{center}
	\includegraphics[scale=0.5]{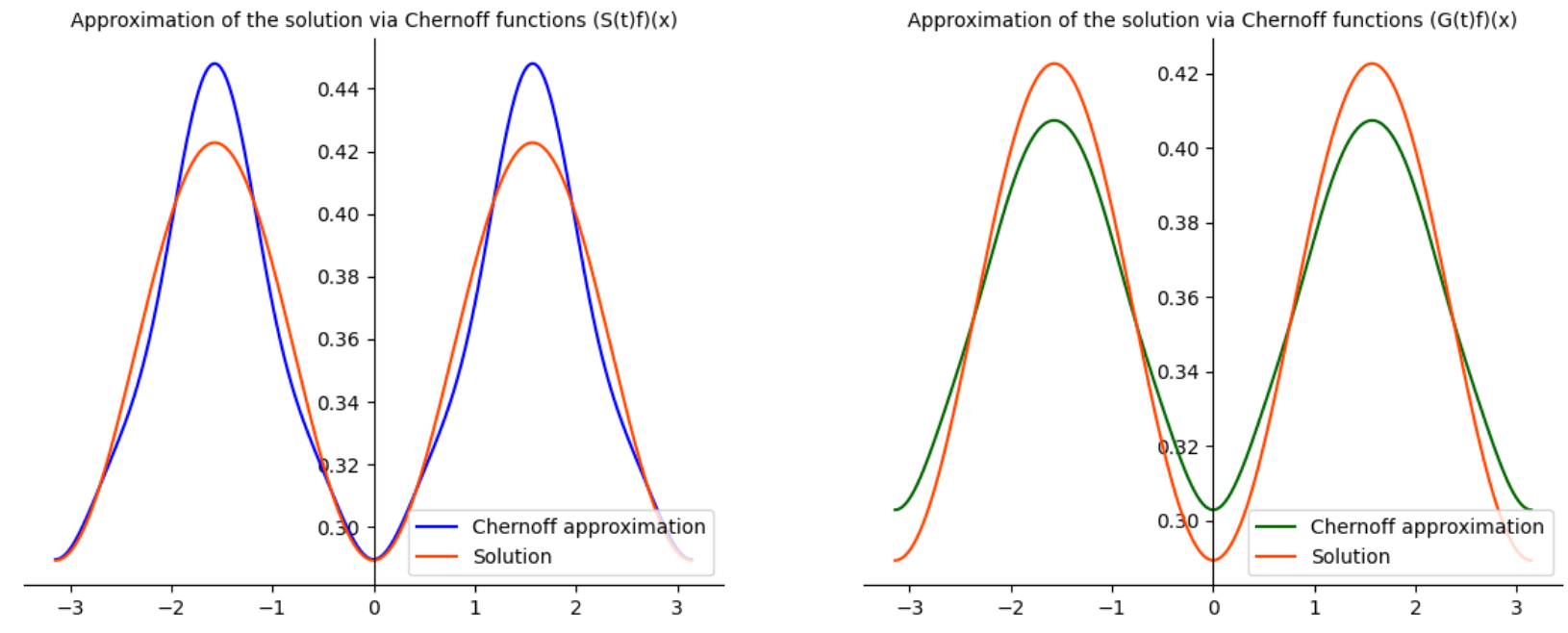}\\
\end{center}
n=4
\begin{center}
	\includegraphics[scale=0.5]{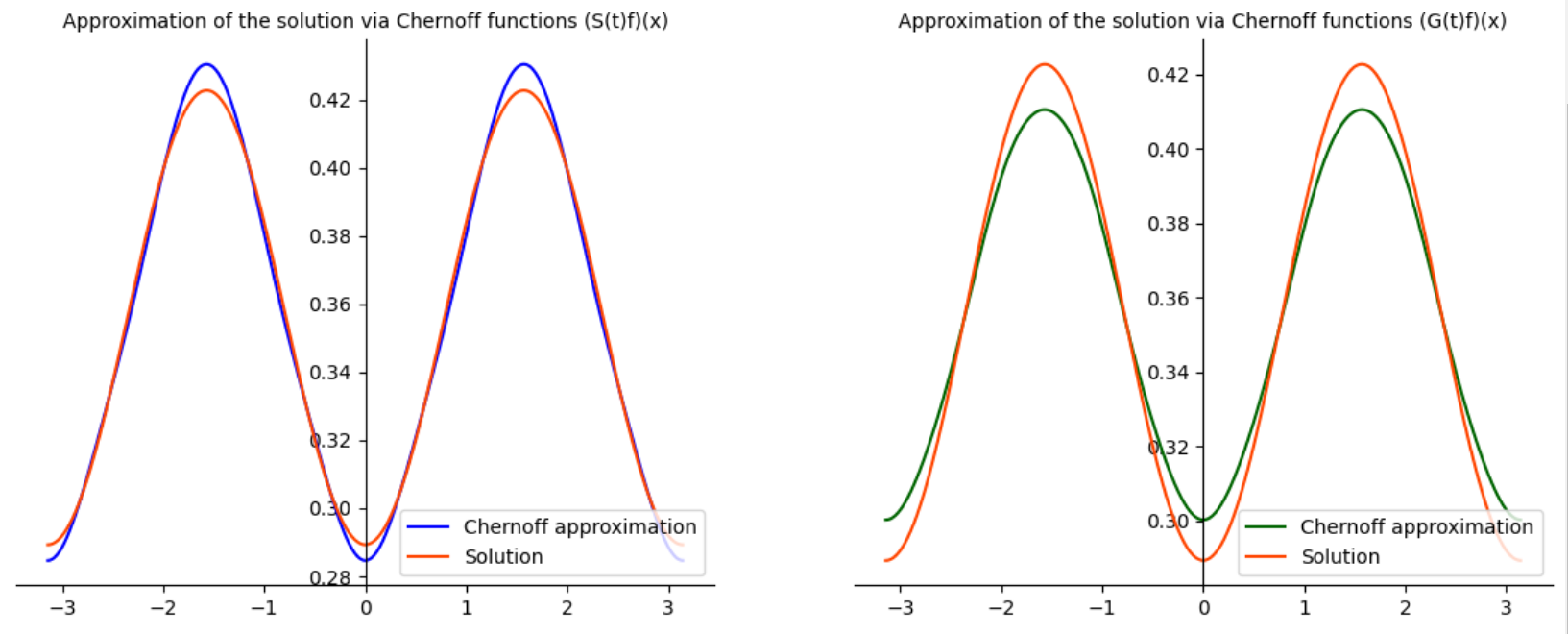}\\
\end{center}
n=5
\begin{center}
	\includegraphics[scale=0.5]{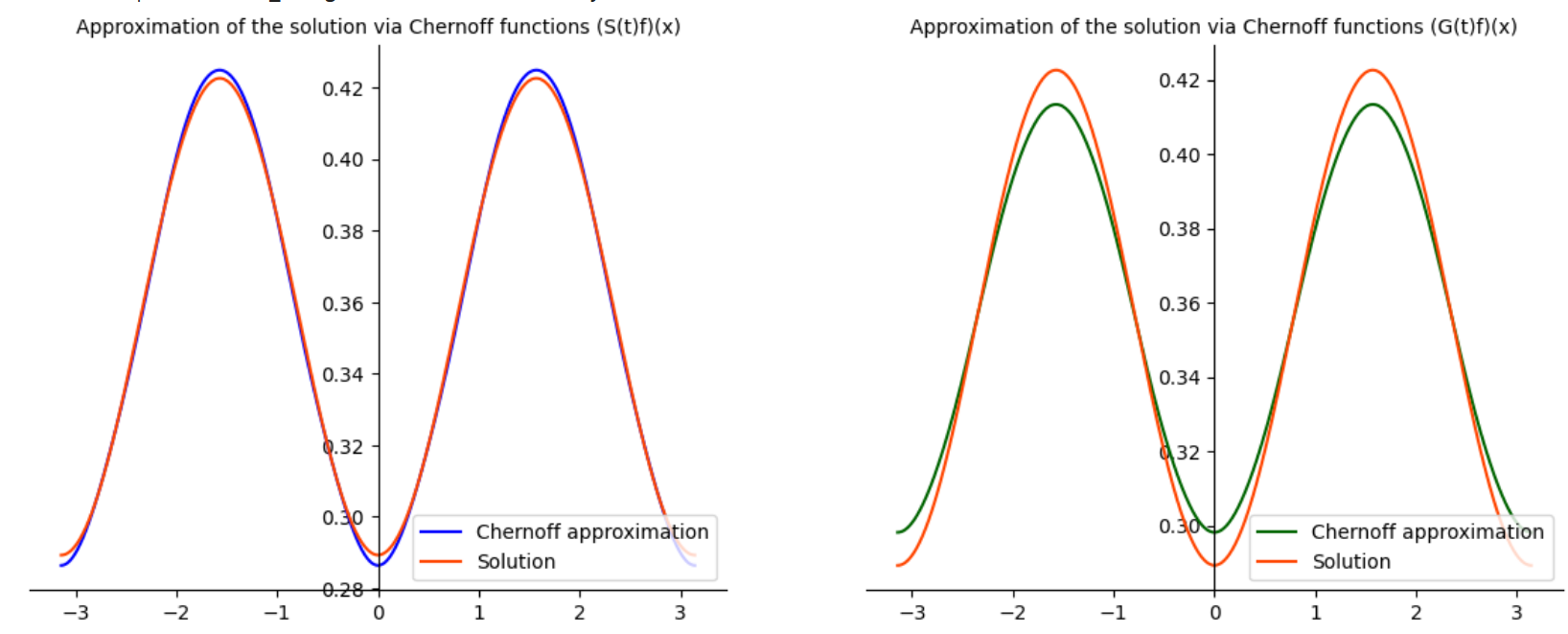}\\
\end{center}
\newpage
n=6
\begin{center}
	\includegraphics[scale=0.5]{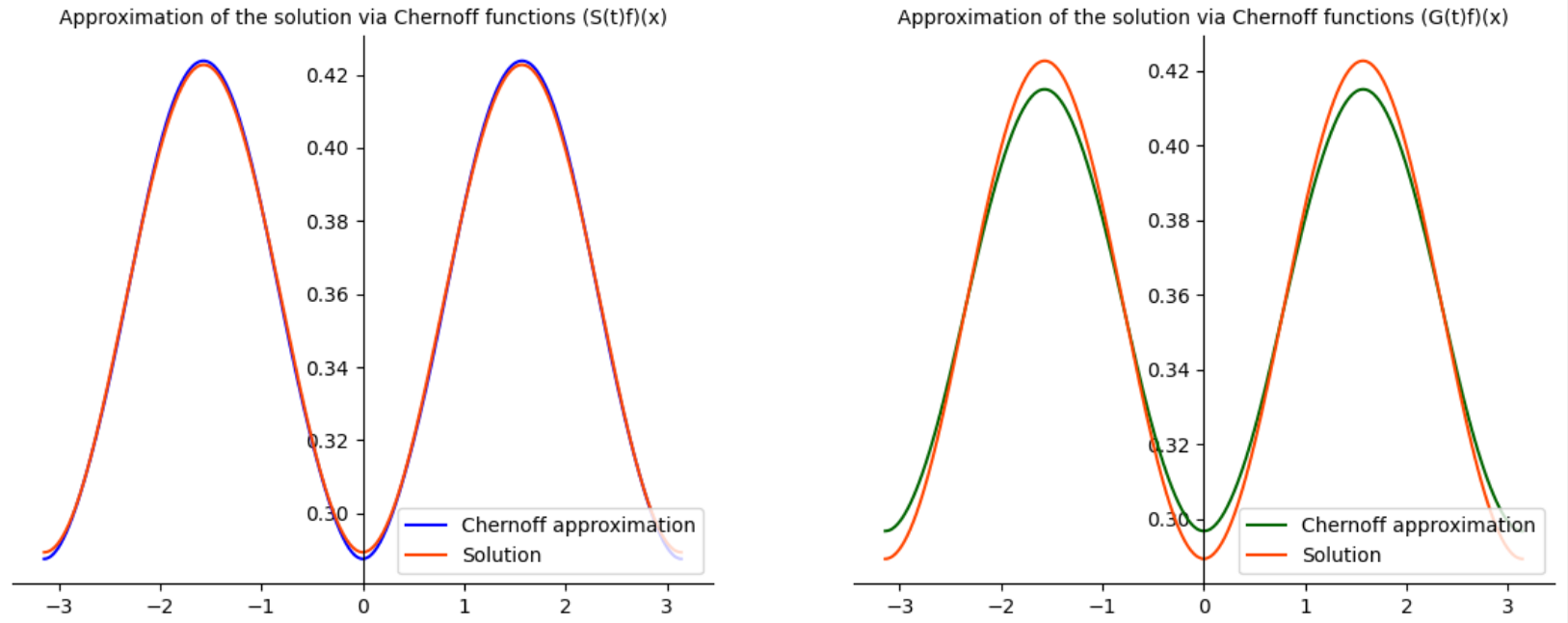}\\
\end{center}
n=7
\begin{center}
	\includegraphics[scale=0.5]{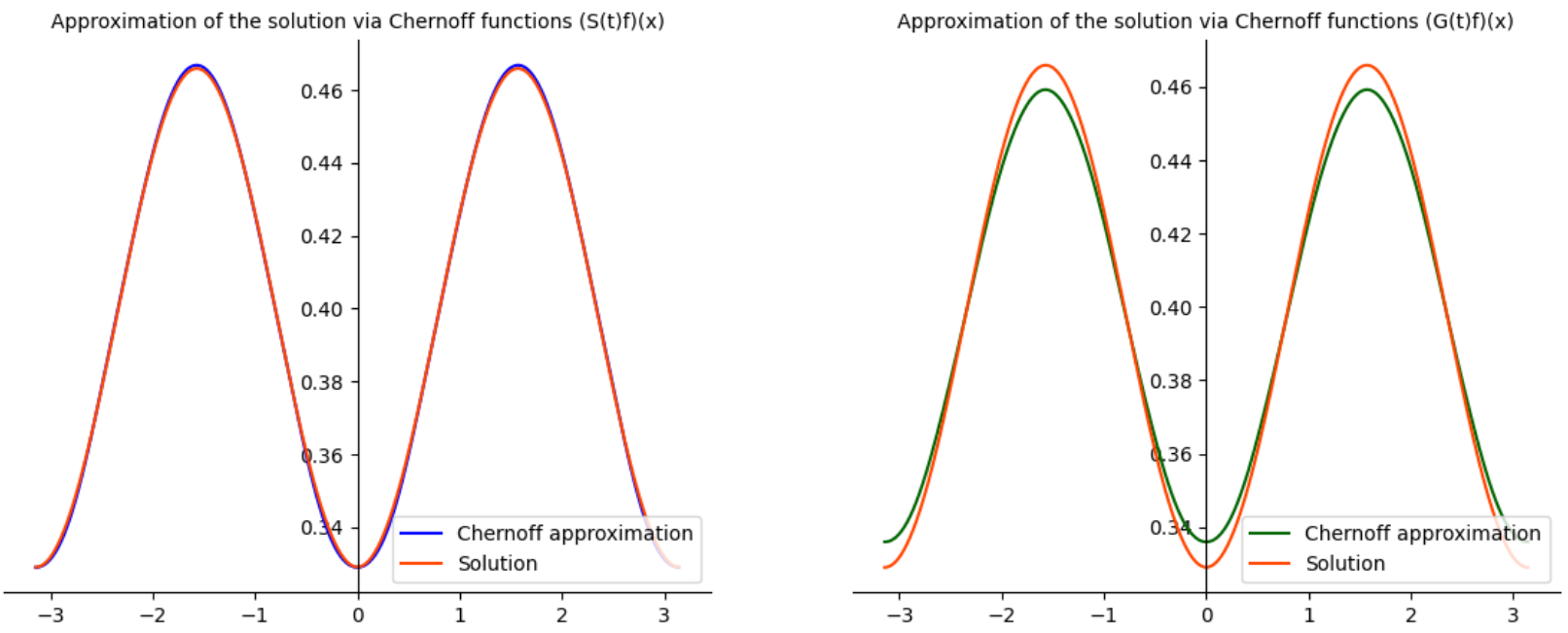}\\
\end{center}
n=8
\begin{center}
	\includegraphics[scale=0.5]{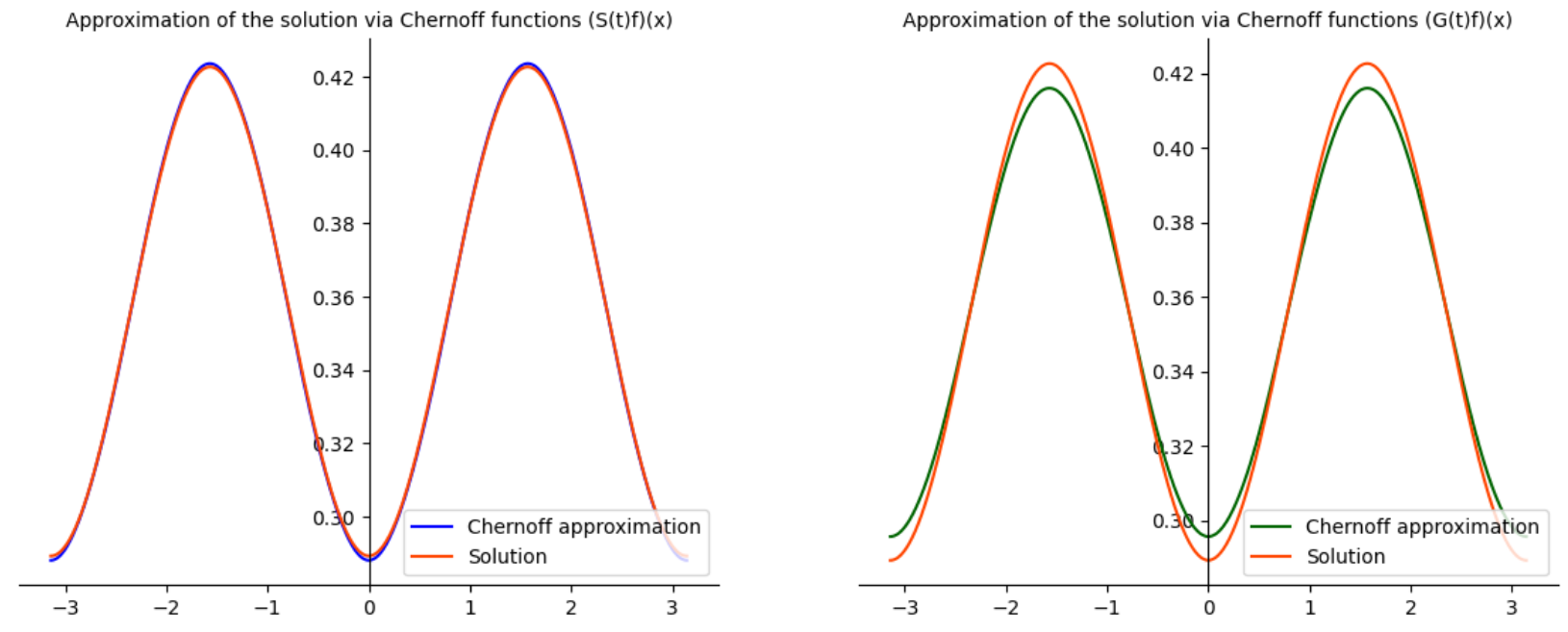}\\
\end{center}
\newpage
n=9
\begin{center}
	\includegraphics[scale=0.5]{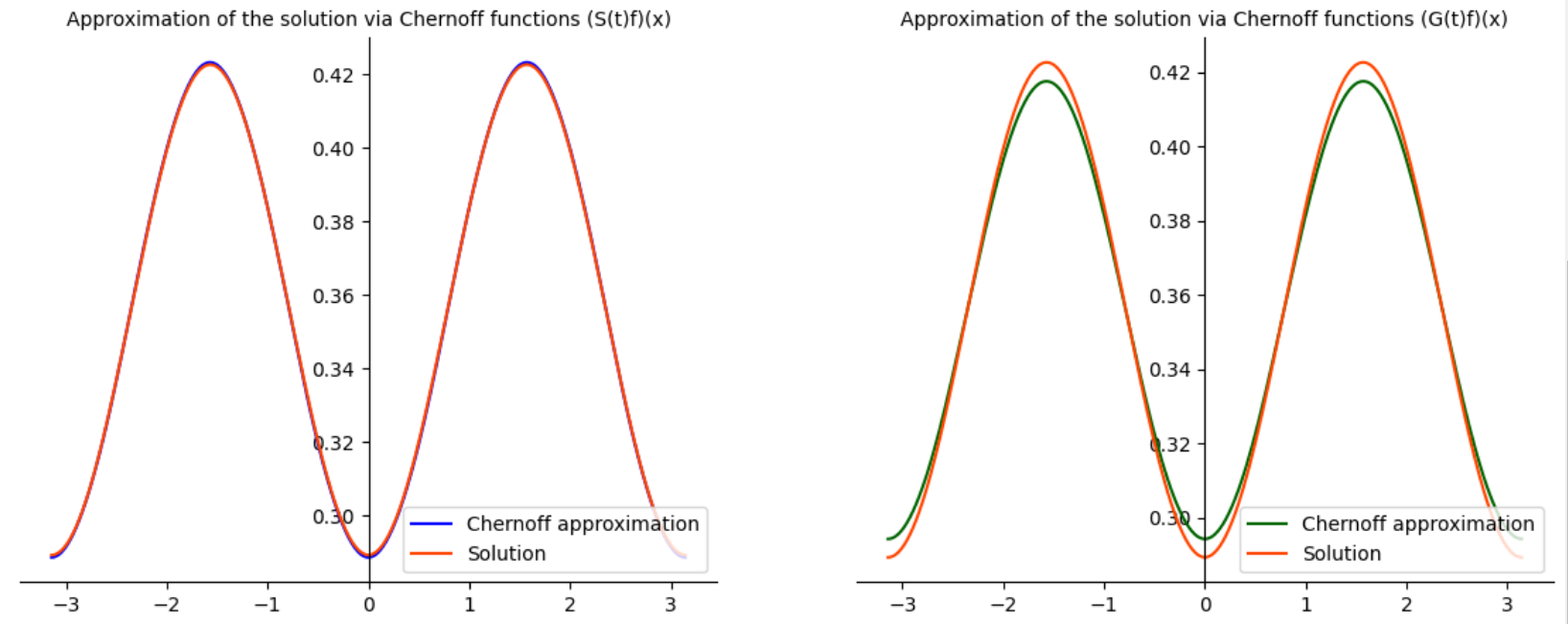}\\
\end{center}
n=10
\begin{center}
	\includegraphics[scale=0.5]{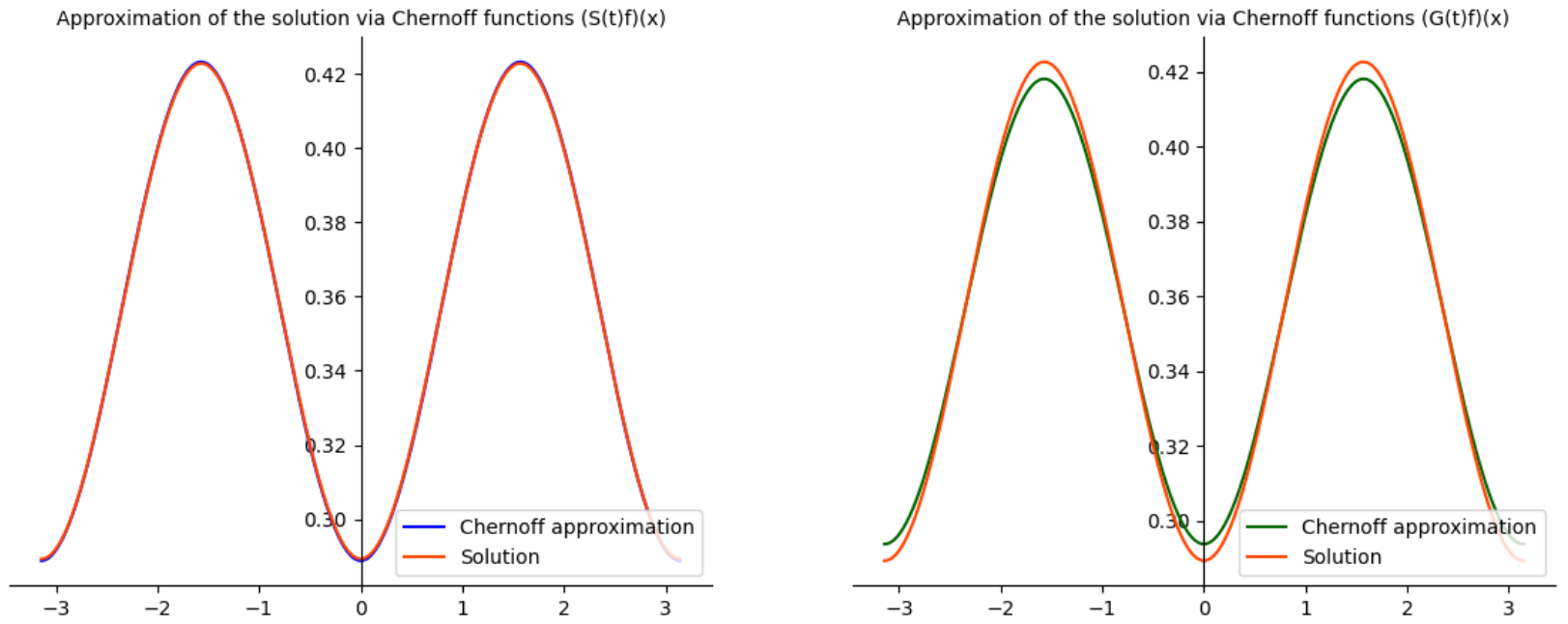}\\
\end{center}

\subsection{$u_0(x)=|\sin(x)|\cdot\sin(x)$}

n=1
\begin{center}
	\includegraphics[scale=0.5]{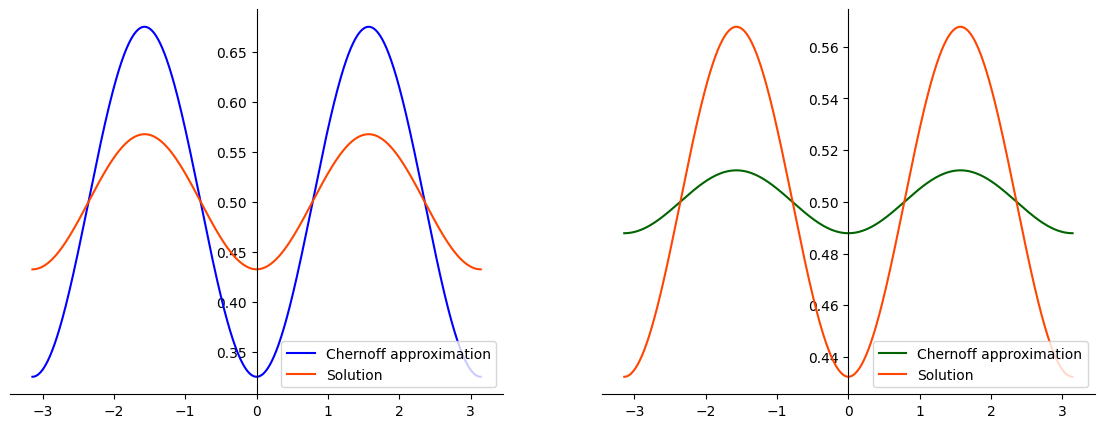}\\
\end{center}
n=2
\begin{center}
	\includegraphics[scale=0.5]{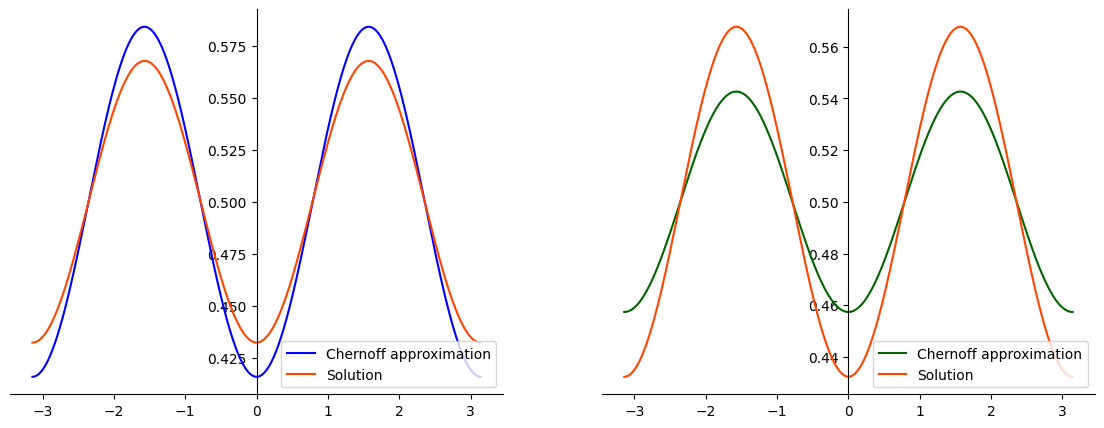}\\
\end{center}

n=3
\begin{center}
	\includegraphics[scale=0.5]{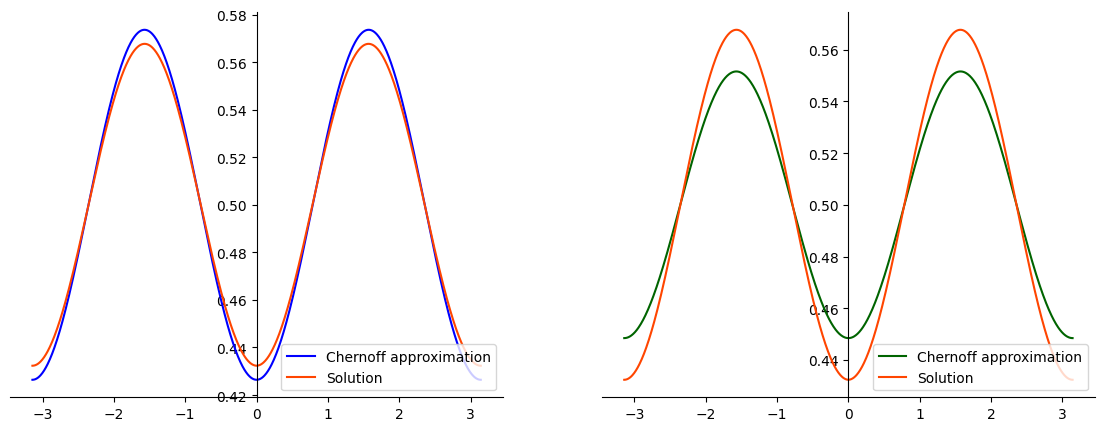}\\
\end{center}
n=4
\begin{center}
	\includegraphics[scale=0.5]{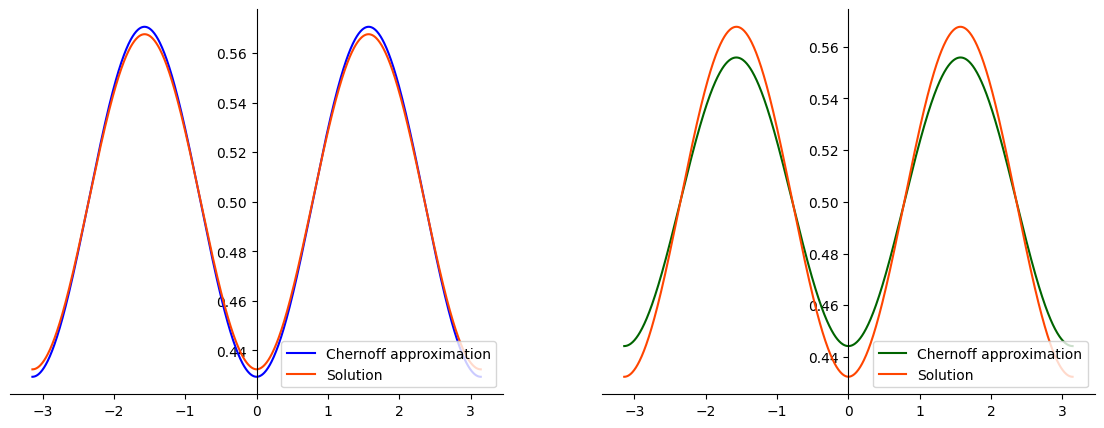}\\
\end{center}
n=5
\begin{center}
	\includegraphics[scale=0.5]{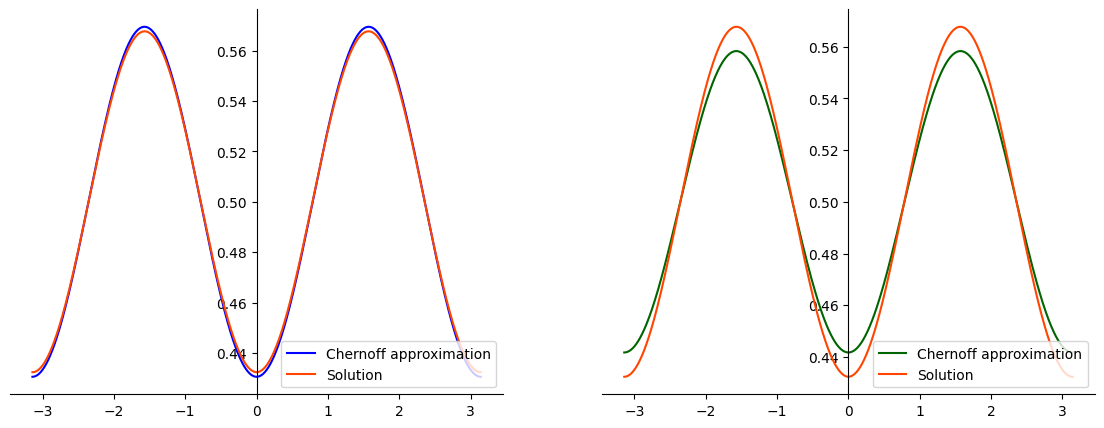}\\
\end{center}

n=6
\begin{center}
	\includegraphics[scale=0.5]{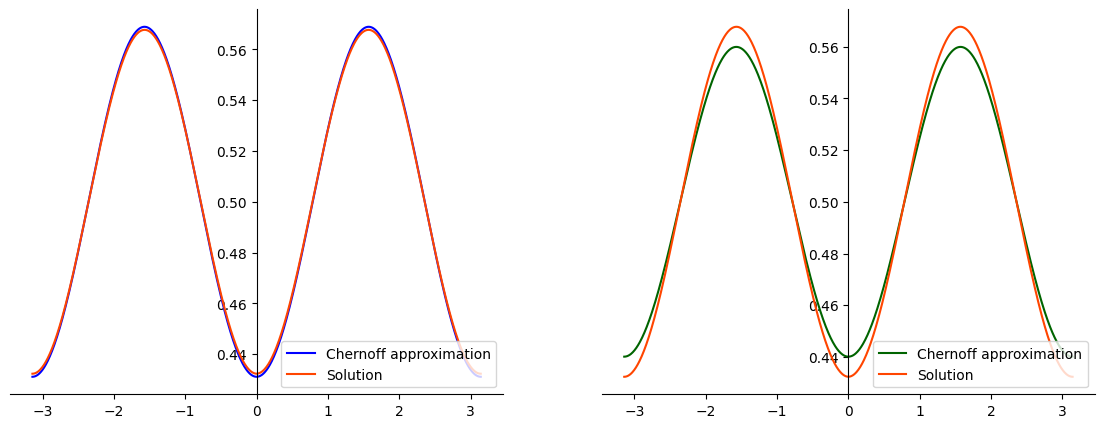}\\
\end{center}
n=7
\begin{center}
	\includegraphics[scale=0.5]{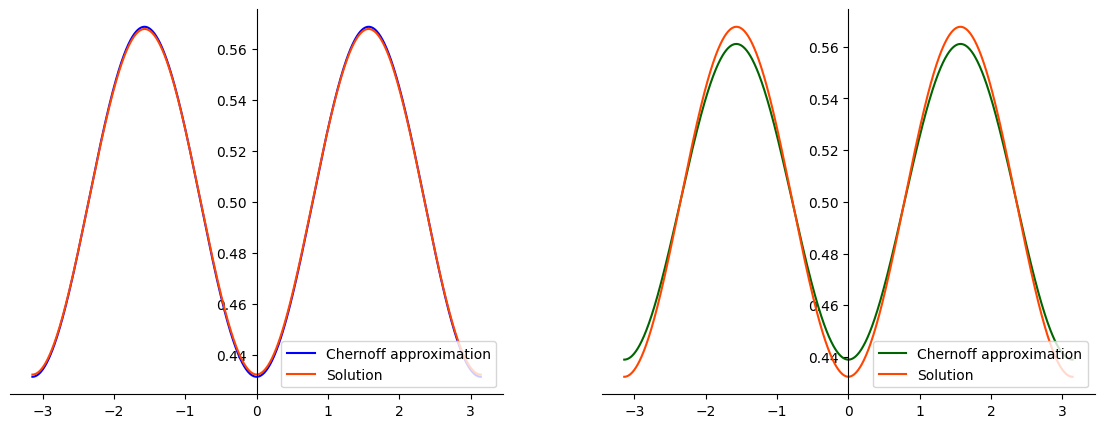}\\
\end{center}
n=8
\begin{center}
	\includegraphics[scale=0.5]{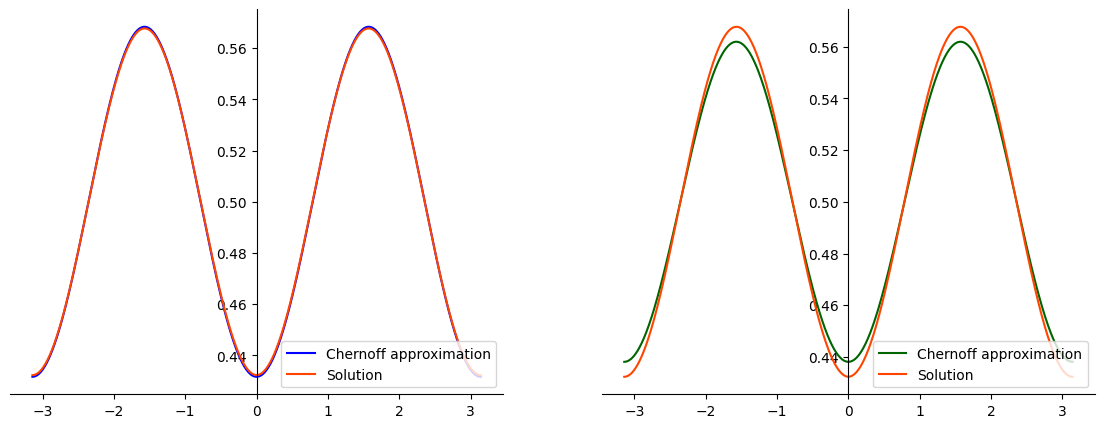}\\
\end{center}

n=9
\begin{center}
	\includegraphics[scale=0.5]{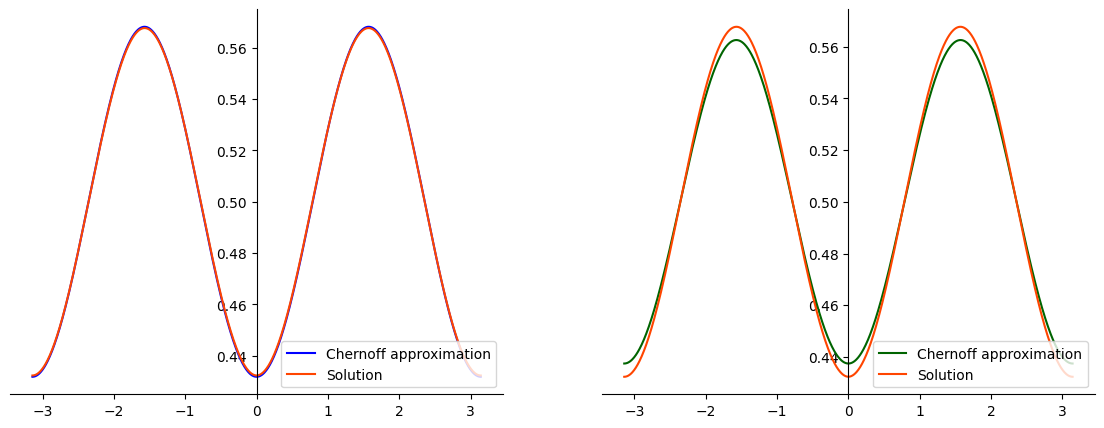}\\
\end{center}
n=10
\begin{center}
	\includegraphics[scale=0.5]{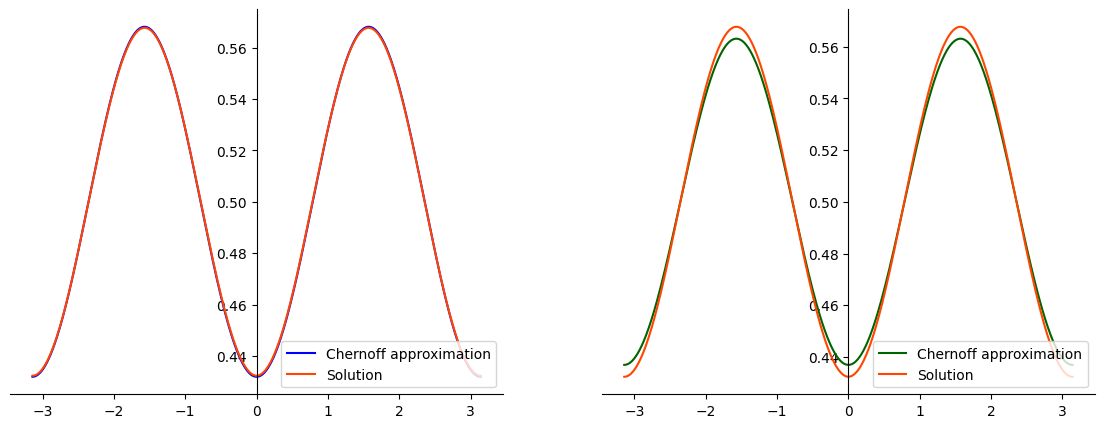}\\
\end{center}

\subsection{$u_0(x)=\sin(x)\cdot|\sin(x)|^{2}$}

n=1
\begin{center}
	\includegraphics[scale=0.5]{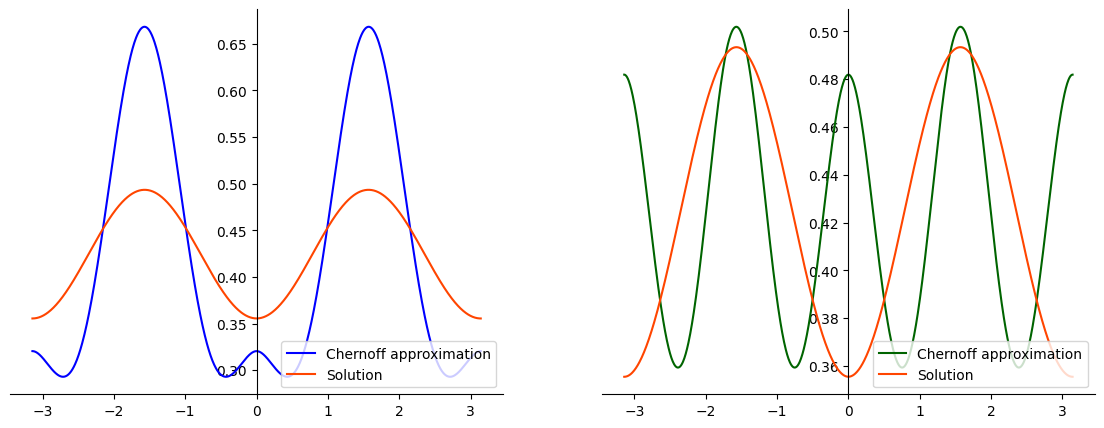}\\
\end{center}
n=2
\begin{center}
	\includegraphics[scale=0.5]{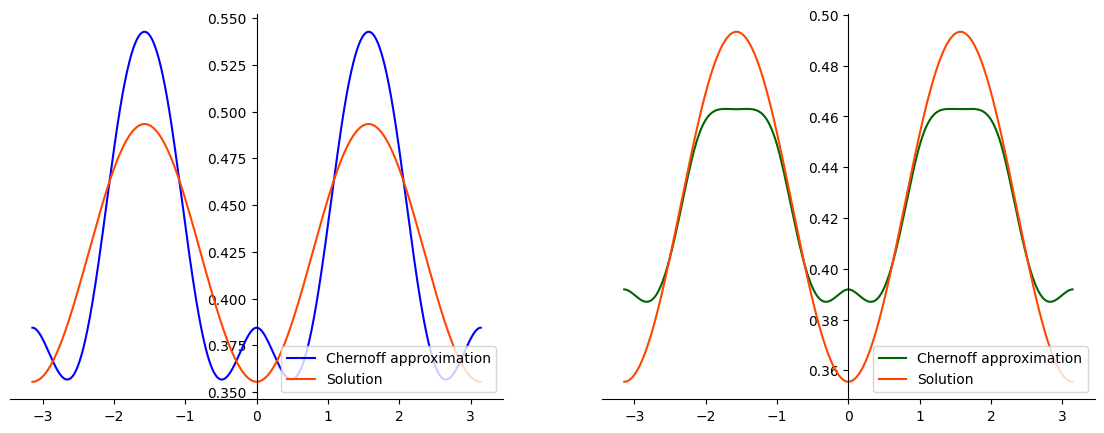}\\
\end{center}

n=3
\begin{center}
	\includegraphics[scale=0.5]{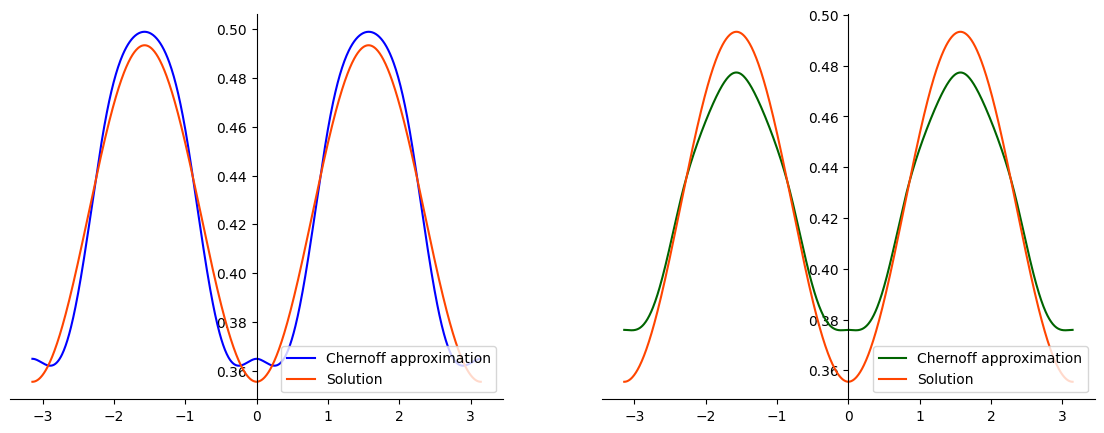}\\
\end{center}
n=4
\begin{center}
	\includegraphics[scale=0.5]{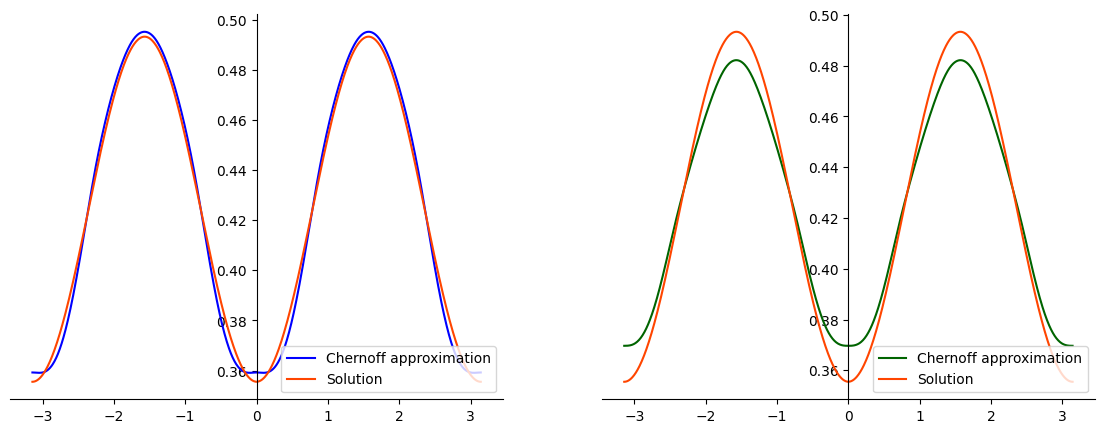}\\
\end{center}
n=5
\begin{center}
	\includegraphics[scale=0.5]{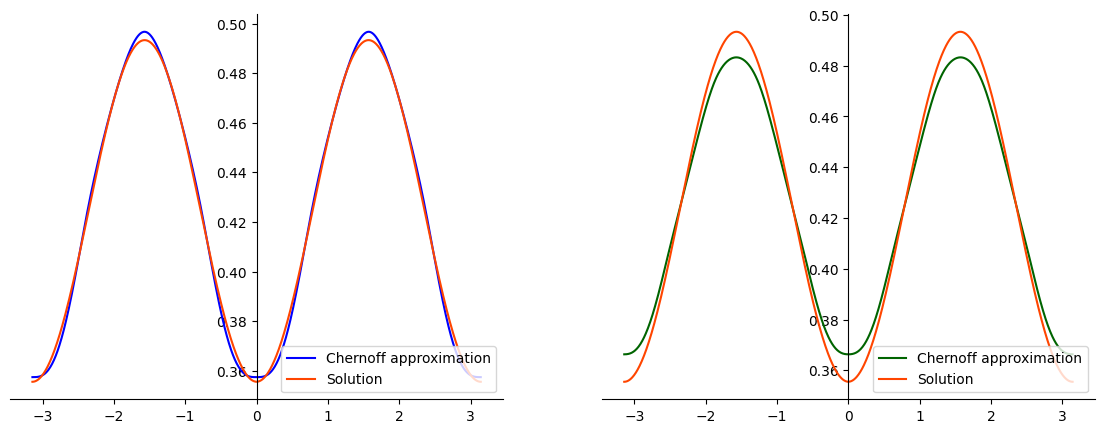}\\
\end{center}

n=6
\begin{center}
	\includegraphics[scale=0.5]{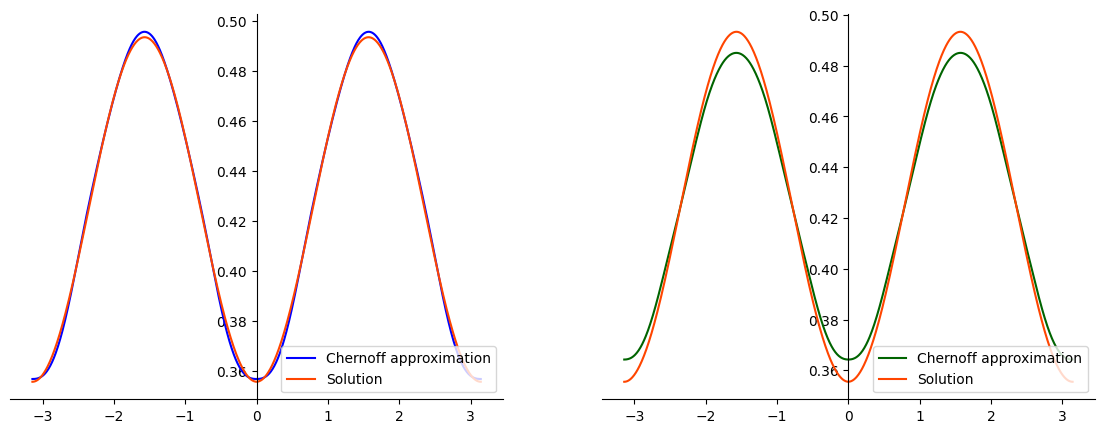}\\
\end{center}
n=7
\begin{center}
	\includegraphics[scale=0.5]{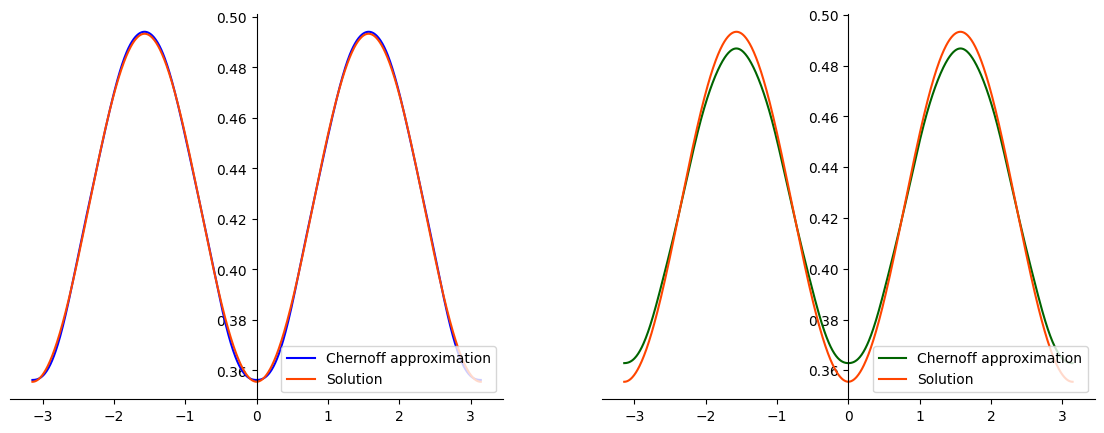}\\
\end{center}
n=8
\begin{center}
	\includegraphics[scale=0.5]{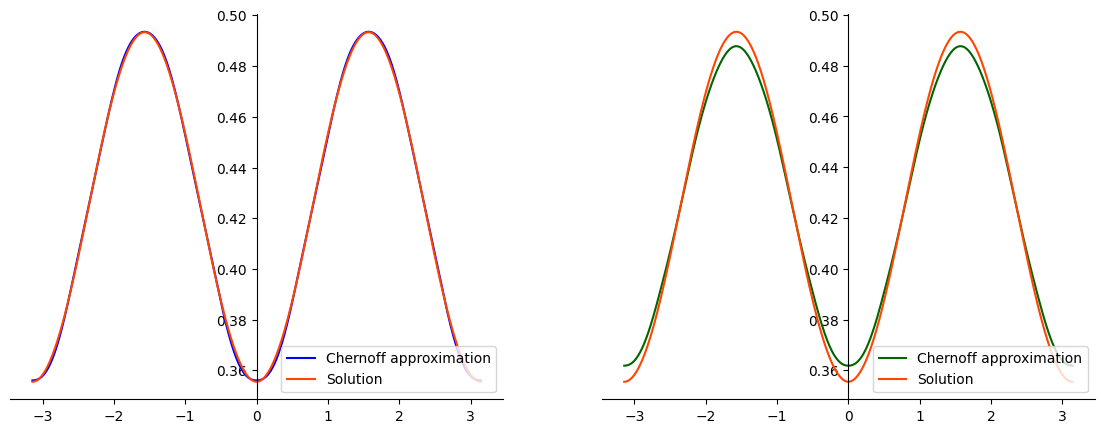}\\
\end{center}

n=9
\begin{center}
	\includegraphics[scale=0.5]{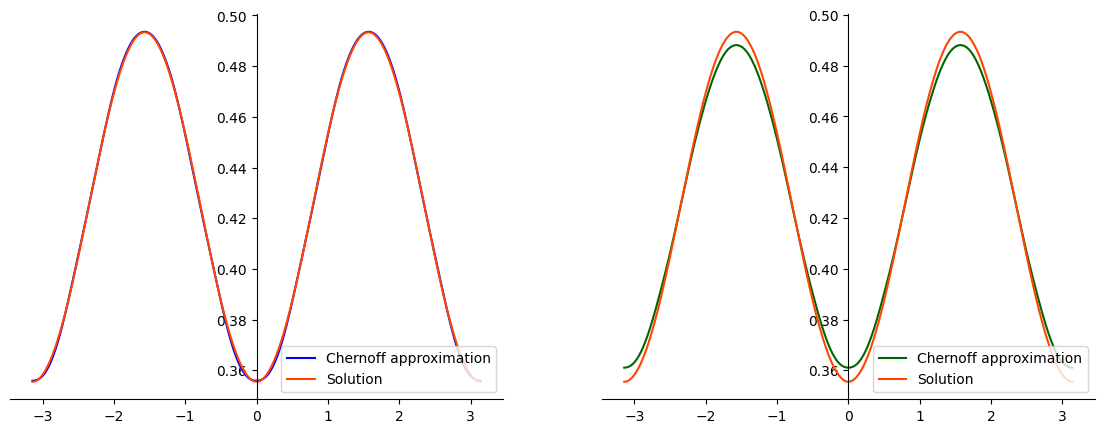}\\
\end{center}
n=10
\begin{center}
	\includegraphics[scale=0.5]{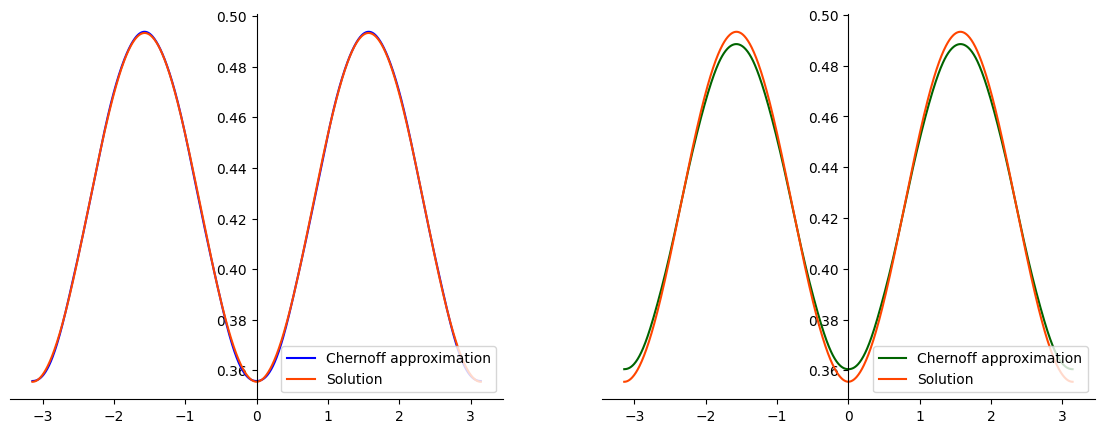}\\
\end{center}

\subsection{$u_0(x)=e^{-|x|}$}

n=1
\begin{center}
	\includegraphics[scale=0.5]{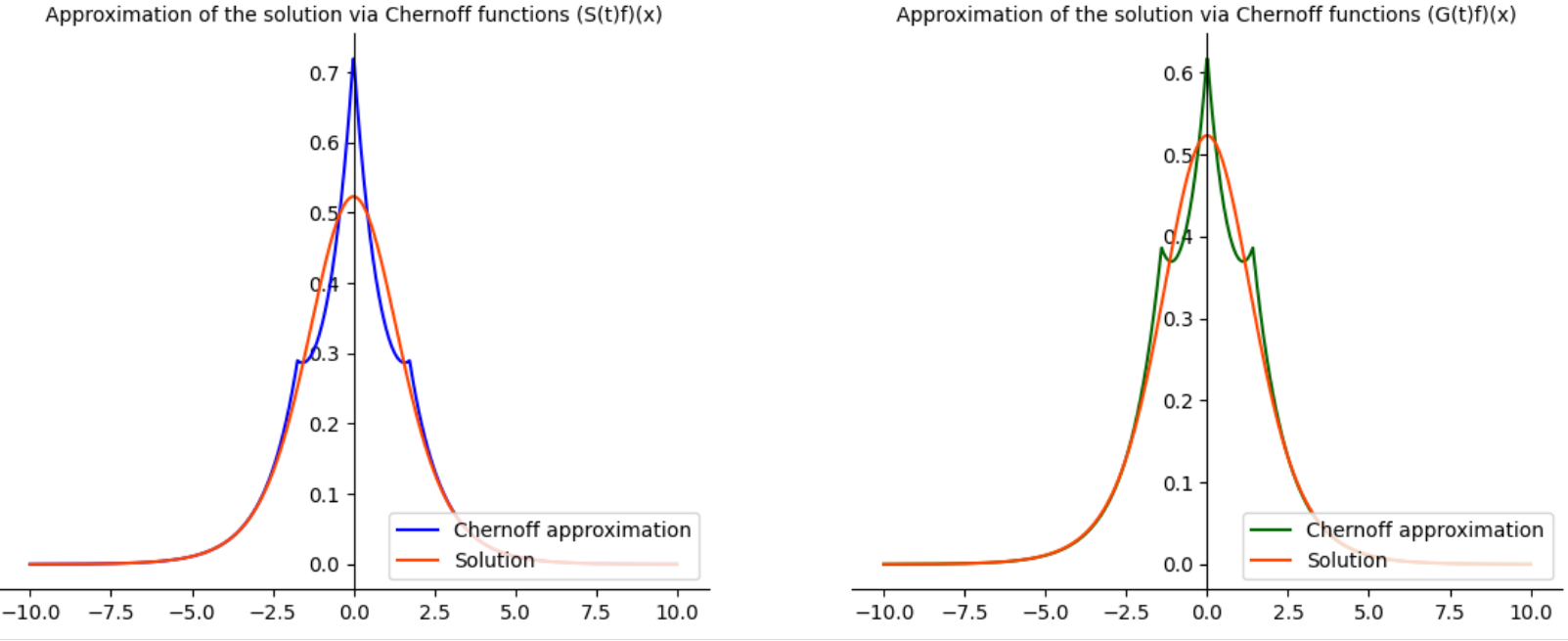}\\
\end{center}
n=2
\begin{center}
	\includegraphics[scale=0.5]{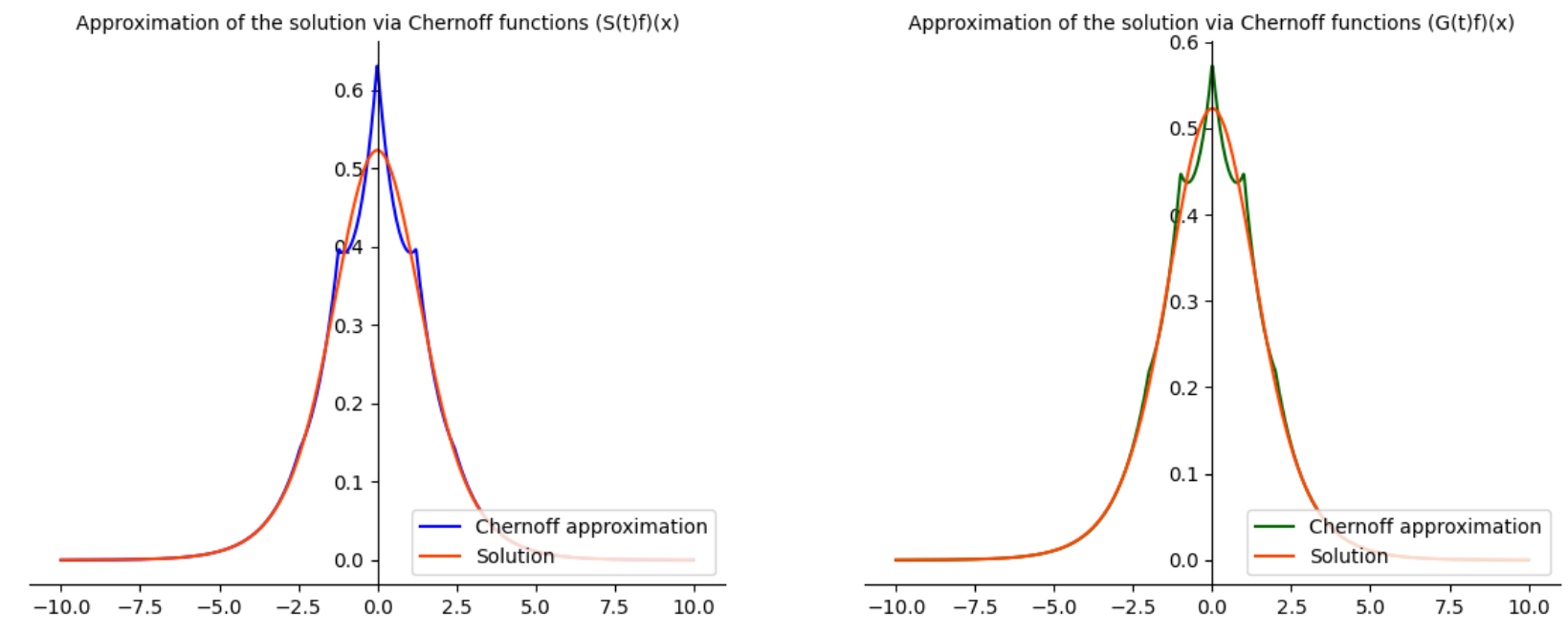}\\
\end{center}

n=3

\begin{center}
	\includegraphics[scale=0.5]{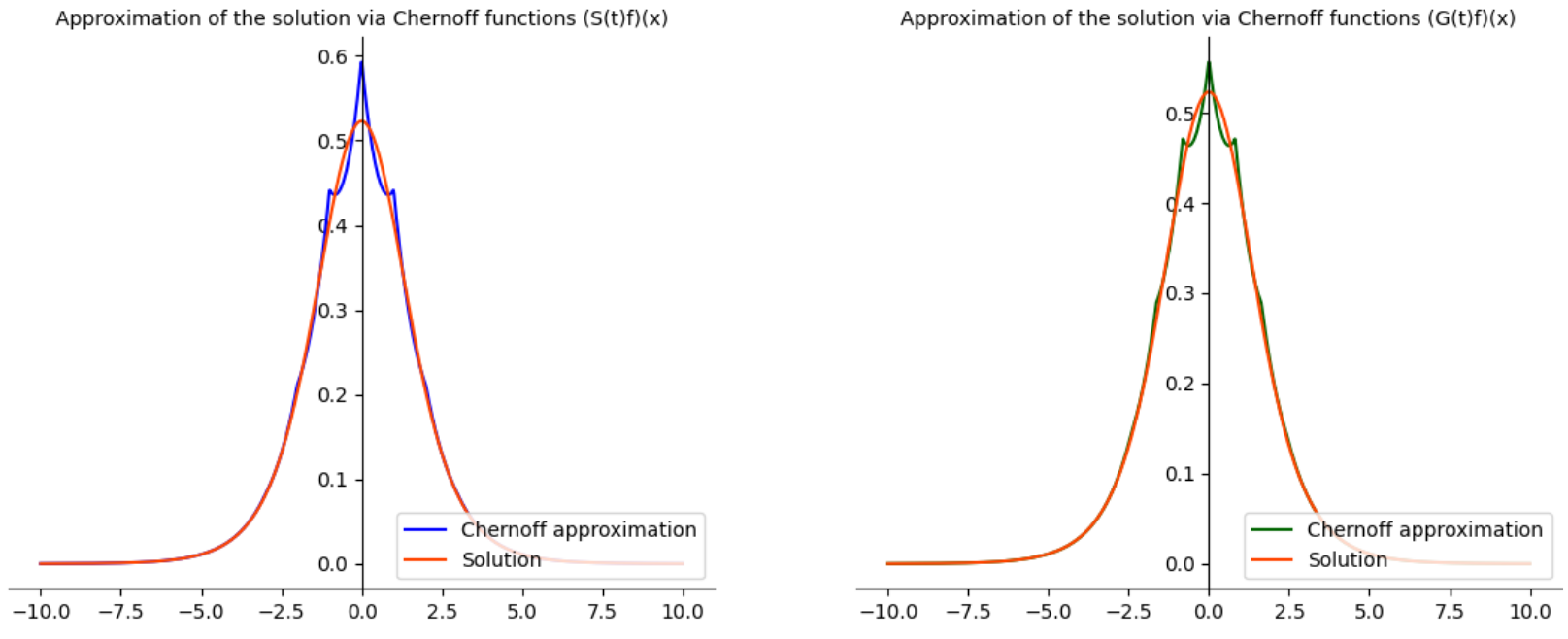}\\
\end{center}
n=4
\begin{center}
	\includegraphics[scale=0.5]{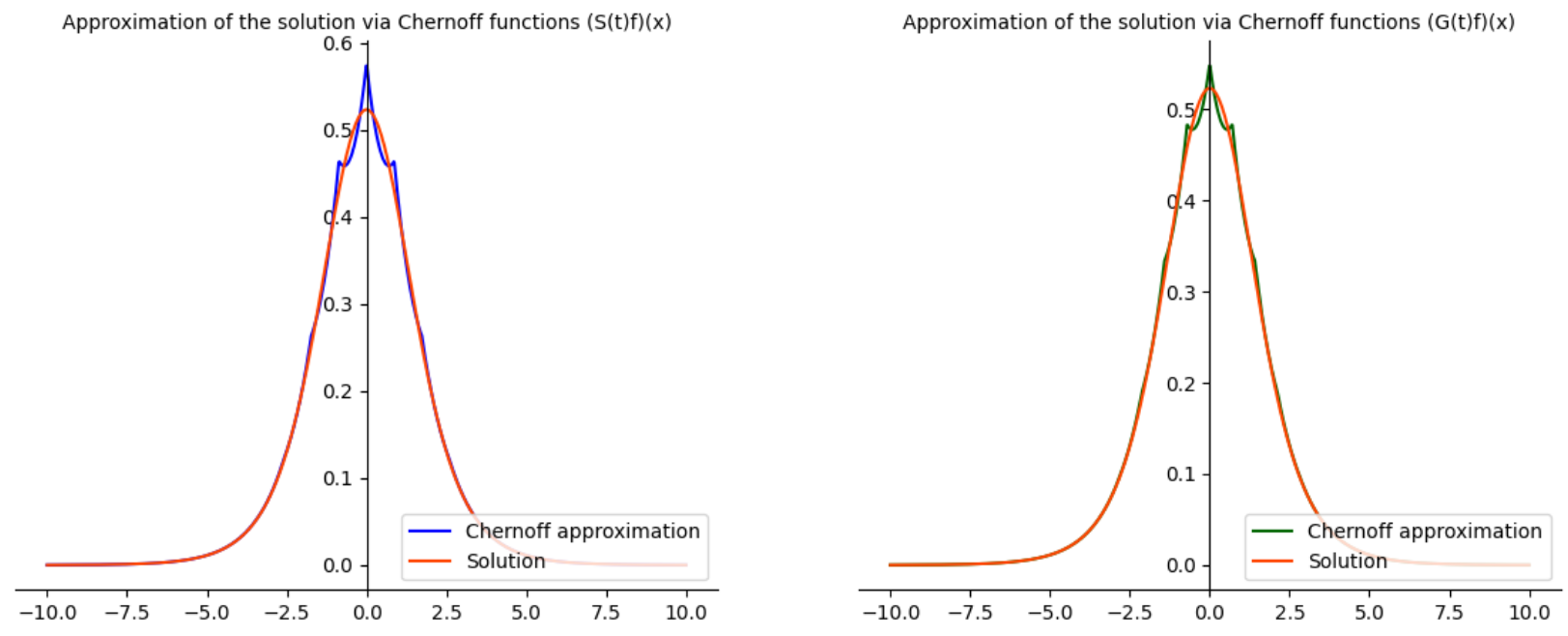}\\
\end{center}
n=5
\begin{center}
	\includegraphics[scale=0.5]{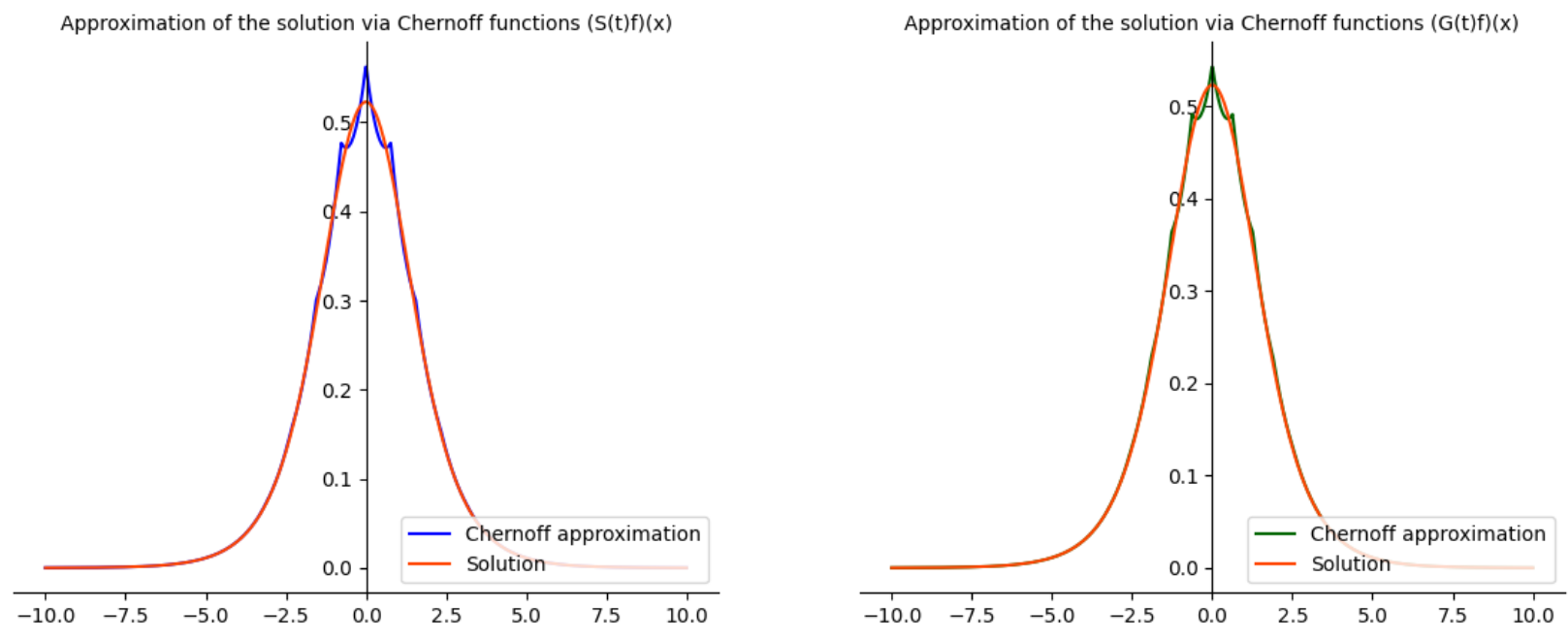}\\
\end{center}

n=6
\begin{center}
	\includegraphics[scale=0.5]{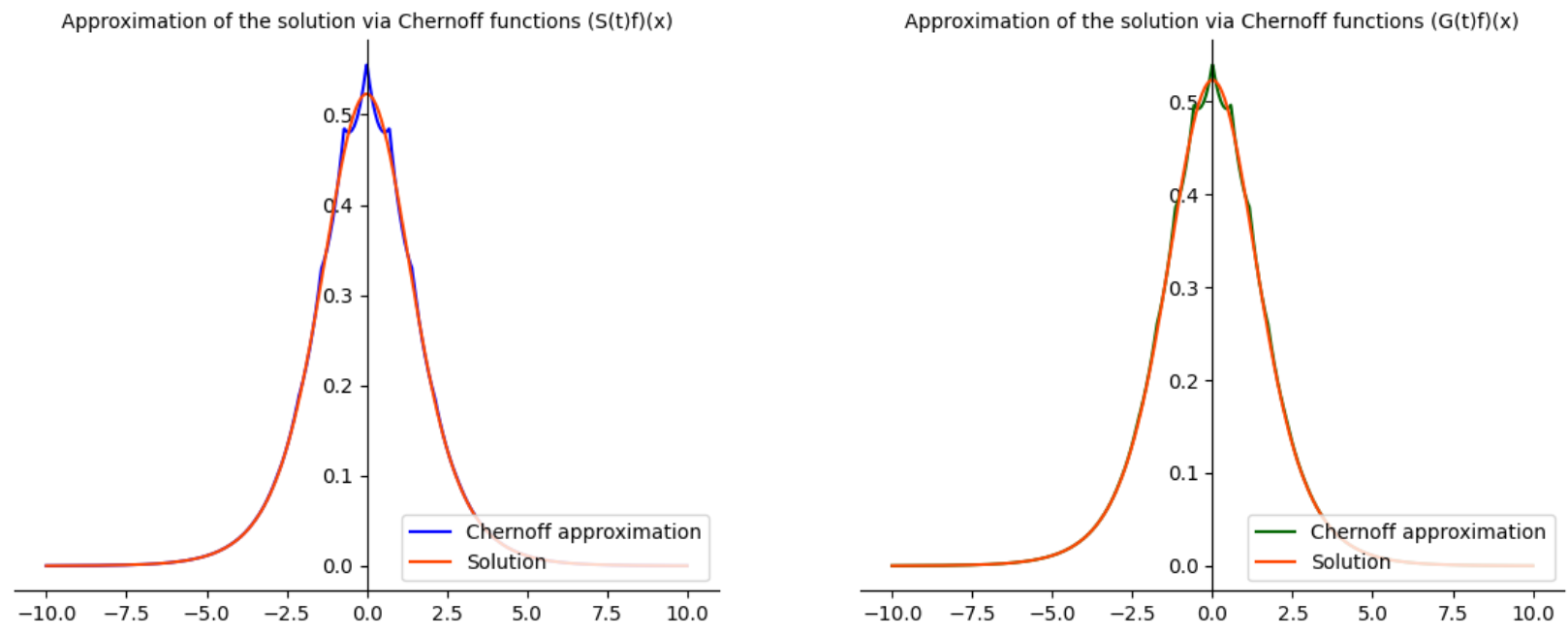}\\
\end{center}
n=7
\begin{center}
	\includegraphics[scale=0.5]{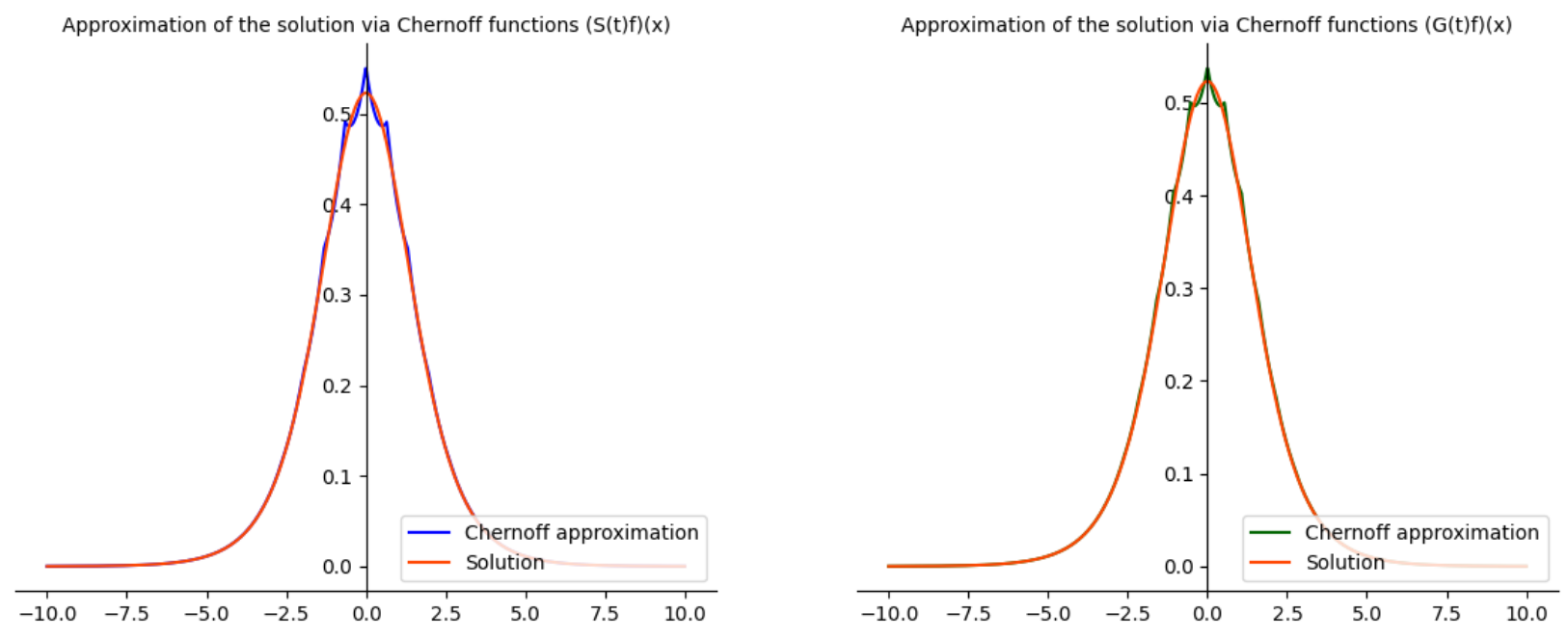}\\
\end{center}
n=8
\begin{center}
	\includegraphics[scale=0.5]{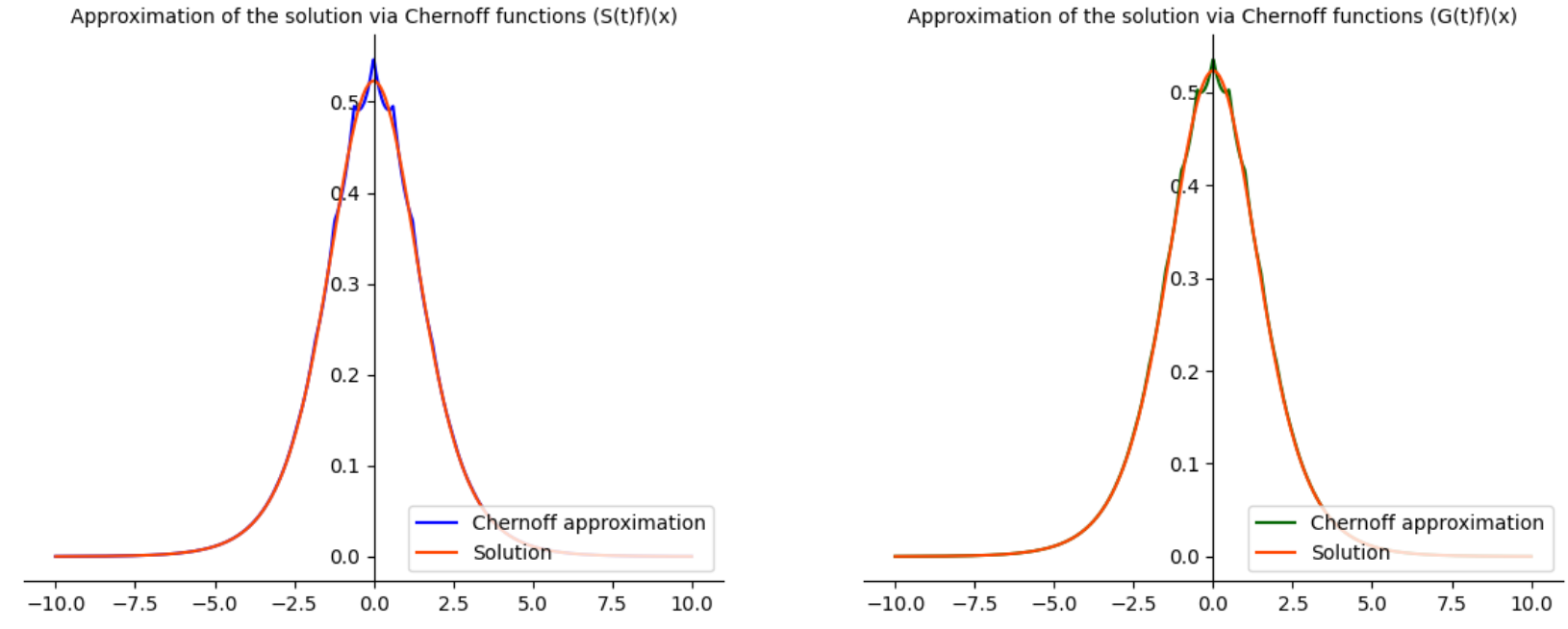}\\
\end{center}

n=9
\begin{center}
	\includegraphics[scale=0.5]{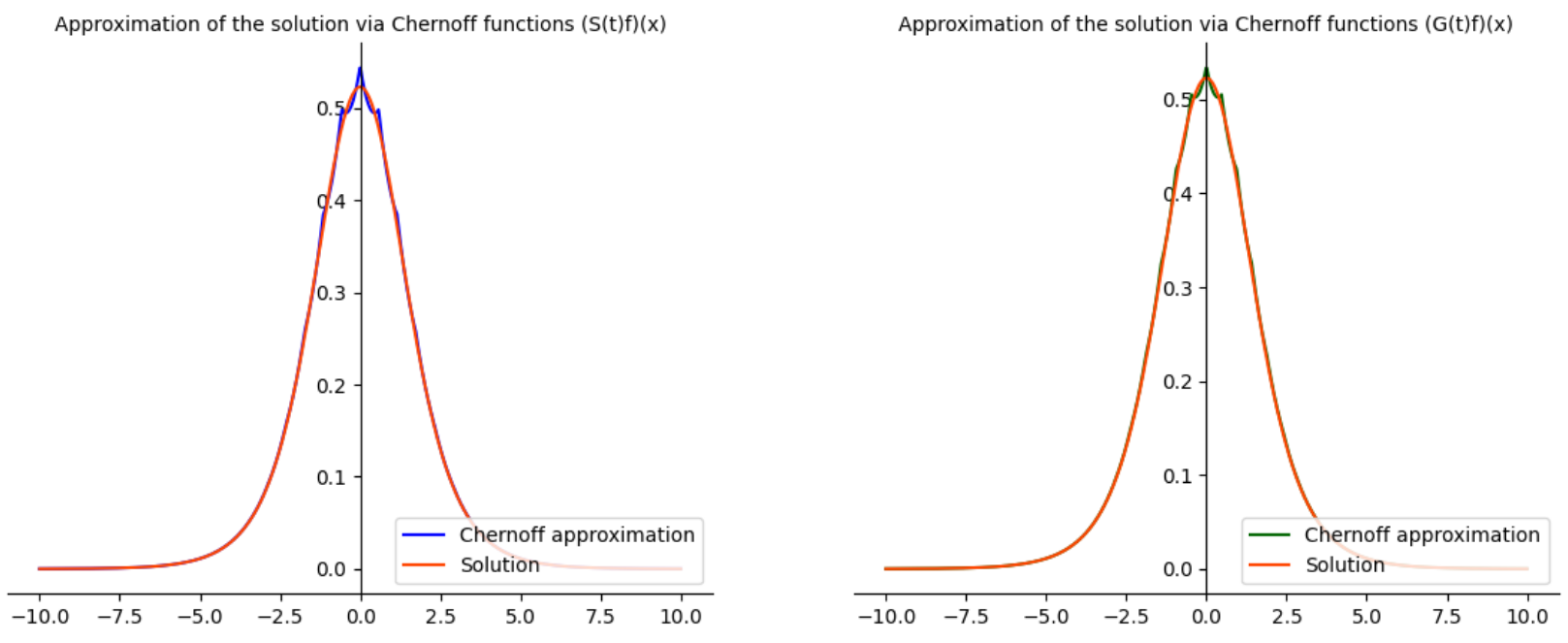}\\
\end{center}
n=10
\begin{center}
	\includegraphics[scale=0.5]{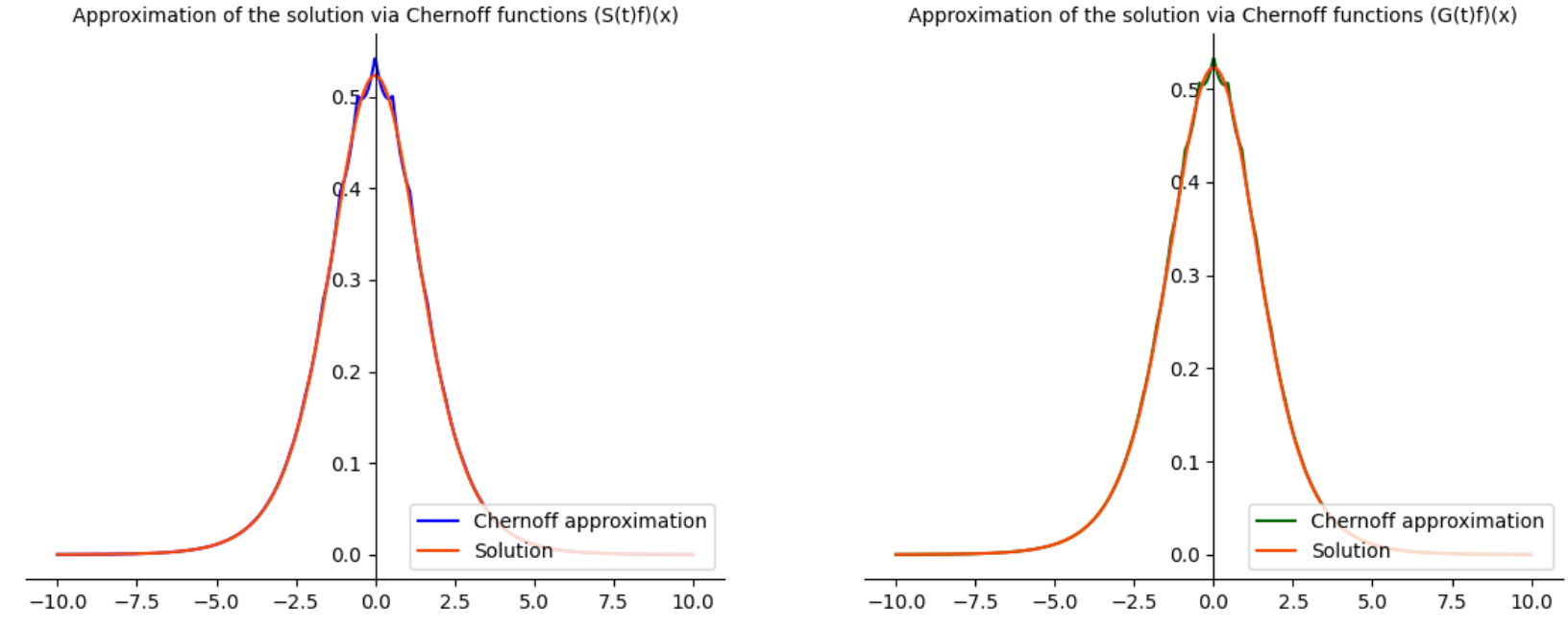}\\
\end{center}

\subsection{$u_0(x)=e^{-|x|^{1/4}}$}

n=1
\begin{center}
	\includegraphics[scale=0.5]{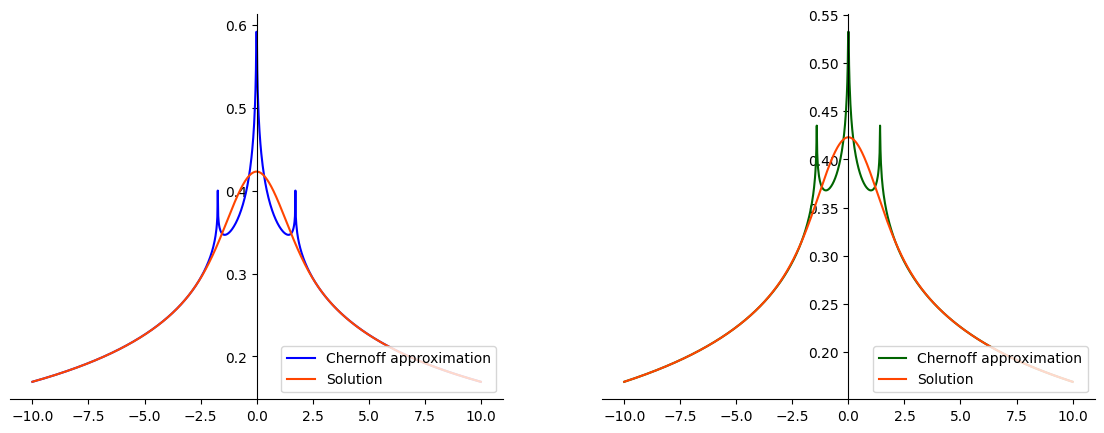}\\
\end{center}
n=2
\begin{center}
	\includegraphics[scale=0.5]{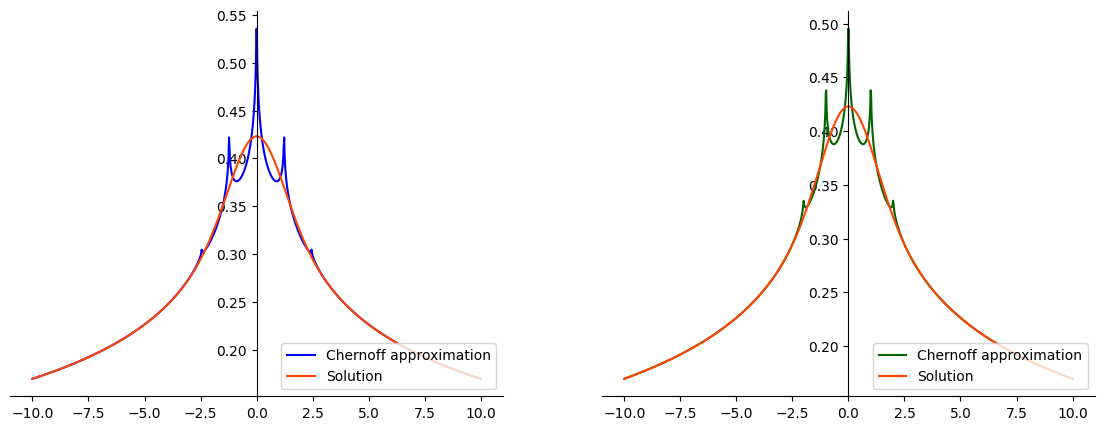}\\
\end{center}

n=3

\begin{center}
	\includegraphics[scale=0.5]{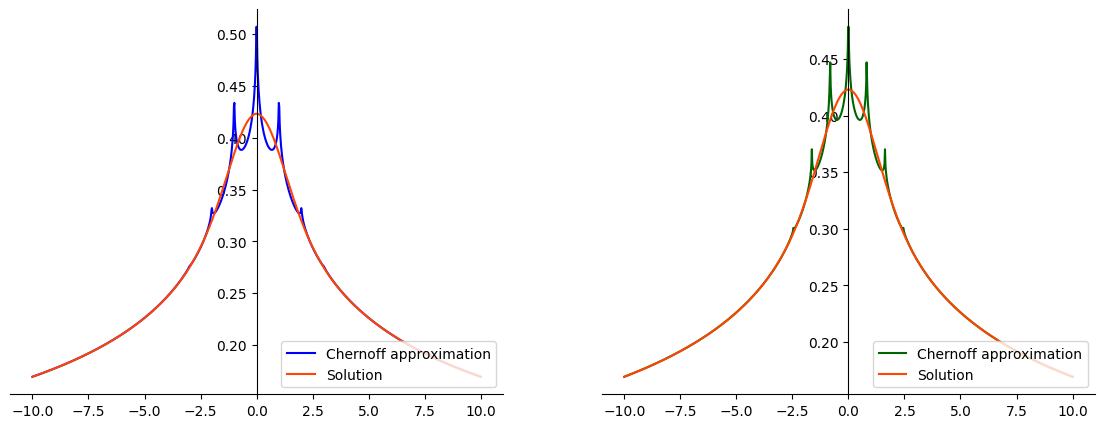}\\
\end{center}
n=4
\begin{center}
	\includegraphics[scale=0.5]{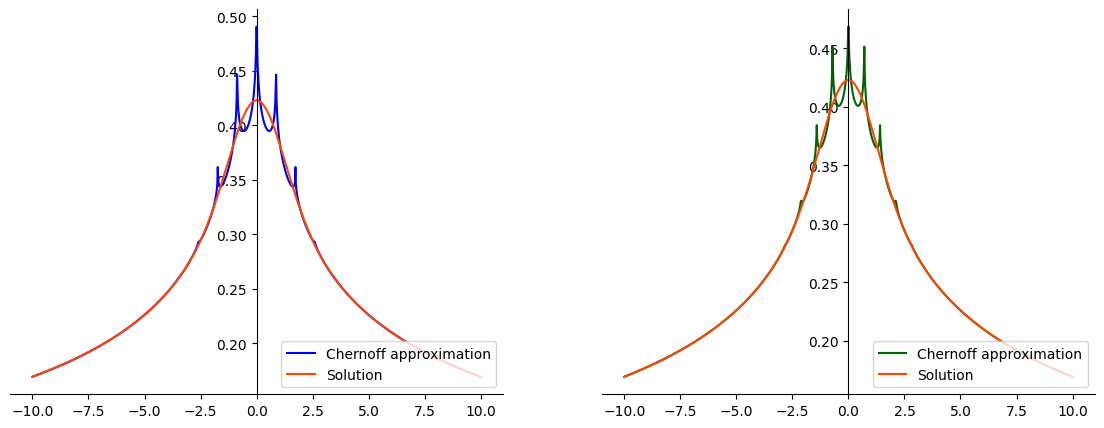}\\
\end{center}
n=5
\begin{center}
	\includegraphics[scale=0.5]{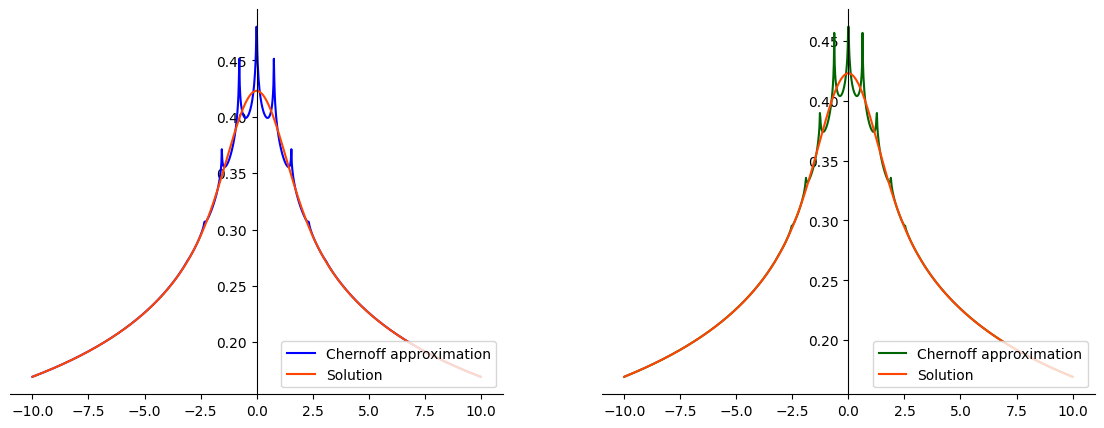}\\
\end{center}

n=6
\begin{center}
	\includegraphics[scale=0.5]{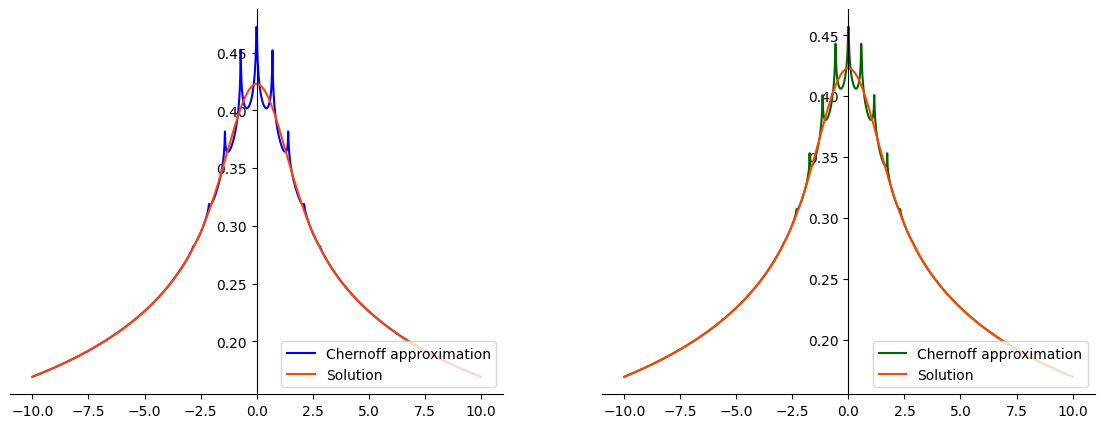}\\
\end{center}
n=7
\begin{center}
	\includegraphics[scale=0.5]{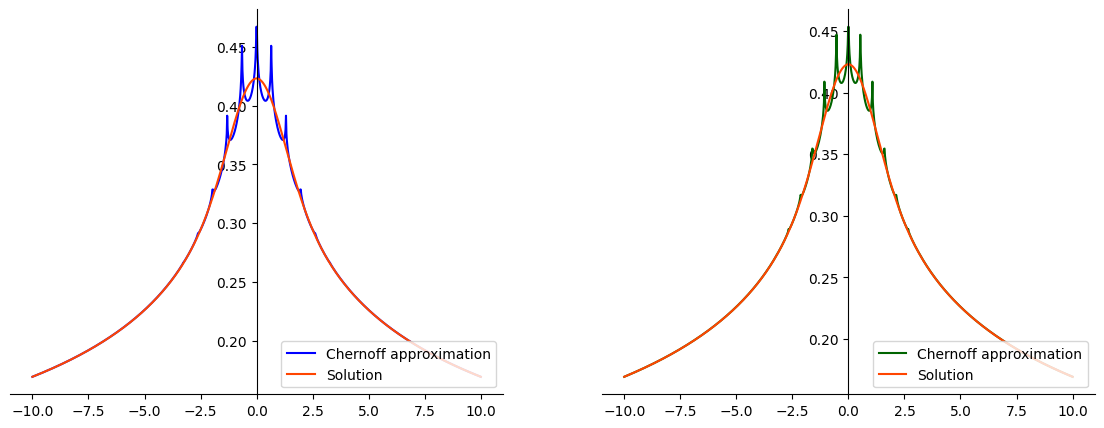}\\
\end{center}
n=8
\begin{center}
	\includegraphics[scale=0.5]{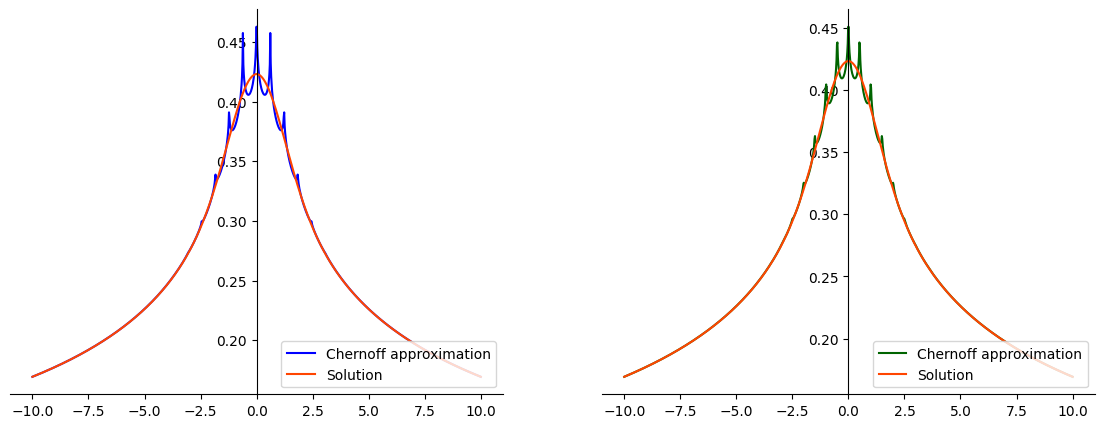}\\
\end{center}

n=9
\begin{center}
	\includegraphics[scale=0.5]{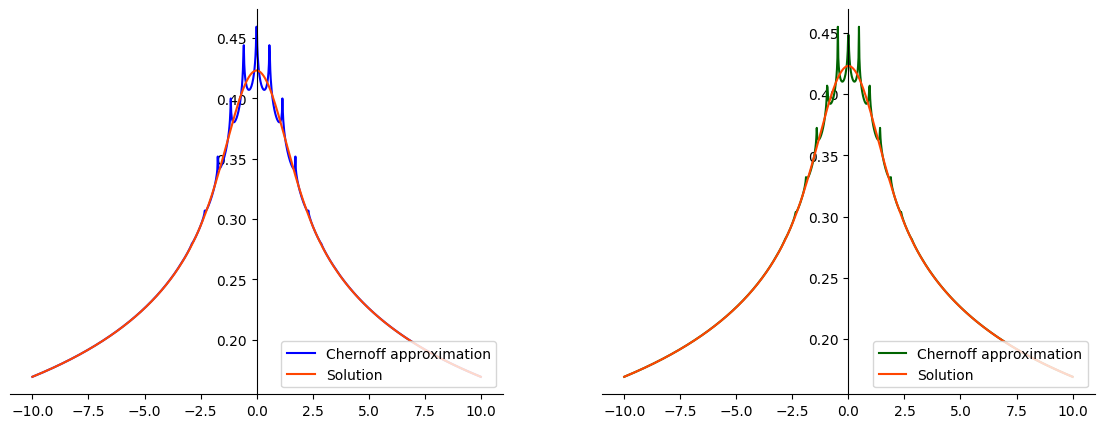}\\
\end{center}
n=10
\begin{center}
	\includegraphics[scale=0.5]{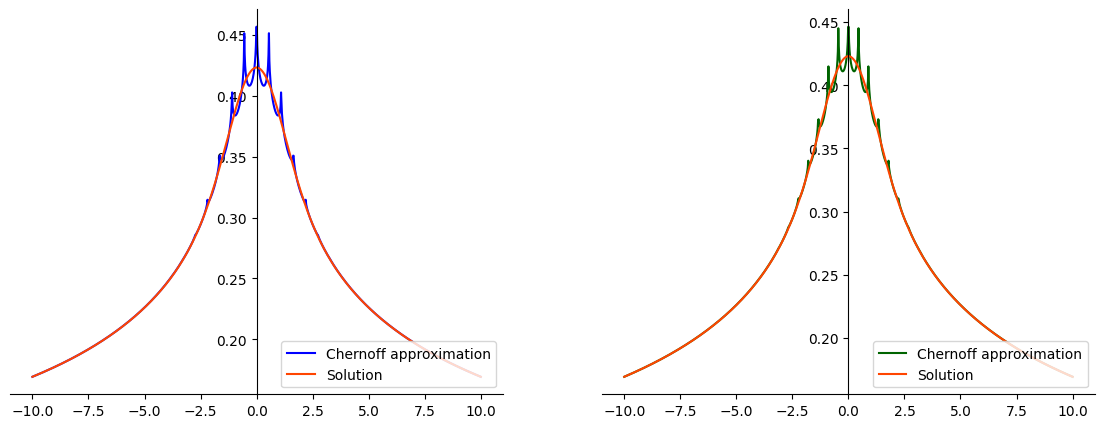}\\
\end{center}

\subsection{$u_0(x)=e^{-|x|^{1/2}}$}

n=1
\begin{center}
	\includegraphics[scale=0.5]{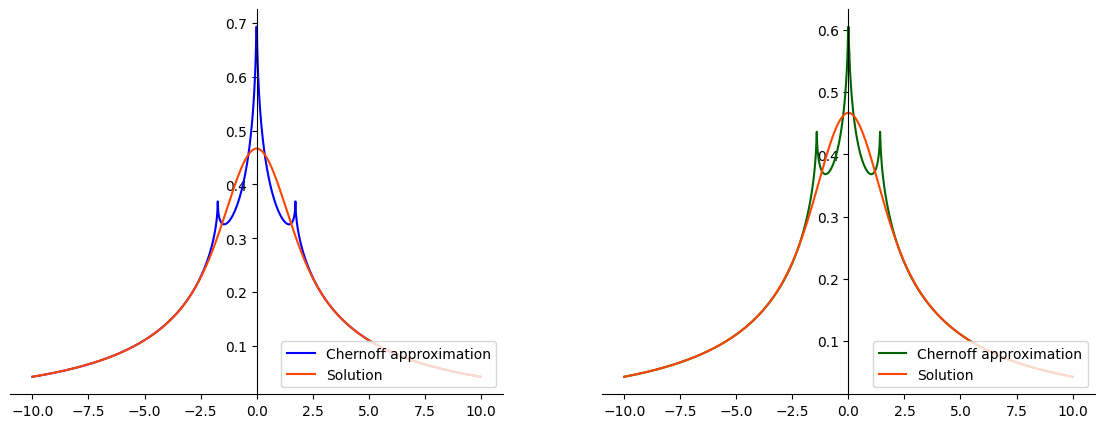}\\
\end{center}
n=2
\begin{center}
	\includegraphics[scale=0.5]{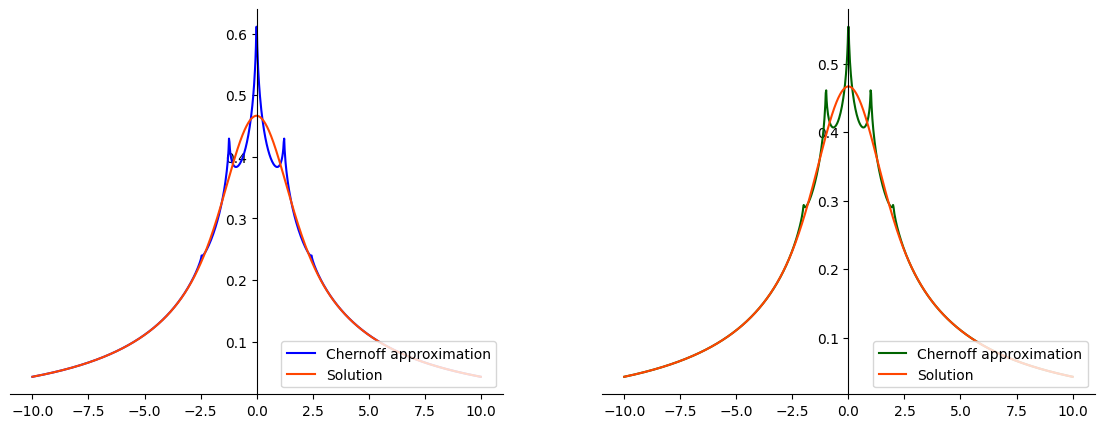}\\
\end{center}

n=3

\begin{center}
	\includegraphics[scale=0.5]{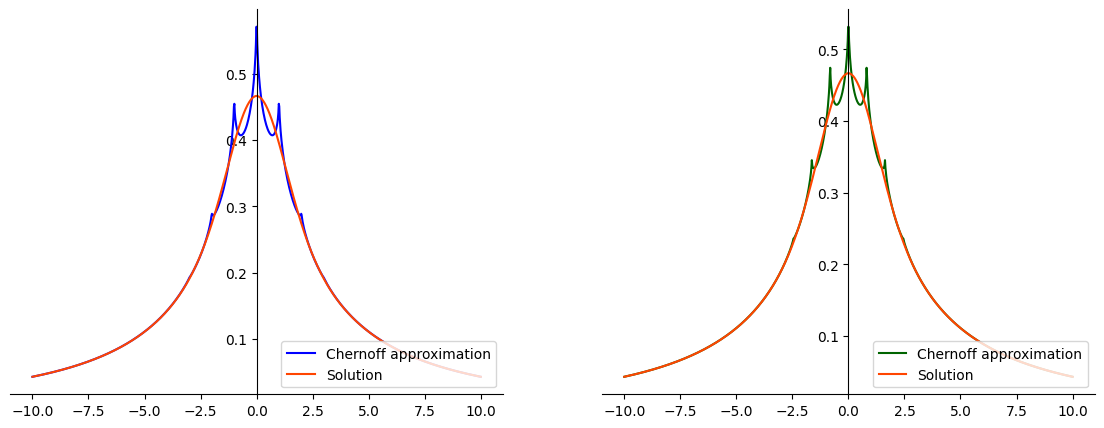}\\
\end{center}
n=4
\begin{center}
	\includegraphics[scale=0.5]{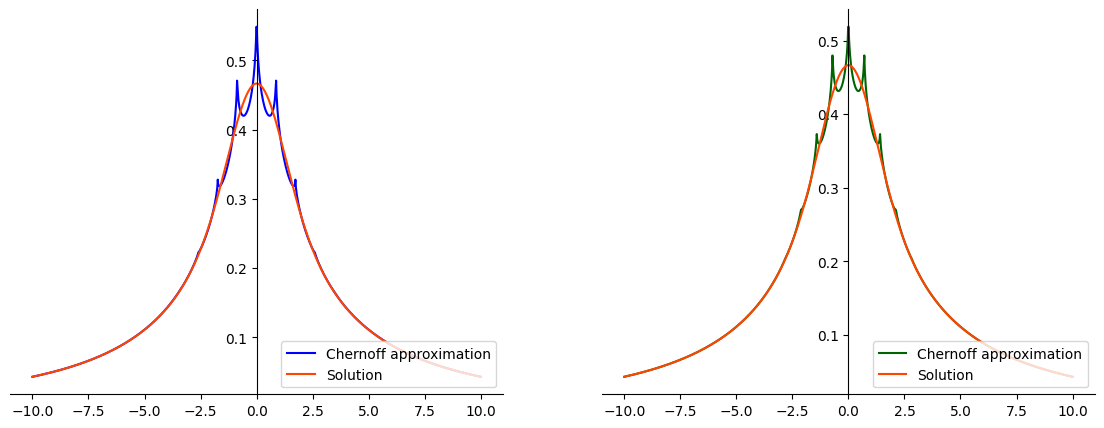}\\
\end{center}
n=5
\begin{center}
	\includegraphics[scale=0.5]{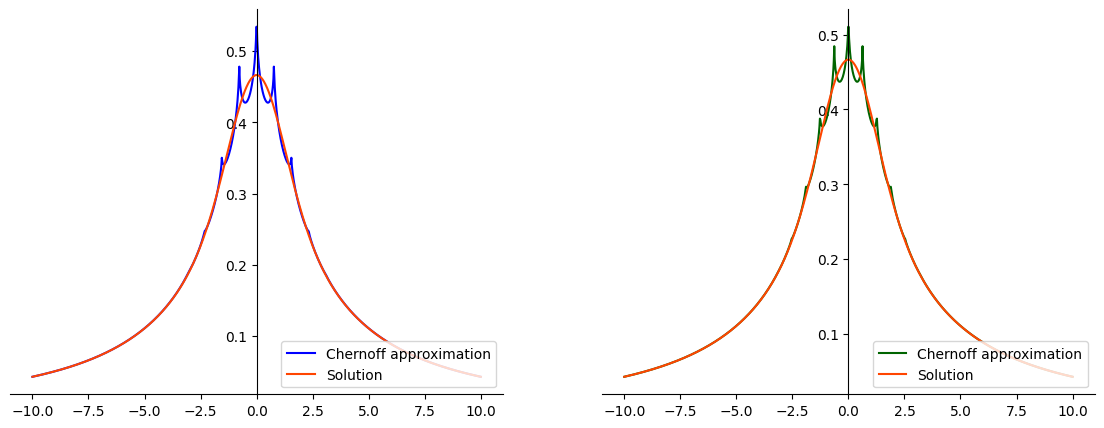}\\
\end{center}

n=6
\begin{center}
	\includegraphics[scale=0.5]{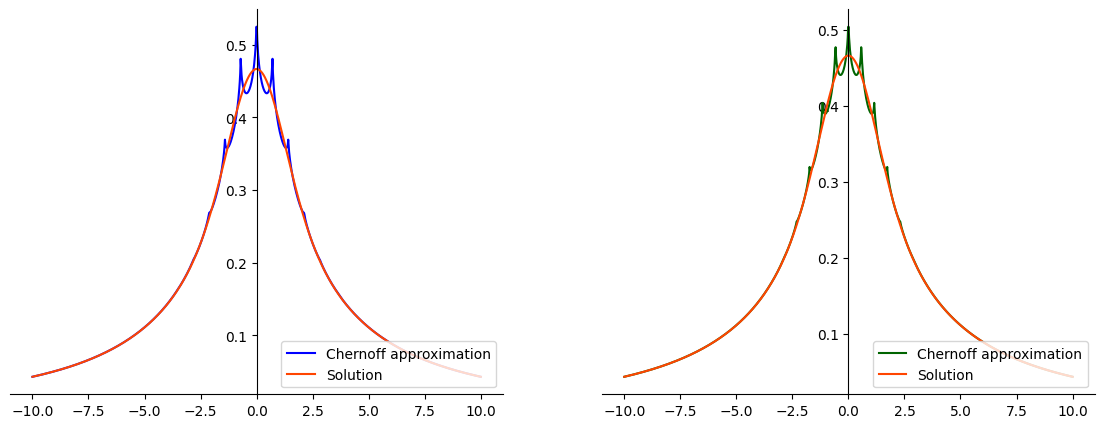}\\
\end{center}
n=7
\begin{center}
	\includegraphics[scale=0.5]{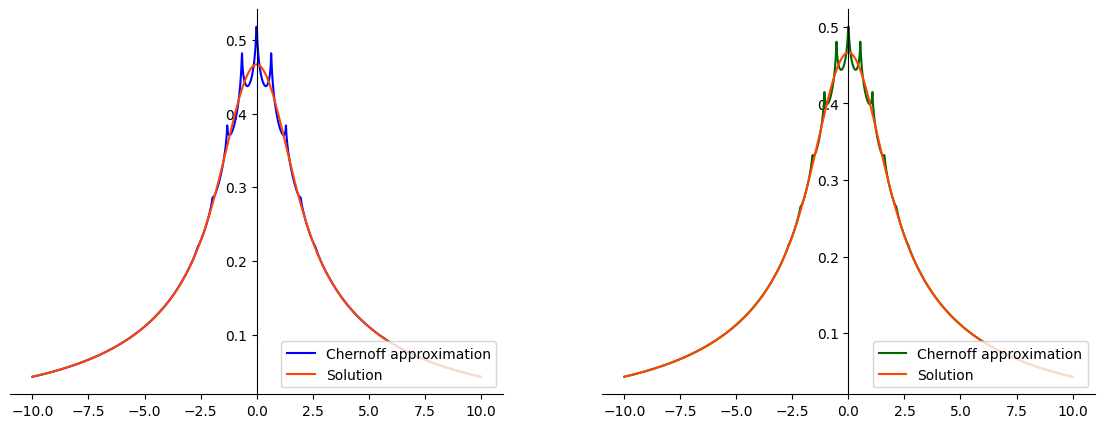}\\
\end{center}
n=8
\begin{center}
	\includegraphics[scale=0.5]{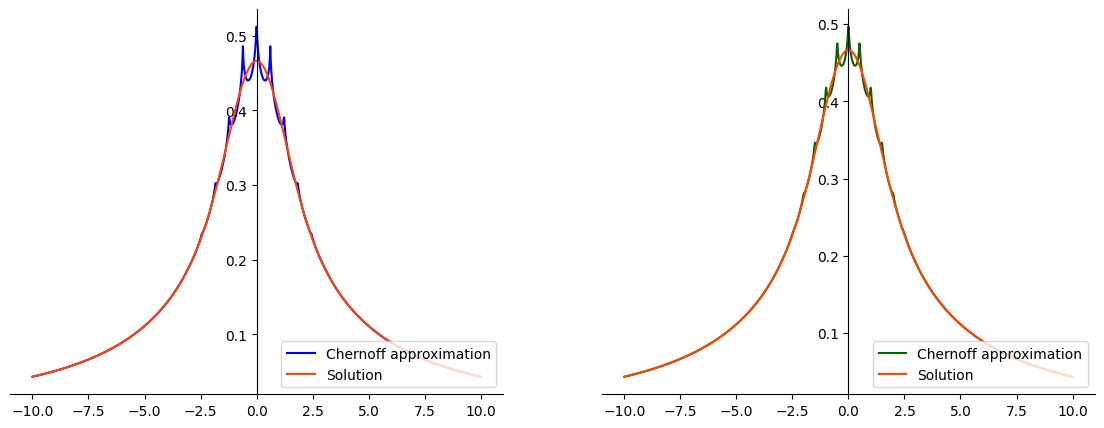}\\
\end{center}

n=9
\begin{center}
	\includegraphics[scale=0.5]{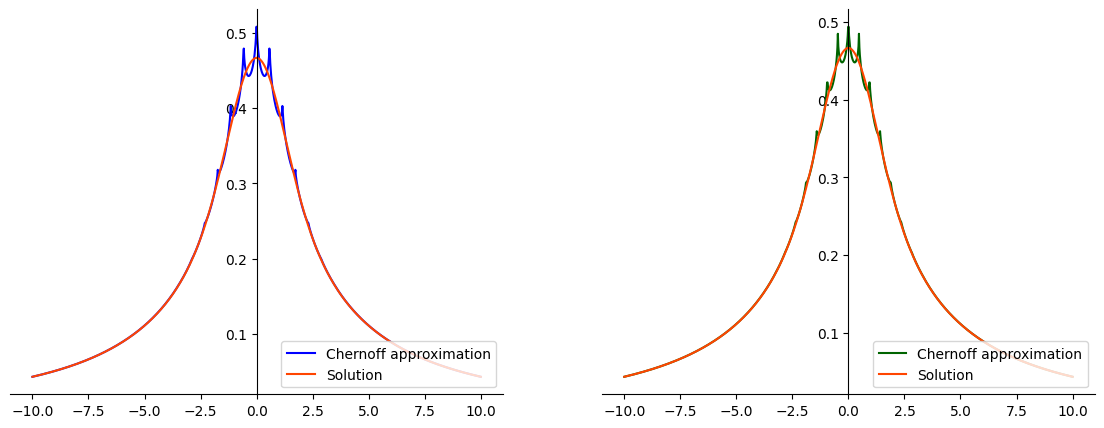}\\
\end{center}
n=10
\begin{center}
	\includegraphics[scale=0.5]{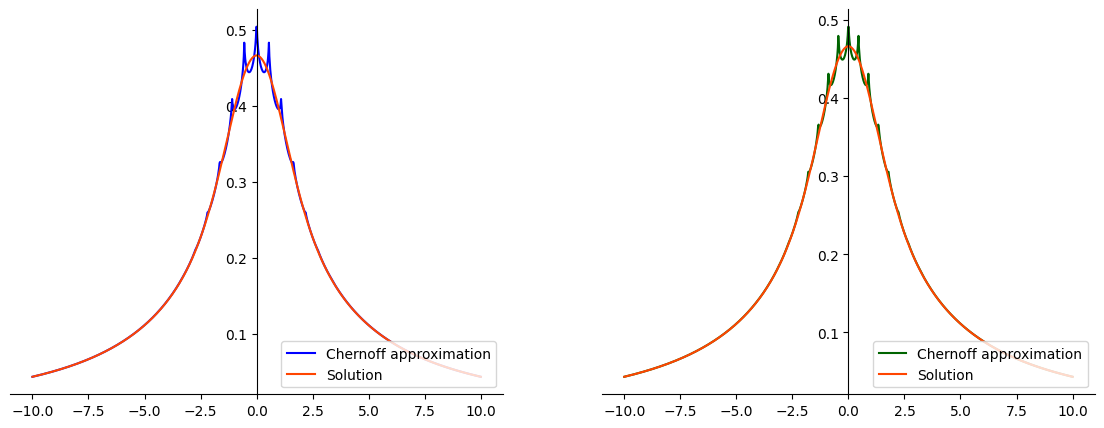}\\
\end{center}

\subsection{$u_0(x)=e^{-|x|^{3/4}}$}

n=1
\begin{center}
	\includegraphics[scale=0.5]{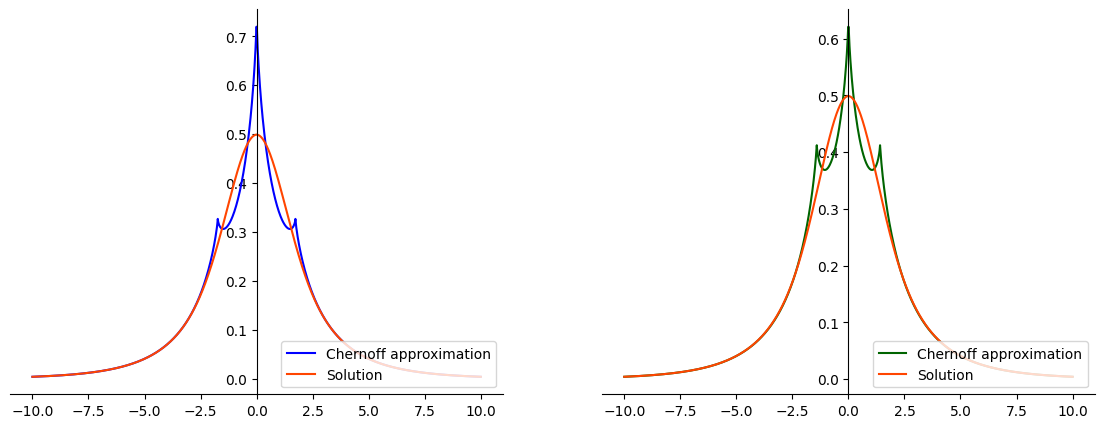}\\
\end{center}
n=2
\begin{center}
	\includegraphics[scale=0.5]{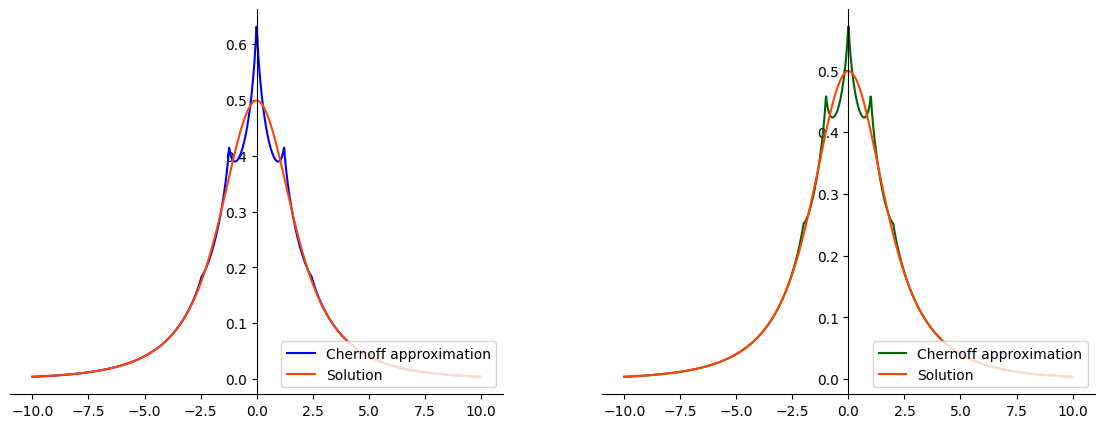}\\
\end{center}

n=3

\begin{center}
	\includegraphics[scale=0.5]{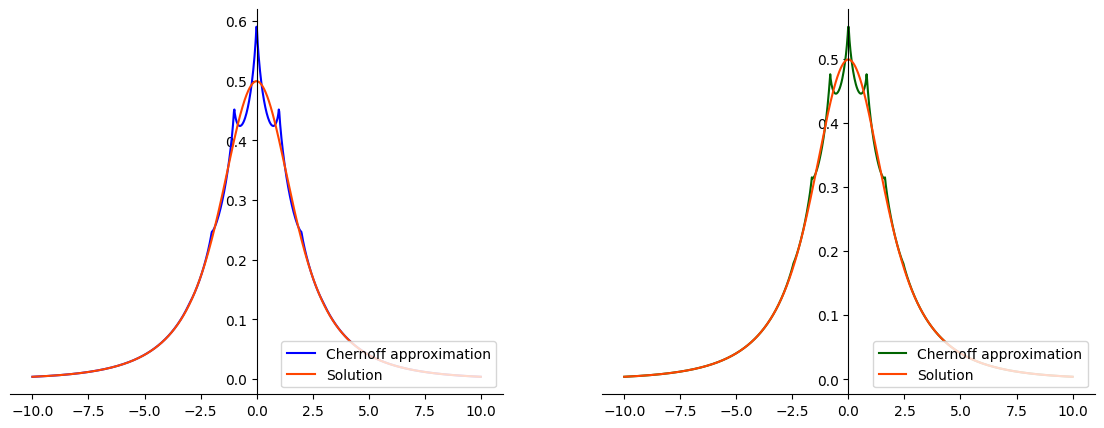}\\
\end{center}
n=4
\begin{center}
	\includegraphics[scale=0.5]{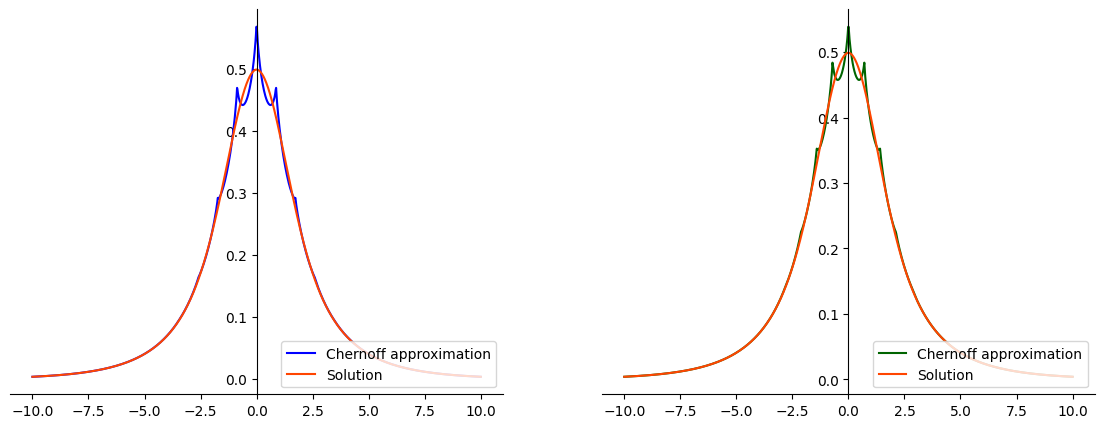}\\
\end{center}
n=5
\begin{center}
	\includegraphics[scale=0.5]{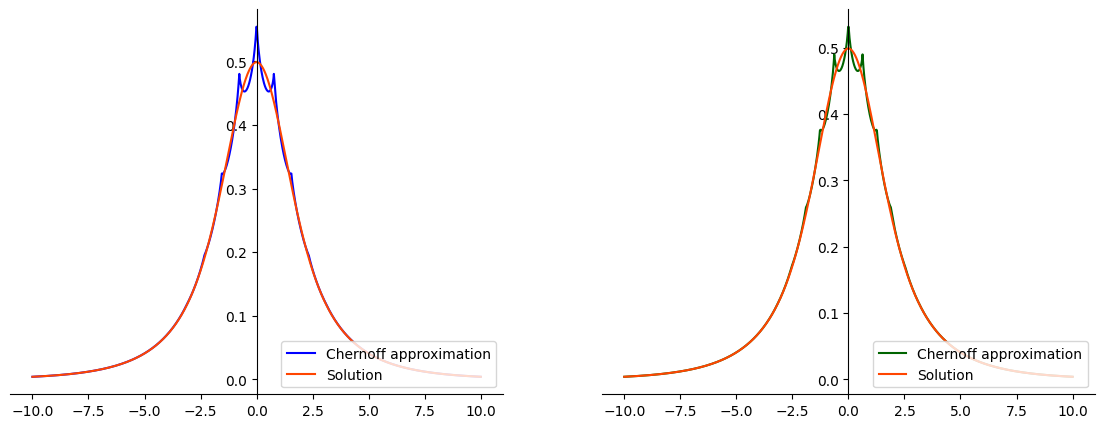}\\
\end{center}
n=6
\begin{center}
	\includegraphics[scale=0.5]{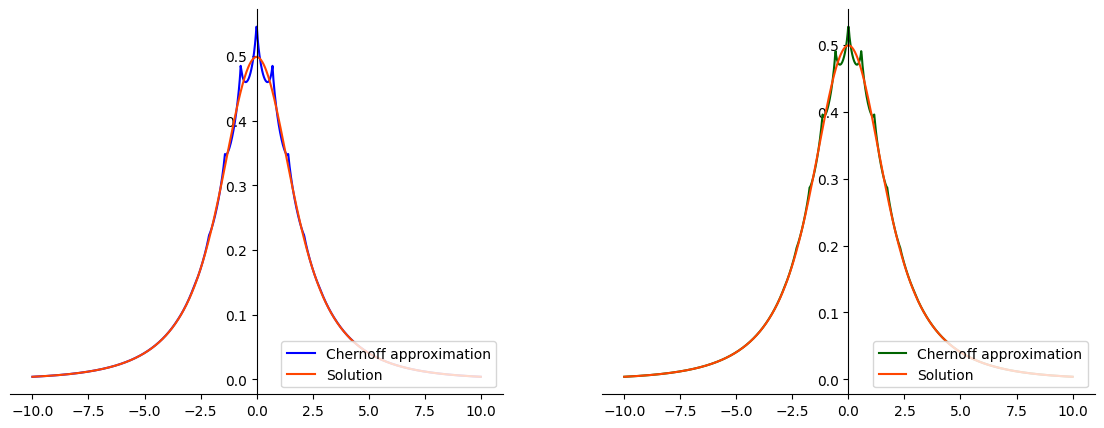}\\
\end{center}
n=7
\begin{center}
	\includegraphics[scale=0.5]{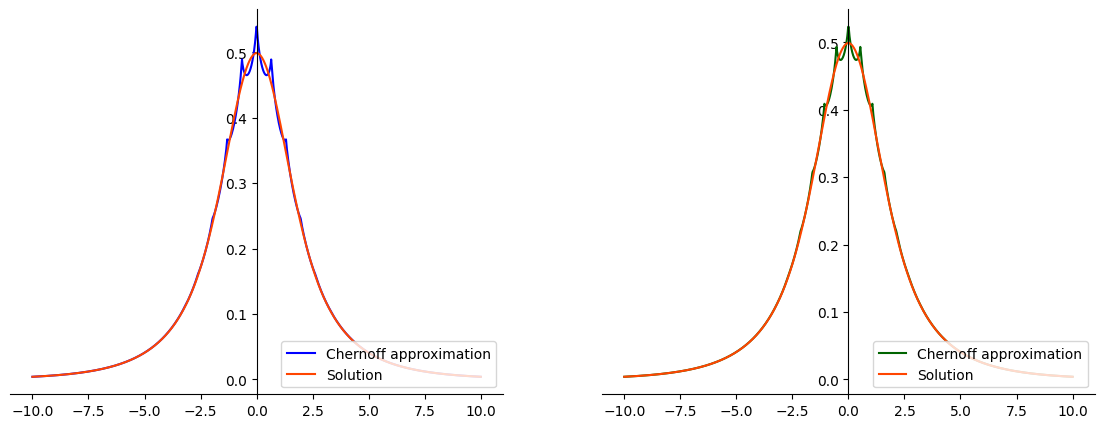}\\
\end{center}
n=8
\begin{center}
	\includegraphics[scale=0.5]{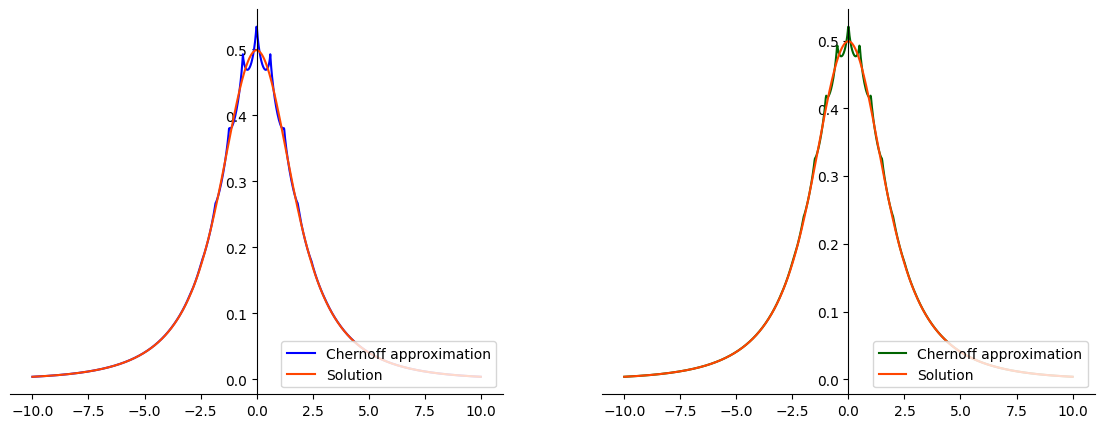}\\
\end{center}

n=9
\begin{center}
	\includegraphics[scale=0.5]{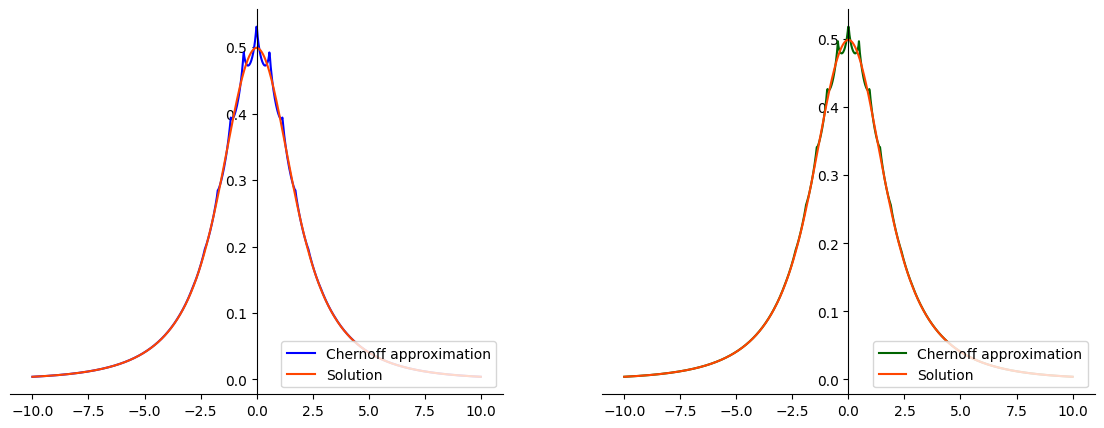}\\
\end{center}
n=10
\begin{center}
	\includegraphics[scale=0.5]{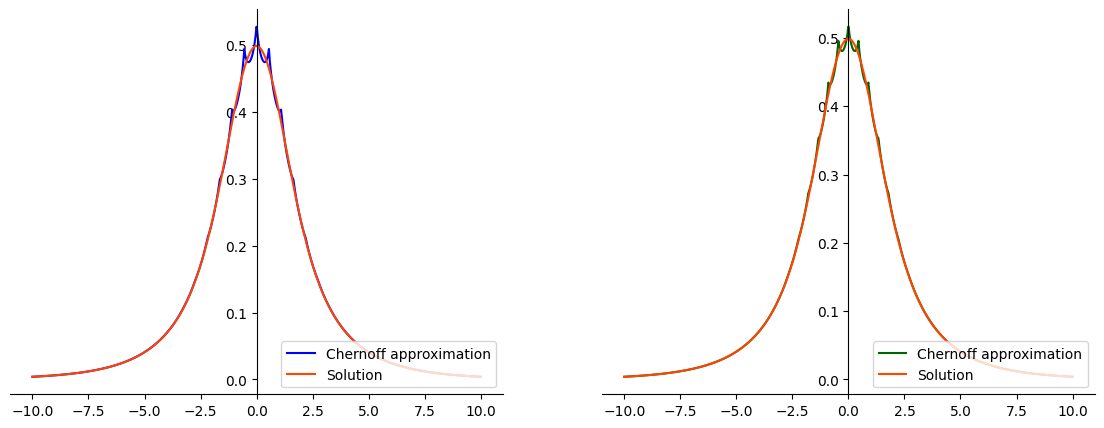}\\
\end{center}

\subsection{$u_0(x)=e^{-|x|^{3/2}}$}

n=1
\begin{center}
	\includegraphics[scale=0.5]{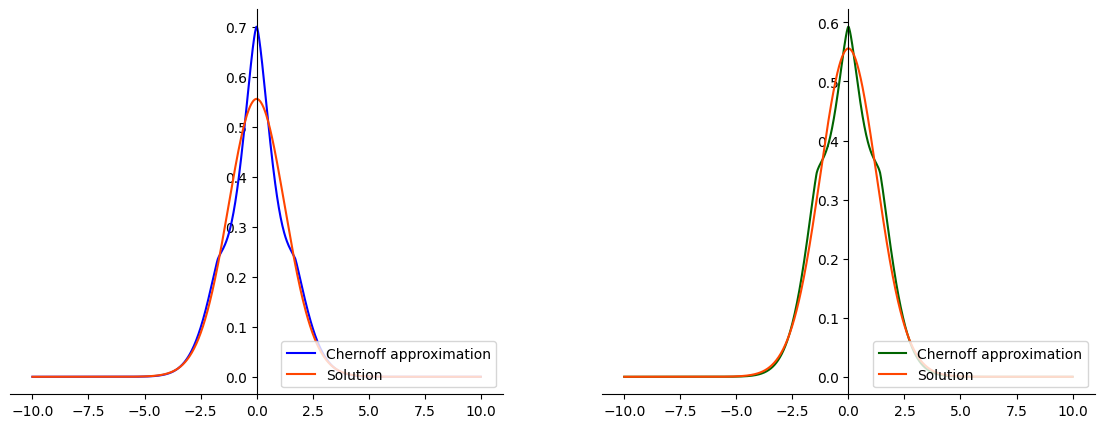}\\
\end{center}
n=2
\begin{center}
	\includegraphics[scale=0.5]{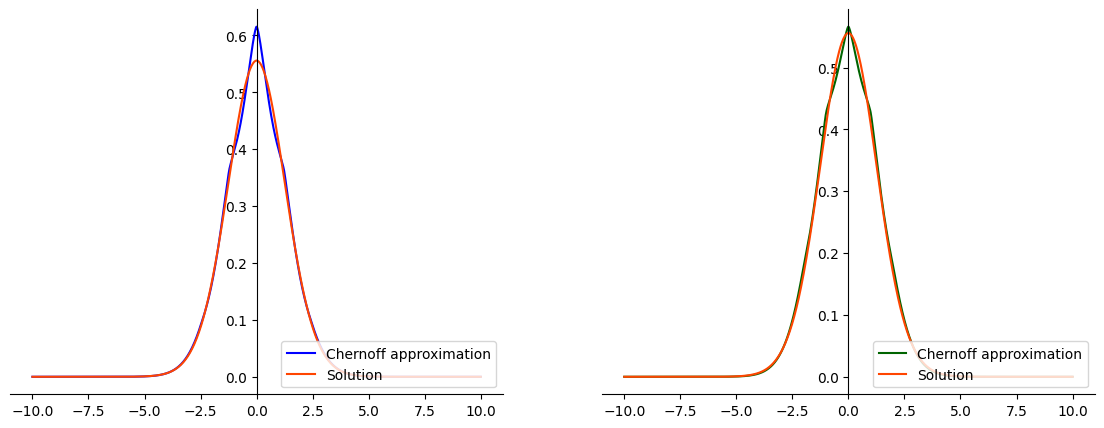}\\
\end{center}

n=3

\begin{center}
	\includegraphics[scale=0.5]{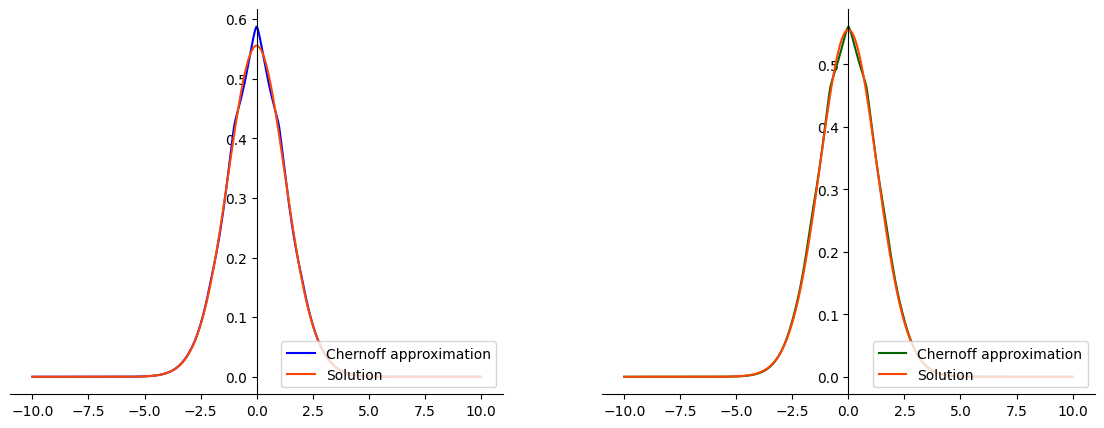}\\
\end{center}
n=4
\begin{center}
	\includegraphics[scale=0.5]{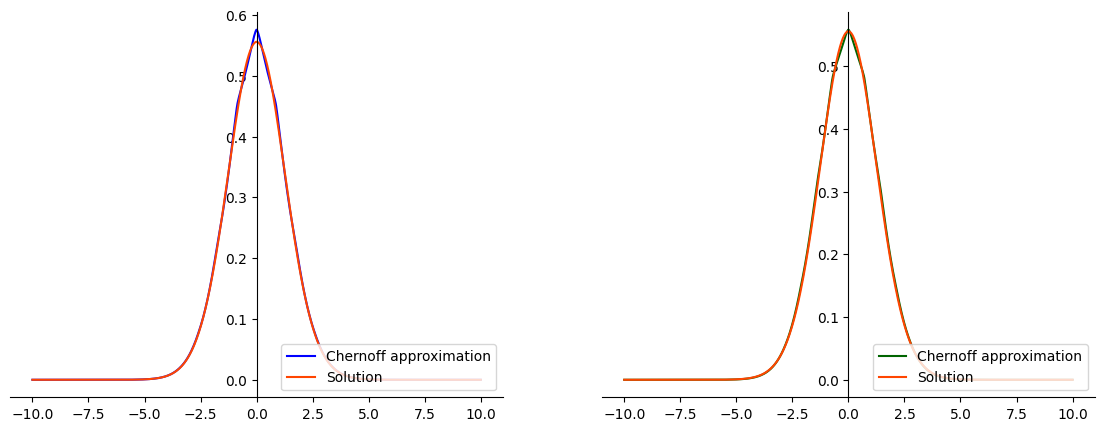}\\
\end{center}
n=5
\begin{center}
	\includegraphics[scale=0.5]{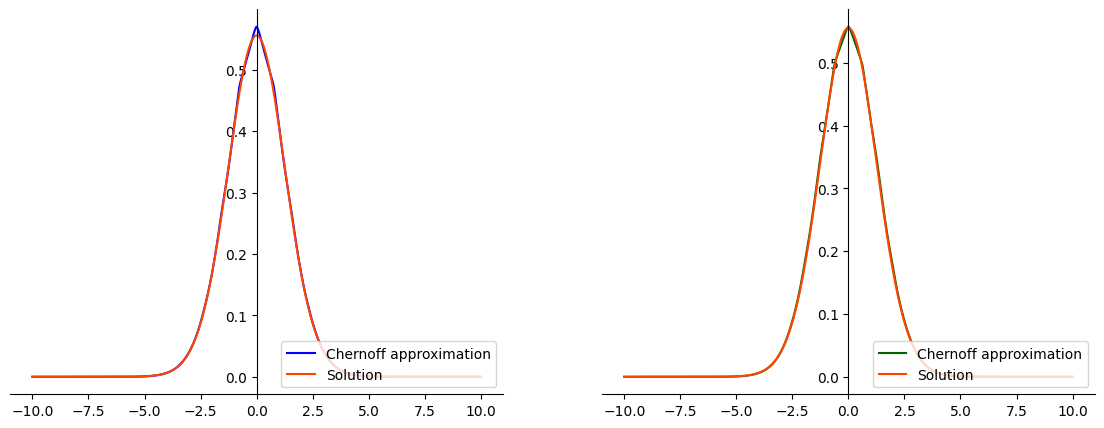}\\
\end{center}

n=6
\begin{center}
	\includegraphics[scale=0.5]{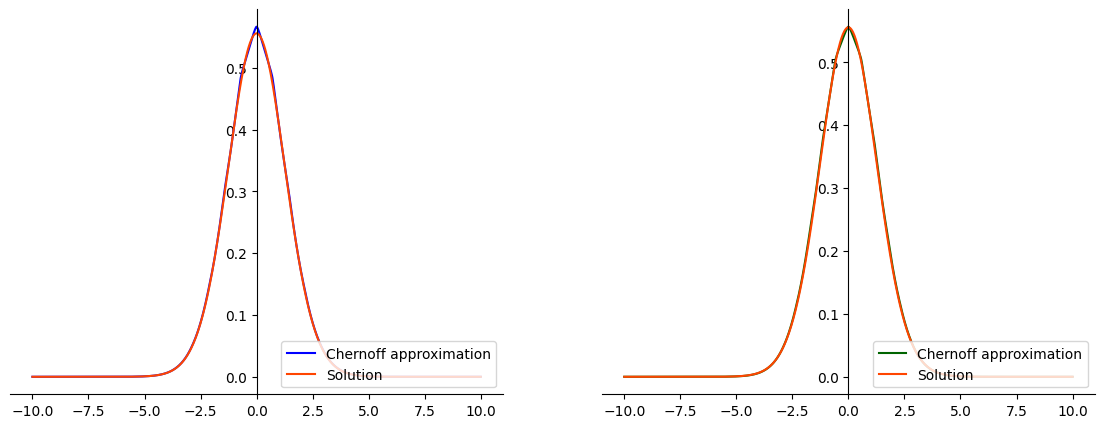}\\
\end{center}
n=7
\begin{center}
	\includegraphics[scale=0.5]{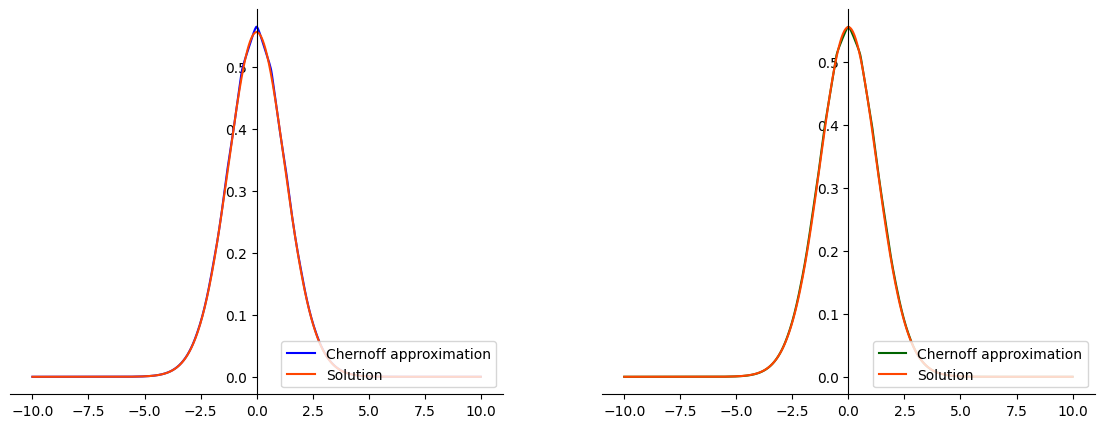}\\
\end{center}
n=8
\begin{center}
	\includegraphics[scale=0.5]{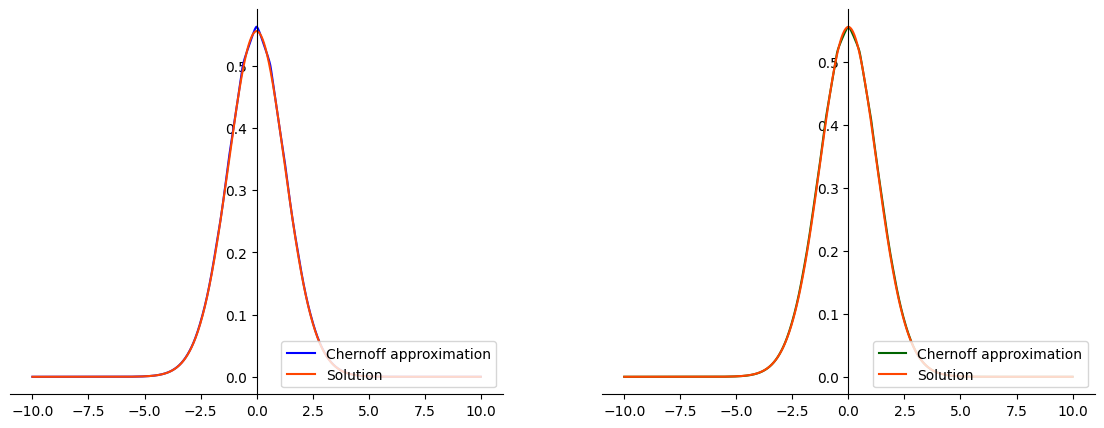}\\
\end{center}

n=9
\begin{center}
	\includegraphics[scale=0.5]{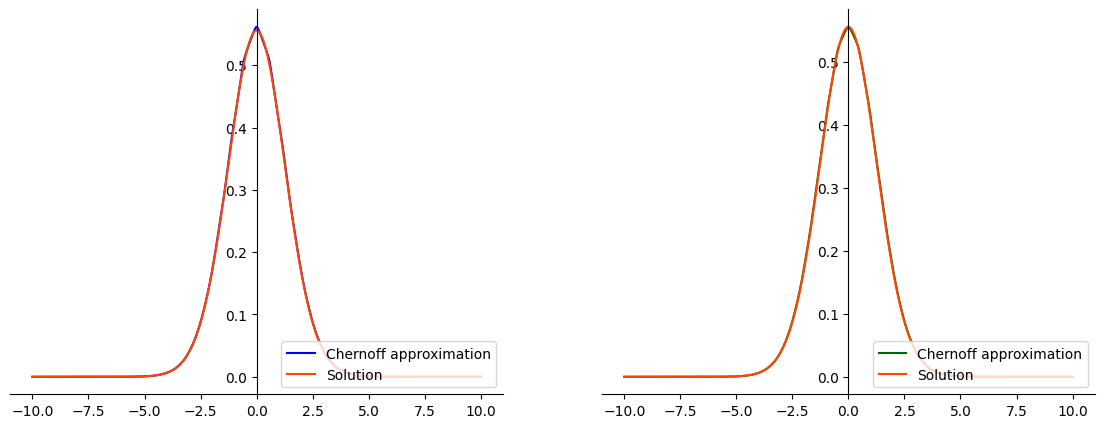}\\
\end{center}
n=10
\begin{center}
	\includegraphics[scale=0.5]{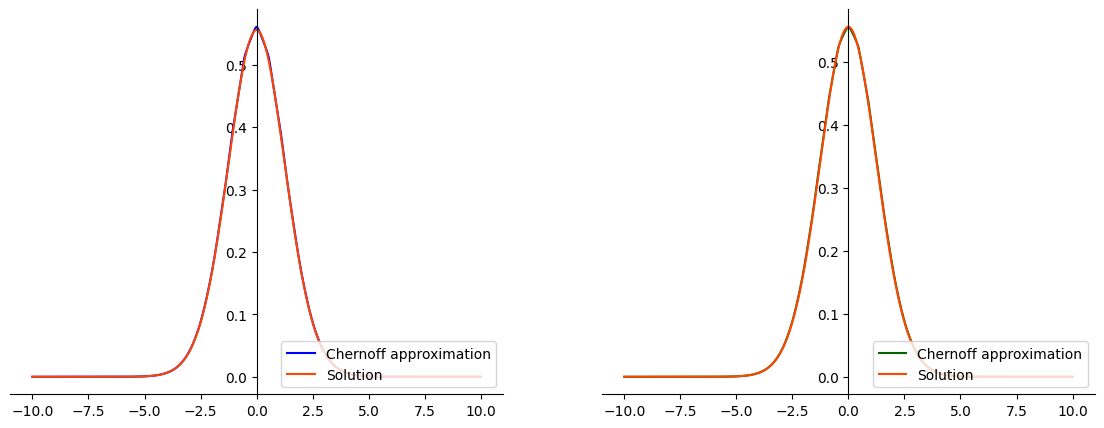}\\
\end{center}

\subsection{$u_0(x)=e^{-|x|^{5/2}}$}

n=1
\begin{center}
	\includegraphics[scale=0.5]{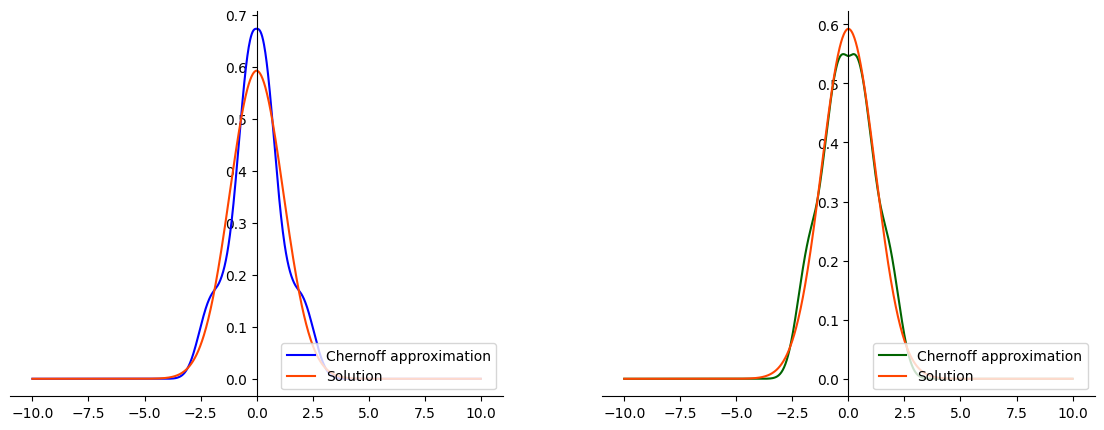}\\
\end{center}
n=2
\begin{center}
	\includegraphics[scale=0.5]{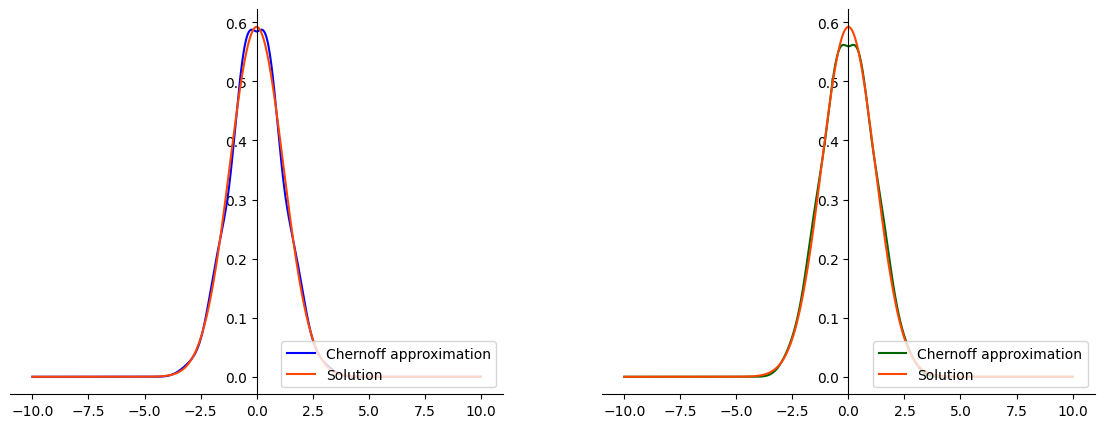}\\
\end{center}

n=3

\begin{center}
	\includegraphics[scale=0.5]{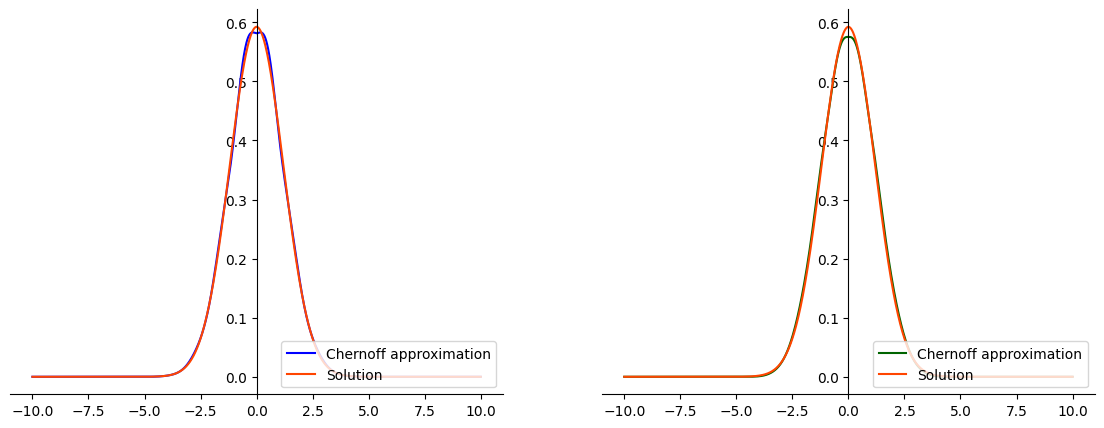}\\
\end{center}
n=4
\begin{center}
	\includegraphics[scale=0.5]{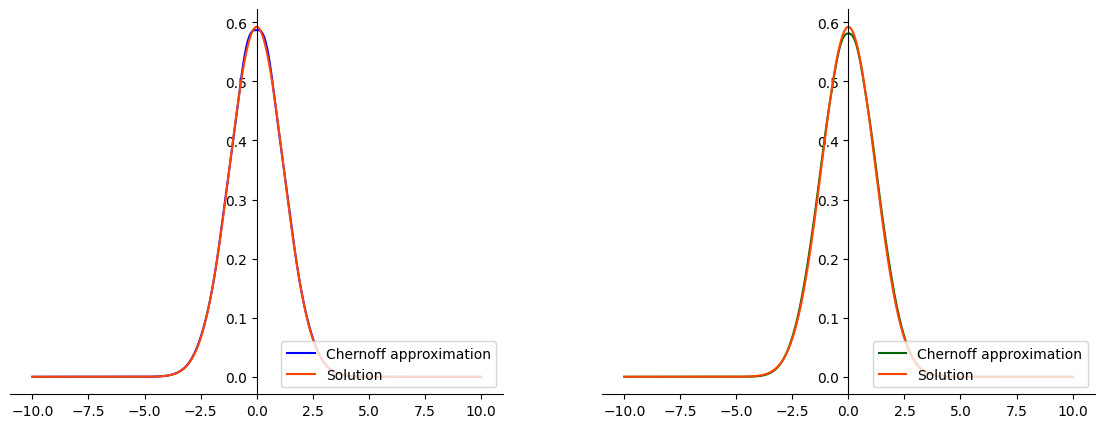}\\
\end{center}
n=5
\begin{center}
	\includegraphics[scale=0.5]{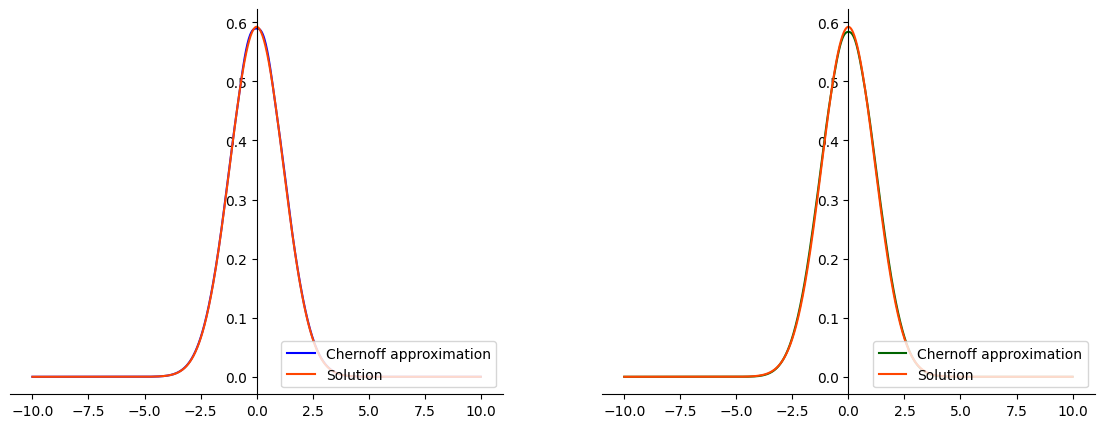}\\
\end{center}

n=6
\begin{center}
	\includegraphics[scale=0.5]{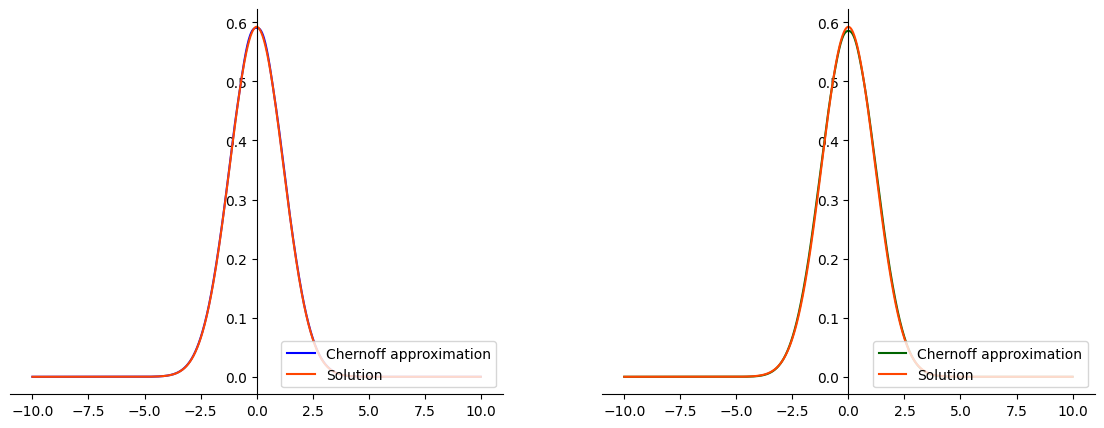}\\
\end{center}
n=7
\begin{center}
	\includegraphics[scale=0.5]{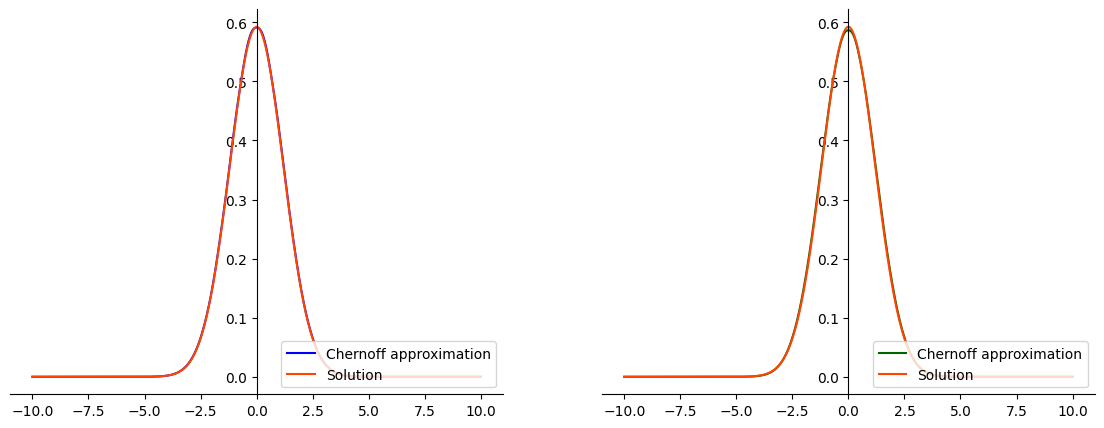}\\
\end{center}
n=8
\begin{center}
	\includegraphics[scale=0.5]{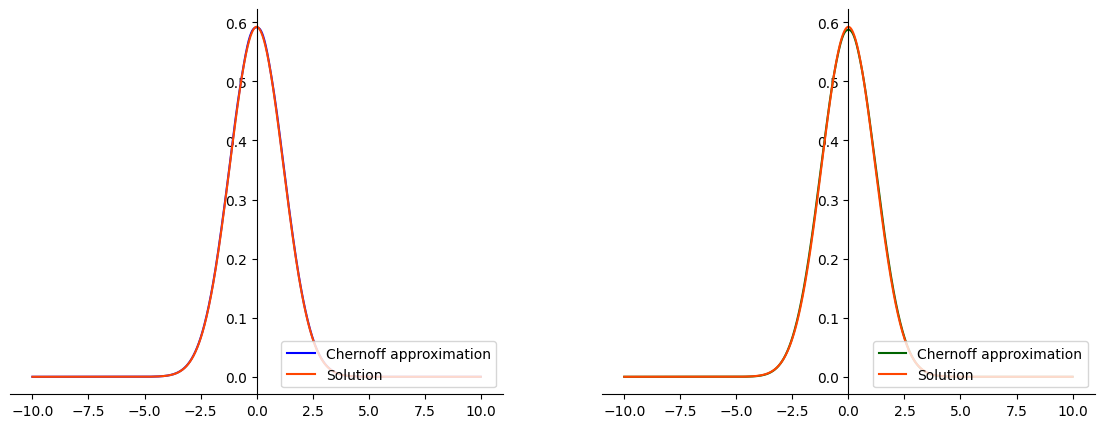}\\
\end{center}

n=9
\begin{center}
	\includegraphics[scale=0.5]{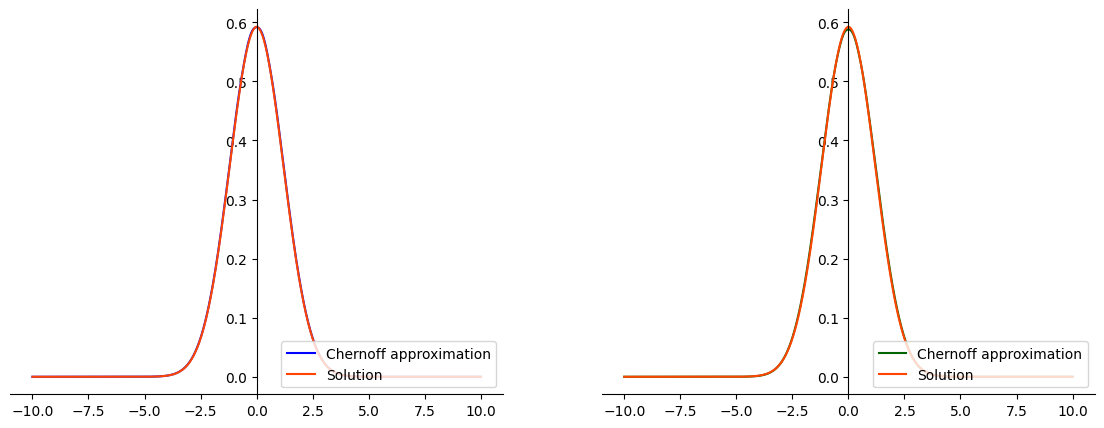}\\
\end{center}
n=10
\begin{center}
	\includegraphics[scale=0.5]{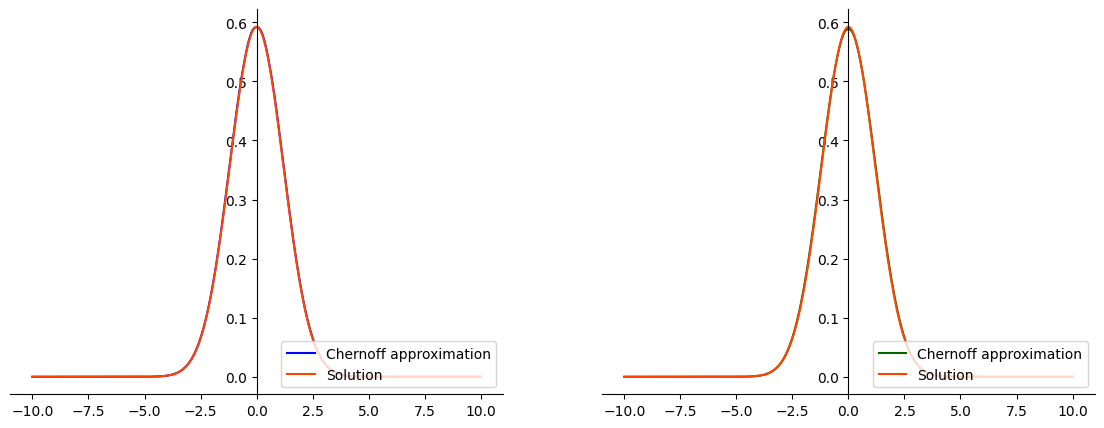}\\
\end{center}

\subsection{$u_0(x)=e^{-|x|^{7/2}}$}

n=1
\begin{center}
	\includegraphics[scale=0.5]{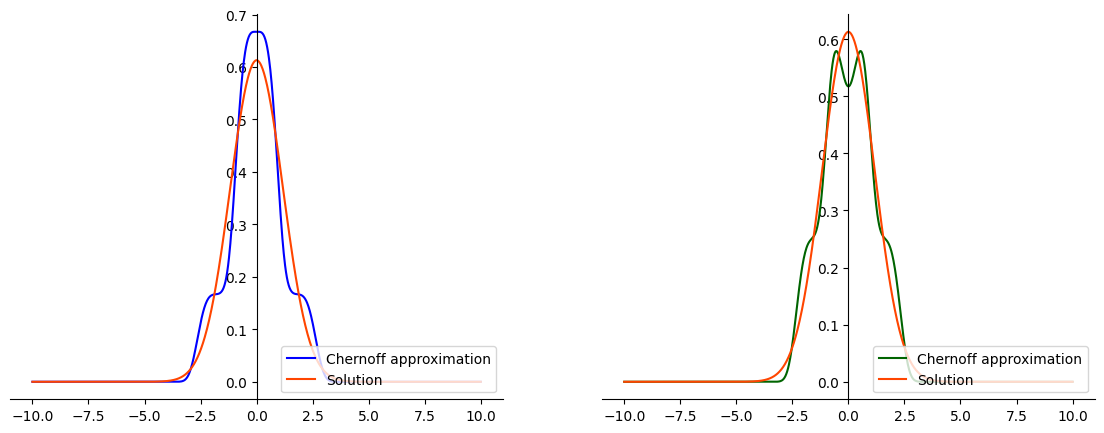}\\
\end{center}
n=2
\begin{center}
	\includegraphics[scale=0.5]{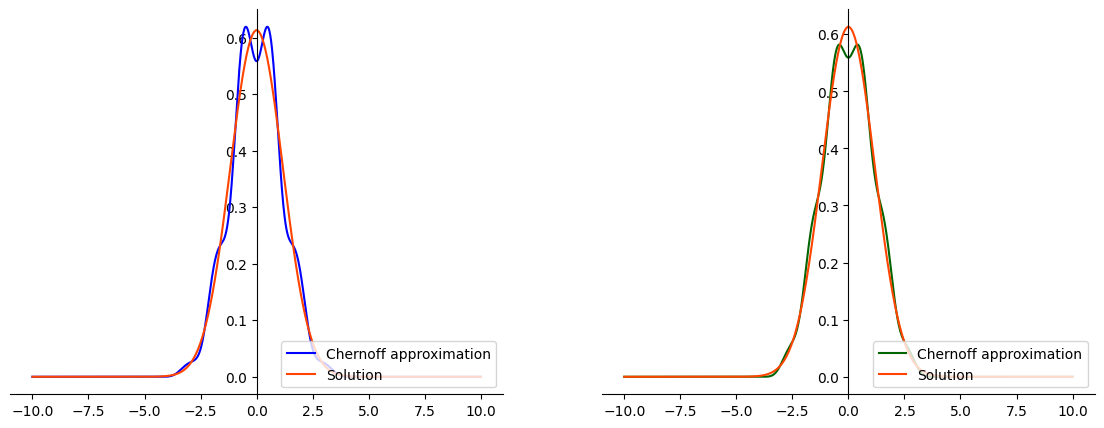}\\
\end{center}

n=3

\begin{center}
	\includegraphics[scale=0.5]{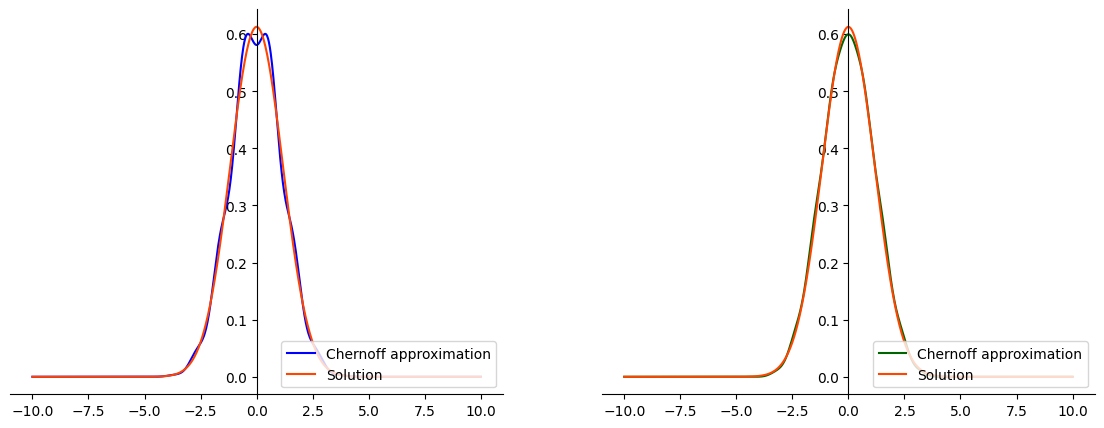}\\
\end{center}
n=4
\begin{center}
	\includegraphics[scale=0.5]{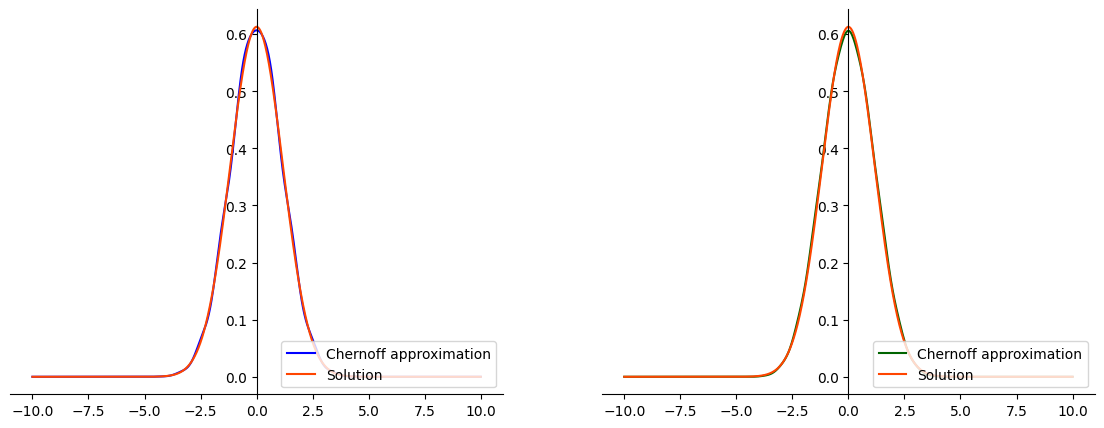}\\
\end{center}
n=5
\begin{center}
	\includegraphics[scale=0.5]{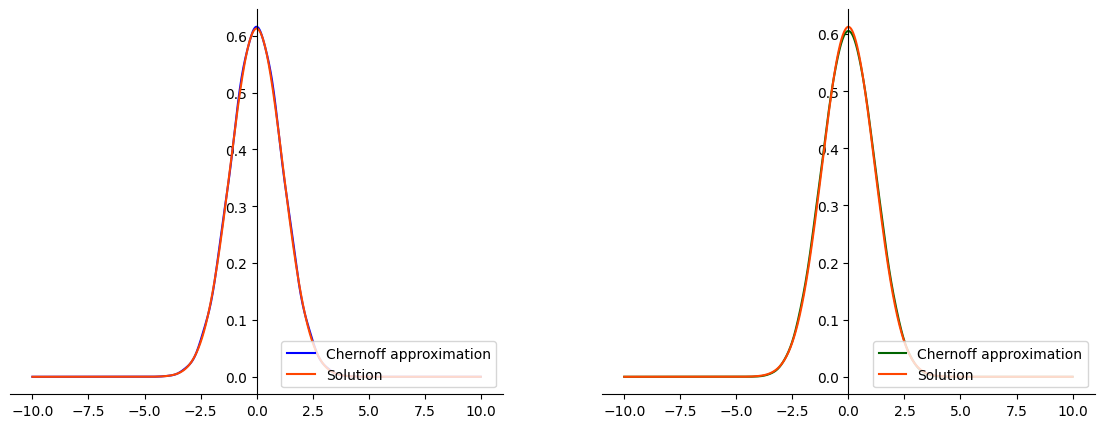}\\
\end{center}

n=6
\begin{center}
	\includegraphics[scale=0.5]{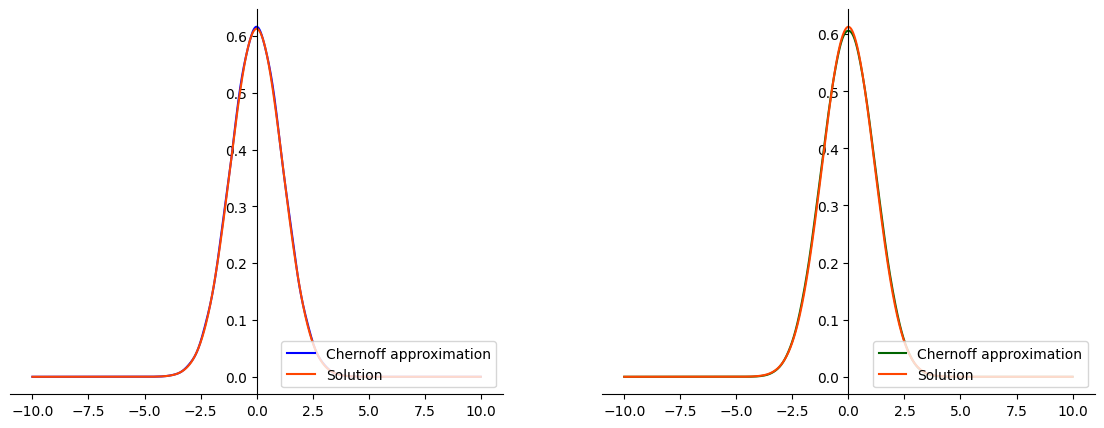}\\
\end{center}
n=7
\begin{center}
	\includegraphics[scale=0.5]{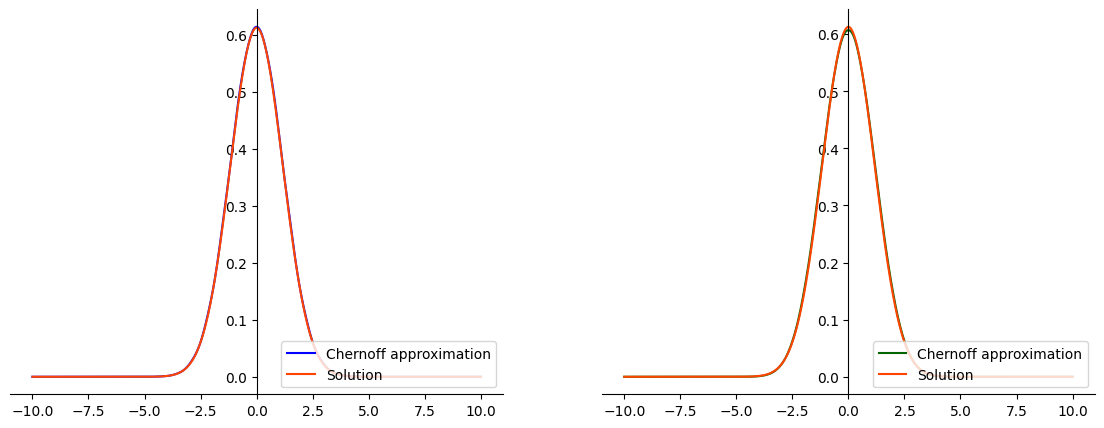}\\
\end{center}
n=8
\begin{center}
	\includegraphics[scale=0.5]{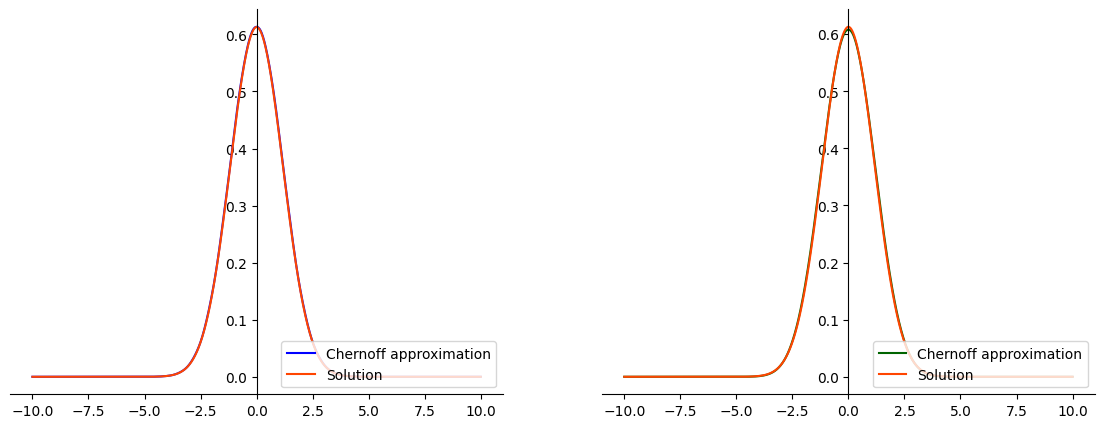}\\
\end{center}

n=9
\begin{center}
	\includegraphics[scale=0.5]{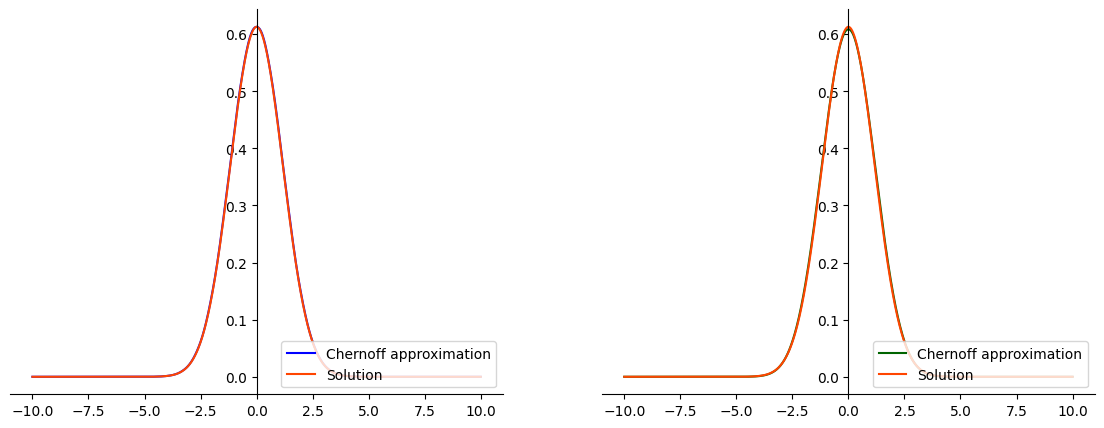}\\
\end{center}
n=10
\begin{center}
	\includegraphics[scale=0.5]{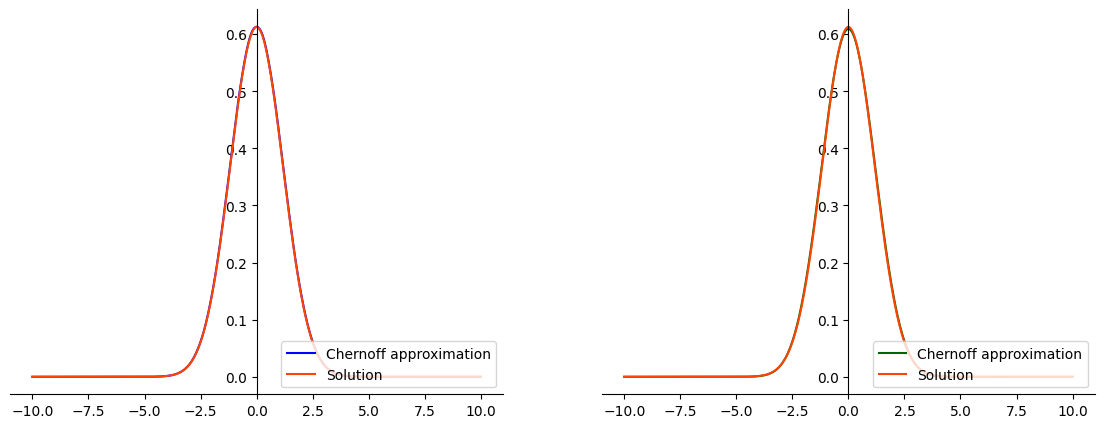}\\
\end{center}

\subsection{$u_0(x)=e^{-|x|^{9/2}}$}

n=1
\begin{center}
	\includegraphics[scale=0.5]{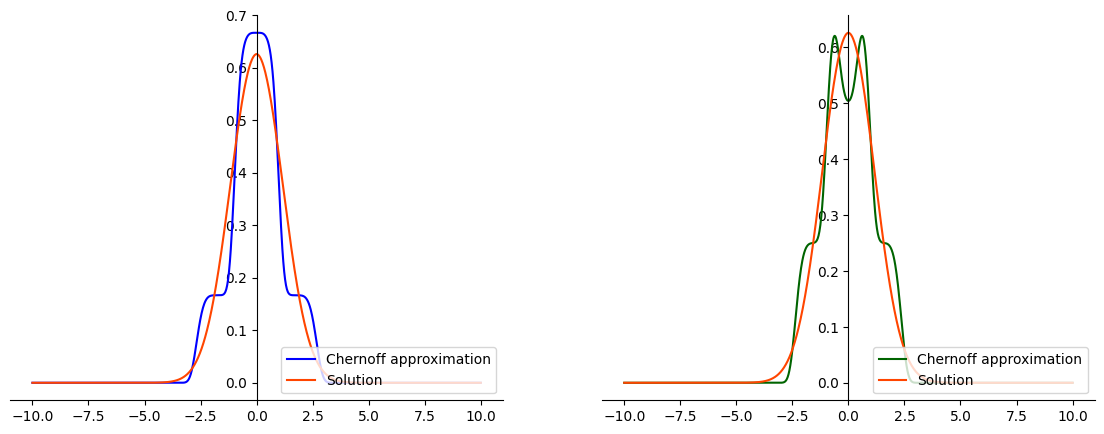}\\
\end{center}
n=2
\begin{center}
	\includegraphics[scale=0.5]{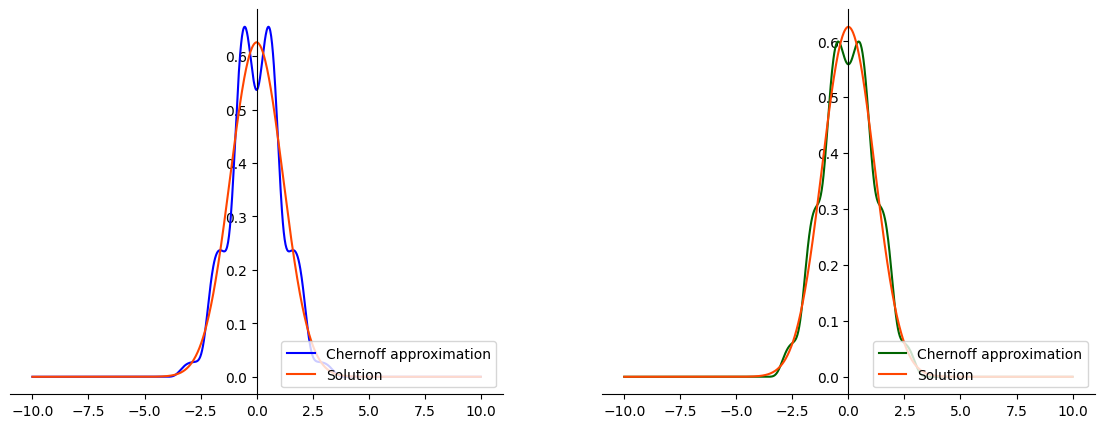}\\
\end{center}

n=3

\begin{center}
	\includegraphics[scale=0.5]{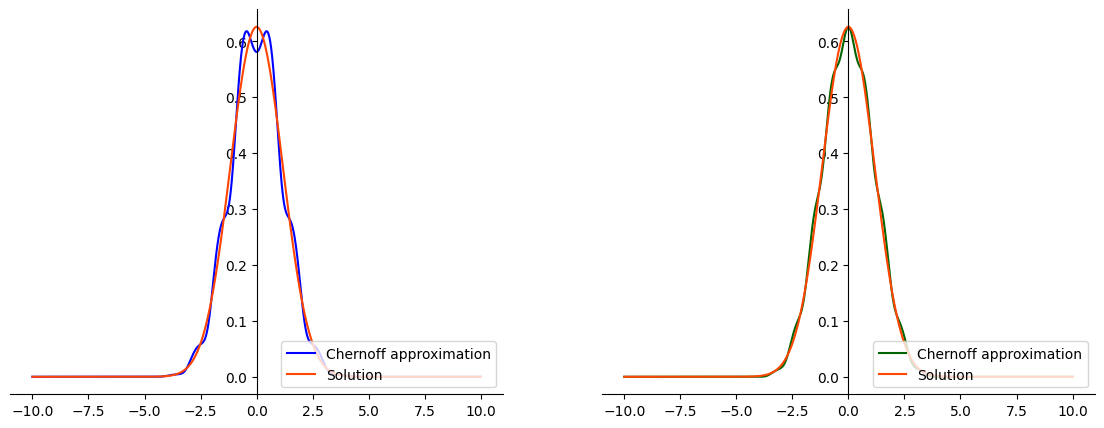}\\
\end{center}
n=4
\begin{center}
	\includegraphics[scale=0.5]{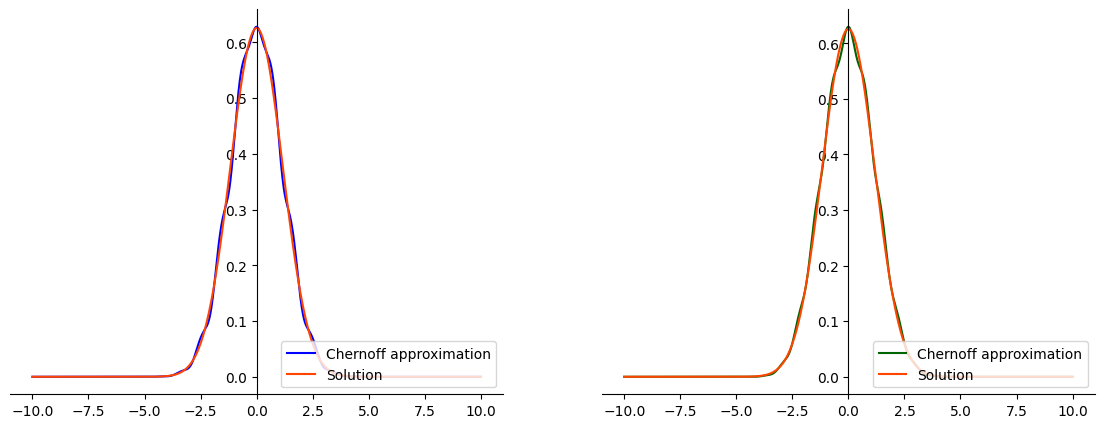}\\
\end{center}
n=5
\begin{center}
	\includegraphics[scale=0.5]{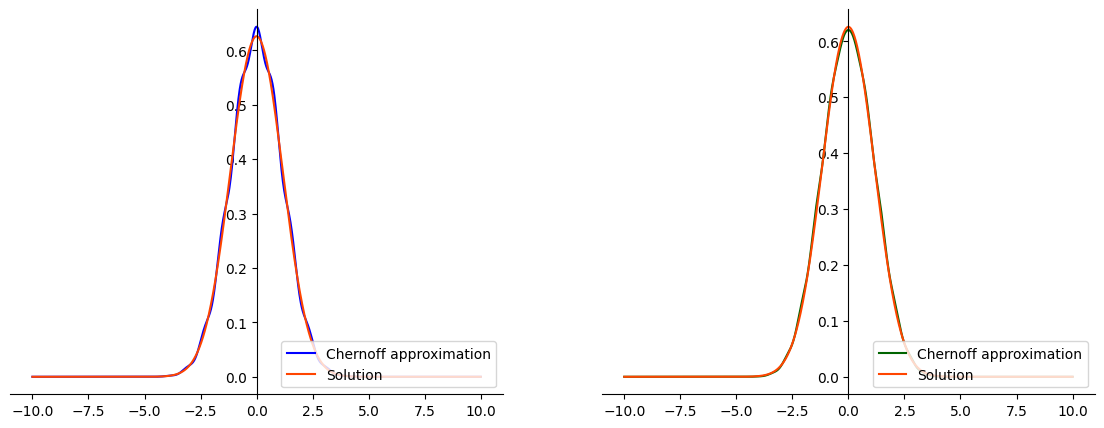}\\
\end{center}

n=6
\begin{center}
	\includegraphics[scale=0.5]{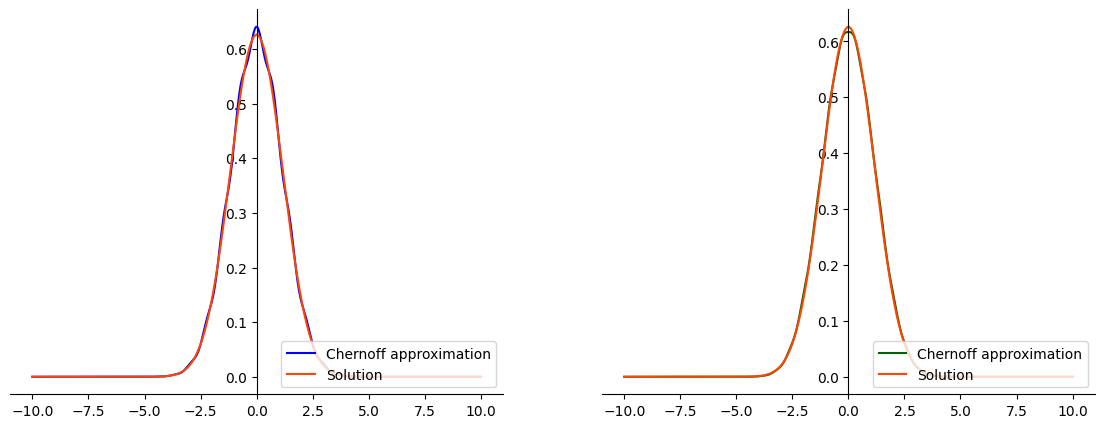}\\
\end{center}
n=7
\begin{center}
	\includegraphics[scale=0.5]{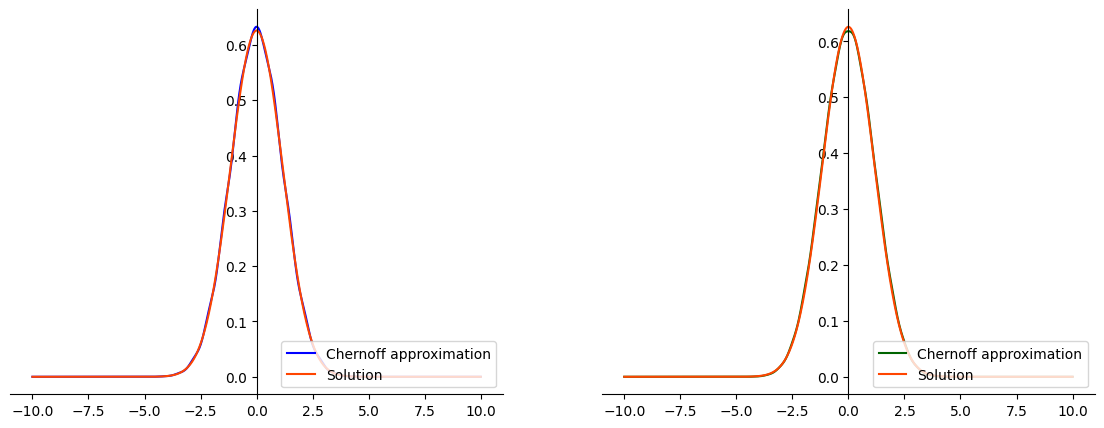}\\
\end{center}
n=8
\begin{center}
	\includegraphics[scale=0.5]{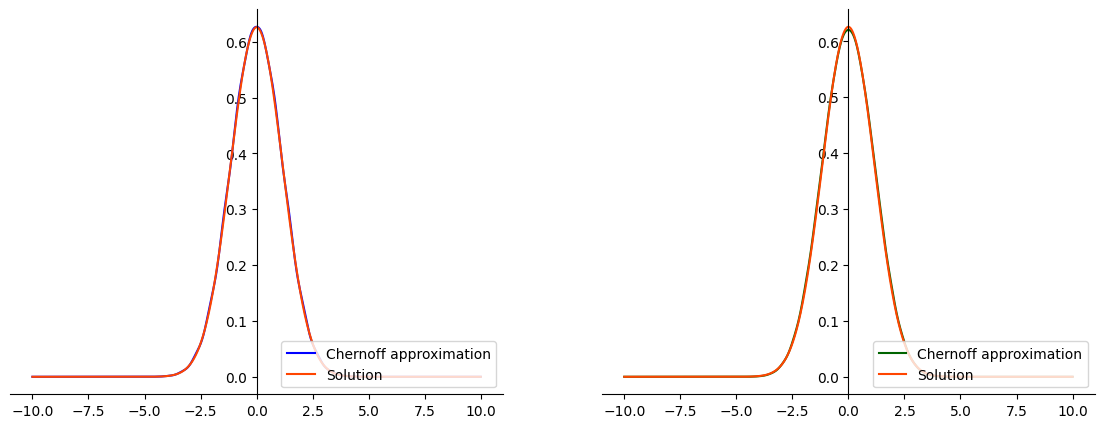}\\
\end{center}

n=9
\begin{center}
	\includegraphics[scale=0.5]{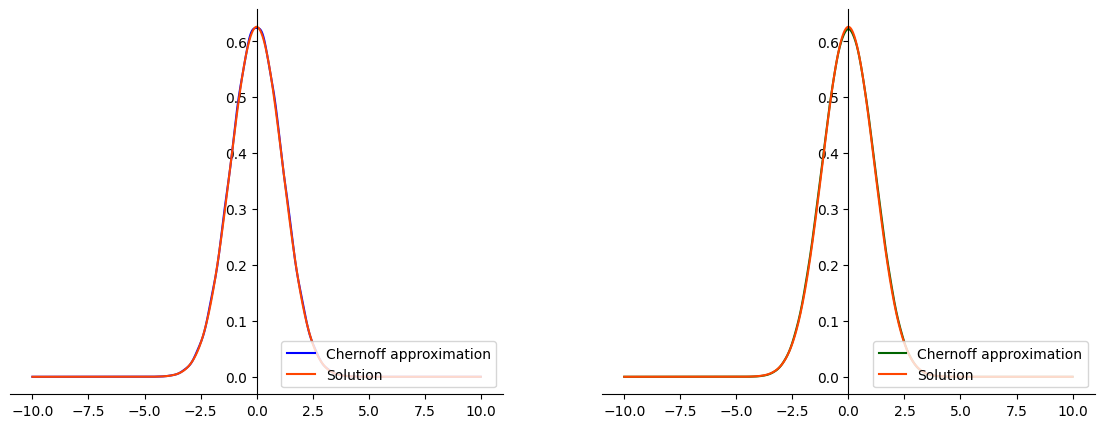}\\
\end{center}
n=10
\begin{center}
	\includegraphics[scale=0.5]{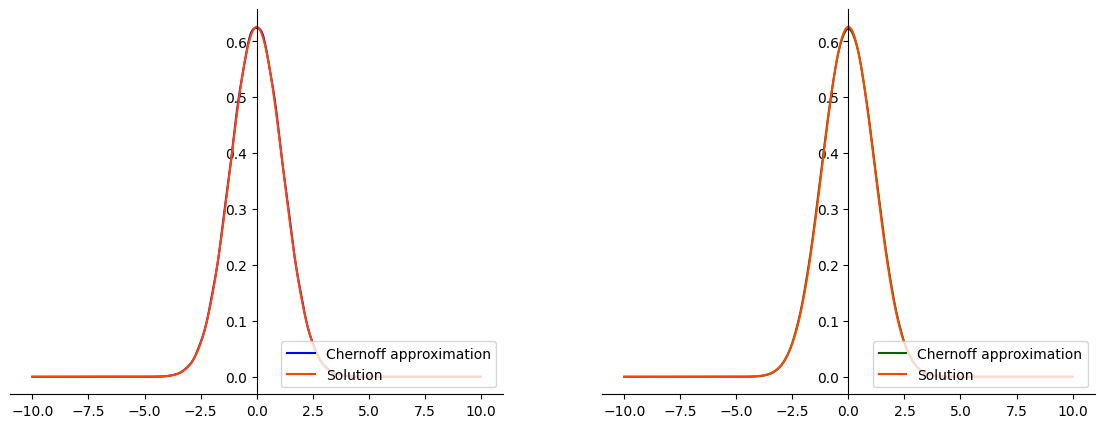}\\
\end{center}

\subsection{$u_0(x)= e^{-|x|\cdot x^{2}}$}

n=1
\begin{center}
	\includegraphics[scale=0.5]{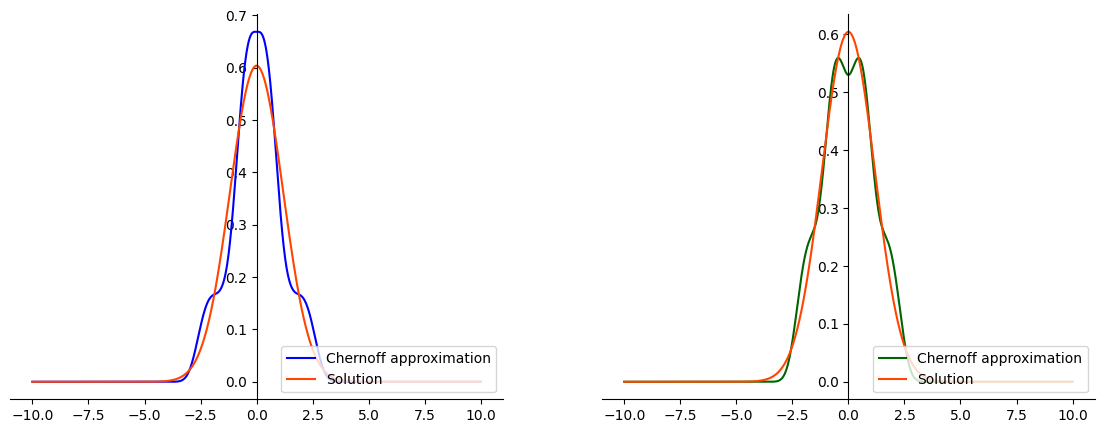}\\
\end{center}
n=2
\begin{center}
	\includegraphics[scale=0.5]{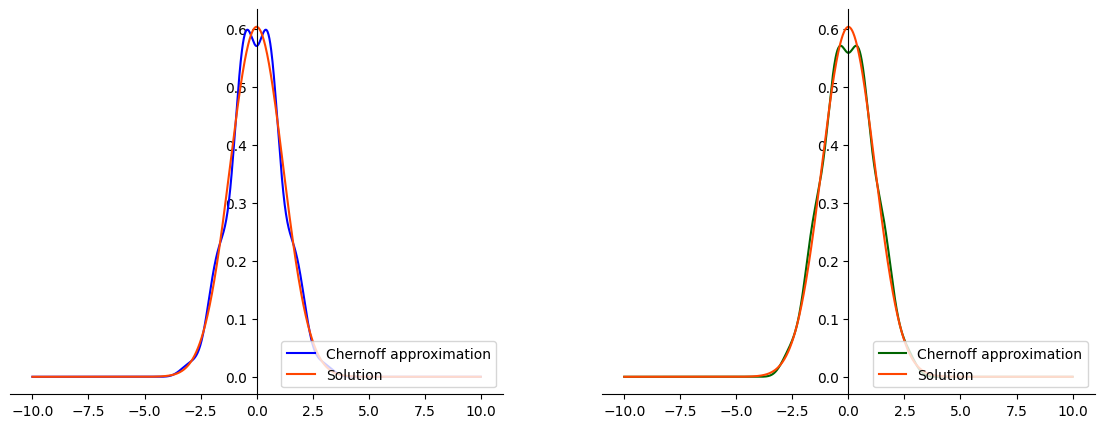}\\
\end{center}

n=3

\begin{center}
	\includegraphics[scale=0.5]{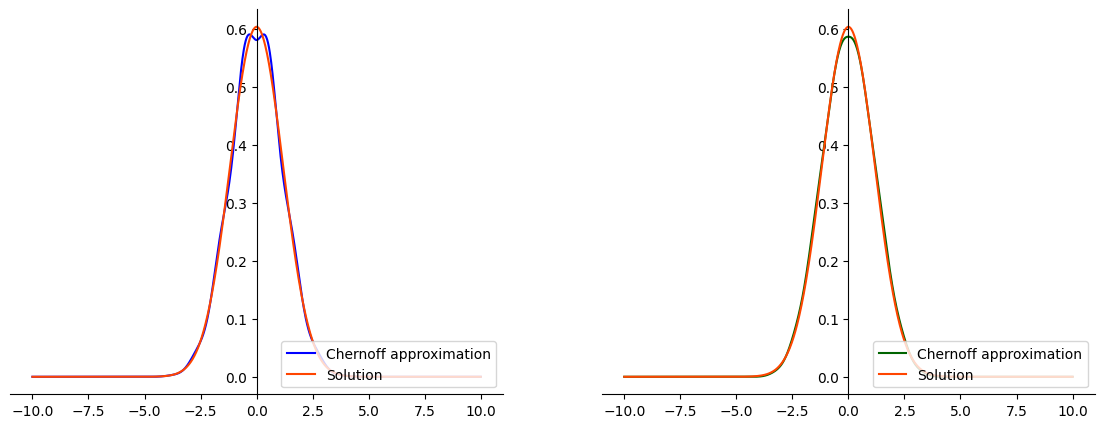}\\
\end{center}
n=4
\begin{center}
	\includegraphics[scale=0.5]{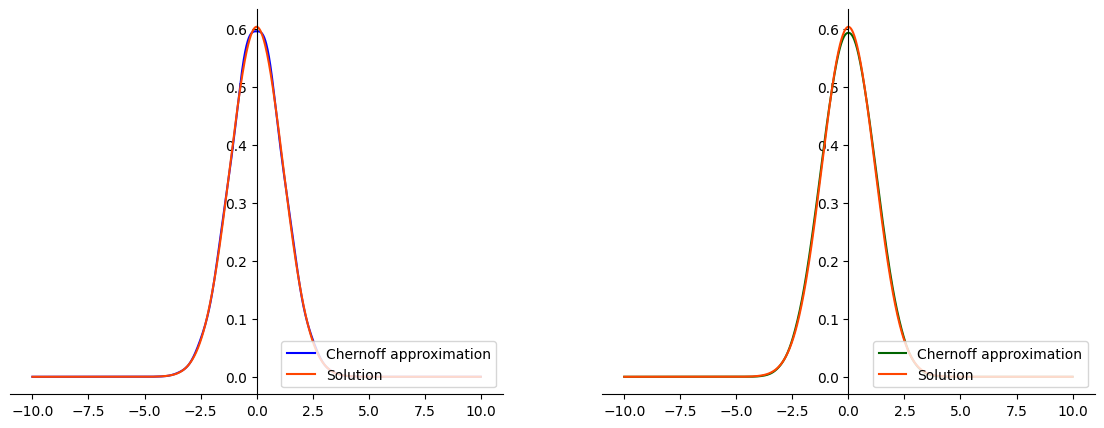}\\
\end{center}
n=5
\begin{center}
	\includegraphics[scale=0.5]{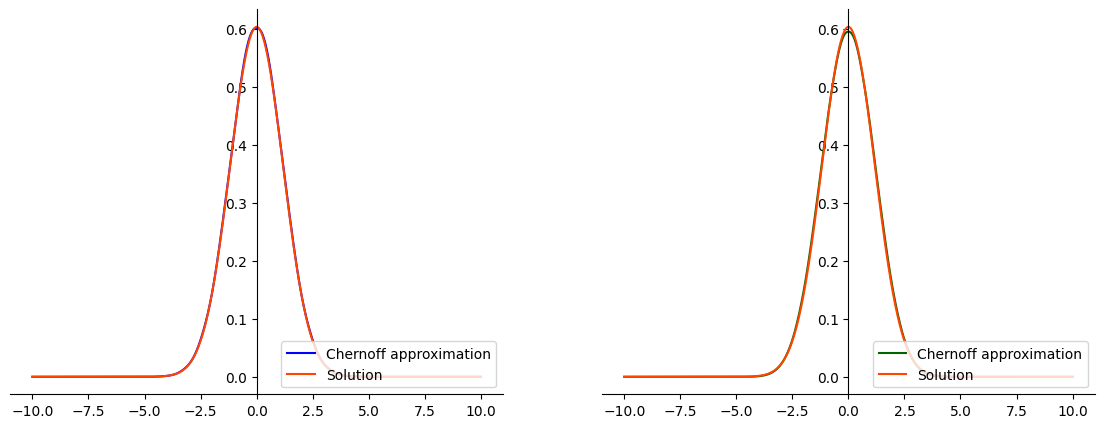}\\
\end{center}

n=6
\begin{center}
	\includegraphics[scale=0.5]{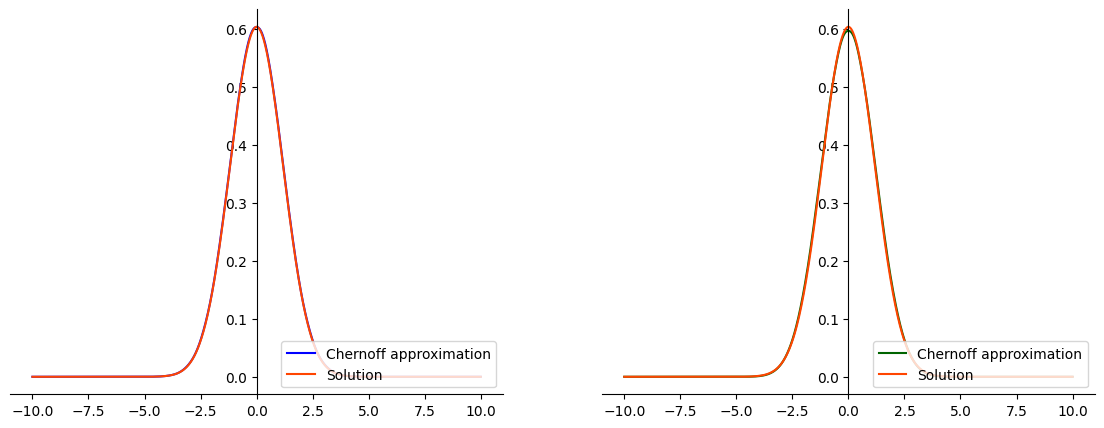}\\
\end{center}
n=7
\begin{center}
	\includegraphics[scale=0.5]{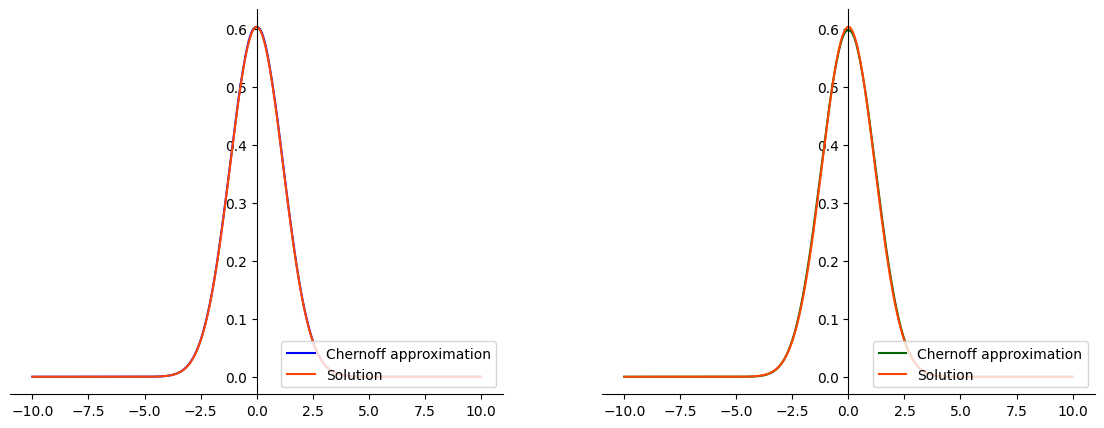}\\
\end{center}
n=8
\begin{center}
	\includegraphics[scale=0.5]{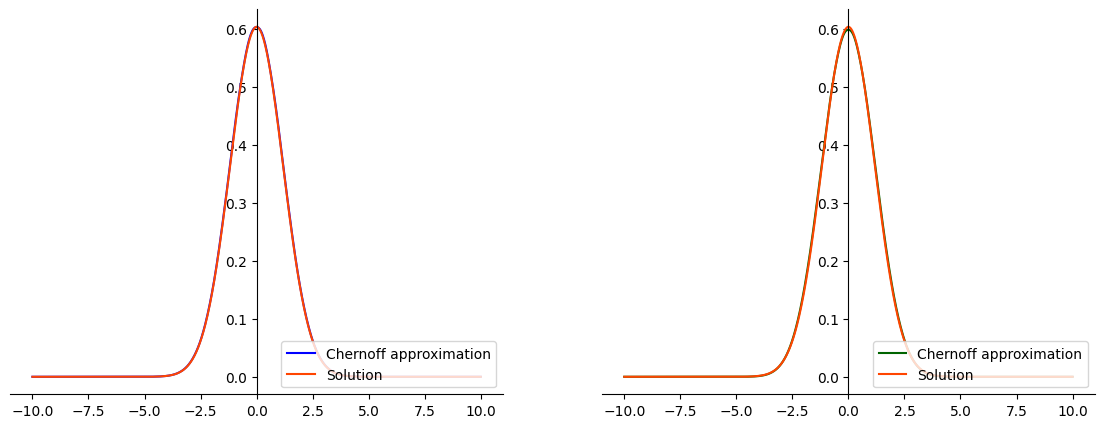}\\
\end{center}

n=9
\begin{center}
	\includegraphics[scale=0.5]{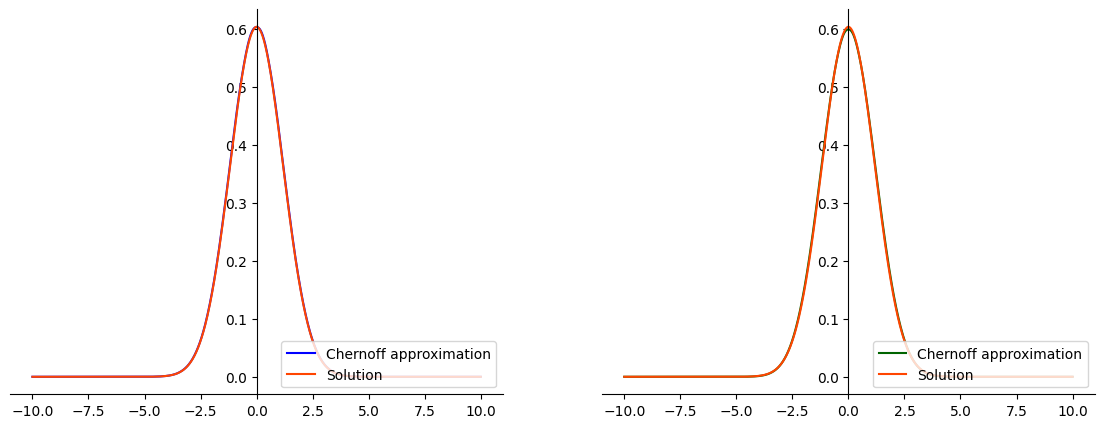}\\
\end{center}
n=10
\begin{center}
	\includegraphics[scale=0.5]{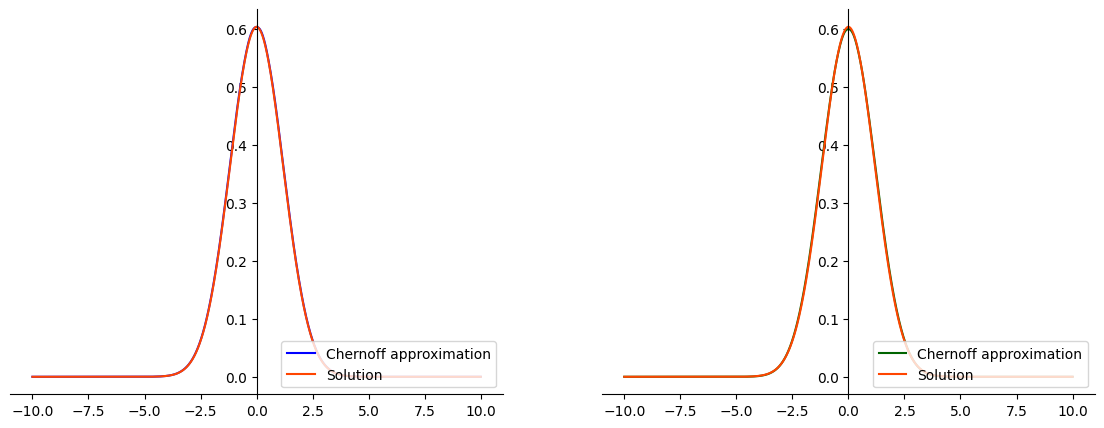}\\
\end{center}

\subsection{$u_0(x)= e^{-x^2}\cdot\sin(x)$}

n=1
\begin{center}
	\includegraphics[scale=0.5]{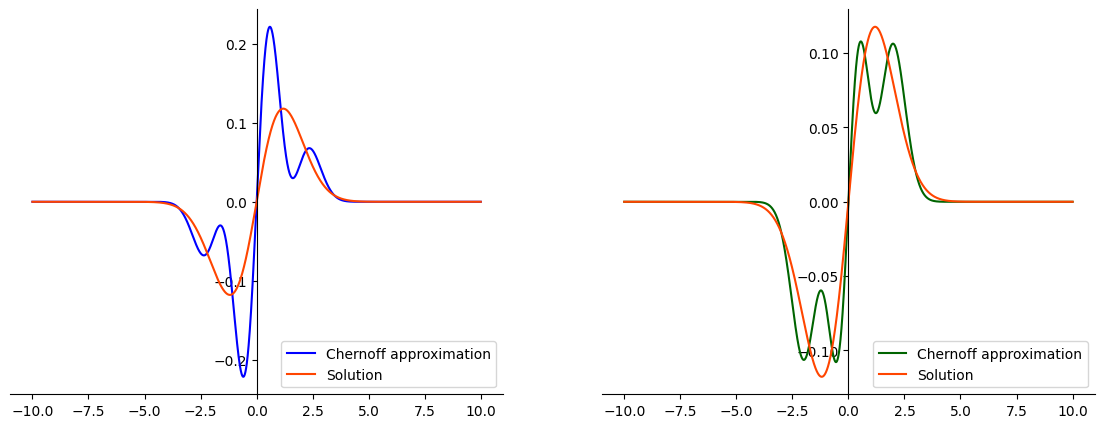}\\
\end{center}
n=2
\begin{center}
	\includegraphics[scale=0.5]{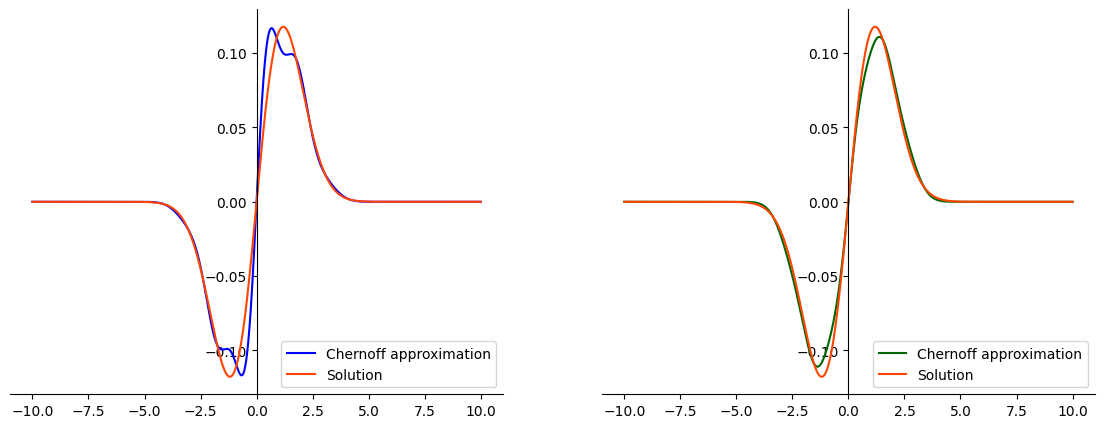}\\
\end{center}

n=3

\begin{center}
	\includegraphics[scale=0.5]{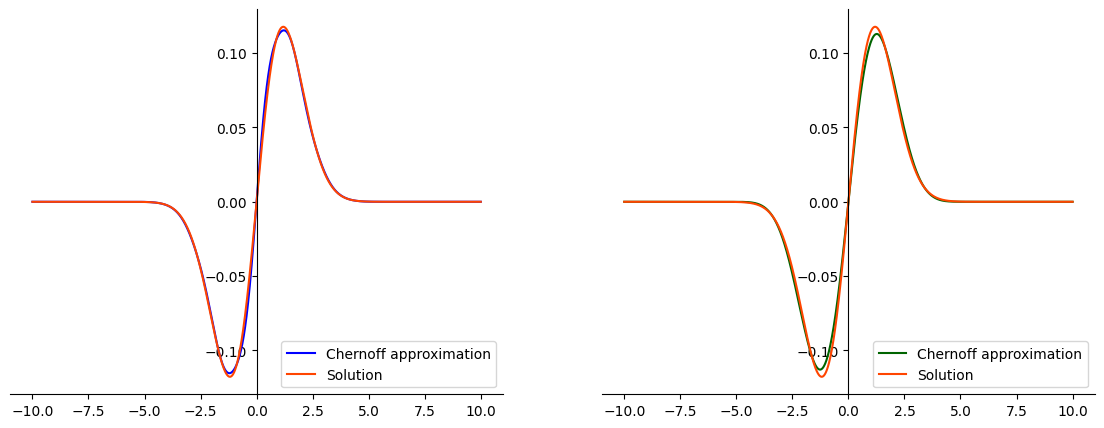}\\
\end{center}
n=4
\begin{center}
	\includegraphics[scale=0.5]{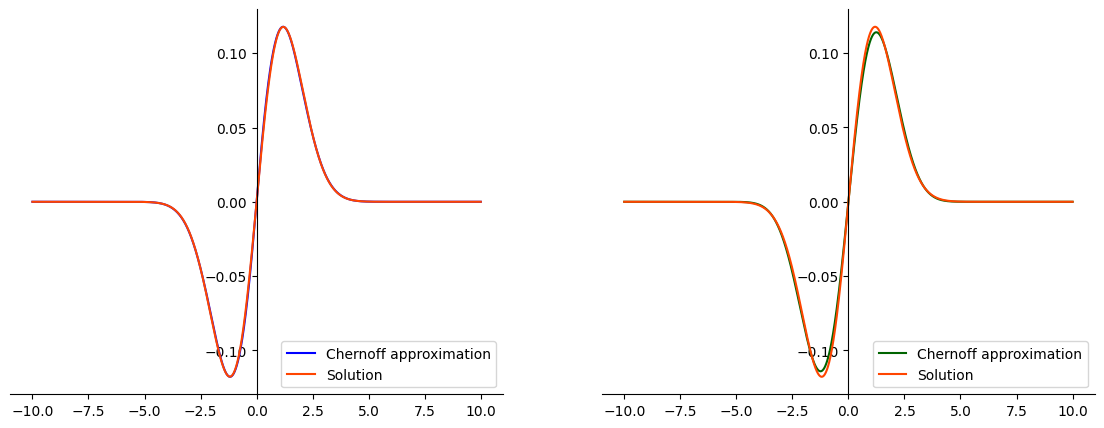}\\
\end{center}
n=5
\begin{center}
	\includegraphics[scale=0.5]{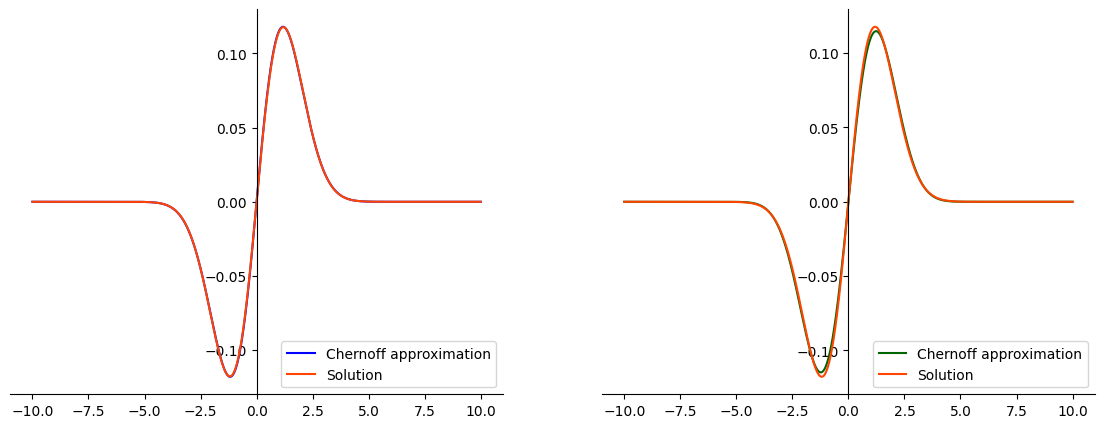}\\
\end{center}

n=6
\begin{center}
	\includegraphics[scale=0.5]{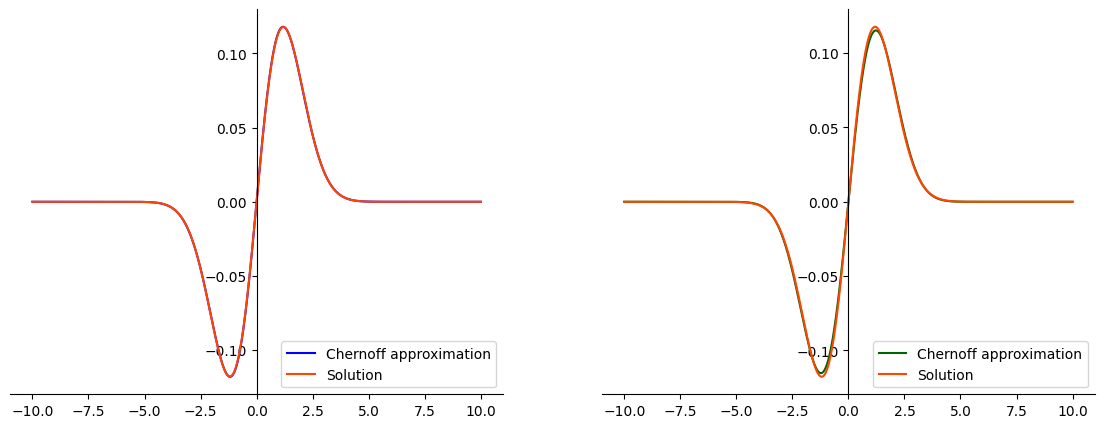}\\
\end{center}
n=7
\begin{center}
	\includegraphics[scale=0.5]{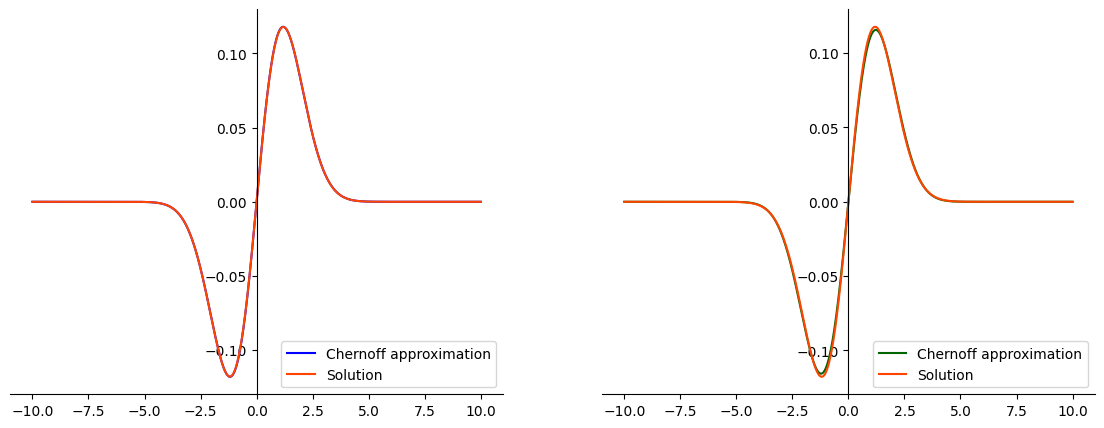}\\
\end{center}
n=8
\begin{center}
	\includegraphics[scale=0.5]{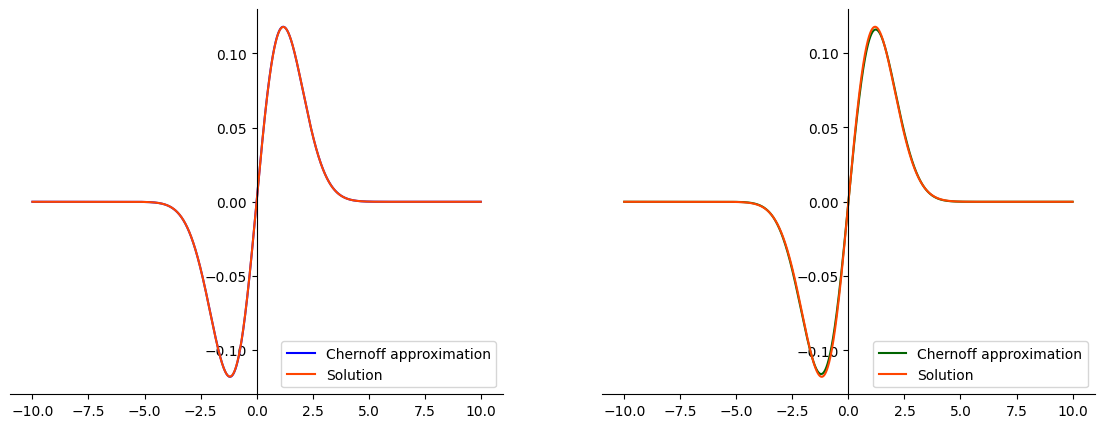}\\
\end{center}

n=9
\begin{center}
	\includegraphics[scale=0.5]{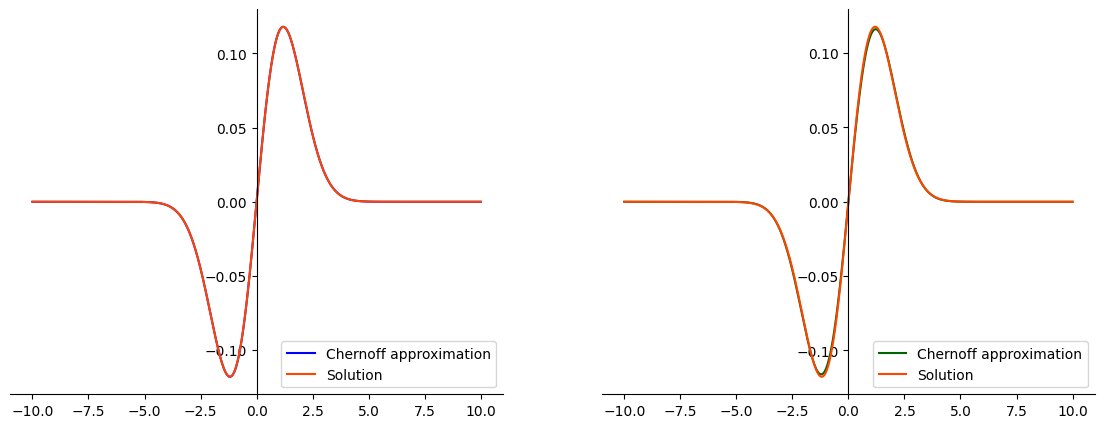}\\
\end{center}
n=10
\begin{center}
	\includegraphics[scale=0.5]{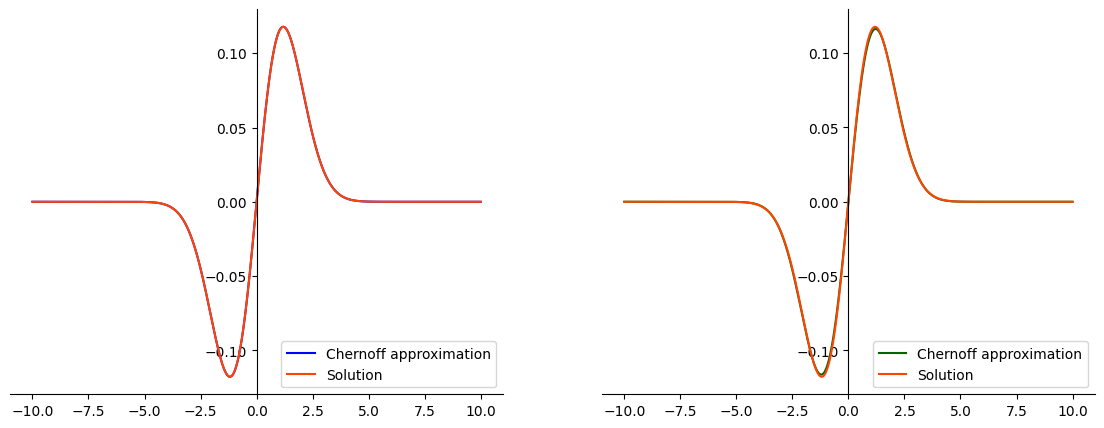}\\
\end{center}

\subsection{$u_0(x)= e^{x|x|}\cdot e^{x^4}$}

n=1
\begin{center}
	\includegraphics[scale=0.5]{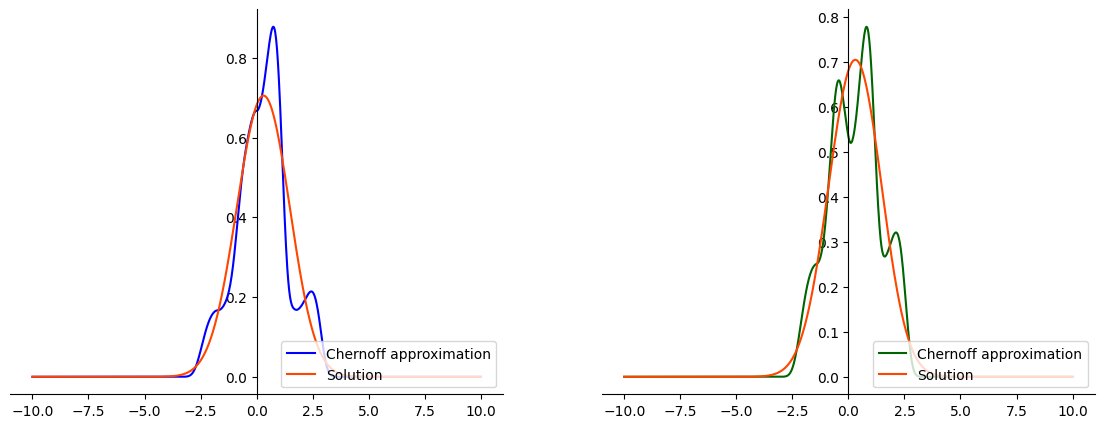}\\
\end{center}
n=2
\begin{center}
	\includegraphics[scale=0.5]{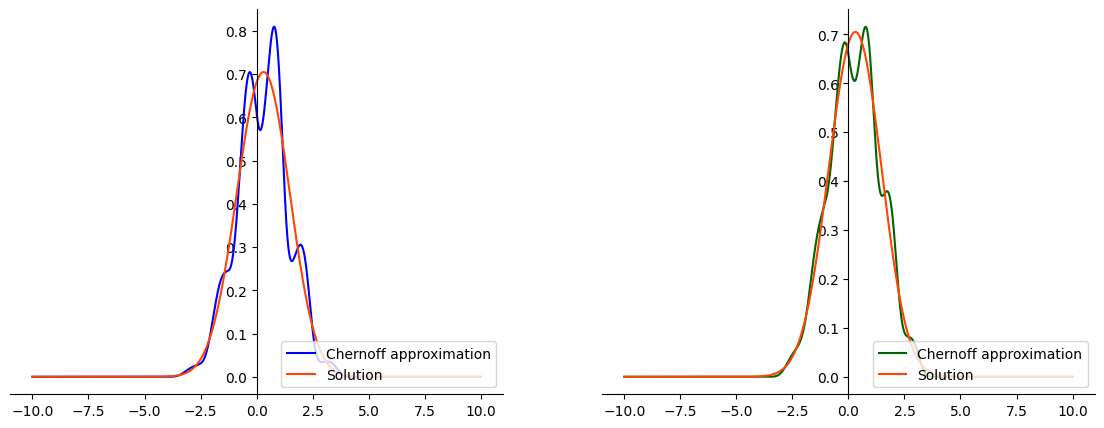}\\
\end{center}

n=3

\begin{center}
	\includegraphics[scale=0.5]{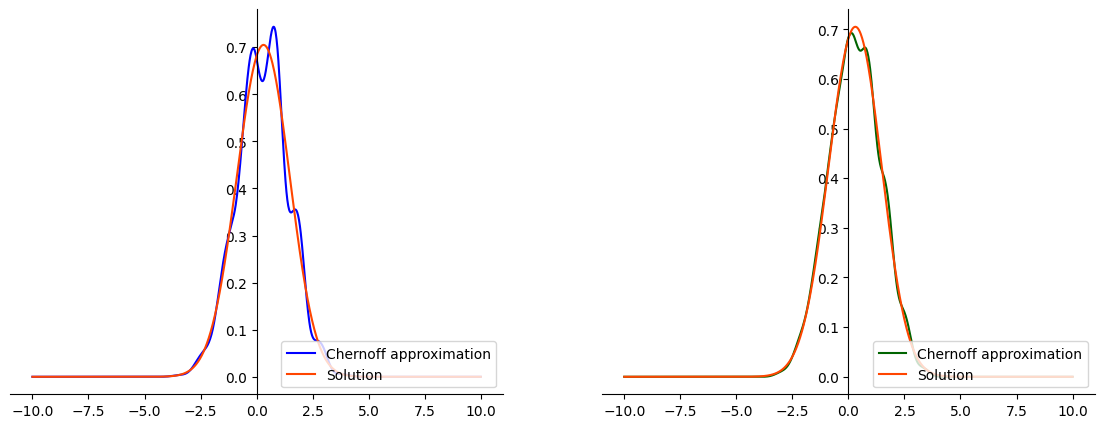}\\
\end{center}
n=4
\begin{center}
	\includegraphics[scale=0.5]{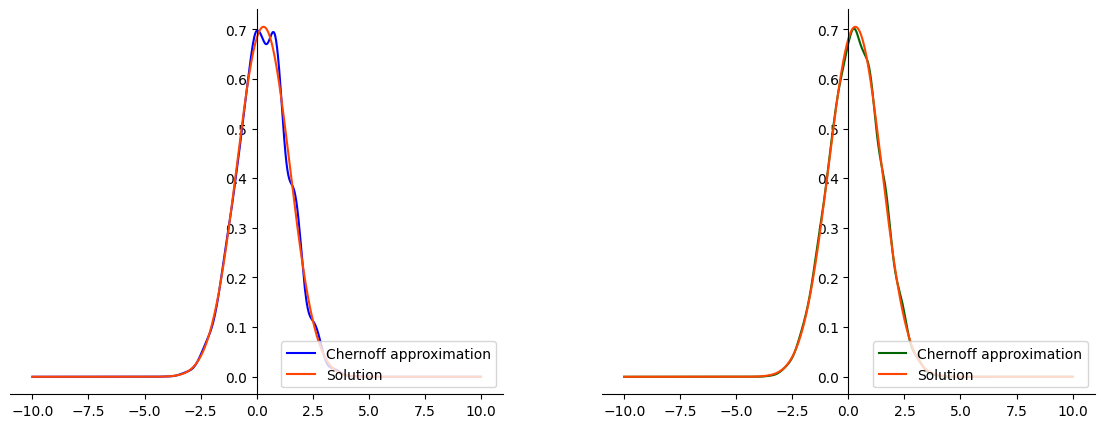}\\
\end{center}
n=5
\begin{center}
	\includegraphics[scale=0.5]{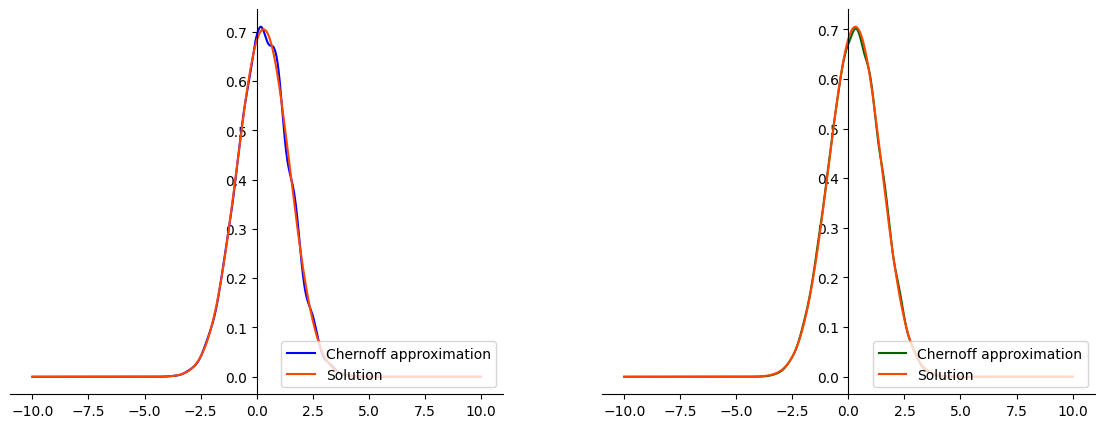}\\
\end{center}

n=6
\begin{center}
	\includegraphics[scale=0.5]{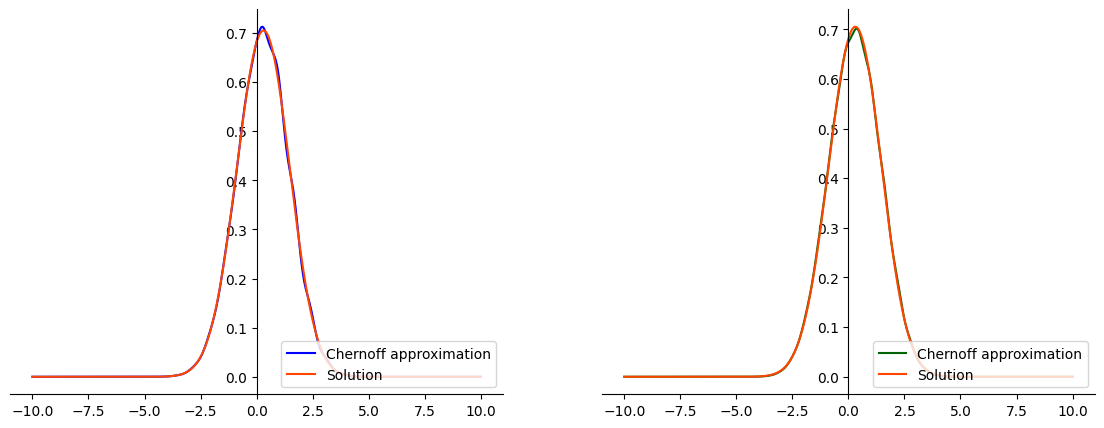}\\
\end{center}
n=7
\begin{center}
	\includegraphics[scale=0.5]{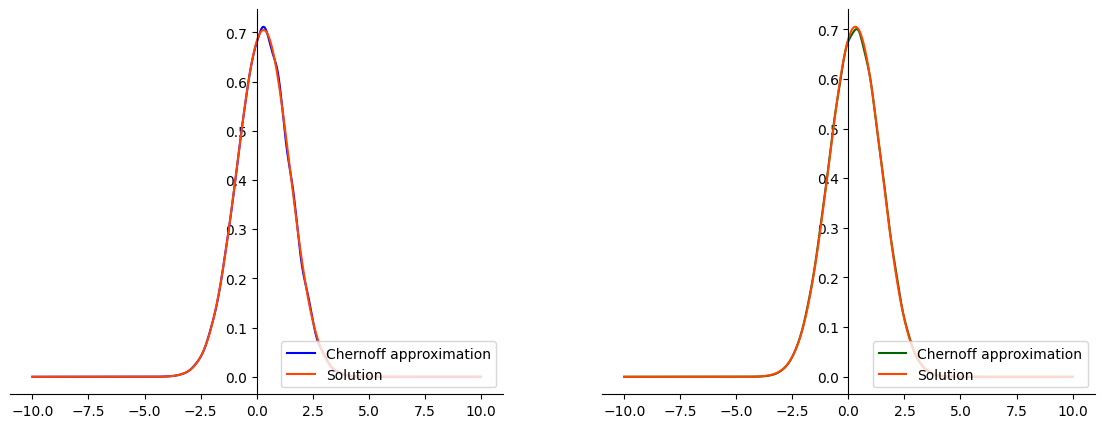}\\
\end{center}
n=8
\begin{center}
	\includegraphics[scale=0.5]{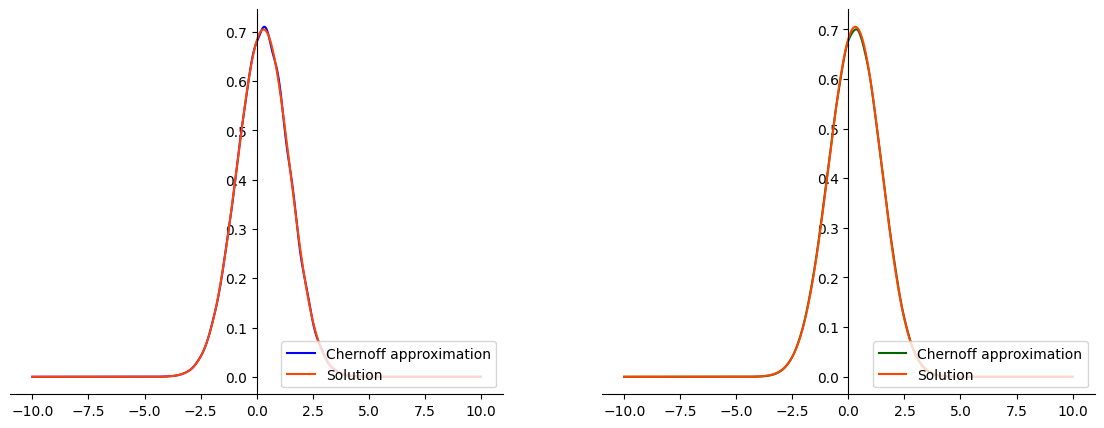}\\
\end{center}

n=9
\begin{center}
	\includegraphics[scale=0.5]{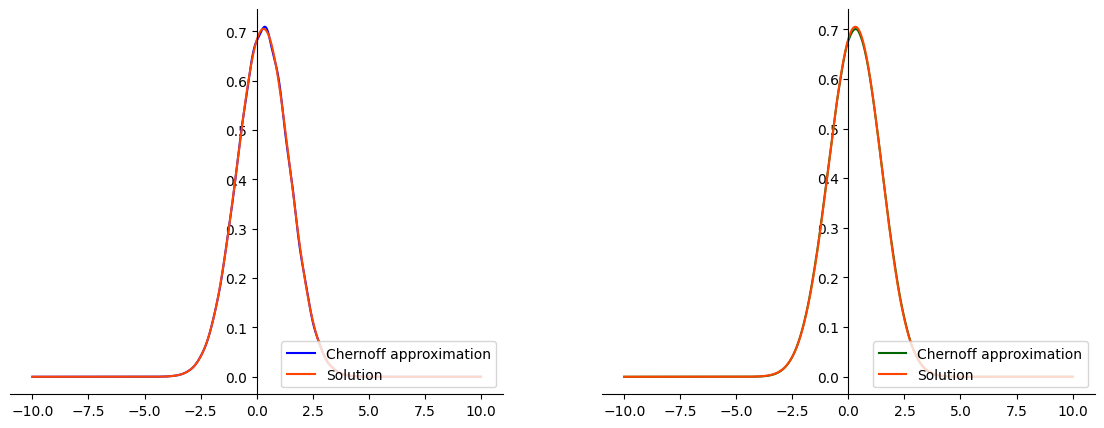}\\
\end{center}
n=10
\begin{center}
	\includegraphics[scale=0.5]{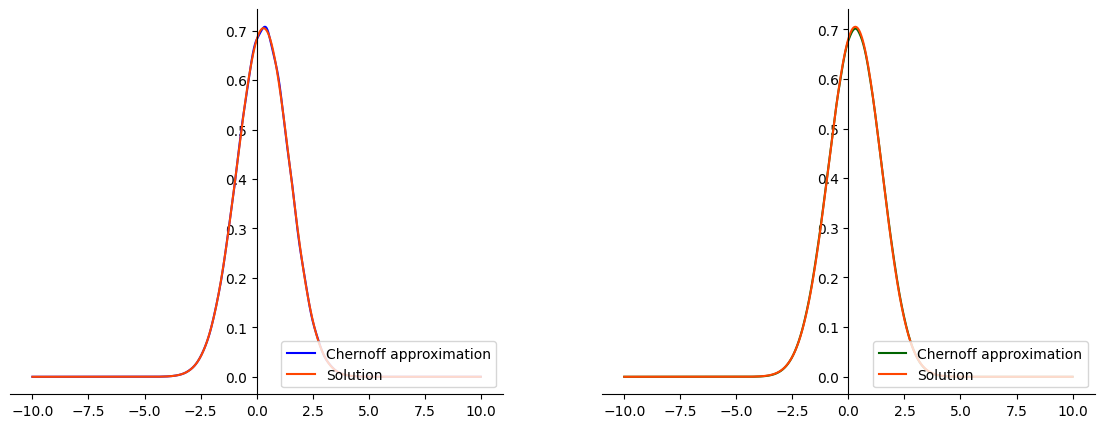}\\
\end{center}
\newpage

\subsection{$u_0(x)= e^{-x^2}$}

n=1
\begin{center}
	\includegraphics[scale=0.5]{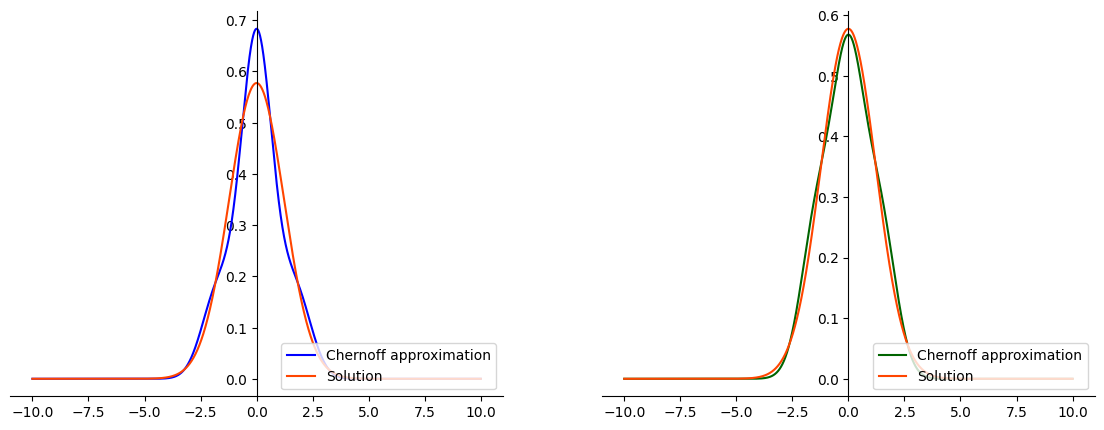}\\
\end{center}
n=2
\begin{center}
	\includegraphics[scale=0.5]{e2_2.png}\\
\end{center}
n=3
\begin{center}
	\includegraphics[scale=0.5]{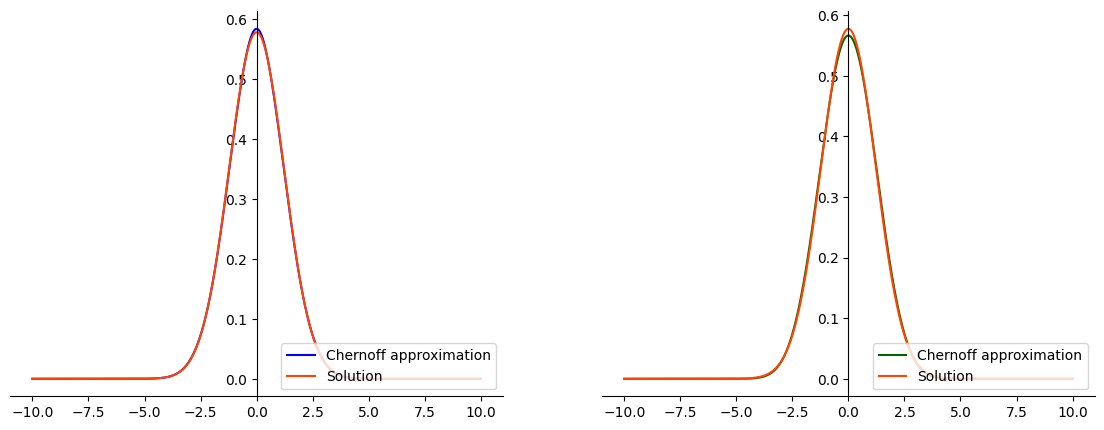}\\
\end{center}
n=4
\begin{center}
	\includegraphics[scale=0.5]{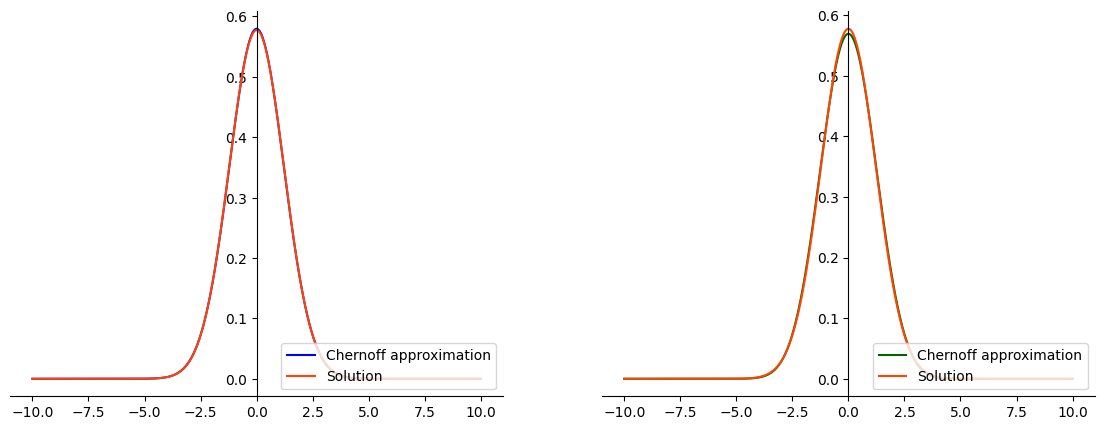}\\
\end{center}
n=5
\begin{center}
	\includegraphics[scale=0.5]{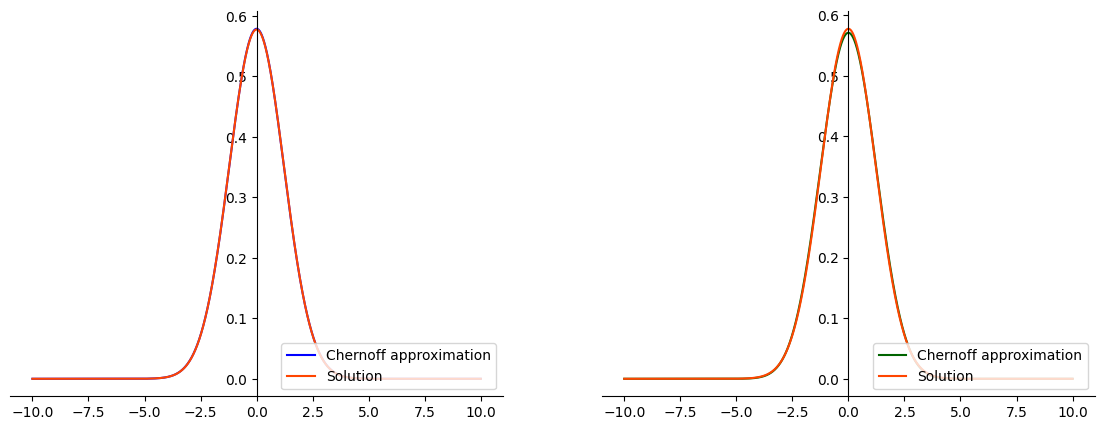}\\
\end{center}
n=6
\begin{center}
	\includegraphics[scale=0.5]{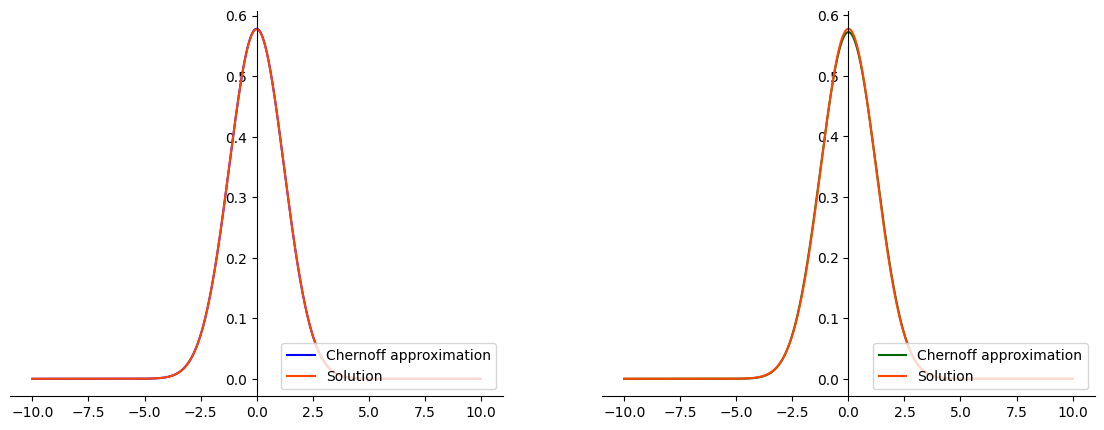}\\
\end{center}
n=7
\begin{center}
	\includegraphics[scale=0.5]{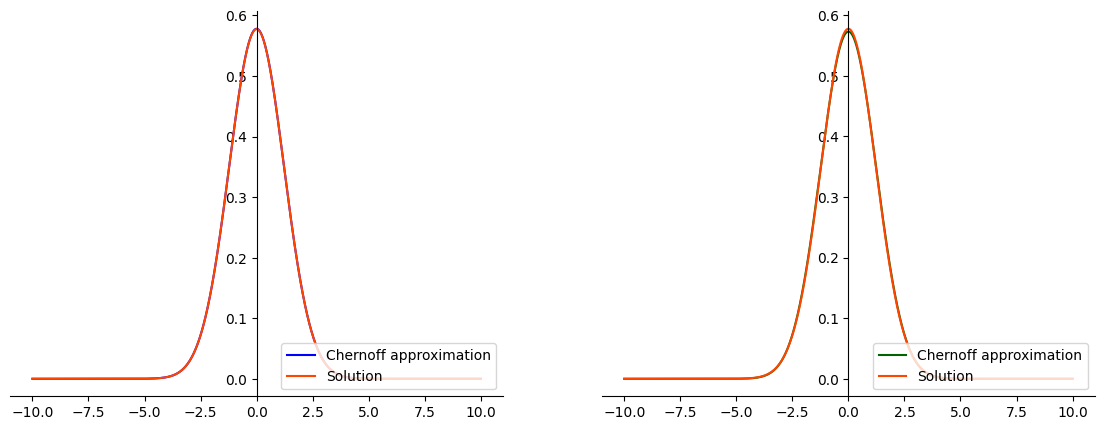}\\
\end{center}
n=8
\begin{center}
	\includegraphics[scale=0.5]{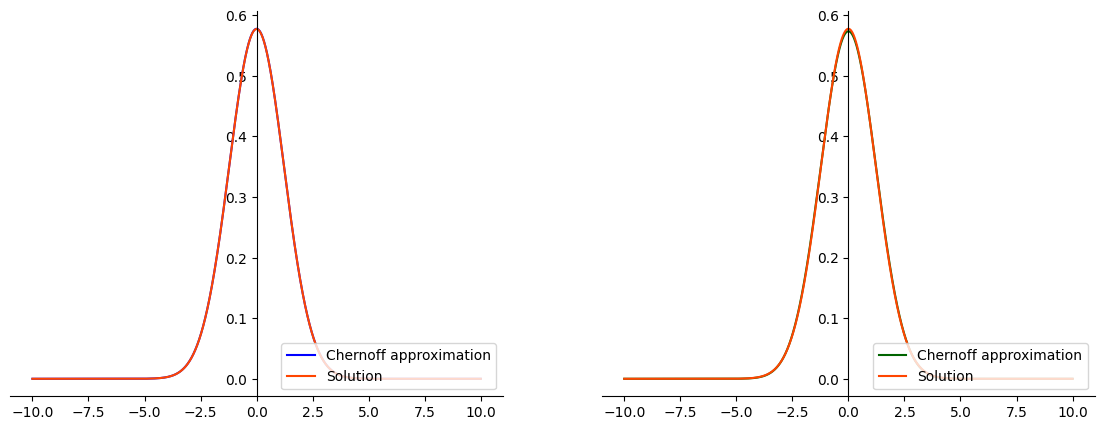}\\
\end{center}
n=9
\begin{center}
	\includegraphics[scale=0.5]{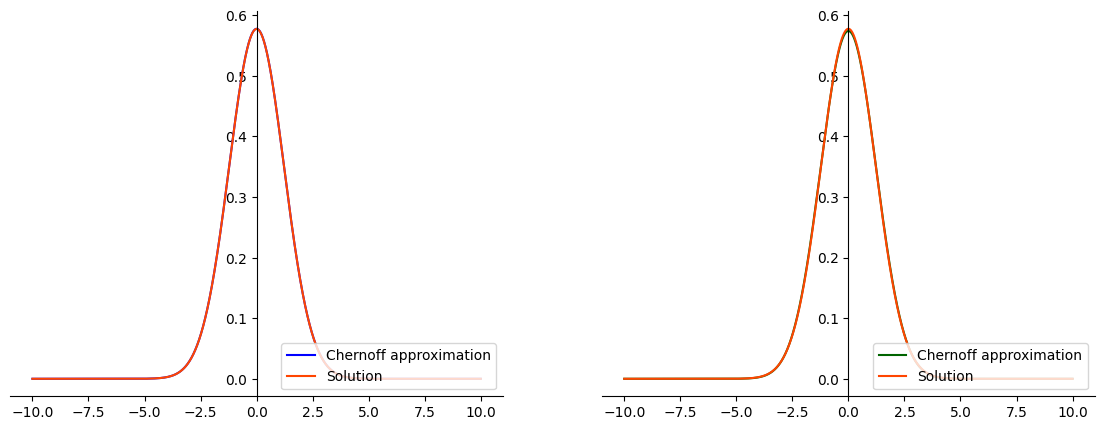}\\
\end{center}
n=10
\begin{center}
	\includegraphics[scale=0.5]{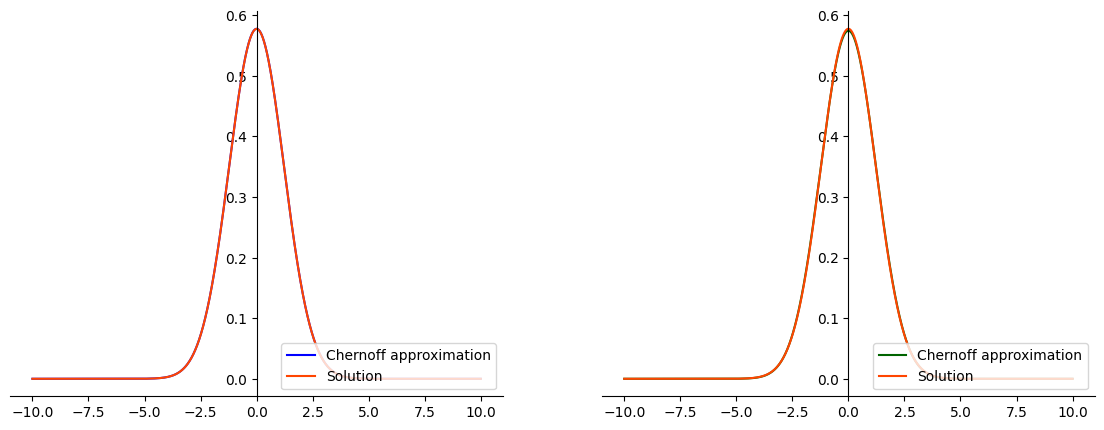}\\
\end{center}

\newpage
\subsection{$u_0(x)= e^{-x^4}$}

n=1
\begin{center}
	\includegraphics[scale=0.5]{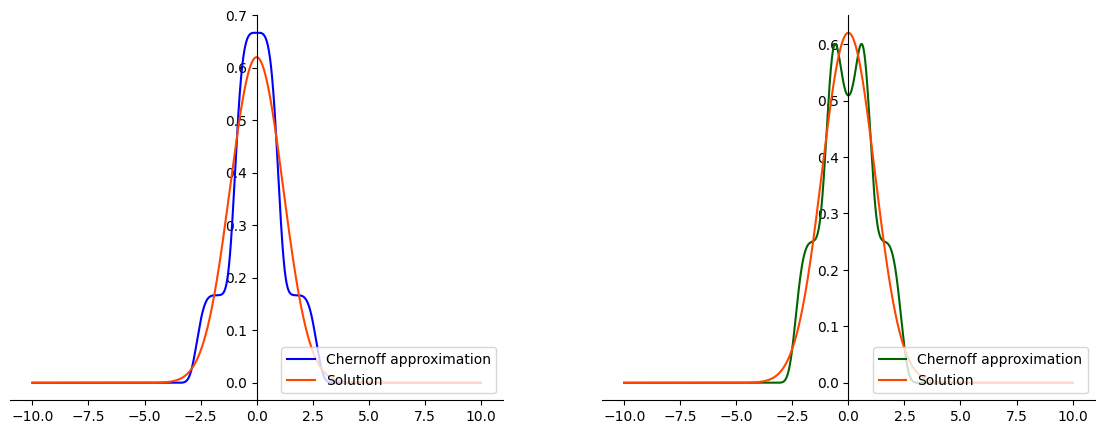}\\
\end{center}
n=2
\begin{center}
	\includegraphics[scale=0.5]{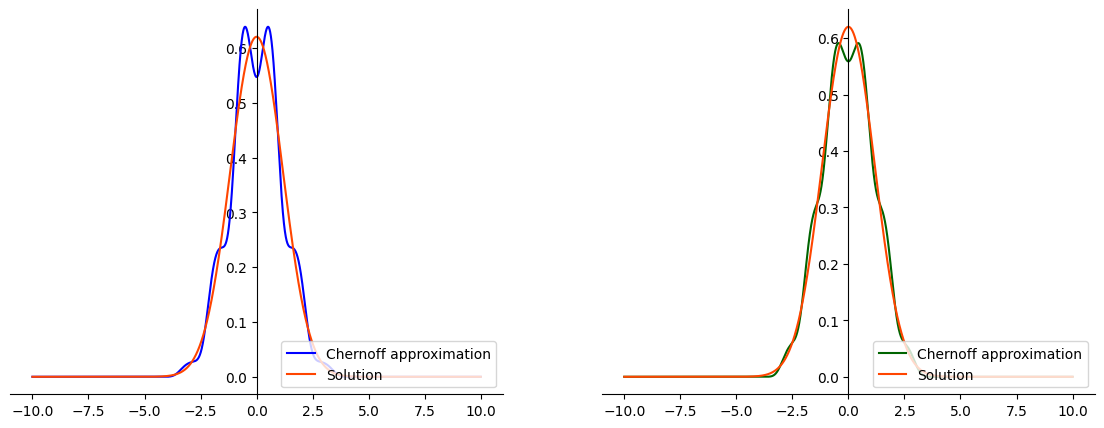}\\
\end{center}
n=3
\begin{center}
	\includegraphics[scale=0.5]{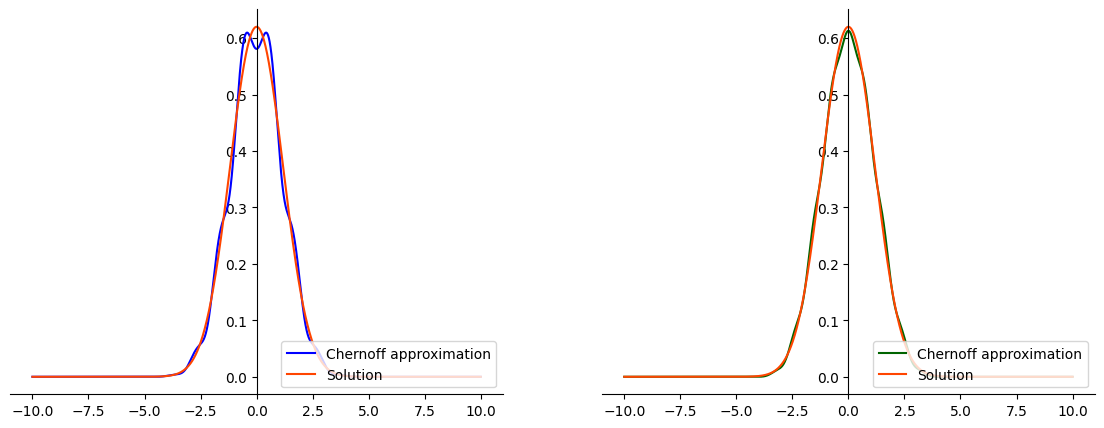}\\
\end{center}
n=4
\begin{center}
	\includegraphics[scale=0.5]{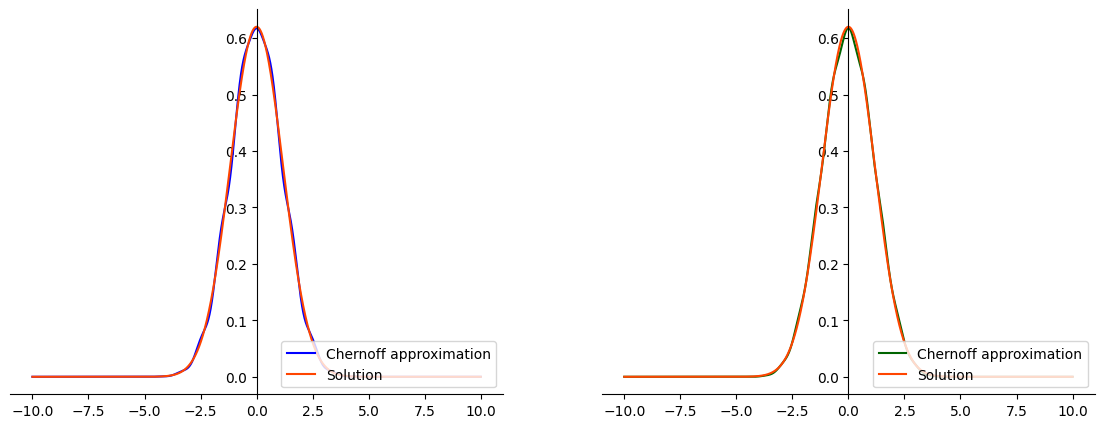}\\
\end{center}
n=5
\begin{center}
	\includegraphics[scale=0.5]{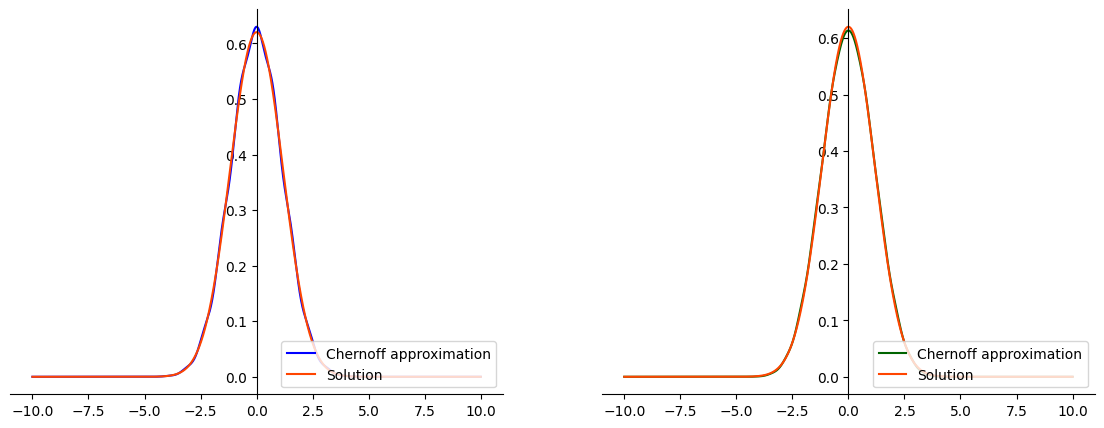}\\
\end{center}
n=6
\begin{center}
	\includegraphics[scale=0.5]{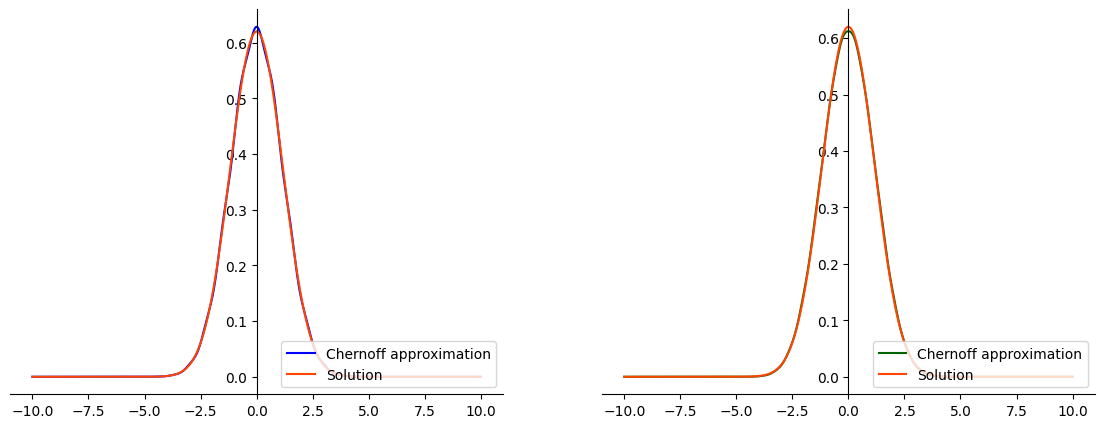}\\
\end{center}
n=7
\begin{center}
	\includegraphics[scale=0.5]{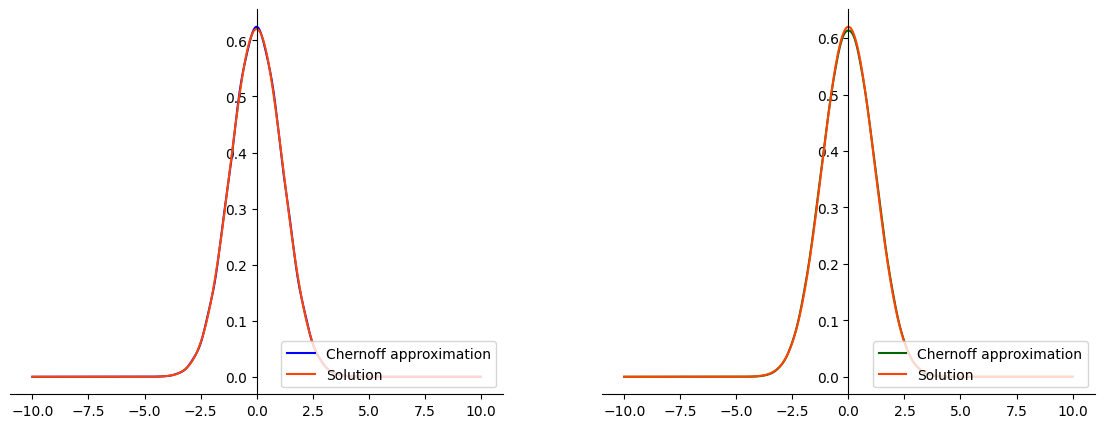}\\
\end{center}
n=8
\begin{center}
	\includegraphics[scale=0.5]{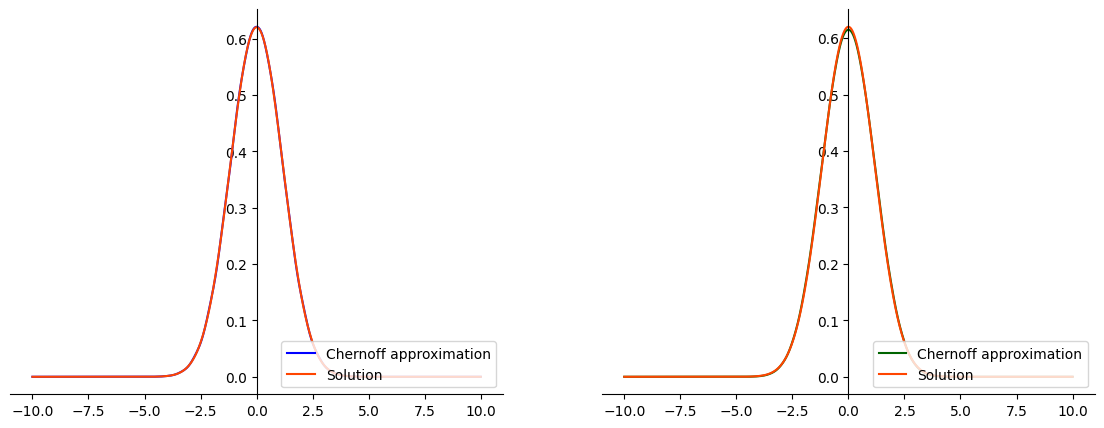}\\
\end{center}

n=9
\begin{center}
	\includegraphics[scale=0.5]{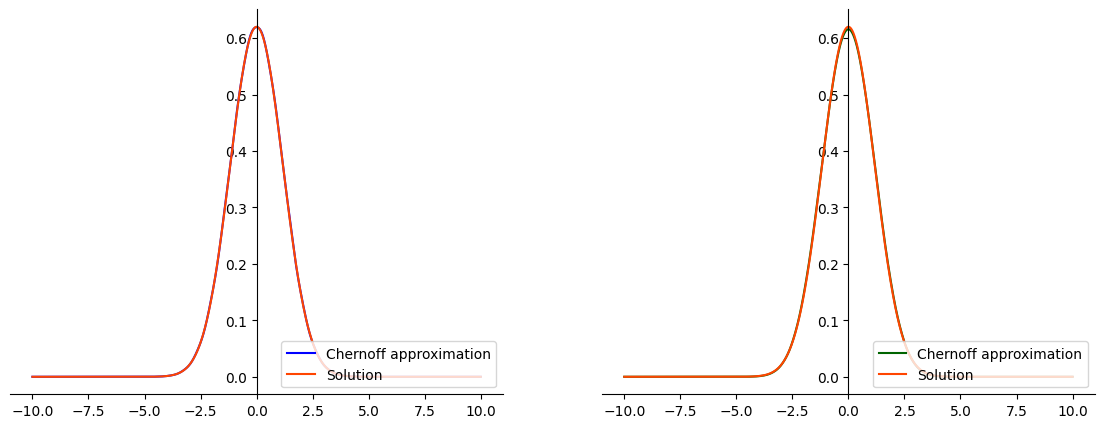}\\
\end{center}
n=10
\begin{center}
	\includegraphics[scale=0.5]{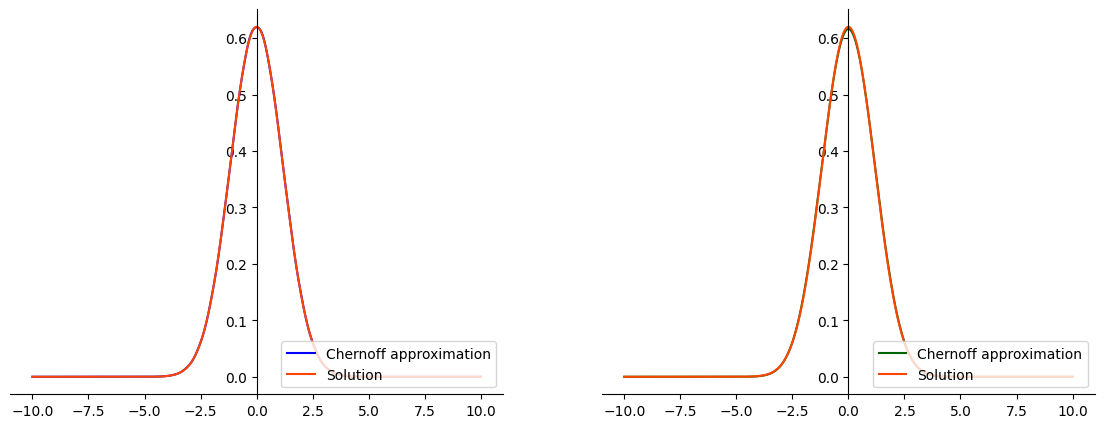}\\
\end{center}

\newpage
\subsection{$u_0(x)= e^{-x^6}$}

n=1
\begin{center}
	\includegraphics[scale=0.5]{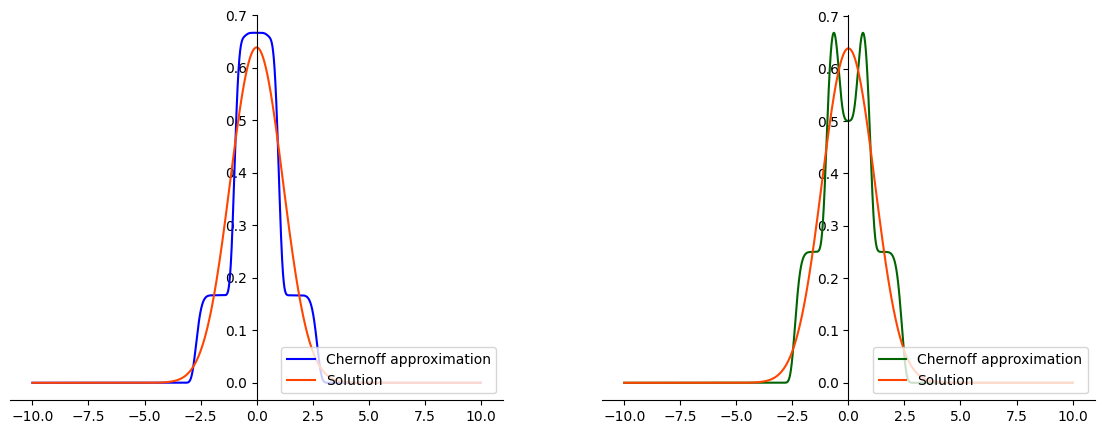}\\
\end{center}
n=2
\begin{center}
	\includegraphics[scale=0.5]{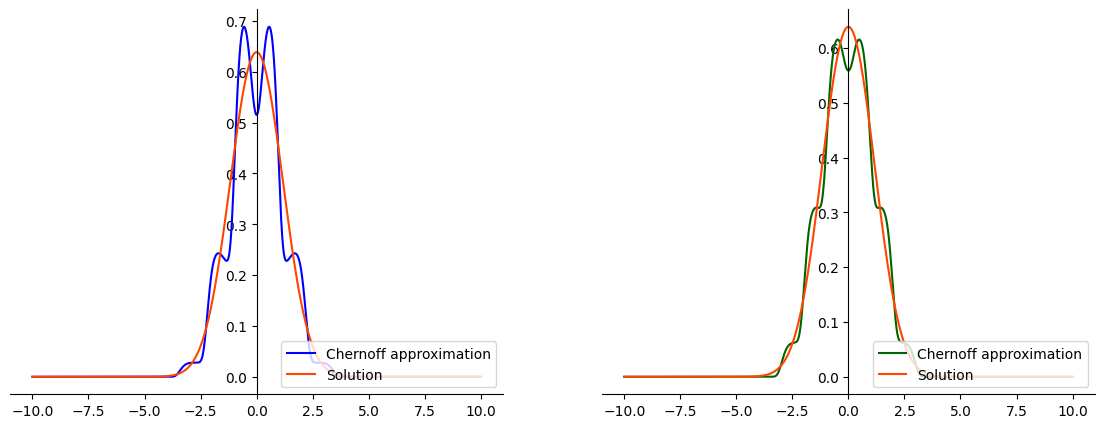}\\
\end{center}
n=3
\begin{center}
	\includegraphics[scale=0.5]{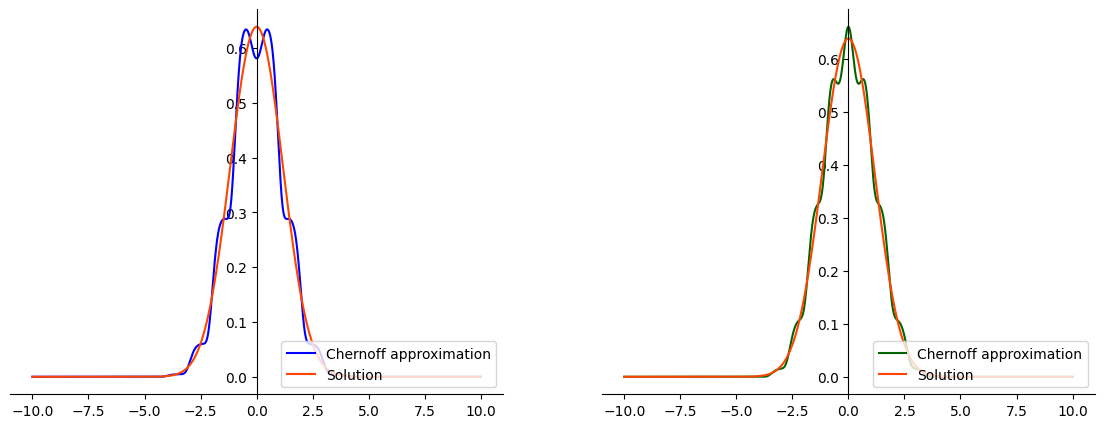}\\
\end{center}
n=4
\begin{center}
	\includegraphics[scale=0.5]{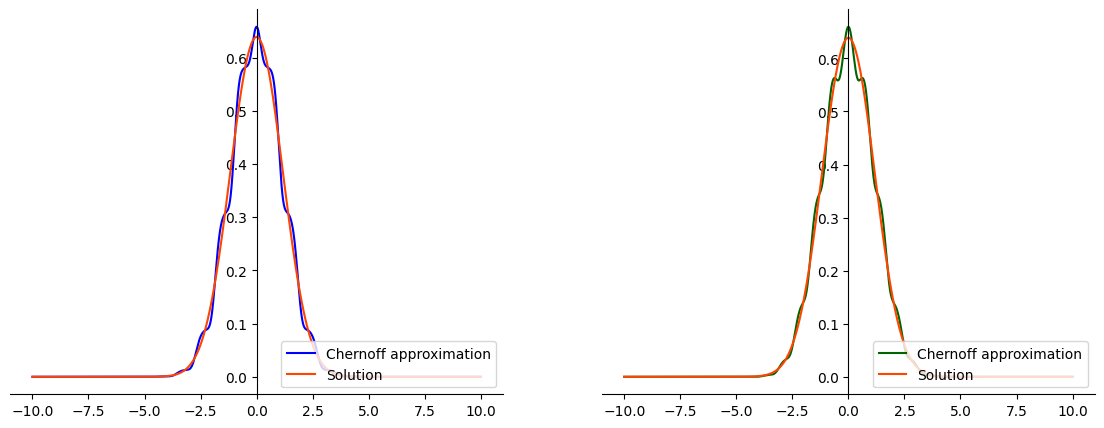}\\
\end{center}
n=5
\begin{center}
	\includegraphics[scale=0.5]{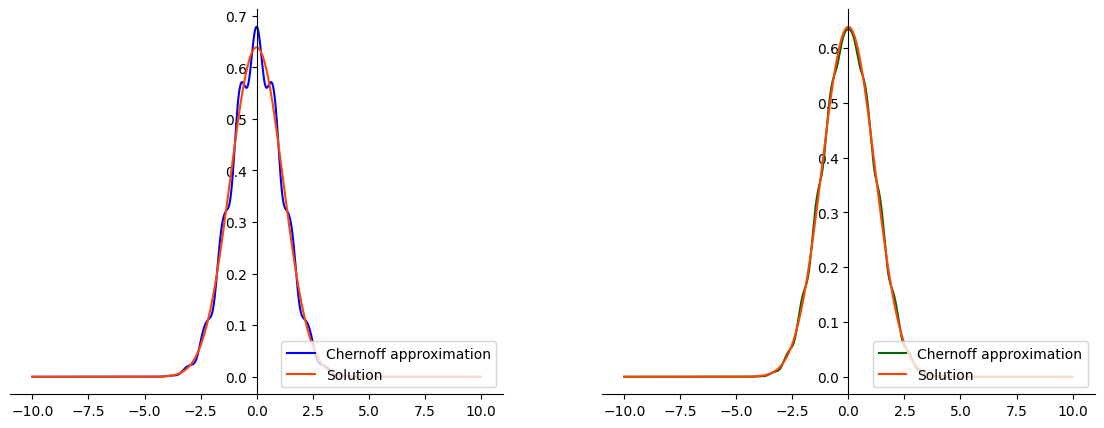}\\
\end{center}
n=6
\begin{center}
	\includegraphics[scale=0.5]{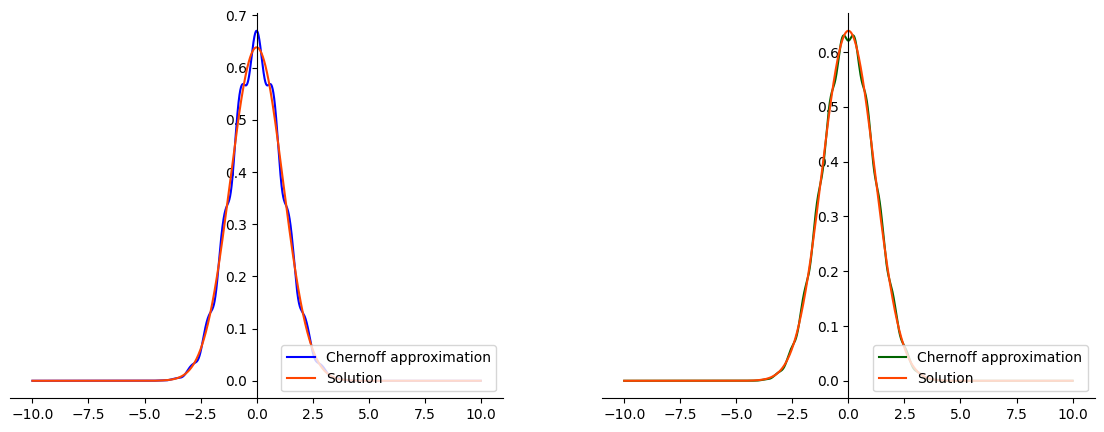}\\
\end{center}
n=7
\begin{center}
	\includegraphics[scale=0.5]{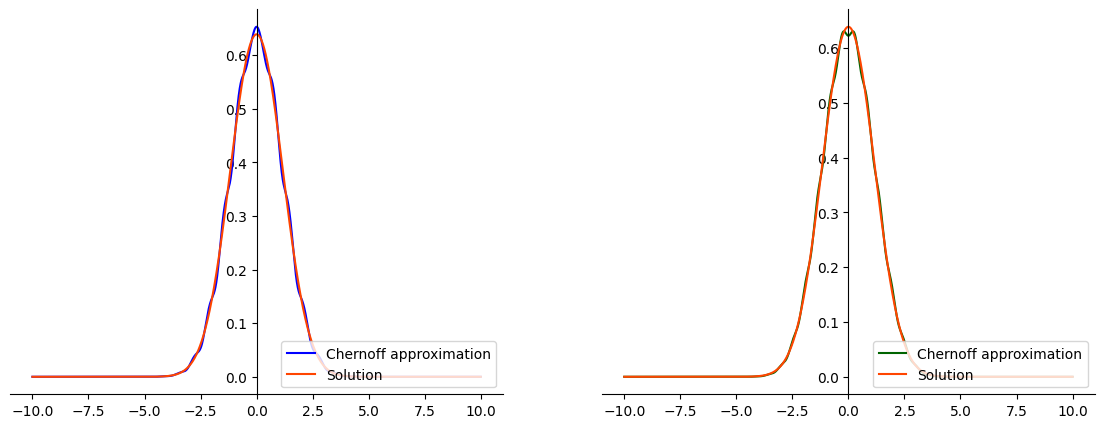}\\
\end{center}
n=8
\begin{center}
	\includegraphics[scale=0.5]{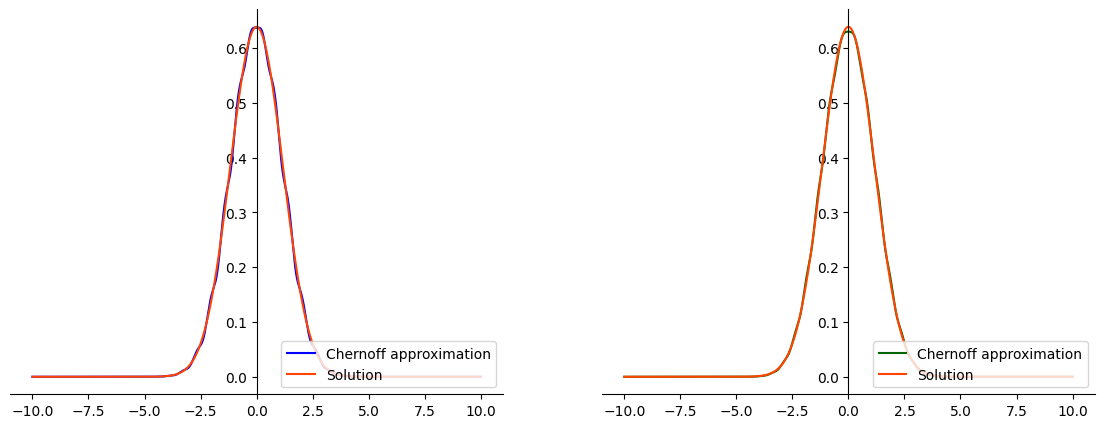}\\
\end{center}

n=9
\begin{center}
	\includegraphics[scale=0.5]{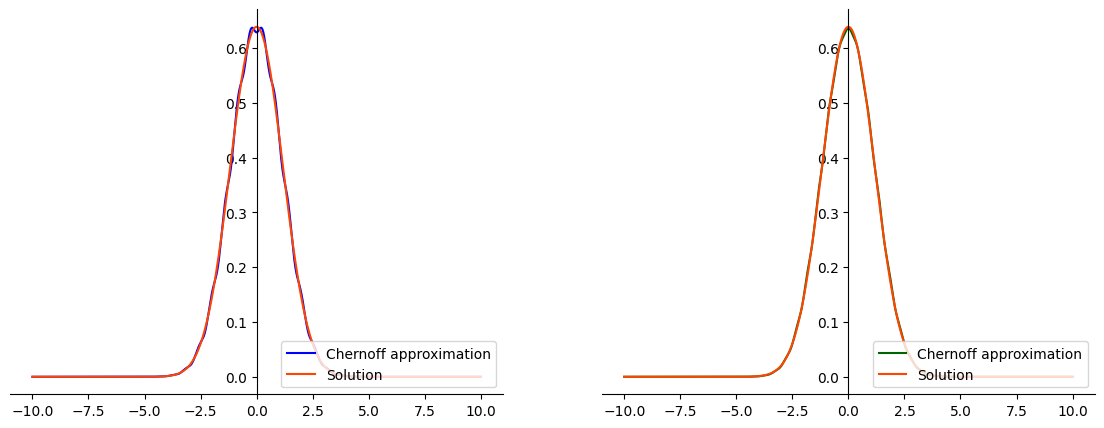}\\
\end{center}
n=10
\begin{center}
	\includegraphics[scale=0.5]{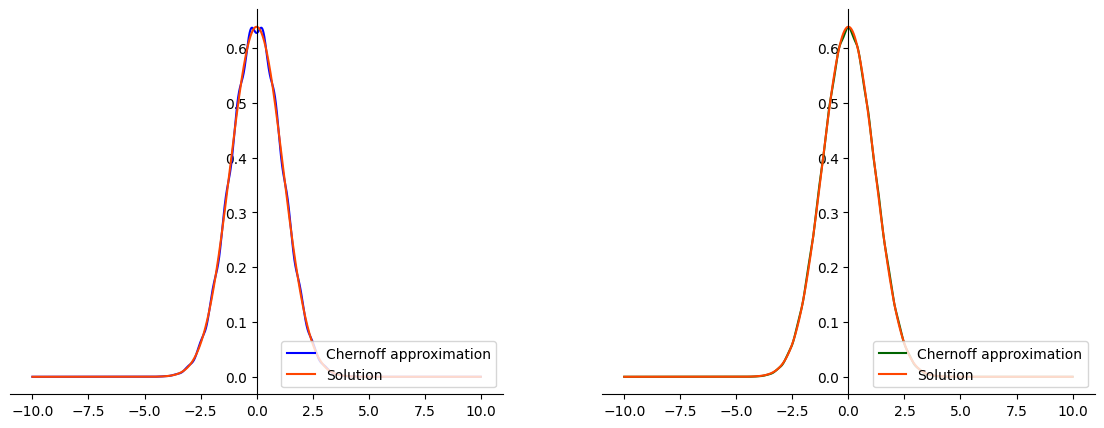}\\
\end{center}

\section*{Appendix B: Python 3 code}
\addcontentsline{toc}{section}{Appendix B: Python 3 code}

\begin{center}
	\includegraphics[scale=0.6]{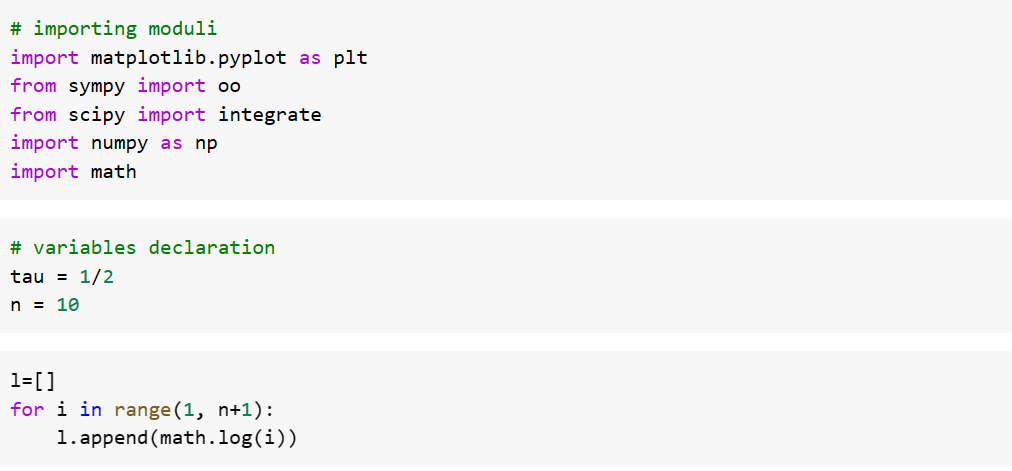}\\
\end{center}
\begin{center}
	\includegraphics[scale=0.6]{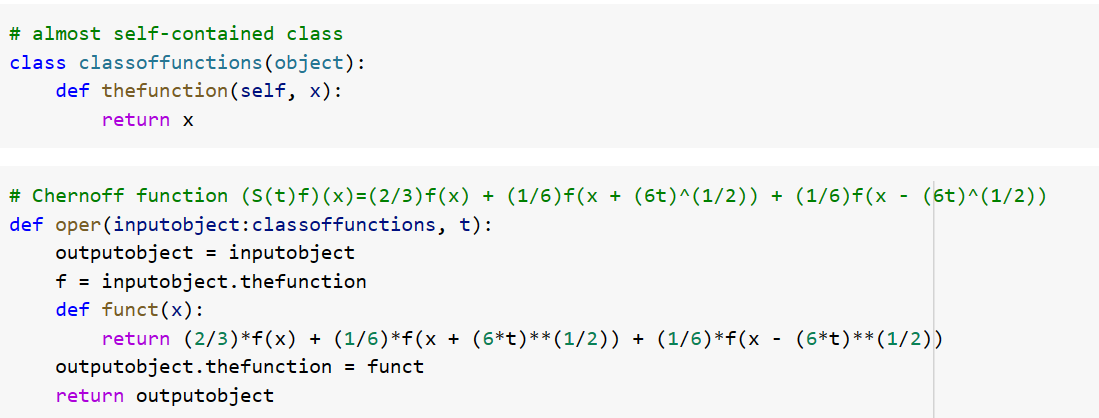}\\
\end{center}\begin{center}
	\includegraphics[scale=0.6]{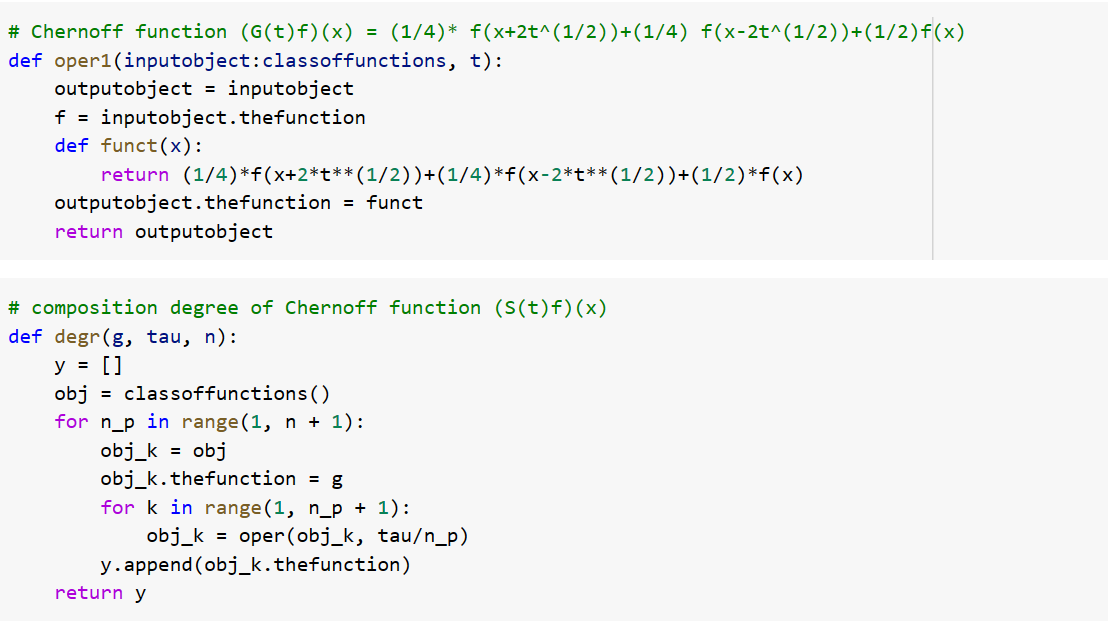}\\
\end{center}\begin{center}
	\includegraphics[scale=0.6]{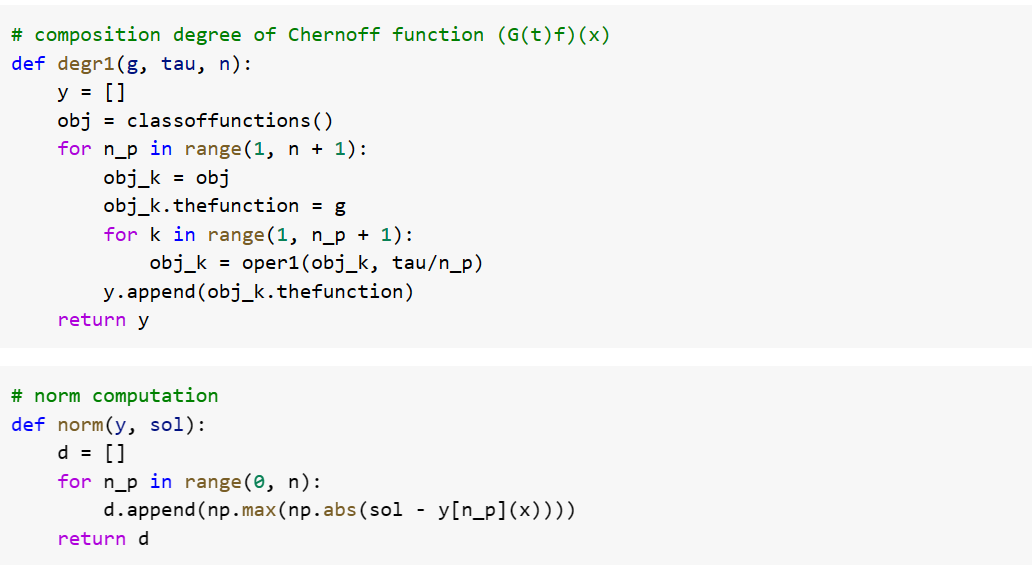}\\
\end{center}

\end{document}